\documentclass[12pt,dvips]{amsart} %% USE THIS FOR DOUBLE SIDED

\usepackage[dvips]{epsfig}
\usepackage{fancyhdr}
\usepackage{afterpage}

\usepackage{chngcntr}
\usepackage{graphicx}
\usepackage{latexsym}
\usepackage{amsfonts}
\usepackage{amscd}
\usepackage{amssymb}
\usepackage{amsmath}
\usepackage{amsthm}
\usepackage{bm}
\usepackage{psfrag}
\usepackage{url}
\usepackage{lscape}
\usepackage[usenames]{color}
\usepackage[all]{xy}
\usepackage{enumitem}
\usepackage{pinlabel}
\usepackage[lmargin=1.25in,rmargin=1.25in,tmargin=1.25in,bmargin=1.25in]{geometry}

\usepackage[unicode, linktocpage, bookmarksnumbered=true,backref=true]{hyperref}

\newtheorem{proposition}{Proposition} [section]
\newtheorem{lemma}[proposition]{Lemma}
\newtheorem{theorem}[proposition]{Theorem}
\newtheorem{corollary}[proposition]{Corollary}

\theoremstyle{definition}

\numberwithin{equation}{section}

\def\I{\mathbb{I}}
\def\N{\mathbb{N}}
\def\R{\mathbb{R}}

\def\Z{\mathbb{Z}}

\def\F{\mathbb{F}}

\def\one{\mathbf{1}}

\def\CF {\widehat{\operatorname{CF}}}
\def\CFK {\widehat{\operatorname{CFK}}}
\def\CFL {\widehat{\operatorname{CFL}}}
\def\CFD {\widehat{\operatorname{CFD}}}
\def\CFA {\widehat{\operatorname{CFA}}}
\def\CFAA {\widehat{\operatorname{CFAA}}}

\def\CFDD {\widehat{\operatorname{CFDD}}}
\def\CFAD {\widehat{\operatorname{CFAD}}}
\def\CFDA {\widehat{\operatorname{CFDA}}}
\def\HFK {\widehat{\operatorname{HFK}}}

\def\HF {\widehat{\operatorname{HF}}}

\def\AA {\mathcal{A}}
\def\BB {\mathcal{B}}

\def\FF {\mathcal{F}}
\def\HH {\mathcal{H}}
\def\II {\mathcal{I}}
\def\MM {\mathcal{M}}

\def\VV {\mathcal{V}}

\def\XX {\mathcal{X}}
\def\YY {\mathcal{Y}}
\def\ZZ {\mathcal{Z}}

\def\s {\mathfrak{s}}

\def\a{\mathbf{a}}
\def\x{\mathbf{x}}
\def\y{\mathbf{y}}
\def\z{\mathbf{z}}

\newcommand{\abs}[1] {\left\lvert #1 \right\rvert}

\newcommand{\gen}[1] {\left\langle #1 \right\rangle}
\def\Th{^{\text{th}}}
\def\minus{\smallsetminus} \def\co{\colon\thinspace}

 \DeclareMathOperator{\Spin}{Spin}
 \DeclareMathOperator{\ind}{ind}
 \DeclareMathOperator{\id}{id}

\newcommand{\twoline}[2]{\txt{ $\scriptstyle #1$ \\ $\scriptstyle
#2$}}
\newcommand{\threeline}[3]{\txt{ $\scriptstyle #1$ \\ $\scriptstyle
#2$ \\ $\scriptstyle #3$}}

\DeclareMathOperator{\hor}{hor} \DeclareMathOperator{\ver}{vert}
\DeclareMathOperator{\unst}{unst}

\definecolor{darkblue}{rgb}{0,0,0.5}
\definecolor{darkred}{rgb}{0.5,0,0}
\definecolor{darkgreen}{rgb}{0,0.5,0}
\definecolor{purple}{rgb}{0.5,0,0.5}
\definecolor{teal}{rgb}{0,0.5,0.5}

\def\Eta{\mathrm{H}}
\def\Kappa{\mathrm{K}}

\def\Oz{{}Ozsv\'ath{}}

\title[Doubling Operators and Bordered Floer Homology]{Knot Doubling Operators and Bordered Heegaard Floer Homology}
\author{Adam Simon Levine}
\address {Mathematics Department \\ Princeton University \\ Fine Hall, Washington Road \\ Princeton, NJ 08544}
\email {asl2@math.princeton.edu}

\thanks{The author was supported by NSF grants DMS-0739392 and DMS-1004622. This paper appears in \emph{Journal of Topology} \textbf{5} (2012) 651-712.}

\begin{document}

\maketitle

\begin{abstract}
We use bordered Heegaard Floer homology to compute the $\tau$ invariant of a
family of satellite knots obtained via twisted infection along two components
of the Borromean rings, a generalization of Whitehead doubling. We show that
$\tau$ of the resulting knot depends only on the two twisting parameters and
the values of $\tau$ for the two companion knots. We also include some notes on
bordered Heegaard Floer homology that may serve as a useful introduction to the
subject.
\end{abstract}

\tableofcontents

\section{Introduction} \label{sec:introduction}
A knot in the $3$-sphere is called \emph{topologically slice} if it bounds a
locally flatly embedded disk in the $4$-ball, and \emph{smoothly slice} if the
disk can be taken to be smoothly embedded. Two knots are called (topologically
or smoothly) \emph{concordant} if they are the ends of an embedded annulus in
$S^3 \times I$; thus, a knot is slice if and only if it is concordant to the
unknot. More generally, a link is (topologically or smoothly) \emph{slice} if
it bounds a disjoint union of appropriately embedded disks. The study of
concordance --- especially the relationship between the notions of topological
and smooth sliceness --- is one of the major areas of active research in knot
theory, and it is closely tied to the perplexing differences between
topological and smooth $4$-manifold theory.

\begin{figure} \centering
\psfrag{Wh+(K)}{$Wh_+(K)$} \psfrag{Wh-(K)}{$Wh_-(K)$} \psfrag{BD(K)}{$BD(K)$}
\psfrag{K}{$K$}
\includegraphics{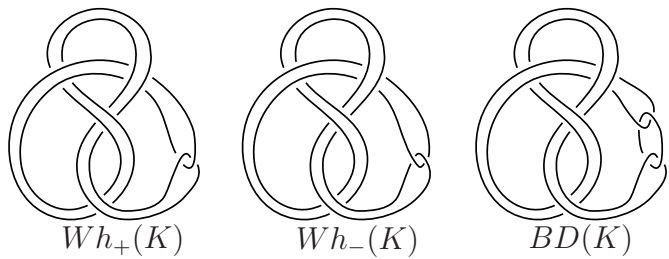}
\caption{The positive and negative Whitehead doubles and the Bing double of the
figure-eight knot.} \label{fig:figure8}
\end{figure}

While all known explicit constructions of slice disks use smooth techniques,
the early obstructions to sliceness --- including the Alexander polynomial, the
signature, J. Levine's algebraic concordance group, and Casson--Gordon
invariants --- arise from the algebraic topology of the complement of a slice
disk, so they only obstruct a knot from being topologically slice. However, in
the 1980s, Freedman \cite{FreedmanQuinn} showed that any knot whose Alexander
polynomial is $1$ is topologically slice, even though it is difficult to
describe the slice disks explicitly. In particular, the \emph{untwisted
positive and negative Whitehead doubles} of any knot $K$, denoted $Wh_\pm(K)$
(Figure \ref{fig:figure8}), are topologically slice. Around the same time, the
advent of Donaldson's gauge theory made it possible to show that some of these
examples are not smoothly slice. Akbulut [unpublished] first proved in 1983
that the positive, untwisted Whitehead double of the right-handed trefoil is
not smoothly slice. Later, using results of Kronheimer and Mrowka on
Seiberg--Witten theory, Rudolph \cite{RudolphObstruction} showed that any
nontrivial knot that is \emph{strongly quasipositive} cannot be smoothly slice.
In particular, the positive, untwisted Whitehead double of a strongly
quasipositive knot is strongly quasipositive; thus, by induction, any iterated
positive Whitehead double of a strongly quasipositive knot is topologically but
not smoothly slice. Bi\v{z}aca \cite{Bizaca} used this result to give explicit
constructions of exotic smooth structures on $\R^4$.

Using Heegaard Floer homology, Ozsv\'ath and Szab\'o \cite{OSz4Genus} defined
an additive, integer-valued knot invariant $\tau(K)$, defined as the minimum
Alexander grading of an element of $\HFK(S^3,K)$ that survives to the
$E^\infty$ page of the spectral sequence from $\HFK(S^3,K)$ to $\HF(S^3)$. The
$\tau$ invariant provides a lower bound on the genus of smooth surfaces in the
four-ball bounded by $K$: $\abs{\tau(K)} \le g_4(K)$. In particular, if $K$ is
smoothly slice, then $\tau(K)=0$. This fact can be used to generalize many of
the previously known results about knots that are topologically but not
smoothly slice. For example, Hedden \cite{HeddenWhitehead} computed the value
of $\tau$ for all \emph{twisted Whitehead doubles} in terms of $\tau$ of the
original knot:
\begin{equation} \label{eq:tau-Wh}
\tau(Wh_+(K,t)) = \begin{cases} 1 & t < 2 \tau(K) \\ 0 & t \ge 2 \tau(K).
\end{cases}
\end{equation}
(An analogous formula for negative Whitehead doubles follows from the fact that
$\tau(\bar K) = - \tau(K)$.) In particular, if $\tau(K)>0$, then
$\tau(Wh_+(K,0))=1$, so $Wh_+(K,0)$ (the untwisted Whitehead double of $K$) is
not smoothly slice. Since $\tau$ of a strongly quasipositive knot is equal to
its genus \cite{LivingstonComputations}, Rudolph's result follows from
Hedden's. There is a famous conjecture (Problem 1.38 on Kirby's problem list
\cite{KirbyList}) that the untwisted Whitehead double of $K$ is smoothly slice
if and only if $K$ is smoothly slice. However, it is not yet known whether, for
instance, the positive Whitehead double of the left-handed trefoil is smoothly
slice. Indeed, it seems that gauge theory invariants have a fundamental
asymmetry that makes them unable to detect such examples, which likely places
the ``only if'' direction of this conjecture beyond the scope of currently
existing techniques.

\begin{figure} \centering
\psfrag{(a)}{(a)} \psfrag{(b)}{(b)} \psfrag{J,s}[cc][cc]{$J,s$}
\psfrag{K,t}[cc][cc]{$K,t$}
\includegraphics{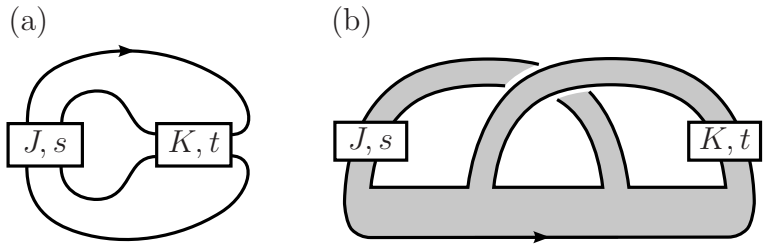}
\caption[The knot $D_{J,s}(K,t)$, along with a Seifert surface.] {(a) The knot
$D_{J,s}(K,t)$. (b) A genus-$1$ Seifert surface for $D_{J,s}(K,t)$.}
\label{fig:djskt}
\end{figure}

We consider the following generalization of Whitehead doubling. For knots $J,K$
and integers $s,t$, let $D_{J,s}(K,t)$ denote the knot shown in Figure
\ref{fig:djskt}(a); the box marked $K,t$ (resp.~$J,s$) indicates that the
strands are tied along $t$-framed (resp.~$s$-framed) parallel copies of the
tangle $K \minus \{\text{pt}\}$ (resp.~$J \minus \{\text{pt}\}$). (We give a
more formal definition below.) If $J$ is the unknot and $s = \pm 1$, then
$D_{J,s}(K,t)$ is the $t$-twisted $\mp$ Whitehead double of $K$.

A genus-$1$ Seifert surface for $D_{J,s}(K,t)$ is shown in Figure
\ref{fig:djskt}(b). From the Seifert form of this surface, we can compute that
the Alexander polynomial of $D_{J,s}(K,t)$ is
\[
\Delta_{D_{J,s}(K,t)}(T) = stT + (1-2st) + st T^{-1}.
\]
In particular, this equals $1$ whenever $s=0$ or $t=0$. By Freedman's theorem,
$D_{J,s}(K,0)$ is therefore topologically slice. Moreover, if $K$ is smoothly
slice, then $D_{J,s}(K,0)$ is smoothly slice for any $(J,s)$. To see this,
perform a ribbon move to eliminate the band that is tied into $J$; the
resulting two-component link, consisting of two parallel copies of $K$ with
linking number $0$, is then the boundary of two parallel copies of a slice disk
for $K$. The conjecture about sliceness of untwisted Whitehead doubles described above has
many potential generalizations in terms of $D_{J,s}(K,0)$ satellites, all
apparently equally difficult.

As a partial result in this direction, we prove the following theorem, which
generalizes Hedden's result:
\begin{theorem} \label{thm:taudjskt}
Let $J$ and $K$ be knots, and let $s,t \in \Z$. Then
\[
\tau(D_{J,s}(K,t)) =
\begin{cases}
1 & s<2\tau(J) \text{ and } t<2 \tau(K) \\
-1 & s>2\tau(J) \text{ and } t>2 \tau(K) \\
0 & \text{otherwise}.
\end{cases}
\]
\end{theorem}
In particular, if $\tau(K)>0$ and $s<2\tau(J)$, or if $\tau(K)<0$ and
$s>2\tau(J)$, then $D_{J,s}(K,0)$ is topologically but not smoothly slice.

We now provide a more rigorous description of $D_{J,s}(K,t)$. Suppose $L$ is a
link in $S^3$, and $\gamma$ is an oriented curve in $S^3 \minus L$ that is
unknotted in $S^3$. For any knot $K \subset S^3$ and $t \in \Z$, we may form a
new link $I_{\gamma,K,t}(L)$, the \emph{$t$-twisted infection of $L$ by $K$
along $\gamma$}, by deleting a neighborhood of $\gamma$ and gluing in a copy of
the exterior of $K$ by a map that takes a Seifert-framed longitude of $K$ to a
meridian of $\gamma$ and a meridian of $K$ to a $t$-framed longitude of
$\gamma$. Since $S^3 \minus \gamma = S^1 \times D^2$, the resulting
$3$-manifold is simply $\infty$ surgery on $K$, i.e. $S^3$; the new link
$I_{\gamma,K,t}(L)$ is defined as the image of $L$. Infecting along the
boundary of a disk perpendicular to a group of strands formalizes the notion of
``tying the strands into a knot.'' Moreover, given an unlink $\gamma_1, \dots,
\gamma_n$ disjoint from $L$, we may infect simultaneously along all the
$\gamma_i$; the result may be denoted $I_{\gamma_1,K_1,t_1; \, \cdots; \,
\gamma_n,K_n,t_n}(L)$, and the order of the tuples $(\gamma_i,K_i,t_i)$ does
not matter.

\begin{figure}
\psfrag{B1}{$B_1$} \psfrag{B2}{$B_2$} \psfrag{B3}{$B_3$} \centering
\includegraphics{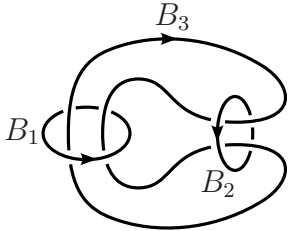}
\caption{The Borromean rings.} \label{fig:borromean}
\end{figure}

Let $B = B_1 \cup B_2 \cup B_3$ denote the Borromean rings, oriented as shown
in Figure \ref{fig:borromean}. Then $D_{J,s}(K,t)$ is the knot obtained from
$B_3$ by $s$-twisted infection by $J$ along $B_1$ and $t$-twisted infection by
$K$ along $B_2$:
\[
D_{J,s}(K,t) = I_{B_1,J,s; \, B_2,K,t}(B_3).
\]

The theory of \emph{bordered Heegaard Floer homology}, developed recently by
Lipshitz, \Oz, and Thurston \cite{LOTBordered, LOTBimodules}, is well-suited to
the problem of computing Heegaard Floer invariants of knots obtained via
infection. Briefly, it associates to a $3$-manifold with boundary a module over
an algebra associated to the boundary, so that if a $3$-manifold $Y$ is
decomposed as $Y = Y_1 \cup_\phi Y_2$, where $\phi:
\partial Y_1 \xrightarrow{\cong}
\partial Y_2$, the chain complex $\CF(Y)$ may be computed as the derived tensor
product of the invariants associated to $Y_1$ and $Y_2$. If a knot $K$ is
contained (nulhomologously) in, say, $Y_1$, then we may obtain the filtration
on $\CF(Y)$ corresponding to $K$ via a filtration on the algebraic invariant of
$Y_1$. The theory also includes bimodules associated to manifolds with two
boundary components.

In our setting, let $Y$ denote the exterior of $B_1 \cup B_2$, and let $X_J$
and $X_K$ denote the exteriors of $J$ and $K$, respectively. For suitable
gluing maps $\phi_1: \partial X_J \to \partial_1 Y$ and $\phi_2:
\partial X_K \to \partial_2 Y$ (where $\partial Y = \partial_1 Y \amalg \partial_2
Y$), the glued-up manifold $(Y \cup_{\phi_1} X_J) \cup_{\phi_2} X_K$ is $S^3$,
and the image of $B_3 \subset Y$ is $D_{J,s}(K,t)$. We shall define suitable
\emph{bordered structures} $\YY$, $\XX_J^s$, and $\XX_K^t$ on $Y$, $X_J$, and
$X_K$, respectively, so as to induce these gluing maps. By the \emph{gluing
theorem} of Lipshitz, \Oz, and Thurston, the filtered chain complex for $(S^3,
D_{J,s}(K,t))$ can then computed as a special tensor product of the modules
associated to $\YY$, $\XX_J^s$, and $\XX_K^t$:
\[
\CF(S^3,D_{J,s}(K,t)) \simeq (\CFAA(\YY,B_3,0) \boxtimes \CFD(\XX_J^s))
\boxtimes \CFD(\XX_K^t).
\]
All of this terminology will be explained in Section
\ref{sec:bordered-background}. Using the formula for $\CFD(\XX_J^s)$ and
$\CFD(\XX_K^t)$ proven by Lipshitz, \Oz, and Thurston \cite{LOTBordered} and a
direct computation of $\CFAA(\YY,B_3,0)$ using holomorphic disks (given in
Section \ref{sec:heegaard}), we shall explicitly evaluate this double tensor
product and compute its homology in Section \ref{sec:tensorproduct}, leading to
the proof of Theorem \ref{thm:taudjskt}. While the proof is fairly technical,
it illustrates the power of the new bordered techniques: using a single
computation involving holomorphic disks (which can in principle be performed
entirely combinatorially) and some lengthy but straightforward algebra, we are
able to obtain a statement about the Floer homology an infinite family of
knots. The proof relies on some computer-assisted computation using \emph{Mathematica}; the details are described in Appendix \ref{sec:appendix}.

In Section \ref{sec:doubling-topology}, we present a few more results
concerning knots of the form $D_{J,s}(K,t)$. Following the approach of
Livingston and Naik \cite{LivingstonNaikDoubled}, we show that if $\nu$ is any
concordance invariant that shares certain formal properties with $\tau$ ---
e.g., Rasmussen's $s$ invariant coming from Khovanov homology --- then
$\nu(D_{J,s}(K,t)) = \tau(D_{J,s}(K,t))$ when $\abs{s}$ and $\abs{t}$ are
sufficiently large. We also provide a family of examples of knots of the form
$D_{J,s}(K,t)$ that are smoothly slice, generalizing a result of Casson about
Whitehead doubles.

Finally, Theorem \ref{thm:taudjskt} has a useful application to the study of
Whitehead doubles of links, which was the author's original motivation for
considering it. Specifically, we consider the Whitehead doubles of links
obtained by \emph{iterated Bing doubling}. Given a knot $K$, the
\emph{(untwisted) Bing double} of $K$ is the satellite link $BD(K) =
I_{B_1,K,0}(B_2 \cup B_3)$, as shown in Figure \ref{fig:figure8}. More
generally, given a link $L$, we may replace a component by its Bing double
(contained in a tubular neighborhood of that component), and iterate this
procedure. Bing doubling one component of the Hopf link yields the Borromean
rings; accordingly, we define the family of \emph{generalized Borromean links}
as the set of all links obtained as iterated Bing doubles of the Hopf link.
Using Theorem \ref{thm:taudjskt}, the author proves in
\cite{LevineBingWhitehead}:
\begin{corollary} \label{cor:bingwhitehead}
Let $L$ be any link obtained by iterated Bing doubling from either:
\begin{enumerate}
\item Any knot $K$ with $\tau(K)>0$, or
\item The Hopf link.
\end{enumerate}
Then the all-positive Whitehead double of $L$, $Wh_+(L)$, is not smoothly
slice.
\end{corollary}
The links in (1) are boundary links, so their Whitehead doubles are all
topologically slice by a result of Freedman \cite{FreedmanWhitehead3}. On the
other hand, it is not yet known whether the Whitehead doubles of iterated Bing
doubles of the Hopf link are topologically slice; indeed, this question is one
of the major unsolved problems in four-dimensional topological surgery theory.
Once again, we see a strong dependence on chirality; our proof breaks down when
clasps of both signs are used. For further details, see
\cite{LevineBingWhitehead}.

\subsubsection*{Acknowledgements}
A version of this paper made up a large portion of the author's thesis at Columbia
University. The author is grateful to his advisor, Peter Ozsv\'ath, and the
other members of his defense committee, Robert Lipshitz, Dylan Thurston, Paul
Melvin, and Denis Auroux, for their suggestions, and to Rumen Zarev, Ina
Petkova, Jen Hom, and Matthew Hedden for many helpful conversations about
bordered Heegaard Floer homology. Additionally, he thanks the Mathematical
Sciences Research Institute for hosting him in Spring 2010, when much of this
research was conducted.

\section{Background on bordered Heegaard Floer homology}
\label{sec:bordered-background}
In this section, we give a brief description of the bordered Heegaard Floer
invariants, with the aim of defining the terms used later in the paper and
illustrating the procedures for computation. We discuss only bordered manifolds
with toroidal boundary components, which has the advantage of greatly
simplifying some of the definitions. All of this material can be found in the
two \emph{magna opera} of Lipshitz, Ozsv\'ath, and Thurston \cite{LOTBordered,
LOTBimodules}.

\subsection{Algebraic objects} \label{subsec:objects}

We first recall the main algebraic constructions used in \cite{LOTBordered,
LOTBimodules}, with the aim of describing how to work with them
computationally. Let $(\AA, d)$ be a unital differential algebra over $\F =
\F_2$, and assume that the set of $\II$ of idempotents in $\AA$ is a commutative subring of $\AA$ and possesses a basis $\{\iota_i\}$ over $\F$ such that $\iota_i \iota_j = \delta_{ij} \iota_i$ and $\sum_i \iota_i = \one$, the identity element of $\AA$. (All of the definitions that follow can be stated in terms of
differential graded algebras, but we suppress all grading information for
brevity.)

\begin{itemize}
\item A \emph{(right) $\AA_\infty$ module} or \emph{type $A$ structure} over $\AA$ is an $\F$-vector space $M$,
equipped with a right action of $\II$ such that $M = \bigoplus_i M \iota_i$ as a vector space, and multiplication maps
\[
m_{k+1}\co M \underset{\II}{\otimes} \underbrace{\AA \underset{\II}{\otimes}
\dots \underset{\II}{\otimes} \AA}_{k \text{ times}} \to M
\]
satisfying the $\AA_\infty$ relations: for any $x \in M$ and $a_1, \dots, a_n \in \AA$,
\begin{equation} \label{eq:a-relation}
\begin{split}
0 & = \sum_{i=0}^n m_{n-i+1}(m_{i+1}(x,a_1, \dots, a_i), a_{i+1}, \dots, a_n) \\
& + \sum_{i=1}^{n} m_{n+1}(x, a_1, \dots, a_{i-1}, d(a_i), a_{i+1}, \dots, a_n) \\
& + \sum_{i=1}^{n-1} m_n(x, a_1, \dots, a_{i-1}, a_i a_{i+1}, a_{i+2}, \dots, a_n).
\end{split}
\end{equation}
We also require that $m_2(x,\one)=x$ and $m_k(x, \dots, \one, \dots)=0$ for
$k>2$.

The module $M$ is called \emph{bounded} if $m_k=0$ for all $k$ sufficiently
large. If $M$ is a bounded type $A$ structure with basis $\{x_1,\dots,x_n\}$,
we encode the multiplications using a matrix whose entries are formal sums of
finite sequences of elements of $\AA$, where having an $(a_1, \dots, a_k)$ term
in the $i,j\Th$ entry means that the coefficient of $x_j$ in $m_{k+1}(x_i, a_1,
\dots, a_k)$ is nonzero. We write $1$ rather than an empty sequence to signify
the $m_1$ multiplication. For brevity, we frequently write $a_1 \cdots a_k$
rather than $(a_1, \dots, a_k)$; in this context, concatenation is \emph{not}
interpreted as multiplication in the algebra $\AA$.

\item A \emph{(left) type $D$ structure} over $\AA$ is an $\F$-vector space
$N$, equipped with a left action of $\II$ such that $N = \bigoplus_i \iota_i N$, and a map
\[
\delta_1\co N \to A \underset{\II}{\otimes} N
\]
satisfying the relation
\begin{equation} \label{eq:d-relation}
(\mu \otimes \id_N) \circ (\id_\AA \otimes \delta_1) \circ \delta_1 + (d \otimes \id_N) \circ \delta_1 = 0,
\end{equation}
where $\mu\co \AA \otimes \AA \to \AA$ denotes the multiplication on $\AA$.

If $N$ is a type $D$ structure, the tensor product $\AA \underset{\II}{\otimes}
N$ is naturally a left differential module over $\AA$, with module structure
given by $a \cdot (b \otimes x) = ab \otimes x$, and differential $\partial(a
\otimes x) = a \cdot \delta_1(x) + d(a) \otimes x$. Condition
\eqref{eq:d-relation} translates to $\partial^2=0$.

Given a type-$D$ module $N$, define maps
\[
\delta_k\co N \to \underbrace{\AA \underset{\II}{\otimes} \dots
\underset{\II}{\otimes} \AA}_{k \text{ times}} \underset{\II}{\otimes} N
\]
by $\delta_0 = \id_N$ and $\delta_k = (\id_{\AA^{\otimes k-1}} \otimes
\delta_1) \circ \delta_{k-1}$. We say $N$ is \emph{bounded} if $\delta_k=0$ for
all $k$ sufficiently large.

Given a basis $\{y_1, \dots, y_n\}$ for $N$, we may encode $\delta_1$ as an $n
\times n$ matrix $(b_{ij})$ with entries in $\AA$, such that $\delta_1 x_i =
\sum_{j=1}^n b_{ij} \otimes x_j$. To encode $\delta_k$ in matrix form, we take
the $k\Th$ power of the matrix for $\delta_1$, except that instead of
evaluating multiplication in $\AA$, we simply concatenate tensor products of
elements.

If $d=0$, \eqref{eq:d-relation} is equivalent to the statement that the square
of the matrix for $\delta_1$ (where now we do evaluate multiplication in $\AA$)
is zero.

\item
If $M$ is a right type $A$ structure, $N$ is a left type $D$ structure, and at
least one of them is bounded, we may form the \emph{box tensor product} $M
\boxtimes N$. As a vector space, this is $M \underset{\II}{\otimes} N$, with
differential
\[
\partial^\boxtimes(x \otimes y) = \sum_{k=0}^\infty (m_{k+1} \otimes \id_N)(x
\otimes \delta_k(y)).
\]
Given matrix representations of the multiplications on $M$ and the $\delta_k$
maps on $N$, it is easy to write down the differential on $M
\boxtimes N$ in terms of the basis $\{x_i \otimes y_j\}$.

\item Now let $(\AA, d_\AA)$ and $(\BB, d_\BB)$ be differential algebras. Lipshitz, Ozsv\'ath, and
Thurston \cite{LOTBimodules} define various types of $(\AA, \BB)$-bimodules. We
do not define these in full detail, but we mention some of the basic notions.

A \emph{type $DD$ structure} is simply a type $D$ structure over the ring $\AA
\underset{\F}{\otimes} \BB$. That is, the map $\delta_1$ outputs terms of the
form $a \otimes b \otimes x$, where $a \in \AA$ and $b \in \BB$.

A \emph{type $AA$ structure} consists of a vector space $M$ with multiplications
\[
m_{1,i,j}\co M \otimes \AA^{\otimes i} \otimes \BB^{\otimes j} \to M,
\]
satisfying a version of the $\AA_\infty$ relation \eqref{eq:a-relation}. As
above, all tensor products are taken over the rings of idempotents, $\II_\rho
\subset \AA$ and $\II_\sigma \subset \BB$. Our notation differs a bit from that
of \cite{LOTBimodules} in that we think of both algebras as acting on the
right.

A \emph{type $DA$ structure} is a vector space $N$ with maps
\[
\delta_1^{1+j}\co N \otimes \BB^{\otimes j} \to \AA \otimes N
\]
satisfying an appropriate relation that combines \eqref{eq:a-relation} and
\eqref{eq:d-relation}. A \emph{type $AD$ structure} is defined similarly,
except that the roles of $\AA$ and $\BB$ are interchanged.

The box tensor product of two bimodules, or of a module and a bimodule, can be
defined assuming at least one of the factors is bounded (in an appropriate
sense). See \cite[Subsection 2.3.2]{LOTBimodules} for details.

\item
A filtration on a type $A$ structure $M$ is a filtration $\dots \subseteq \FF_p
\subseteq \FF_{p+1} \subseteq \dots$ of $M$ as a vector space, such that
$m_{k+1}(\FF_p \otimes \AA^{\otimes k}) \subseteq \FF_p$. Similarly, a filtration on a type $D$ structure $N$ is a filtration of $N$ such that $\delta_1(\FF_p) \subseteq \AA \otimes \FF_p$. It is easy to
extend these definitions to the various types of bimodules. A filtration on $M$
or $N$ naturally induces a filtration on $M \boxtimes N$.

\end{itemize}

\subsection{The torus algebra}
\label{subsec:torus}

The \emph{pointed matched circle for the torus}, $\ZZ$, consists of an oriented
circle $Z$, equipped with a basepoint $z \in Z$, a tuple $\a =
(a_1,a_2,a_3,a_4)$ of points in $Z \minus \{z\}$ (ordered according to the
orientation on $Z \minus \{z\}$), and the equivalence relation $a_1 \sim a_3$,
$a_2 \sim a_4$. The genus-$1$, one-boundary-component surface $F^\circ(\ZZ)$ is
obtained by identifying $Z$ with the boundary of a disk $D$ and attaching
$1$-handles $h_1$ and $h_2$ that connect $a_1$ to $a_3$ and $a_2$ to $a_4$,
respectively. By attaching a $2$-handle along $\partial F^\circ(\ZZ)$, we
obtain the closed surface $F(\ZZ)$. There is an orientation-reversing
involution $r\co Z \to Z$ that fixes $z$, interchanges $a_1$ and $a_4$, and
interchanges $a_2$ and $a_3$, which extends to a diffeomorphism $r\co F(\ZZ)
\to -F(\ZZ)$ that interchanges $h_1$ and $h_2$.

The algebra $\AA = \AA(\ZZ,0)$ is generated as a vector space over $\F$ by two
\emph{idempotents} $\iota_0, \iota_1$ and six \emph{Reeb elements} $\rho_1,
\rho_2, \rho_3, \rho_{12}, \rho_{23}, \rho_{123}$. For each sequence of
consecutive integers $I = (i_1,\dots,i_k) \subset \{1,2,3\}$, we have
$\iota_{[i_1-1]} \rho_I = \rho_I \iota_{[i_k]} = \rho_I$, where $[j]$ denotes
the residue of $j$ modulo $2$. The nonzero multiplications among the Reeb
elements are: $\rho_1 \rho_2 = \rho_{12}$, $\rho_2 \rho_3 = \rho_{23}$, $\rho_1
\rho_{23} = \rho_{12} \rho_3 = \rho_{123}$. All other products are zero, as is
the differential. Let $\II$ denote the subring of idempotents of $\AA$; it is
generated as a vector space by $\iota_0$ and $\iota_1$. The identity element is
$\one = \iota_0 + \iota_1$.

By abuse of notation, we identify $\rho_1$ with the oriented arc of $Z$ from
$a_1$ to $a_2$, $\rho_2$ with the arc from $a_2$ to $a_3$, $\rho_3$ with the
arc from $a_3$ to $a_4$, and $\rho_{12}$, $\rho_{23}$, and $\rho_{123}$ with
the appropriate concatenations.

\subsection{Bordered 3-manifolds and their invariants}
\label{subsec:bordered}

A \emph{bordered $3$-manifold with boundary $F(\ZZ)$} consists of the data $\YY
= (Y,\Delta,z',\phi)$, where $Y$ is an oriented $3$-manifold with a single
boundary component, $\Delta$ is a disk in $\partial Y$, $z' \in
\partial \Delta$, and $\phi: F(\ZZ) \to \partial(Y)$ is a diffeomorphism
taking $D$ to $\Delta$ and $z$ to $z'$. The map $\phi$ is specified (up to
isotopy fixing $\Delta$ pointwise) by the images of the core arcs of the two
one-handles in $F^\circ(\ZZ)$. We may analogously define a \emph{bordered
$3$-manifold with boundary $-F(\ZZ)$}. The diffeomorphism $r\co F(\ZZ) \to
-F(\ZZ)$ provides a one-to-one correspondence between these two types of
bordered manifolds.

A bordered $3$-manifold $\YY$ may be presented by a bordered Heegaard diagram
\[\HH = (\Sigma, \{\alpha_1^c, \dots, \alpha_{g-1}^c, \alpha_1^a, \alpha_2^a\}, \{\beta_1, \dots,
\beta_g\}, z),\] where $\Sigma$ is a surface of genus $g$ with one boundary
components, $\{\alpha_1^c, \dots, \alpha_{g-1}^c\}$ and $\{\beta_1, \dots,
\beta_g\}$ are tuples of homologically linearly independent, disjoint circles
in $\Sigma$, and $\alpha_1^a$ and $\alpha_2^a$ are properly embedded arcs that
are disjoint from the $\alpha$ circles and linearly independent from them in
$H_1(\Sigma,\partial\Sigma)$. If we identify $(\partial \Sigma, z,
\partial \Sigma \cap (\alpha_1^a \cup \alpha_2^a))$ with $\ZZ$ --- where
$\partial \Sigma$ is given the boundary orientation --- we obtain a bordered
$3$-manifold with boundary parametrized by $F(\ZZ)$ by attaching handles along
the $\alpha$ and $\beta$ circles. If instead we identify $\partial \Sigma$ with
$-\ZZ$, we obtain a bordered $3$-manifold with boundary parametrized by
$-F(\ZZ)$.

Let $\mathfrak{S}(\HH)$ denote the set of unordered $g$-tuples of points $\x =
\{x_1, \dots, x_g\}$ such that each $\alpha$ circle and each $\beta$ circle
contains exactly one point of $\x$ and each $\alpha$ arc contains at most one
point of $\x$. Let $X(\HH)$ denote the $\F_2$-vector space spanned by
$\mathfrak{S}(H)$.

For generators $\x, \y \in \mathfrak{S}(\HH)$, let $\pi_2(\x,\y)$ denote the
set of homology classes of maps $u\co S \to \Sigma \times [0,1] \times [-2,2]$,
where $S$ is a surface with boundary, taking $\partial S$ to
\begin{multline*}
\left( (\bm{\alpha} \times \{1\} \, \cup \, \bm \beta \times \{0\} \, \cup \,
(\partial\Sigma \minus z) \times [0,1] ) \times [-2,2] \right) \, \cup \\
\left(\x \times [0,1] \times \{-2\} \right) \, \cup \, \left( \y \times [0,1]
\times \{2\}\right)
\end{multline*}
and mapping to the relative fundamental homology class of ${\left(\x \times
[0,1] \times \{-2\} \right)} \cup {\left( \y \times [0,1] \times
\{2\}\right)}$. Each element $B \in \pi_2(\x,\y)$ is determined by its
\emph{domain}, the projection of $B$ to $H_2(\Sigma, \bm{\alpha} \cup
\bm{\beta} \cup \partial \Sigma; \Z)$. The group $H_2(\Sigma, \bm{\alpha} \cup
\bm{\beta} \cup \partial \Sigma; \Z)$ is freely generated by the closures of
the components of $\Sigma \minus (\bm{\alpha} \cup \bm{\beta})$, which we call
\emph{regions}. The domain of any $B \in \pi_2(\x,\y)$ satisfies the following
conditions:
\begin{itemize}
\item The multiplicity of the region containing the basepoint $z$ is
$0$.\footnote{In classical Heegaard Floer homology, the definition of
$\pi_2(\x,\y)$ does not include this requirement.}
\item For each point $p \in \bm{\alpha} \cap \bm{\beta}$, if we identify an oriented
neighborhood of $p$ with $\R^2$, taking $p$ to the origin and the $\alpha$ and
$\beta$ segments containing $p$ to the $x$- and $y$-axes, respectively, and let
$n_1(p)$, $n_2(p)$, $n_3(p)$, and $n_4(p)$ denote the multiplicities in $D$ of
the regions in the four quadrants, then
\begin{equation} \label{eq:pointconditions}
n_1(p) - n_2(p) + n_3(p) -n_4(p) =
\begin{cases}
1 & p \in \x \minus \y \\
-1 & p \in \y \minus \x \\
0 & \text{otherwise}.
\end{cases}
\end{equation}
\end{itemize}
Conversely, any such domain represents some $B \in \pi_2(\x,\y)$. Thus, finding the elements of $\pi_2(\x,\y)$ is a simple matter of linear algebra. A homology class $B \in \pi_2(\x,\y)$ is called \emph{positive} if the
regions in its domain all have non-negative multiplicity; only positive classes
can support holomorphic representatives.

We shall describe only the invariant $\CFD$ here, since we do not compute
$\CFA$ explicitly from a Heegaard diagram in this paper.

We identify the boundary of $\Sigma$ with $-\ZZ$. Assume that the $\alpha$ arcs
are labeled so that $\alpha_1^a \cap \partial \Sigma = \{a_1,a_3\}$ and
$\alpha_2^a \cap \partial \Sigma = \{a_2, a_4\}$.

Define a function $I_D\co \mathfrak{S}(\HH) \to \{\iota_0,\iota_1\}$ by
\begin{equation} \label{eq:idempotents}
I_D(\x) = \begin{cases} \iota_0 & \x \cap \alpha_2^a \ne \varnothing \\ \iota_1
& \x \cap \alpha_1^a \ne \varnothing. \end{cases}
\end{equation}
Define a left action of $\II$ on $X(\HH)$ by $\iota_i \cdot \x= \delta(\iota_i,
I_D(\x)) \x$, where $\delta$ is the Kronecker delta.

For each of the oriented arcs $\rho_I \subset \ZZ$, let $-\rho_I$ denote
$\rho_I$ with its opposite orientation. (That is, $-\rho_1$ goes from $a_2$ to
$a_1$, etc.) Given $\x \in \mathfrak{S}(\HH)$ and a sequence $\vec \rho =
(-\rho_{I_1}, \dots, -\rho_{I_k})$, the pair $(\x, \vec\rho)$ is called
\emph{strongly boundary monotonic} if the initial point of $-\rho_{I_1}$ is on
the same $\alpha$ circle as $\x$, and for each $i>1$, the initial point of
$-\rho_{I_i}$ and the final point of $-\rho_{I_{i-1}}$ are paired in $\ZZ$.

If $B \in \pi_2(\x,\y)$ is a positive class, then $\partial^\partial B$ (the
intersection of the domain of $B$ with the boundary of $\Sigma$) may be
expressed (non-uniquely) as a sum of arcs $-\rho_{I_i}$. Specifically, we say
that the pair $(B,\vec\rho)$ is \emph{compatible} if $(\x, \vec\rho)$ is
strongly boundary monotonic and $\partial^\partial B = \sum_{i=1}^k
(-\rho_{I_i})$. If $(B, \vec\rho)$ is compatible, the \emph{index} of $(B,
\vec\rho)$ is defined in \cite[Definition~5.61]{LOTBordered} as
\begin{equation} \label{eq:index}
\ind(B, \vec\rho) = e(B) + n_\x(B) + n_\y(B) + \abs{\vec\rho} +
\iota(\vec\rho),
\end{equation}
where $e(B)$ is the Euler measure of $B$; $n_\x(B)$ (resp.~$n_\y(B)$) is the
sum over points $x_i \in \x$ (resp.~$y_i \in \y$) of the average of the
multiplicities of the regions incident to $x_i$ (resp.~$y_i$), $\abs{\vec\rho}$
is the number of entries in $\vec\rho$, and $\iota(\vec\rho)$ is a
combinatorially defined quantity \cite[Equation~5.58]{LOTBordered} that
measures the overlapping of the arcs $\rho_{I_i}$. The index $\ind(B,
\vec\rho)$ is equal to one plus the expected dimension of a certain moduli
space $\MM^B(\x,\y,\vec\rho)$ of $J$-holomorphic curves in $\Sigma \times [0,1]
\times \R$ in the homology class $B$ whose asymptotics near $\partial \Sigma
\times [0,1] \times \R$ are specified by $\vec \rho$. In particular, if
$\ind(B, \vec \rho)=1$, then this moduli space contains finitely many points.
We do not give the full definition here; see \cite[Chapter~5]{LOTBordered} for
the details.

For each $\x,\y \in \mathfrak{S}(\x)$ and $B \in \pi_2(\x,\y)$, define
\[
a^B_{\x,\y} = \sum_{ \substack{ \{ \vec\rho =
(-\rho_{I_1}, \dots, -\rho_{I_k}) \, \mid \\ (B,\vec\rho) \text{ compatible}, \\
\ind(B, \vec\rho) = 1\} }} \#(\MM^B(\x,\y,\vec\rho)) \, \rho_{I_1} \dots
\rho_{I_k} \in \AA,
\]
where the count of points in $\MM^B(\x,\y,\vec\rho)$ is taken modulo $2$. We
define $\delta_1\co  X(\HH) \to \AA \underset{\II}{\otimes} X(\HH)$ by
\begin{equation} \label{eq:delta1}
\delta_1 (\x) = \sum_{\y \in \mathfrak{S}(\HH)} \sum_{B \in \pi_2(\x,\y)}
a^B_{\x,\y} \otimes \y.
\end{equation}
This defines a type $D$ structure, which we denote $\CFD(\HH)$. The
verification of \eqref{eq:d-relation} is a version of the standard
$\partial^2=0$ argument in Floer theory. (Henceforth, if $\x$, $\y$, and $\vec \rho$ are understood from the context, we shall write $\MM(B)$ in place of $\MM^B(\x,\y, \vec\rho)$. If we need to be explicit about the choice of complex structure $J$ on $\Sigma$, we shall write $\MM_J(B)$ or $\MM^B_J(\x,\y, \vec\rho)$.)

\begin{proposition} \label{prop:chordsallowed}
$\ $
\begin{enumerate}
\item The only sequences of chords that can contribute nonzero terms to
$\delta_1$ are the empty sequence, $(-\rho_1)$, $(-\rho_2)$, $(-\rho_3)$,
$(-\rho_1, -\rho_2)$, $(-\rho_2, -\rho_3)$, $(-\rho_{123})$, and $(-\rho_1,
-\rho_2, -\rho_3)$. Therefore, only classes whose multiplicities in the
boundary regions of $\Sigma$ are $0$ or $1$ can count for $\delta_1$.

\item If $B \in \pi_2(\x,\y)$ is a positive class whose domain has multiplicity $1$
in the regions abutting $\rho_1$ and $\rho_2$ (resp.~$\rho_2$ and $\rho_3$) and
$0$ in the region abutting $\rho_3$ (resp.~$\rho_1$), then $B$ may count for
the differential only if $\x$ and $\y$ contain points of $\alpha^a_1$
(resp.~$\alpha^a_2$).
\end{enumerate}
\end{proposition}

\begin{proof}
For the first statement, the only other sequences of chords for which the
product of algebra elements in the definition of $a_{\x,\y}^B$ is nonzero are
$(-\rho_{12})$, $(-\rho_{23})$, $(-\rho_1, -\rho_{23})$, and $(-\rho_{12},
-\rho_3)$. The two latter sequences are not strongly boundary monotonic. If $B
\in \pi_2(\x,\y)$ is a positive class compatible with $(-\rho_{12})$, then $\x$
and $\y$ both contain points on $\alpha_1^a$, since otherwise $B$ would have a
boundary component without a $\beta$ segment. Therefore, $I_D(\y) = \iota_1$.
Since the tensor product is taken over the ring of idempotents,
\[
\rho_{12} \otimes \y = \rho_{12} \otimes \iota_1 \y = \rho_{12} \iota_1 \otimes
\y = 0,
\]
so the contribution of $B$ to $\delta_1(\x)$ is zero. A similar argument
applies for the sequence $(-\rho_{23})$. The second statement follows
immediately from the same argument.
\end{proof}

The invariant $\CFA$ is a type $A$ structure associated to a bordered Heegaard
diagram whose boundary is identified with $\ZZ$. We do not give all the details
here. The analogue of Proposition \ref{prop:chordsallowed} does not hold for
$\CFA$; one must consider domains with arbitrary multiplicities on the boundary
and a much larger family of sequences of chords. Therefore, it is generally
easier to compute $\CFD$.

We conclude this section with the \emph{gluing theorem}:
\begin{theorem}[Lipshitz--Ozsv\'ath--Thurston \cite{LOTBordered}]
Suppose $\YY_1$ and $\YY_2$ are bordered $3$-manifolds, and $Y = Y_1 \cup_\phi
Y_2$ is the manifold obtained by gluing them together along their boundaries,
where $\phi\co -\partial Y_1 \to \partial Y_2$ is the map induced by the
bordered structures. Then
\[
\CF(Y) \simeq \CFA(\YY_1) \boxtimes \CFD(\YY_2),
\]
provided that at least one of the modules is bounded (so that the box tensor
product is defined).
\end{theorem}

\subsection{Bimodules} \label{subsec:bimodules}

In \cite{LOTBimodules}, Lipshitz, Ozsv\'ath, and Thurston also define
invariants for a \emph{bordered manifold with two boundary components}.
Essentially, this consists of a manifold $Y$ with two boundary components
$\partial_L Y$ and $\partial_R Y$, with parametrizations of the two boundary
components just like in the single-component case, and a framed arc connecting
the two boundary components. Here, we assume that both boundary components are tori; see \cite[Chapter 5]{LOTBimodules} for the full definition.

If both boundary components are parametrized by $-F(\ZZ)$, the associated
invariant is a type $DD$ structure over two copies of $\AA$, denoted
$\CFDD(\YY)$; if both are parametrized by $F(\ZZ)$, the invariant is a type
$AA$ structure, denoted $\CFAA(\YY)$; and similarly there are invariants
$\CFAD(\YY)$ and $\CFDA(\YY)$. For simplicity, we denote the two copies of $\AA$ by $\AA_\rho$
and $\AA_\sigma$; in the latter, the Reeb elements are written $\sigma_1$,
$\sigma_2$, etc.

In fact, we consider only a direct summand of each bimodule, denoted
$\CFDD(\YY,0)$, $\CFAA(\YY,0)$, etc., which is all that is necessary to compute
the Floer complex of a manifold obtained by gluing two separate
one-boundary-component manifolds to the two boundary components of $Y$. The
other summands are only necessary if one wishes to glue together the two
boundary components of $Y$.

As in the previous discussion, we describe only the construction of $\CFDD$. A
bordered manifold with two toroidal boundary components may be presented by an
\emph{arced bordered Heegaard diagram}
\[\HH = (\Sigma, \{\alpha_1^c, \dots, \alpha_{g-2}^c, \alpha_1^L, \alpha_2^L, \alpha_1^R, \alpha_2^R \}, \{\beta_1, \dots,
\beta_g\}, \z),\] where now $\partial\Sigma$ has two components $\partial_L
\Sigma$ and $\partial_R \Sigma$, on which the arcs $\alpha_i^L$ and
$\alpha_i^R$ have their respective boundaries, and $\z$ is an arc in the
complement of all the $\alpha$ and $\beta$ circles and $\alpha$ arcs connecting
the two boundary components.

We define $\mathfrak{S}(\HH)$ and $X(\HH)$ just in the
single-boundary-component case. Let $\mathfrak{S}(\HH,0)$ be the subset of
$\mathfrak{S}(\HH)$ consisting of $g$-tuples $\x$ containing one point in
$\alpha_1^L \cup \alpha_2^L$ and one point in $\alpha_1^R \cup \alpha_2^R$, and
let $X(\HH,0)$ be the $\F$-vector space generated by $\mathfrak{S}(\HH,0)$.
This is the underlying vector space for the invariants $\CFDD(\HH,0)$,
$\CFAA(\HH,0)$, etc.

To define $\CFDD(\HH,0)$, identify both boundary components of $\Sigma$ with
$-\ZZ$. Each generator of $\CFDD(\HH,0)$ has associated idempotents in
$\AA_\rho$ and $\AA_\sigma$, as in \eqref{eq:idempotents}. The differential
\[
\delta_1\co X(\HH,0) \to (\AA_\rho \otimes \AA_\sigma) \underset{\II_\rho
\otimes \II_\sigma}{\otimes} X(\HH,0)
\]
is then defined essentially the same way as with $\CFD$ of a
single-boundary-component diagram. Specifically, for a homology class $B \in
\pi_2(\x,\y)$ and sequences of chords $\vec\rho = (-\rho_{I_1}, \dots,
-\rho_{I_k})$ and $\vec\sigma = (-\sigma_{J_1}, \dots, -\sigma_{J_l})$ on the
two boundary components, the definitions of compatibility and of the index
$\ind(B,\vec\rho,\vec\sigma)$ are as above. Define
\[
a^B_{\x,\y} = \sum_{ \substack{\{ (\vec \rho, \vec \sigma)  \, \mid \\ (B,\vec\rho, \vec\sigma) \text{ compatible}, \\
\ind(B, \vec\rho, \vec\sigma) = 1\} }} \#(\MM^B(\x,\y,\vec\rho, \vec\sigma)) \,
\rho_{I_1} \dots \rho_{I_k} \otimes \sigma_{J_1} \dots \sigma_{J_l} \in
\AA_\rho \otimes \AA_\sigma.
\]
The map $\delta_1$ is then given by \eqref{eq:delta1} just as above. An
analogue of Proposition \ref{prop:chordsallowed} also holds in this setting.
For further details, see \cite[Section 6]{LOTBimodules}.

The gluing theorem generalizes naturally to bimodules. For instance, if $Y_1$
has a single boundary component parametrized by $F(\ZZ)$, $Y_2$ has two
boundary components parametrized by $-F(\ZZ)$, and $\phi\co -\partial Y_1 \to
\partial_L Y_2$ is the map induced by the parametrizations, then
\[
\CFD(\YY_1 \cup_\phi \YY_2) \simeq \CFA(\YY_1) \underset{\AA_\rho}{\boxtimes}
\CFDD(\YY_2,0).
\]
The remaining generalizations are found in \cite[Theorems 11,
12]{LOTBimodules}.

Finally, we mention the \emph{identity $AA$ bimodule} \cite[Subsection
10.1]{LOTBimodules}. Consider the manifold $\I= F(\ZZ) \times I$. Parametrize
$\partial_R Y = F(\ZZ) \times \{1\}$ by inclusion and $\partial_L Y = F(\ZZ)
\times \{0\}$ (whose boundary-induced orientation is opposite to the standard
orientation of $F(\ZZ)$) by the composition $F(\ZZ) \xrightarrow{r} -F(\ZZ)
\hookrightarrow F(\ZZ) \times \{0\}$; thus, both boundary components are
parametrized by $F(\ZZ)$ as opposed to $-F(\ZZ)$. The bijection between
bordered manifolds with boundary $-F(\ZZ)$ and bordered manifolds with boundary
$F(\ZZ)$ may be given by $Y \mapsto Y \cup \I$. Thus, if $\HH$ is any bordered
Heegaard diagram with one boundary component, then the type $A$ module
$\CFA(\HH)$ (where we identify $\partial \Sigma$ with $\ZZ$) is chain homotopy
equivalent to $\CFAA(\I,0) \boxtimes \CFD(\HH)$ (where, in the second factor,
we identify $\partial \Sigma$ with $-\ZZ$).\footnote{Our presentation here is a bit different from that of Lipshitz, Ozsv\'ath, and Thurston, who describe $\CFAA(\I,0)$ as a bimodule over two separate algebras, $\AA(\ZZ)$ and $\AA(-\ZZ)= \AA(\ZZ)^{\operatorname{op}}$. The latter happens to be isomorphic to $\AA(\ZZ)$ because of the involution $r$, so the two boundary components of $\II$ are effectively interchangeable. For the purposes of this introduction, we find it clearer to suppress the distinction between $\AA(\ZZ)$ and $\AA(-\ZZ)$, at the cost of being more explicit about $r$.} As mentioned above, it is easier to compute $\CFD$ explicitly from a Heegaard diagram than $\CFA$; by taking a
tensor product with $\CFAA(\I,0)$, we can always avoid the latter.

\begin{theorem}[Lipshitz--Ozsv\'ath--Thurston] \label{thm:identityAA}
The type $AA$ module $\CFAA(\I,0)$ has generators $w_1, w_2, x, y, z_1, z_2$,
with $\AA_\infty$ multiplications as illustrated in Figure
\ref{fig:identityAA}. That is, $m_{1,0,0}(w_1) = w_2$, $m_{1,0,1}(w_1, \sigma_1) = y$, $m_{1,1,1}(w_1, \sigma_{12},\rho_{23}) = w_2$, and so on. (See below for more on this notation.)
\end{theorem}

\begin{figure} \centering
\[
\xy
 (0,50)*{w_1}="w1";
 (0,0)*{w_2}="w2";
 (50,50)*{z_1}="z1";
 (50,0)*{z_2}="z2";
 (25,35)*{y}="y";
 (25,15)*{x}="x";
 {\ar|{1+\sigma_{12} \rho_{23}} "w1";"w2"};
 {\ar|{1+\sigma_{23} \rho_{12}} "z1";"z2"};
 {\ar|{\sigma_1} "w1";"y"};
 {\ar|{\rho_1} "z1";"y"};
 {\ar|{\sigma_{12} \rho_2} "w1";"x"};
 {\ar|{\sigma_2 \rho_{12}} "z1";"x"};
 {\ar|{\sigma_{23} \rho_2} "y";"z2"};
 {\ar|{\sigma_2 \rho_{23}} "y";"w2"};
 {\ar|{\sigma_2 \rho_2} "y";"x"};
 {\ar|{\sigma_3} "x";"z2"};
 {\ar|{\rho_3} "x";"w2"};
 {\ar@/_12.5pc/|{\twoline{\sigma_{123}\rho_2 +}{\sigma_3 \sigma_2 \sigma_1 \rho_2}} "w1";"z2"};
 {\ar@/^12.5pc/|{\sigma_2 \rho_{123}} "z1";"w2"};
\endxy
\]
\caption{The identity $AA$ bimodule, $\CFAA(\II,0)$.} \label{fig:identityAA}
\end{figure}

\subsection{Knots in bordered manifolds} \label{subsec:knots}

A \emph{doubly-pointed bordered Heegaard diagram} consists of a bordered
Heegaard diagram $\HH = (\Sigma, \bm{\alpha}, \bm{\beta},z)$ along with an
additional basepoint $w \in \Sigma \minus (\bm{\alpha} \cup \bm{\beta})$. As
explained in \cite[Section 11.4]{LOTBordered}, a doubly-pointed diagram
determines a knot $K \subset Y$ with a single point of $K$ meeting the
basepoint on $\partial Y$; the isotopy class of $K$ relative to this point is invariant under
Heegaard moves missing $w$. Lipshitz, Ozsv\'ath, and Thurston define invariants
$\operatorname{CFD}^-(Y,K)$ and $\operatorname{CFA}^-(Y,K)$ by working over the
algebra $\AA \otimes \F[U]$, where the $U$ powers record the multiplicity of
$w$ in each domain that counts for the differential or multiplications.

If the knot $K$ is nulhomologous in $Y$, we prefer the following alternate
perspective. Push $K$ slightly into the interior of $Y$ (so that it now misses the boundary), and let $F$ be a Seifert surface for $K$. Just as in ordinary knot
Floer homology \cite{OSzKnot, RasmussenThesis}, each generator $\x \in
\mathfrak{S}(\HH)$ has an associated relative spin$^c$ structure $\s_{z,w}(\x)
\in \Spin^c(Y,K)$, and we may define an \emph{Alexander grading} on
$\mathfrak{S}(\HH)$ by
\begin{equation} \label{eq:absalex}
A(x) = \frac12 \gen{c_1(\s_{z,w}(\x)), [F]},
\end{equation}
where $c_1(\s_{z,w}(\x)) \in H^2(Y,K)$ and $[F] \in H_2(Y,K)$. The grading
difference between two generators is given by
\begin{equation} \label{eq:relalex}
A(x) - A(y) = n_w(B)
\end{equation}
where $B \in \pi_2(\x,\y)$ is any domain from $\x$ to $\y$. To verify that the
right-hand side of \eqref{eq:relalex} is well-defined, note that for any
periodic class $P \in \pi_2(\x,\x)$, $n_w(P)$ equals the intersection number of
$K$ with the homology class in $H_2(Y,\partial Y)$ corresponding to $P$, which
must be zero since $K$ is nulhomologous. Further details are completely
analogous to \cite{OSzKnot, RasmussenThesis}.

The Alexander grading on $X(\HH)$ determines a filtration on $\CFA(\HH)$ or
$\CFD(\HH)$, since any domain that counts for the differential or $\AA_\infty$
multiplications has non-negative multiplicity at $w$. We denote the filtered
chain homotopy type by $\CFA(\YY,K)$ or $\CFD(\YY,K)$.

When we evaluate a tensor product $\CFA(\YY_1) \boxtimes \CFD(\YY_2)$, a knot
filtration on one factor extends naturally to a filtration on the whole
complex, which agrees with the filtration that the knot induces on $\CF(Y_1
\cup Y_2)$.

A nulhomologous knot in a bordered manifold with two boundary components may be
handled similarly. For invariance, one point of the knot must be constrained to
lie on the arc connecting the two boundary components, and isotopies must be
fixed in a neighborhood of that point.

\subsection{The edge reduction algorithm} \label{subsec:edge}

We now describe the well-known ``edge reduction'' procedure for chain complexes
and its extension to $\AA_\infty$ modules.

Suppose $(C,\partial)$ is a free chain complex with basis $\{x_1, \dots, x_n\}$
over a ring $R$. For each $i,j$, let $a_{ij}$ be the coefficient of $x_j$ in
$\partial x_i$ with respect to this basis. If some $a_{ij}$ is invertible in
$R$, define a new basis $\{y_1, \dots, y_n\}$ by setting $y_i = x_i$, $y_j =
\partial x_i$, and for each $k \ne i,j$, $y_k = x_k - a_{kj} a_{ij}^{-1} x_i$,
where $a_{kj}$ is the coefficient of $x_j$ in $\partial x_k$. With respect to
the new basis, the coefficient of $y_j$ in $\partial y_k$ is zero, so the
subspace spanned by $y_i$ and $y_j$ is a direct summand subcomplex with trivial
homology. Thus, $C$ is chain homotopy equivalent to the subcomplex $C'$ spanned
by $\{y_k \mid k \ne i,j\}$, in which the coefficient of $y_l$ in $\partial
y_k$ is $a_{kl} - a_{kj} a_{ij}^{-1} a_{il}$.

When $R=\F_2$, a convenient way to represent a chain complex $(C,\partial)$
with basis $\{x_i\}$ is a directed graph $\Gamma_{C,\partial,\{x_i\}}$ with
vertices corresponding to basis elements and an edge from $x_i$ to $x_j$
whenever $a_{ij}=1$. To obtain $\Gamma_{C',\partial,\{y_k\}}$ from
$\Gamma_{C,\partial,\{x_i\}}$ as above, we delete the vertices $x_i$ and $x_j$
and any edges going into or out of them. For each $k$ and $l$ with edges $x_k
\to x_j$ and $x_i \to x_l$, we either add an edge from $x_k$ to $x_l$ (if there
was not one previously) or eliminate the edge from $x_k$ to $x_l$ (if there was
one). We call this procedure \emph{canceling the edge from $x_i$ to $x_j$}. The
vertices of the resulting graph should be labeled with $\{y_k \mid k \ne
i,j\}$, but by abuse of notation we frequently continue to refer to them with
$\{x_k \mid k \ne i,j\}$ instead.

By iterating this procedure until no more edges remain, we compute the homology
of $C$. If the matrix $(a_{ij})$ is sparse, this tends to be a very efficient
algorithm for computing homology. If $C$ is a graded complex and the basis
$\{x_1, \dots, x_n\}$ consists of homogeneous elements, then $y_k$ is clearly
homogeneous with the same grading as $x_k$, so we can compute the homology as a
graded group.

If $C$ has a filtration $\cdots \subseteq F_p \subseteq F_{p+1} \subseteq
\cdots$, the \emph{filtration level} of an element of $C$ is the unique $p$ for
which that element is in $F_p \minus F_{p-1}$. To compute the spectral sequence
associated to the filtration, we cancel edges in increasing order of the amount
by which they decrease filtration level. At each stage, this guarantees that
the filtration level of $y_k$ equals that of $x_k$. The complex that remains
after we delete all edges that decrease filtration level by $k$ is the
$E^{k+1}$ page in the spectral sequence, and the vertices that remain after all
edges are deleted is the $E^\infty$ page. In particular, when $C = \CF(S^3,K)$,
the filtered complex associated to a knot $K \subset S^3$, the total homology
of $C$ is $\HF(S^3;\F) \cong \F$, so a unique vertex survives after all
cancellations are complete. The filtration level of this vertex is, by
definition, the invariant $\tau(K)$.

More generally, over an arbitrary ring $R$, we may represent $(C,\partial)$ by
a labeled, directed graph, where now we label an edge from $x_i$ to $x_j$ by
$a_{ij}$, often omitting the label when $a_{ij}=1$. When we cancel an unlabeled
edge from $x_i$ to $x_j$, we replace a zigzag
\[
x_k \xrightarrow{a_{kj}} x_j \longleftarrow x_i \xrightarrow{a_{il}} x_l
\]
with an edge
\[
x_k \xrightarrow {-a_{kl} a_{il}} x_l
\]
if no such edge existed previously, and either relabel or delete such an edge
if it did exist. Of course, when $R$ is not a field, this procedure is not
guaranteed to eliminate all edges or to yield a result that is independent of
the choice of the order in which the edges are deleted, but it is still often a
useful way to simplify a chain complex.

The same procedure works for type $D$ structures over the torus algebra $\AA$,
as can be seen by looking at the ordinary differential module obtained by
taking the tensor product with $\AA$ as above.

Edge cancellation for type $A$ structures is slightly more complicated. We work
only with bounded modules for simplicity. Suppose $M$ is a bounded type $A$
structure over $\AA$ with a basis $\{x_1, \dots, x_n\}$. As above, we may
describe the multiplications using a matrix of formal sums of finite sequences
of elements of $\AA$, and we may represent the nonzero entries using a labeled graph.
The empty sequence will be denoted by $1$, and we often omit the label of an edge whose label is $1$. If there is an unlabeled edge from $x_i$ to $x_j$ then we may cancel $x_i$ and $x_j$, replacing a zigzag
\[
x_k \xrightarrow{(a_1, \dots, a_p)} x_j \longleftarrow x_i \xrightarrow{(b_1,
\dots, b_q)} x_l
\]
by an edge
\[
x_k \xrightarrow {(a_1, \dots, a_p, b_1, \dots, b_q)} x_l
\]
(or eliminating such an edge if one already exists). The $\AA_\infty$ module
$M'$ described by the resulting graph is then $\AA_\infty$ chain homotopic to
$M$. If $M$ is a filtered $\AA_\infty$-module and the edge being canceled is
filtration-preserving (i.e., $x_i$ and $x_j$ have the same filtration level),
then $M'$ is filtered $\AA_\infty$ chain homotopic to $M$. Similar techniques
may also be used for bimodules. (As in Figure \ref{fig:identityAA}, we frequently omit the parentheses and commas on the edge labels for conciseness; with this notation, concatenation does not indicate multiplication in $\AA$.)

\subsection{\texorpdfstring{$\widehat{\mathrm{CFD}}$}{CFD} of knot complements}

For any knot $K$, let $X_K$ denote the exterior of $K$. For $t\in \Z$, let
$\XX_K^t$ denote the bordered structure on $X_K$ determined by a map $\phi:
-F(\ZZ) \to \partial X_K$ sending $h_1$ to a $t$-framed longitude (relative to
the Seifert framing) and $h_2$ to a meridian of $K$. Lipshitz, Ozsv\'ath, and
Thurston \cite{LOTBordered} give a complete computation of $\CFD(\XX_K^t)$ in
terms of the knot Floer complex of $K$, which we now describe.

In the computation that follows, we will need to work with two different framed
knot complements, $\XX_J^s$ and $\XX_K^t$. We first state the results for
$\CFD(\XX_J^s)$ and then indicate how to modify the notation for
$\CFD(\XX_K^t)$. Define $r = \abs{2\tau(J)-s}$, and say that $\dim \HFK(S^3,J)
= 2n+1$.

We may find two distinguished bases for $\operatorname{CFK}^-(S^3,J)$: a
``vertically reduced'' basis $\{\tilde\xi_0, \dots, \tilde\xi_{2n}\}$, with
``vertical arrows'' $\tilde\xi_{2j-1} \to \tilde\xi_{2j}$ of length $k_j \in
\N$, and a ``horizontally reduced'' basis $\{\tilde\eta_0, \dots,
\tilde\eta_{2n}\}$, with ``horizontal arrows'' $\tilde\xi_{2j-1} \to
\tilde\xi_{2j}$ of length $l_j \in \N$. (See \cite[Section 11.5]{LOTBordered} for
the definitions.) Denote the change-of-basis matrices by $(x_{p,q})$ and
$(y_{p,q})$, so that
\begin{equation} \label{eq:xpq-ypq}
\tilde\xi_p = \sum_{q=0}^{2n} x_{p,q} \tilde\eta_q \quad \text{and} \quad
\tilde\eta_p = \sum_{q=0}^{2n} y_{p,q} \tilde\xi_q.
\end{equation}
In all known instances, the two bases may be taken to be equal as sets (up to a
permutation), but it has not been proven that this holds in general.

According to \cite[Theorems 11.27, A.11]{LOTBordered}, the structure of
$\CFD(\XX_J^t)$ is as follows. The part in idempotent $\iota_0$ (i.e., $\iota_0
\CFD(\XX_J^s)$) has dimension $2n+1$, with designated bases $\{\xi_0, \dots,
\xi_{2n}\}$ and $\{\eta_0, \dots, \eta_{2n}\}$ related by \eqref{eq:xpq-ypq}
without the tildes. The part in idempotent $\iota_1$ (i.e., $\iota_1
\CFD(\XX_J^s)$) has dimension $r+\sum_{j=1}^n (k_j + l_j)$, with basis
\[
\{\gamma_1, \dots, \gamma_r\} \cup \bigcup_{j=1}^n \{\kappa^j_1, \dots,
\kappa^j_{k_j}\} \cup \bigcup_{j=1}^n \{\lambda^j_1, \dots, \lambda^j_{l_j}\}.
\]

For $j=1,\dots, n$, corresponding to the vertical arrow $\tilde\eta_{2j-1} \to
\tilde\eta_{2j}$, there are differentials
\begin{equation} \label{eq:vertchain}
\xi_{2j} \xrightarrow{\rho_{123}} \kappa^j_1 \xrightarrow{\rho_{23}} \cdots
\xrightarrow{\rho_{23}} \kappa^j_{k_j} \xleftarrow{\rho_1} \xi_{2j-1}.
\end{equation}
(In other words, $\delta_1 (\xi_{2j})$ has a $\rho_{123} \otimes \kappa^j_1$
term, and so on.) We refer to the subspace of $\CFD(\XX_J^s)$ spanned by the
generators in \eqref{eq:vertchain} as a \emph{vertical stable chain}.
Similarly, corresponding to the horizontal arrow $\eta_{2j-1} \to \eta_{2j}$ of
length $l_j$, there are differentials
\begin{equation} \label{eq:horizchain}
\eta_{2j-1} \xrightarrow{\rho_3} \lambda^j_1 \xrightarrow{\rho_{23}} \cdots
\xrightarrow{\rho_{23}} \lambda^j_{l_j} \xrightarrow{\rho_2} \eta_{2j},
\end{equation}
and the generators here span a \emph{horizontal stable chain}. Finally, the
generators $\xi_0, \eta_0, \gamma_1, \dots, \gamma_r$ span the \emph{unstable
chain}, with differentials depending on $s$ and $\tau(J)$:
\begin{equation} \label{eq:unstchain}
\begin{cases}
\eta_0 \xrightarrow{\rho_3} \gamma_1 \xrightarrow{\rho_{23}} \cdots
\xrightarrow{\rho_{23}} \gamma_r \xleftarrow{\rho_1} \xi_0 & s <
2\tau(J) \\
\xi_0 \xrightarrow{\rho_{12}} \eta_0 & s=2\tau(J) \\
\xi_0 \xrightarrow{\rho_{123}} \gamma_1 \xrightarrow{\rho_{23}} \cdots
\xrightarrow{\rho_{23}} \gamma_r \xrightarrow{\rho_2} \eta_0 & s > 2\tau(J).
\end{cases}
\end{equation}
In some instances, as with the unknot and the figure-eight knot, we may have
$\xi_0=\eta_0$.

For $\CFD(\XX_K^t)$, we modify the preceding two paragraphs by replacing all
lower-case letters with capital letters. Specifically, $\iota_0 \CFD(\XX_K^t)$
has bases $\{\Xi_0, \dots, \Xi_{2N}\}$ and $\{\Eta_0, \dots, \Eta_{2N}\}$
related by change-of-basis matrices $(X_{P,Q})$ and $(Y_{P,Q})$ as in
\eqref{eq:xpq-ypq}; $\iota_1 \CFD(\XX_K^t)$ has basis
\[
\{\Gamma_1, \dots, \Gamma_R\} \cup \bigcup_{J=1}^N \{\Kappa^J_1, \dots,
\Kappa^J_{K_J}\} \cup \bigcup_{J=1}^N \{\Lambda^J_1, \dots, \Lambda^J_{L_J}\};
\]
and the differentials are just as in \eqref{eq:vertchain},
\eqref{eq:horizchain}, and \eqref{eq:unstchain}, suitably
modified.\footnote{The reader should take care to distinguish capital eta
($\Eta$) and kappa ($\Kappa$) from the Roman letters $H$ and $K$. We find that
the mnemonic advantage of using parallel notation for the generators of
$\CFD(\XX_J^s)$ and $\CFD(\XX_K^t)$ outweighs any confusion that may arise.} In
the discussion below, we shall treat $\CFD(\XX_K^t)$ as a type $D$ structure
over a copy of $\AA_\sigma$ in which the elements are referred to as
$\sigma_1$, $\sigma_2$, etc., to facilitate taking the double tensor product.

In Section \ref{sec:tensorproduct}, we shall frequently use the following
proposition to simplify computations:

\begin{proposition} \label{prop:norho1rho2}
In the matrix entries for the higher maps $\delta_k$ for $\CFD(\XX_J^s)$, there
are no sequences of elements containing $\rho_1 \otimes \rho_2$, $\rho_1
\otimes \rho_{23}$, $\rho_2 \otimes \rho_3$, or $\rho_{12} \otimes \rho_3$.
\end{proposition}

\begin{proof}
The only instances of $\rho_1$ in $\CFK(\XX_J^s)$ are $\xi_{2j-1}
\xrightarrow{\rho_1} \kappa^j_{k_j}$ in the vertical chains and $\xi_0
\xrightarrow{\rho_1} \gamma_r$ in the unstable chain when $s<2\tau(J)$, and
$\delta_1(\kappa^j_{k_j}) = \delta_1(\gamma_r) = 0$. Thus, $\rho_1 \otimes
\rho_2$ and $\rho_1 \otimes \rho_{23}$ may not occur in $\delta_k$. Similarly,
the only instances of $\rho_2$ and $\rho_{12}$ are $\lambda^j_{l_j}
\xrightarrow{\rho_2} \eta_{2j}$ in the horizontal chains, $\gamma_r
\xrightarrow{\rho_2} \eta_0$ in the unstable chain when $s>2\tau(J)$, and
$\xi_0 \xrightarrow{\rho_{12}} \eta_0$ when $s=2\tau(J)$, and the only
instances of $\rho_3$ are $\eta_{2j-1} \xrightarrow{\rho_3} \lambda^j_1$ in the
horizontal chains and $\eta_0 \xrightarrow{\rho_3} \gamma_1$ in the unstable
chain when $s<2\tau(J)$. Thus, no element that is at the head of a $\rho_2$ or
$\rho_{12}$ arrow is also at the tail of a $\rho_3$ arrow.
\end{proof}

%\section{Direct computation of $\widehat{\mathrm{CFAA}}(Y,B_3)$} \label{sec:heegaard}
\section{Direct computation of \texorpdfstring{$\widehat{\mathrm{CFAA}}(Y,B_3)$}{CFAA(Y,B\textunderscore3)}} \label{sec:heegaard}
As above, let $B = B_1 \cup B_2 \cup B_3 \subset S^3$ denote the Borromean
rings. Let $Y$ denote the complement of a neighborhood of $B_1 \cup B_2$; then
$B_3$ is a nulhomologous knot in $Y$. Let $\partial_L Y$ and $\partial_R Y$ be
the boundary components coming from $B_1$ and $B_2$, respectively. We define a
strongly bordered structure $\YY$ on $Y$ (in the sense of \cite[Definition
5.1]{LOTBimodules}) so that the map $\phi_L: F(\ZZ) \to \partial_L Y$
(resp.~$\phi_R: F(\ZZ) \to \partial_R Y$) takes $h_1$ to a meridian of $B_1$
(resp.~$B_2$) and $h_2$ to a Seifert-framed longitude of $B_1$ (resp.~$B_2$).
It follows that the glued manifold $(\YY \cup_{\partial_L Y} \XX_J^s)
\cup_{\partial_R Y} \XX_K^t$, is $S^3$, and the image of $B_3$ is the knot
$D_{J,s}(K,t)$.\footnote{Because we are gluing the two boundary components of
$\YY$ to separate single-boundary-component bordered manifolds, the choice of
framed arc connecting $\partial_L Y$ and $\partial_R Y$ does not affect the
final computation of the tensor product, so we suppress all reference to it.}
Thus, we must compute the filtered type $AA$ bimodule $\CFAA(\YY,B_3,0)$.

\subsection{A Heegaard diagram for \texorpdfstring{$(\YY,B_3)$}{(Y,B_3)}} \label{subsec:heegaarddiagram}

\begin{figure} \centering
\psfrag{w}{$w$} \psfrag{zz}{$\mathbf{z}$}
\includegraphics{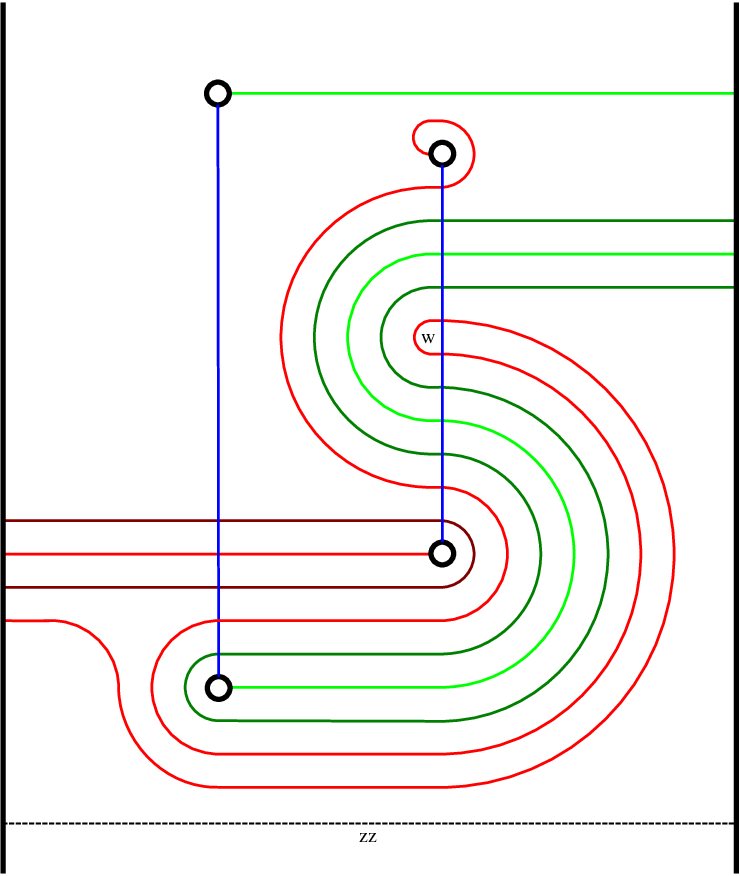}
\caption{The arced Heegaard diagram $\HH$.} \label{fig:heegaard1}
\end{figure}

\begin{figure} \centering
\psfrag{w}{$w$} \psfrag{z}{$z$} \psfrag{A}{{\color{purple} $A$}}
\psfrag{C}{{\color{teal} $C$}}
\includegraphics{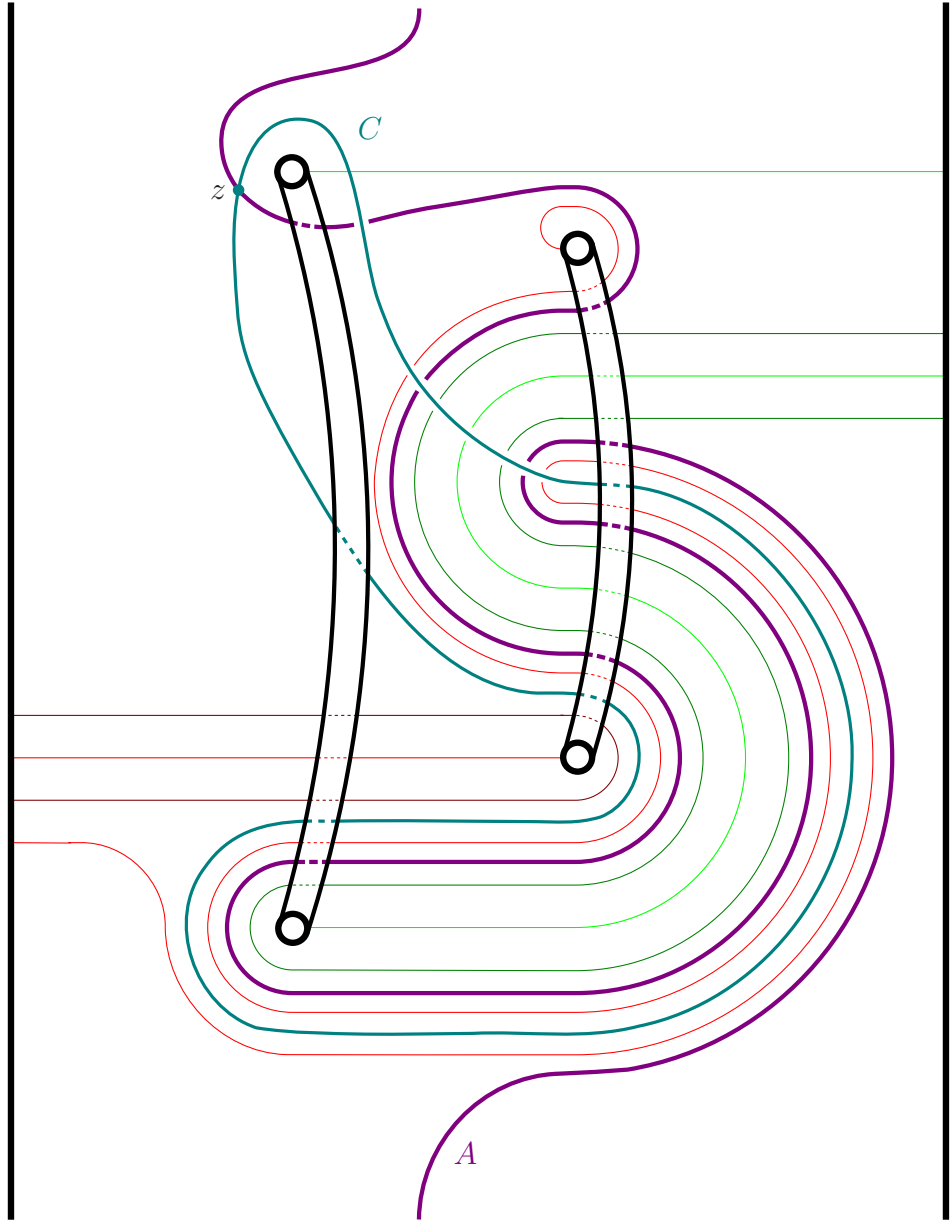}
\caption[The manifold $Y(\HH_{dr})$.] {The manifold $Y(\HH_{dr})$. The $\alpha$
arcs from $\HH$ (the thin red and green curves) and the circle $A$ (purple) sit
in the plane of the page, while the knot $C$ (turquoise) is above the plane of the page (i.e., in the interior of
$Y(\HH_{dr})$) except at the point $z$.} \label{fig:heegaardknot1}
\end{figure}

\begin{figure} \centering
\psfrag{w}{$w$} \psfrag{z}{$z$} \psfrag{A}{{\color{purple} $A$}}
\psfrag{C}{{\color{teal} $C$}} \psfrag{a1}{$a_1$} \psfrag{a2}{$a_2$}
\psfrag{a3}{$a_3$} \psfrag{a4}{$a_4$}
\includegraphics{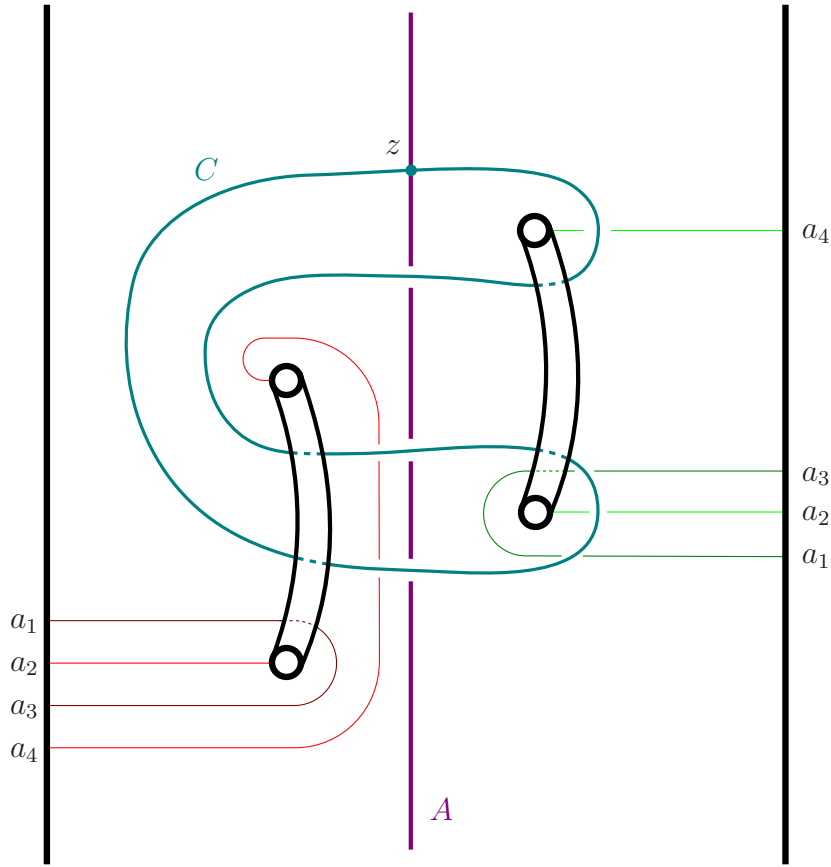}
\caption[The result of isotoping Figure \ref{fig:heegaardknot1}.] {The result
of isotoping Figure \ref{fig:heegaardknot1}. Each boundary component is
identified with $\ZZ$.} \label{fig:heegaardknot2}
\end{figure}

\begin{proposition}
The arced Heegaard diagram $\HH$ (with extra basepoint $w$) shown in Figure
\ref{fig:heegaard1} determines the pair $(\YY,B_3)$.
\end{proposition}

\begin{proof}
As in \cite[Construction 5.4]{LOTBimodules}, by cutting along the arc
$\mathbf{z}$, we obtain a bordered Heegaard diagram with a single boundary
component, $\HH_{dr}$, which we view as rectangle with two tunnels attached.
After attaching $2$-handles to $\HH_{dr} \times [0,1]$ along $\bm{\beta} \times
\{1\}$ and attaching a single $3$-handle, we may view the resulting manifold
$Y(\HH_{dr})$ as $[-1,1] \times \R \times [0,\infty) \subset \R^3$ plus a point at infinity, minus two
tunnels as shown in Figure \ref{fig:heegaardknot1}. The boundary of $Y(\HH_{dr})$ is the union of two
embedded copies of $F^\circ(\ZZ)$ that are determined by the $\alpha$ arcs on
each side; they intersect along a circle $A$. The extra basepoint $w$
determines a knot $C$ in $Y(\HH_{dr})$ with a single point on the boundary: the
union of an arc connecting $w$ to $z$ in the complement of the $\alpha$ arcs
and an arc connecting $z$ to $w$ in the complement of the $\beta$ circles,
pushed into the interior of $Y(\HH_{dr})$ except at $z$. The curves $A$ and $C$
are both shown in Figure \ref{fig:heegaardknot1}.

We obtain Figure \ref{fig:heegaardknot2} from Figure \ref{fig:heegaardknot1} by
an isotopy that slides the tunnel on the right underneath the tunnel on the
left. The circle $A$ can then be identified with the $y$-axis plus the point at
infinity. To obtain $Y(\HH)$, we attach a three-dimensional two-handle along
$A$, which can be seen as $[-\epsilon,\epsilon] \times \R \times (-\infty,0]$
plus the point at infinity. Then $Y(\HH)$ is the complement of a two-component
unlink $(B_1 \cup B_2)$ in $S^3$, and the knot $C$ inside $Y(\HH)$ is $B_3$.
When we identify each component of $\partial \Sigma$ with $\ZZ$, we see that
the $\alpha$ arc connecting the points $a_1$ and $a_3$ is a meridian, and the
$\alpha$ arc connecting $a_2$ and $a_4$ is a $0$-framed longitude, as in the
definition of $\YY$.
\end{proof}

\begin{figure} \centering
 \psfrag{r1}{$\rho_1$} \psfrag{r2}{$\rho_2$} \psfrag{r3}{$\rho_3$}
 \psfrag{s1}{$\sigma_1$} \psfrag{s2}{$\sigma_2$} \psfrag{s3}{$\sigma_3$}
 \psfrag{b1}{{\color{blue} $\beta_1$}}
 \psfrag{b2}{{\color{darkblue} $\beta_2$}}
 \psfrag{aL1}[tl][tl]{{\color{red} $\alpha^L_1$}}
 \psfrag{aL2}{{\color{darkred} $\alpha^L_2$}}
 \psfrag{aR1}{{\color{green} $\alpha^R_1$}}
 \psfrag{aR2}[tl][tl]{{\color{darkgreen} $\alpha^R_2$}}
 \psfrag{w}[cc][cc]{$w$}
 \psfrag{zz}{$\mathbf{z}$}
\includegraphics[scale=0.8]{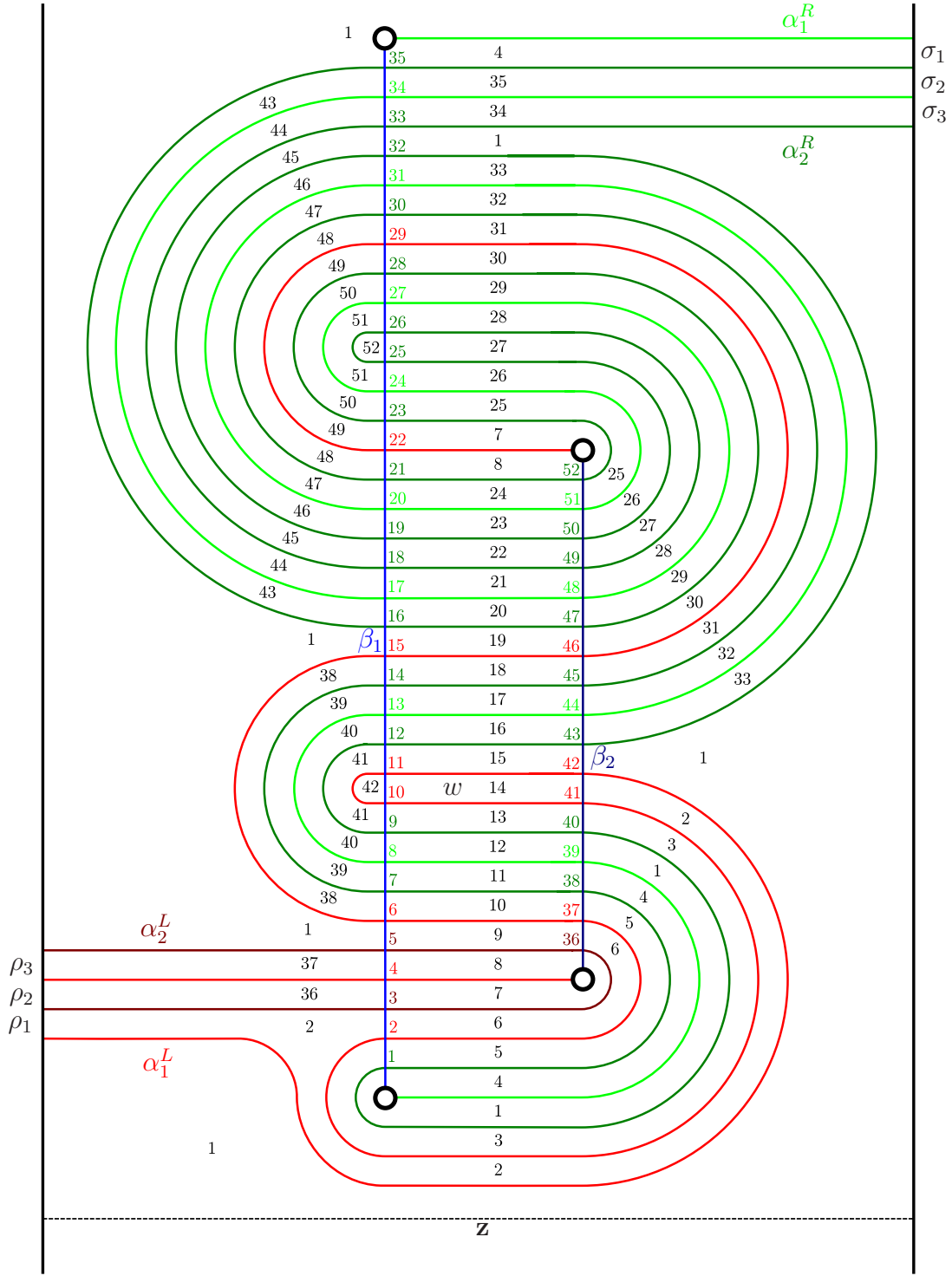}
\caption [The Heegaard diagram $\HH'$.] {The Heegaard diagram $\HH'$, with the
boundary labeled consistently with the conventions for type $D$ structures.}
\label{fig:heegaard2}
\end{figure}

If we try to compute $\CFAA(\HH,0)$ directly, we run into difficulties counting
the holomorphic curves, largely because there is a $14$-sided region that runs
over both handles and shares edges with itself. Instead, it is easier to
perform a sequence of isotopies on the $\alpha$ arcs to obtain the diagram
$\HH'$ shown in Figure \ref{fig:heegaard2}. While $\HH'$ is not a nice diagram
in the sense of Sarkar and Wang \cite{SarkarWang}, the analysis needed to count
the relevant holomorphic curves is vastly simpler. Of course, the drawback is
that the number of generators is much larger.

By Theorem \ref{thm:identityAA}, it suffices to compute $\CFDD(\HH',0)$, as
described previously. Thus, we identify each component of $\partial \Sigma$
with $-\ZZ$. We now describe this computation.

The bimodule $\CFDD(\HH',0)$ is a type $DD$ structure over two copies of the torus algebra $\AA$. We denote these copies by $\AA_\rho$ and $\AA_\sigma$, corresponding to the left and right boundary components of $\HH'$. In $\AA_\sigma$, the Reeb elements are denoted $\sigma_1$, $\sigma_2$, etc. The idempotents in $\AA_\rho$ are denoted $\iota_0^\rho$ and $\iota_1^\rho$, and those in $\AA_\sigma$ are denoted $\iota_0^\sigma$ and $\iota_1^\sigma$. The idempotent maps $I_D^\rho: \mathfrak{S}(\HH',0) \to \{\iota_0^\rho, \iota_1^\rho\}$ and $I_D^\sigma: \mathfrak{S}(\HH',0) \to \{\iota_0^\sigma, \iota_1^\sigma\}$ are defined just as in \eqref{eq:idempotents}.

We denote the regions in $\HH'$ by $R_1, \dots, R_{52}$, as indicated by the black numbers in Figure \ref{fig:heegaard2}.
We label the intersection points of the $\alpha$ and $\beta$ curves $x_1, \dots, x_{52}$, as indicated by the colored numbers. These points are distributed among the various $\alpha$ and
$\beta$ circles as follows:
\[
\begin{array}{|c|c|c|} \hline
 & \beta_1 & \beta_2 \\ \hline
 \alpha_1^L & x_{2},x_{4},x_{6},x_{10},x_{11},x_{15},x_{22},x_{29} & x_{37},x_{41},x_{42},x_{46} \\ \hline
 \alpha_2^L & x_{3},x_{5} & x_{36} \\ \hline
 \alpha_1^R & x_{8},x_{13},x_{17},x_{20},x_{24},x_{27},x_{31},x_{34} & x_{39},x_{44},x_{48},x_{51} \\ \hline
 \alpha_2^R &
  \begin{array}{c} x_{1},x_{7},x_{9},x_{12},x_{14},x_{16},x_{18},x_{19} \\ x_{21},x_{23}, x_{25},x_{26},x_{28},x_{30},x_{32},x_{33},x_{35} \end{array} &
  \begin{array}{c} x_{38},x_{40},x_{43},x_{45} \\ x_{47},x_{49},x_{50},x_{52} \end{array} \\ \hline
\end{array}
\]
The underlying vector space for $\CFDD(\HH',0)$ is generated by the set
$\mathfrak{S}(\HH',0)$, consisting pairs of intersection points with one point
on each $\beta$ circle, one point on either $\alpha_1^L$ or $\alpha_2^L$, and
one point on either $\alpha_1^R$ or $\alpha_2^R$. A simple count shows that
there are $245$ generators.

\subsection{Enumerating index-1 positive domains} \label{subsec:domains}

In order to find all index-1 positive domains in $\HH'$, we begin with the following lemma:

\begin{lemma} \label{lemma:uniquesolution}
For any generators $\x$ and $\y$, the set $\pi_2(\x,\y)$ is nonempty, and there is at most one domain in $\pi_2(\x,\y)$ with any prescribed multiplicities in the six boundary regions ($R_2$, $R_4$, $R_{34}$, $R_{35}$, $R_{36}$, and $R_{37}$).
\end{lemma}

\begin{proof}
For the first statement, the obstruction to $\pi_2(\x,\y)$ being nonempty is an element $\epsilon(\x,\y) \in H_1(Y,\partial Y)$ that is in the image of $H_1(Y) \to H_1(Y, \partial Y)$, and this image is trivial since $H_1(\partial Y) \to H_1(Y)$ is surjective.

The group of periodic domains in $\HH'$ is isomorphic to $H_2(Y,\partial Y) \cong \Z^2$; it is freely generated by
\[
\begin{aligned}
P_1 &= R_2 + R_6 + R_7 - R_{10} - \dots - R_{13} - R_{15} - \dots - R_{18} + R_{25} + \dots + R_{30} + R_{36} \\
& \qquad - R_{38} - \dots - R_{41} + R_{49} + \dots + R_{52} \\
P_2 &= R_4 + \dots + R_{11} + R_{17} + \dots + R_{20} + R_{24} + R_{25} + R_{29} + \dots + R_{32} + R_{35} \\
& \qquad - R_{40} - R_{41} - R_{42} - R_{44} - R_{45} - R_{46} - R_{51} - R_{52}.
\end{aligned}
\]
Thus, any nonzero periodic domain has a nonzero multiplicity at either $R_2$ or $R_4$, so there are no nonzero provincial periodic domains. The uniqueness statement then follows immediately.
\end{proof}

We may algorithmically find all the positive domains with index $1$ by the following procedure. By Proposition \ref{prop:chordsallowed}, the multiplicity of each of the six boundary regions must be $0$ or $1$. For each of the $2^6$ choices of boundary multiplicities and each pair of generators $\x, \y$ (subject to the idempotent restrictions of Proposition \ref{prop:chordsallowed}), we may solve \eqref{eq:pointconditions} to find the unique domain in $\pi_2(\x,\y)$ with the prescribed boundary multiplicities, if one exists. We may then list only those solutions which represent positive classes and have index $1$ for some compatible $\vec\rho$, where the index is computed using \eqref{eq:index}. Specifically, note that if $B$ is domain representing a class in $\pi_2(\x,\y)$ with boundary multiplicities all $0$ or $1$, the quantity $\abs{\vec \rho} + \iota(\vec \rho)$ in \eqref{eq:index} equals $0$ if $B$ is provincial, $\frac12$ if it abuts one component of $\partial \Sigma$, and $1$ if it abuts both components. \label{pageref:iota} The Euler measure of $B$ equals the sum of of the Euler measures of its regions (namely $1-\frac{k}{2}$ for a $2k$-gon), weighted by their multiplicities. Using a \emph{Mathematica} computation, we find that there are $1,013$ positive index-1 domains satisfying the restrictions of Proposition \ref{prop:chordsallowed}.\footnote{More precisely, we mean that there are $1,013$ tuples $(\x,\y,B)$, where $\x$ and $\y$ are generators and $B \in \pi_2(\x,\y)$ is a positive class with index $1$. In some cases, the same domain $B$ may be used for different pairs of generators, such as when $B$ is a bigon. We shall use this abuse of terminology throughout this section.}

We now partition the $1,013$ positive domains with index $1$ into classes that share the same holomorphic geometry and discuss each case that arises. The results are summarized in Table \ref{table:domains}.

\begin{table}
\begin{tabular}{|l|l|l|l|}
  \hline
  Type of domain & Examples & Quantity & Count for differential? \\ \hline
  Bigons & $R_{36}$, $R_{42}$  & 488 & Yes \\
  Quadrilaterals & $R_4$, $R_5$ & 167 & Yes \\
  Domains with a boundary cut & $D_1$, $D_2$, $D_3$, $D_4$ & 52 & Yes \\
  Domains without a boundary cut & $D_5$, $D_6$, $D_{20}$ & 171 & No \\
  Disconnected domains & $D_7$ & 37 & No \\
  Indecomposable annuli & $D_8$, $D_9$, $D_{10}$ & 35 & Yes \\
  Singly decomposable annuli & $D_{11}$, $D_{12}$, $D_{13}$ & 18 & No* \\
  Doubly decomposable annuli & $D_{14}$, $D_{15}$, $D_{16}$ & 7 & No* \\
  Good tori & $D_{17}$, $D_{19}$ & 29 & Yes \\
  Conditional tori & $D_{18}$ & 9 & Yes* \\
  \hline
\end{tabular}
\caption{Summary of the different types of domains, along with whether or not they count for the differential on $\CFDD(\HH',0)$. The starred entries in the fourth column hold when the complex structure on $\Sigma$ is sufficiently stretched.} \label{table:domains}
\end{table}

\subsubsection*{Bigons and quadrilaterals}

In the context of closed Heegaard diagrams, Sarkar and Wang \cite{SarkarWang} showed that in a Heegaard diagram in which
every non-basepointed region is either a bigon or a quadrilateral, the domains with Maslov index $1$ are precisely the embedded bigons and quadrilaterals that are embedded in the Heegaard diagram, and these all support support a unique
holomorphic representatives. (Such a Heegaard diagram is called \emph{nice}.) Lipshitz, Ozsv\'ath, and Thurston proved an
analogous result for bordered diagrams \cite[Proposition 8.4]{LOTBordered}, where now we extend the definition of ``quadrilateral'' to include a region with boundary consisting of one segment of a $\beta$ circle, two segments of $\alpha$ arcs, and one segment of $\partial \Sigma$. The only non-basepointed regions in $\HH'$ that are not bigons or quadrilaterals are $R_2$, $R_4$, $R_7$, and $R_8$, which are hexagons. Therefore, any index-1 domain on our list that does not use one of these four regions automatically supports a unique holomorphic representative.

We may easily find several additional families of domains that are embedded bigons or quadrilaterals, perhaps with one or more boundary punctures, which use at least one of the regions $R_2$, $R_4$, $R_7$, or $R_8$. For instance, $R_2 + R_6 + R_{14} + R_{42}$ is a boundary-punctured bigon from $x_3 x_i$ to $x_2 x_i$ (for any $i \in \{38, 39, 40, 43, 44, 45, 47, \dots, 52\}$) with chord marked $\rho_1$, and $R_4$ is a boundary-punctured rectangle from $x_{35} x_{39}$ to $x_1 x_{38}$ with a chord marked $\sigma_1$.

In total, we find some 488 bigons and 167 quadrilaterals.

\begin{figure}
\begin{center}
\labellist
 \pinlabel (a) at 18 145
 \pinlabel (b) at 203 145
 \pinlabel $x_{15}$ at 74 141
 \pinlabel $x_{22}$ at 80 10
 \pinlabel {{\color{red} $\alpha_1^L$}} at 8 128
 \pinlabel {{\color{blue} $\beta_1$}} at 148 128
 \pinlabel {{\color{darkred} $\alpha_2^L$}} at 80 68
 \pinlabel $\rho_2$ at 89 93
 \pinlabel $\rho_3$ at 72 93
 \pinlabel $x_{23}$ at 192 88
 \pinlabel $x_{15}$ at 274 15
 \pinlabel $x_{46}$ at 358 88
 \pinlabel $x_{52}$ at 271 144
 \pinlabel {{\color{red} $\alpha_1^L$}} at 210 128
 \pinlabel {{\color{blue} $\beta_1$}} at 334 128
 \pinlabel {{\color{darkblue} $\beta_2$}} at 223 58
 \pinlabel {{\color{darkgreen} $\alpha_2^R$}} at 326 58
 \pinlabel {{\color{darkred} $\alpha_2^L$}} at 275 69
 \pinlabel $\rho_2$ at 285 94
 \pinlabel $\rho_3$ at 267 94
\endlabellist
\includegraphics{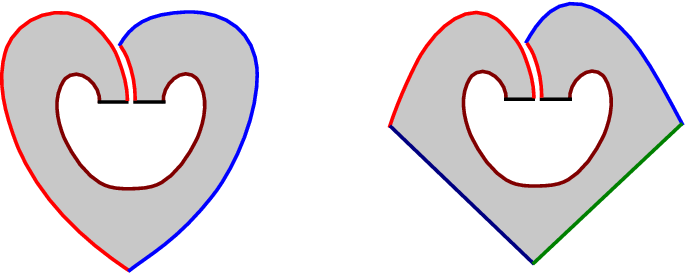}
\end{center}
\caption[The domains $D_1$ and $D_3$.] {The domains $D_1$ (a) and $D_3$ (b).} \label{fig:rho23}
\end{figure}
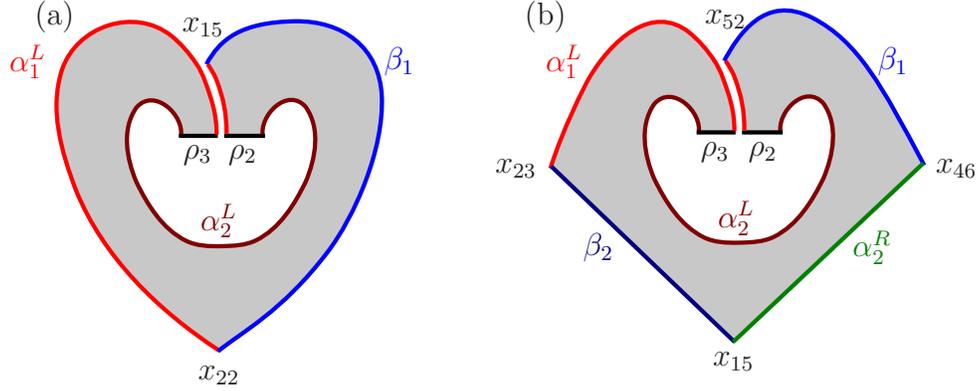

\subsubsection*{Domains with a boundary cut}
Let
\[
D_1 = R_7+R_8+ R_{19} + \dots + R_{30} + R_{36} + R_{37} + R_{49} + R_{50} + R_{51} + R_{52},
\]
For any $i \in \{38,39,40, 43,44,45\}$ $D_1$ represents a class in $\pi_2(x_{15} x_i, x_{22} x_i)$ and has index $1$ with respect to the sequence $(-\rho_2, -\rho_3)$. (If $i \in \{47, \dots, 52\}$, the index is $3$ rather than $1$.) To obtain
a holomorphic representative of $D_1$ compatible with $(-\rho_2, -\rho_3)$, we cut along $\alpha_1^L$ all the way to the boundary, as shown in Figure \ref{fig:rho23}. Thus, we parametrize $D_1$ as a bigon with two separate boundary punctures  rather than as an annulus with a single puncture (which is prohibited by Proposition \ref{prop:chordsallowed}). It is then straightforward to see that $D_1$ supports a unique holomorphic representative and thus provides a differential $x_{15} x_i \xrightarrow{\rho_{23} \otimes \one} x_{22} x_i$ for each $i$ as above. Likewise, for each $j \in \{36, 37, 41, 42, 46\}$, the domain
\[
D_2 = R_4 + R_{11} + R_{17} + R_{20} + R_{24} + R_{25} + R_{29} + R_{32} + R_{35} + R_{39} + R_{43} + R_{47} + R_{50},
\]
representing a class in $\pi_2(x_{35} x_j, x_1 x_j)$, contributes a differential $x_{35} x_j \xrightarrow{\one \otimes \sigma_{12}} x_1 x_j$. In fact, $D_1$ and $D_2$ are the only domains of this form (so they account for $11$ of the $1,013$ total classes).

Similarly, the domains
\begin{align*}
D_3 &= R_7 +R_8 + R_{31} + R_{36}+R_{37} + R_{48} \\
D_4 &= D_2 + R_4 + R_5 + R_{10} + R_{11} + R_{17} + R_{18} + R_{24} + R_{25} +
R_{29} + R_{30} + R_{32} \\
& \qquad + R_{36} + R_{37} + R_{38} + R_{39} + R_{47} +R_{49} + R_{50},
\end{align*}
respectively represent index-$1$ classes in $\pi_2(x_{22} x_{45}, x_{23} x_{46})$, and $\pi_2(x_{35} x_{46}, x_{2}
x_{48})$. The source curve for $D_3$ or $D_4$ is a quadrilateral, with two boundary punctures on one $\alpha$
edge mapping to $-\rho_2$ and $-\rho_3$ and (for $D_4$) a boundary puncture on the other $\alpha$ edge mapping to $\sigma_1$. It is easy to see that these classes all support holomorphic representatives. Thus, we have differentials
$x_{15} x_i \xrightarrow{\rho_{23} \otimes \one} x_{22} x_i$, $x_{22} x_{45} \xrightarrow{\rho_{23} \otimes \one} x_{23} x_{46}$, and $x_{35} x_{46} \xrightarrow{\rho_{23} \otimes \sigma_1} x_2 x_{48}$. We find $41$ domains of this form.

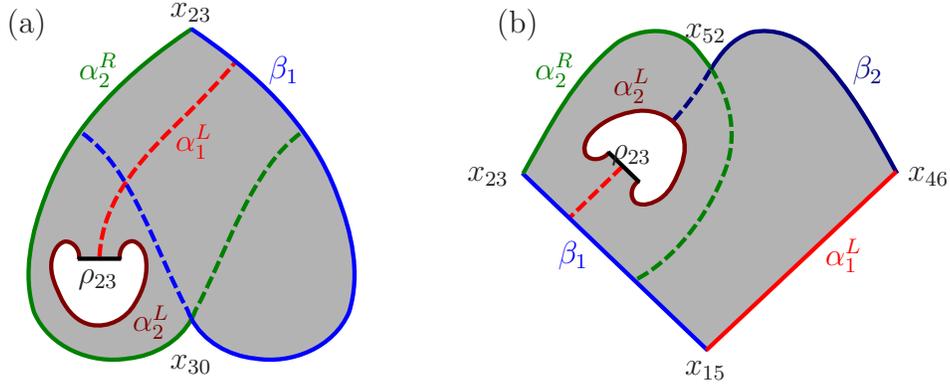
\begin{figure}
\begin{center}
\labellist
 \pinlabel (a) at 18 145
 \pinlabel (b) at 203 145
 \pinlabel $x_{23}$ at 80 150
 \pinlabel $x_{30}$ at 80 17
 \pinlabel {{\color{darkgreen} $\alpha_2^R$}} at 45 128
 \pinlabel {{\color{blue} $\beta_1$}} at 115 128
 \pinlabel {{\color{darkred} $\alpha_2^L$}} at 65 33
 \pinlabel {{\color{red} $\alpha_1^L$}} at 81 101
 \pinlabel $\rho_{23}$ at 45 49
 \pinlabel $x_{23}$ at 192 88
 \pinlabel $x_{15}$ at 274 15
 \pinlabel $x_{46}$ at 358 88
 \pinlabel $x_{52}$ at 274 142
 \pinlabel {{\color{darkgreen} $\alpha_2^R$}} at 217 128
 \pinlabel {{\color{darkblue} $\beta_2$}} at 335 128
 \pinlabel {{\color{blue} $\beta_1$}} at 224 58
 \pinlabel {{\color{red} $\alpha_1^L$}} at 326 58
 \pinlabel {{\color{darkred} $\alpha_2^L$}} at 246 122
 \pinlabel $\rho_{23}$ at 246 96
\endlabellist
\includegraphics{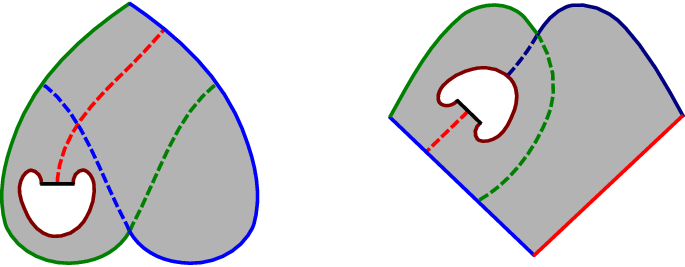}
\end{center}
\caption[The domains $D_5$ and $D_6$.] {The domains $D_5$ (a) and $D_6$ (b).} \label{fig:rho23cantcut}
\end{figure}

\subsubsection*{Domains without a boundary cut}

Let
\[
\begin{aligned}
D_5 &= R_7 + R_8 + R_{36} + R_{37} + R_{48} + R_{49} + R_{50} + R_{51} + R_{52}\\
D_6 &= R_7 + R_8 + R_{19} + R_{20} + R_{21} + R_{22} + R_{23} + R_{24} + R_{36} + R_{37},
\end{aligned}
\]
illustrated in Figure \ref{fig:rho23cantcut}. $D_5$ represents a class in $\pi_2(x_{30} x_j ,x_{23} x_j)$ for each $j \in \{37,41,42,46\}$, and $D_6$ represents a class in $\pi_2(x_{15}x_{52}, x_{23}x_{46})$. We cannot cut these domains along $\alpha_1^L$ as we did with $D_3$ and $D_4$, since in each case, as we travel along $\alpha_1^L$ from the intersection point of $\rho_2 \cap \rho_3$, we reach the boundary of $D_5$ (resp. $D_6$) before reaching $x_j$ (resp. $x_{46}$). Thus, $D_5$ and $D_6$ cannot admit holomorphic representatives. We find $106$ domains like $D_5$ and $62$ domains like $D_6$. (Some of these domains have additional $\sigma_2$ or $\sigma_3$ punctures on their boundaries, but these do not affect the argument above.)

We also find domains such as
\[
D_7 = R_7 + R_8 + R_{14} + R_{36} + R_{37} + R_{42}
\]
which are the disjoint union of an annulus, one of whose boundary component equals $\alpha_2^L \cup \rho_2 \cup \rho_3$, and a bigon. Since the $\alpha_2^L \cup \rho_2 \cup \rho_3$ boundary component of the annulus does not contain a point of either the source or the target generator, there is no way to find a source surface representing $D_7$. We find $37$ domains of this form.

\subsubsection*{Indecomposable annuli}

\begin{figure}
\begin{center}
\labellist
 \pinlabel (a) at 18 158
 \pinlabel (b) at 176 158
 \pinlabel (c) at 348 158
 \pinlabel $x_{4}$ at 80 155
 \pinlabel $x_{6}$ at 80 22
 \pinlabel $x_{38}$ at 80 72
 \pinlabel $x_{52}$ at 80 92
 \pinlabel {{\color{red} $\alpha_1^L$}} at 43 128
 \pinlabel {{\color{blue} $\beta_1$}} at 117 128
 \pinlabel {{\color{darkgreen} $\alpha_2^R$}} at 63 81
 \pinlabel {{\color{darkblue} $\beta_2$}} at 99 81
 \pinlabel $x_6$ at 237 147
 \pinlabel $x_{37}$ at 237 13
 \pinlabel $x_{23}$ at 218 87
 \pinlabel $x_{21}$ at 218 108
 \pinlabel {{\color{darkgreen} $\alpha_2^R$}} at 235 96
 \pinlabel {{\color{darkblue} $\beta_2$}} at 307 128
 \pinlabel {{\color{blue} $\beta_1$}} at 200 96
 \pinlabel {{\color{red} $\alpha_1^L$}} at 214 29
 \pinlabel {{\color{darkred} $\alpha_2^L$}} at 169 128
 \pinlabel $\rho_{123}$ at 180 50
 \pinlabel $x_4$ at 395 155
 \pinlabel $x_{22}$ at 395 11
 \pinlabel $x_j$ at 399 83
 \pinlabel {{\color{red} $\alpha_1^L$}} at 356 128
 \pinlabel {{\color{blue} $\beta_1$}} at 433 128
 \pinlabel {{\color{darkgreen} $\alpha_2^R$}} at 423 92
 \pinlabel {{\color{darkblue} $\beta_2$}} at 385 94
\endlabellist
\includegraphics{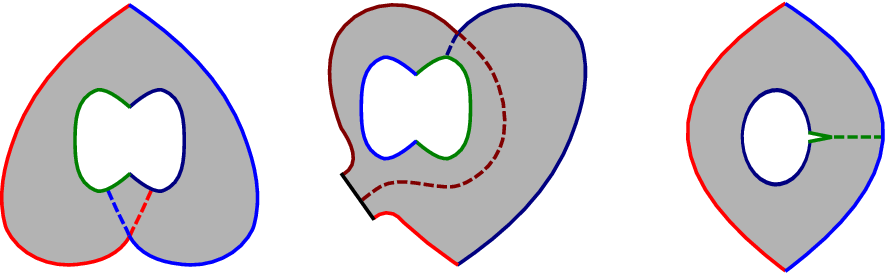}
\end{center}
\caption{The indecomposable annuli $D_8$ (a), $D_9$ (b), and $D_{10}$ (c).} \label{fig:indecomposable}
\end{figure}
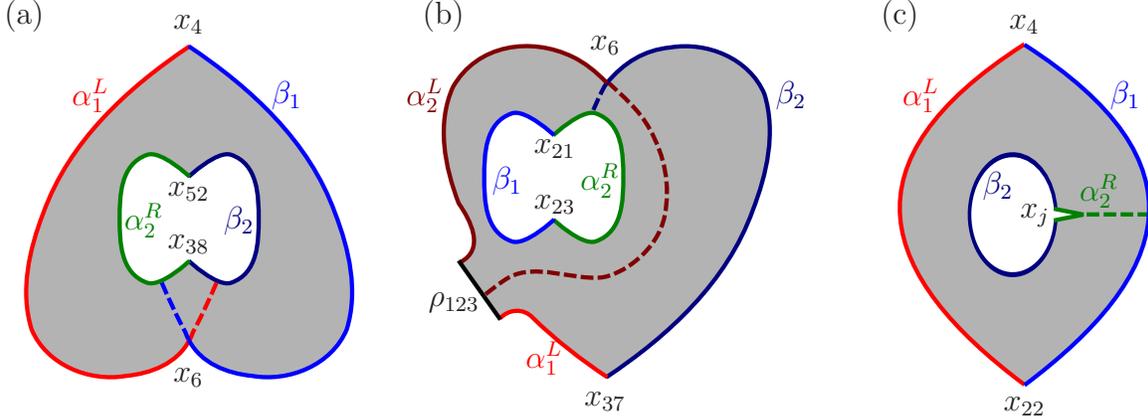

Consider the domains
\[
\begin{aligned}
D_8 &= R_8 + R_9 + R_{10} + R_{18} + R_{31} + R_{38} + R_{48} \in \pi_2(x_4x_{38}, x_6x_{52}) \\
D_9 &= R_2 + R_6 + R_7 + R_8 + R_{14} + R_{36} + R_{37} + R_{42} \in \pi_2(x_{21}x_{36}, x_{23} x_{37})\\
\end{aligned}
\]
shown in Figure \ref{fig:indecomposable}. Each of these domains is topologically an annulus with three convex corners and one concave corner and cannot be decomposed as the composition of an index-0 annulus and a bigon (in contrast to another family of annuli considered below). As illustrated in Figure \ref{fig:indecomposable}, we call the boundary component containing the convex corner the \emph{outer boundary} and the other component the \emph{inner boundary}. There are a total of eight domains with the geometry of $D_8$ and fifteen with the geometry of $D_9$.

\begin{lemma}
For any complex structure on $\Sigma$, the moduli spaces $\MM(D_8)$ and $\MM(D_9)$ each contain an odd number of points. Thus, $D_8$ and $D_9$ always count for the differential.
\end{lemma}

\begin{proof}
Given a choice of complex structure on $\Sigma$, each domain admits a one-dimensional family of conformal structures, depending on the value of a cut parameter $c \in \R$, where $c < 0$ corresponds to cutting along the $\alpha$ curve and $c>0$ corresponds to cutting along the $\beta$ curve. For each value of the cut parameter $c$, let $\theta_0(c)$ (resp.~$\theta_1(c)$) denote the ratio of the conformal length of the $\alpha$ arc of the outer (resp.~inner) boundary to the conformal length of the $\beta$ arc of the outer (resp.~inner) boundary. A given conformal structure admits a holomorphic involution if and only if $\theta_0(c) = \theta_1(c)$, so the number of points in the moduli space of each domain (modulo 2) equals the number of zeros of the function $f(c) = \theta_0(c) - \theta_1(c)$, which for generic choices of complex structure on $\Sigma$ can be assumed to be transverse to $0$. This number is determined by the limiting behavior of $f(c)$, as follows.

For $D_8$, the cut in the $\alpha$ direction approaches the $\beta$ arc of the inner boundary and the cut in the $\beta$ direction approaches the $\alpha$ arc of the inner boundary. Thus, in the limit as we cut in the $\alpha$ direction, $\theta_0(c)$ becomes arbitrarily large and $\theta_1(c)$ approaches $0$, so $\lim_{c \to -\infty} f(c) = \infty$. Similarly, as we cut in the $\beta$ direction, $\theta_0(c)$ approaches $0$ and $\theta_1(c)$ becomes arbitrarily large, so $\lim_{c \to \infty} f(c) = - \infty$. By transversality and the intermediate value theorem, $f(c)$ has an odd number of zeros.

For $D_9$, there is a Reeb chord marked $\rho_{123}$ on the outer boundary. The cut in the $\beta$ direction approaches the $\alpha$ arc of the inner boundary, while the cut in the $\alpha$ direction approaches this boundary puncture. Thus, as we cut in the $\alpha$ direction, $\theta_0(c)$ becomes arbitrarily large, while $\theta_1(c)$ approaches a finite value. Thus, $\lim_{c \to -\infty} f(c) = \infty$, while $\lim_{c \to \infty} f(c) = -\infty$ just as with $D_8$.
\end{proof}

Similarly, the annular domain
\[
D_{10} = \sum_{i=8}^{24} R_i
\]
represents an index-1 class in $\pi_2(x_4 x_j, x_{22} x_j)$ for any of the twelve points $x_j \in (\alpha_1^R \cup \alpha_2^R) \cap \beta_2$. This domain similarly admits a 1-dimensional family of conformal structures given by a cut parameter by $c \in (0, \infty)$. As we increase the cut parameter, the ratio of the length of the $\alpha$ segment of the boundary component containing $x_j$ to the $\beta$ segment of the same tends from $0$ to infinity, while the same ratio on the opposite boundary component tends from a finite value to $0$. Thus, $D_{10}$ counts for the differential for each choice of $x_j$.

\subsubsection*{Decomposable annuli}

We next consider domains whose moduli spaces may depend nontrivially on the choice of complex structure. As a preliminary, let $\eta \subset \Sigma$ be a simple closed curve passing through the regions $R_{1}$, and $R_4, \dots, R_{12}$. For a given complex structure $J$ on $\Sigma$ and $T \in [0, \infty)$, let $J_T$ denote the complex structure obtained by ``stretching the neck'' along $\eta$ by inserting an annulus of width $T$.

Consider the following domains:
\begin{equation} \label{eq:index0annuli}
\begin{aligned}
A_1 &= R_7 + R_8 + R_{48} + R_{49} + R_{30} + R_{31} \\
A_2 &= A_1 + R_{24} + R_{25} + R_{47} + R_{50} + R_{29} + R_{32} \\
A_3 &= A_2 + R_{23} + R_{26} + R_{46} + R_{51} + R_{28} + R_{33} \\
A_4 &= A_1 + R_6 + R_9 \\
A_5 &= A_2 + R_6 + R_9 \\
A_6 &= A_3 + R_6 + R_9 \\
\end{aligned}
\end{equation}
Each of these domains is an index-$0$ annulus, with one boundary component consisting of a segment of $\beta_1$ and a segment of $\alpha_1^L$ or $\alpha_2^L$ and the other consisting of a segment of $\beta_2$ and a segment of $\alpha_1^R$ or $\alpha_2^R$. We call these the two boundary components the \emph{outer boundary} and \emph{inner boundary}, respectively. A choice of complex structure $J$ on $\Sigma$ completely determines a conformal structure on each $A_i$. Let $\Theta_0^i(J)$ (resp. $(\Theta_1^i(J))$) denote the ratio of the conformal length of the $\alpha$ segment of the inner (resp.~outer) boundary of $A_i$ to the conformal length of the $\beta$ segment of the inner (resp.~outer) boundary. We say that $J$ is \emph{sufficiently stretched} if $\Theta_0^i(J) > \Theta_1^i(J)$ for each $i=1, \dots, 6$.

\begin{lemma} \label{lemma:stretched}
For any complex structure $J$ on $\Sigma$ there exists a number $T_0 = T_0(J)$ such that for any $T>T_0$, the complex structure $J_T$ is sufficiently stretched.
\end{lemma}

\begin{proof}
For each $i$, the only intersections of the curve $\eta$ with the boundary of $A_i$ are on the $\alpha$ segment of the outer boundary, so stretching the neck along $\eta$ increases the conformal length of that segment relative to the $\beta$ segment of the outer boundary. Therefore, for large values of $T$, $\Theta_0^i(J_T)$ can be made arbitrarily large, while $\Theta_1^i(J_T)$ approaches some finite value.
\end{proof}

\begin{figure}
\begin{center}
\labellist
 \pinlabel (a) at 18 158
 \pinlabel (b) at 169 158
 \pinlabel (c) at 318 158
 \pinlabel $x_{4}$ at 80 143
 \pinlabel $x_{5}$ at 80 11
 \pinlabel $x_{45}$ at 80 56
 \pinlabel $x_{47}$ at 80 95
 \pinlabel $\rho_2$ at 49 156
 \pinlabel {{\color{red} $\alpha_1^L$}} at 66 150
 \pinlabel {{\color{darkred} $\alpha_2^L$}} at 32 50
 \pinlabel {{\color{blue} $\beta_1$}} at 126 150
 \pinlabel {{\color{darkgreen} $\alpha_2^R$}} at 63 83
 \pinlabel {{\color{darkblue} $\beta_2$}} at 98 83
 \pinlabel $x_3$ at 231 143
 \pinlabel $x_2$ at 231 11
 \pinlabel $x_{45}$ at 199 91
 \pinlabel $x_{47}$ at 199 135
 \pinlabel {{\color{darkgreen} $\alpha_2^R$}} at 191 113
 \pinlabel {{\color{darkblue} $\beta_2$}} at 210 113
 \pinlabel {{\color{blue} $\beta_1$}} at 278 50
 \pinlabel {{\color{red} $\alpha_1^L$}} at 205 27
 \pinlabel {{\color{darkred} $\alpha_2^L$}} at 162 93
 \pinlabel $\rho_{123}$ at 175 47
 \pinlabel $x_{3}$ at 380 143
 \pinlabel $x_{6}$ at 380 11
 \pinlabel $x_{45}$ at 409 76
 \pinlabel $x_{47}$ at 409 126
 \pinlabel $\rho_1$ at 351 156
 \pinlabel {{\color{darkred} $\alpha_2^L$}} at 366 150
 \pinlabel {{\color{red} $\alpha_1^L$}} at 333 50
 \pinlabel {{\color{blue} $\beta_1$}} at 426 50
 \pinlabel {{\color{darkgreen} $\alpha_2^R$}} at 398 100
 \pinlabel {{\color{darkblue} $\beta_2$}} at 420 100
\endlabellist
\includegraphics{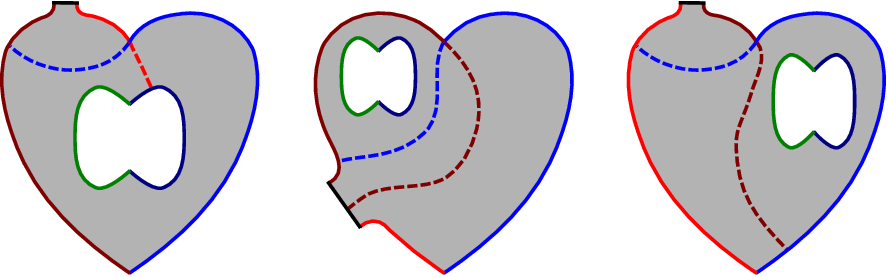}
\end{center}
\caption{The decomposable annuli $D_{11}$ (a), $D_{14}$ (b), and $D_{15}$ (c).} \label{fig:decomposable}
\end{figure}
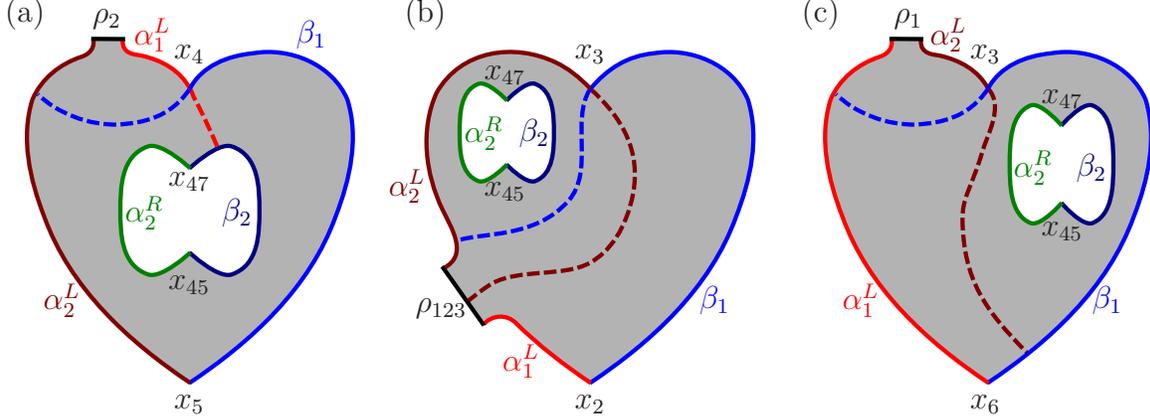

Consider the index-$1$ annuli
\begin{equation} \label{eq:decomposable}
\begin{aligned}
D_{11} &= A_1 + R_{36} \in \pi_2(x_4 x_{45}, x_5 x_{47}) \\
D_{12} &= A_1 + R_{37}\in \pi_2(x_3 x_{45}, x_4 x_{47}) \\
D_{13} &= A_1 + R_{25} + R_{26} + R_{27} + R_{28} + R_{29} + R_{50} + R_{51} + R_{52} \\
& \qquad \qquad \qquad \in \pi_2(x_3 x_{45}, x_5 x_{52})  \\
D_{14} &= A_1 + R_2 + R_{14} + R_{36} + R_{37} + R_{42} \in \pi_2(x_3 x_{45}, x_2 x_{47}) \\
D_{15} &= A_1 + R_2 + R_6 + R_9 + R_{14} + R_{42} \in \pi_2(x_3 x_{45}, x_6 x_{47})\\
D_{16} &= A_1 + R_{20} + R_{24} + R_{25} + R_{29} + R_{32} +R_{35} + R_{43} + R_{47} + R_{50} \\
& \qquad \qquad \qquad \in \pi_2(x_3 x_{44}, x_5 x_{47}),
\end{aligned}
\end{equation}
some of which are shown in Figure \ref{fig:decomposable}. Each of these annuli can be written as a sum of an index-0 annulus and an adjacent bigon, so we call these domains \emph{decomposable}. It is easy to find eighteen other domains of this form, where we take $A_2, \dots, A_6$ in place of $A_1$ in \eqref{eq:decomposable} as applicable. Note that $D_{15}$ and $D_{16}$ can each be decomposed into the sum of an index-0 annulus and an adjacent bigon in a second way as well:

\begin{equation} \label{eq:doublydecomposable}
\begin{aligned}
D_{15} &= A_4 + R_{2} + R_{14} + R_{42} \\
D_{16} &= A_2 + R_{20} + R_{35} + R_{43},
\end{aligned}
\end{equation}
We therefore call these domains \emph{doubly decomposable}.

\begin{lemma} \label{lemma:decomposable}
If $J$ is sufficiently stretched, the moduli spaces of all of the decomposable annuli each contain an even number of points. Thus, these domains do not count for the differential.
\end{lemma}

\begin{proof}
We begin with $D_{11}$. Just as with the indecomposable annuli discussed above, there is a $1$-dimensional family of conformal structures on $D_{11}$ given by a cut parameter at $x_4$. As we cut along $\alpha_1^L$, the cut approaches the $\beta$ arc of the inner boundary, $\theta_0(c)$ becomes arbitrarily large while $\theta_1(c)$ approaches $0$, so $\lim_{c \to -\infty} f(c) = \infty$. On the other hand, cutting along $\beta_1$ degenerates $D_{11}$ into $A_1$ and a bigon (with a Reeb chord). By Gromov compactness, in the limit as $c \to \infty$, the ratios $\theta_0(c)$ and $\theta_1(c)$ approach the corresponding parameters for $A_1$, namely $\Theta_0^1(J_T)$ and $\Theta_1^1(J_T)$. By Lemma \ref{lemma:stretched}, if we choose $T$ large enough that $\Theta_0^1(J_T) > \Theta_1^1(J_T)$, we see that $\lim_{c \to \infty} f(c) >0$, so $f(c)$ has an even number of zeroes, as required.

The arguments for $D_{12}$, $D_{13}$, and $D_{14}$ are very similar. The one modification for $D_{14}$ is that as we cut along $\alpha_2^L$ at $x_3$ out to the boundary puncture, $\theta_1(c)$ approaches a finite value that is not necessarily zero, just as we saw with $D_8$ above. However, $\theta_0(c)$ still approaches $\infty$, so the remainder of the argument carries through unchanged.

For $D_{15}$, cutting along $\alpha_2^L$ at $x_3$ decomposes the domain as in \eqref{eq:decomposable}, while cutting along $\beta_1$ decomposes it as in \eqref{eq:doublydecomposable}. Therefore, $\lim_{c \to -\infty} f(c) = \Theta_0^1(J_T) - \Theta_1^1(J_T)$ and $\lim_{c \to \infty} f(c) = \Theta_0^4(J_T) - \Theta_1^4(J_T)$. By Lemma \ref{lemma:stretched}, we may choose $T$ large enough such that both of these limits are positive numbers, which implies that $f(c)$ has an even number of zeroes. The same analysis goes through for $D_{16}$.
\end{proof}

\subsubsection*{Genus-1 classes}
Having analyzed all the classes represented by planar surfaces, we now turn to classes that are represented by surfaces of genus $1$. It is difficult to determine whether these classes support holomorphic representatives using direct conformal geometry arguments as above. Instead, we will look at how these domains arise in the broken flowlines that are the ends of 1-dimensional moduli spaces --- specifically, the fact the relation $\partial^2 = 0$ and its more complicated $\AA_\infty$ analogues --- to deduce the behavior of these domains indirectly. We shall see that knowledge of the planar classes completely determines which of the genus-1 classes count for the differential.

\begin{figure}
\begin{center}
\labellist
 \pinlabel (a) at 10 150
 \pinlabel (b) at 134 150
 \pinlabel $x_3$ at 13 14
 \pinlabel $x_{23}$ at 13 135
 \pinlabel $x_{36}$ at 89 57
 \pinlabel $x_{52}$ at 89 92
 \pinlabel {{\color{darkred} $\alpha_2^L$}} at 104 35
 \pinlabel {{\color{blue} $\beta_1$}} at 11 76
 \pinlabel {{\color{darkgreen} $\alpha_2^R$}} at 104 114
 \pinlabel {{\color{darkblue} $\beta_2$}} at 88 76
 \pinlabel $x_3$ at 135 14
 \pinlabel $x_{16}$ at 135 135
 \pinlabel $x_{36}$ at 211 57
 \pinlabel $x_{45}$ at 211 92
 \pinlabel {{\color{darkred} $\alpha_2^L$}} at 226 35
 \pinlabel {{\color{blue} $\beta_1$}} at 133 76
 \pinlabel {{\color{darkgreen} $\alpha_2^R$}} at 244 114
 \pinlabel {{\color{darkblue} $\beta_2$}} at 210 76
\endlabellist
\includegraphics{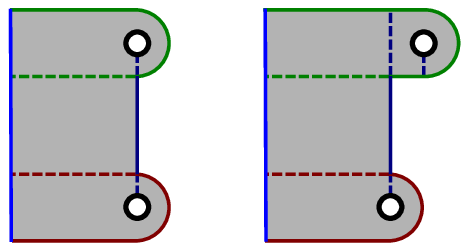}
\end{center}
\caption{The embedded genus-1 domains $D_{17}$ (a) and $D_{18}$ (b).} \label{fig:genus1}
\end{figure}

Consider the domains
\[
\begin{aligned}
D_{17} &= R_7 + \dots + R_{24} \in \pi_2(x_3 x_{52}, x_{23} x_{36}) \\
D_{18} &= R_7  + \dots + R_{19} + R_{30} + R_{31} + R_{48} + R_{49} \in \pi_2(x_3 x_{45}, x_{16} x_{36}) \\
\end{aligned}
\]
shown in Figure \ref{fig:genus1}. Each of these is represented by an embedded or immersed genus-$1$ surface with one boundary component. Any domain of the form $\sum_{i=a}^b R_i$, where $a \in \{4,5,6,7\}$ and $b \in \{24, \dots, 33\}$ are chosen such that the two $\alpha$ segments of the boundary do not lie on the same $\alpha$ curve, has the same holomorphic geometry as $D_{17}$, for a total of 21 domains. Likewise, there are a total of 9 domains with the same geometry as $D_{18}$.

\begin{lemma} $ \ $
\begin{enumerate}
\item For any complex structure $J$ on $\Sigma$, the moduli space $\MM_J(D_{17})$ contains an odd number of points, so $D_{17}$ counts for the differential.

\item If the complex structure $J$ on $\Sigma$ is sufficiently stretched, the moduli space $\MM_J(D_{18})$ contains an odd number of points, so $D_{18}$ counts for the differential.
\end{enumerate}
\end{lemma}

\begin{proof}
Let
\[
\begin{aligned}
E_1 &= R_{36} \in \pi_2(x_4 x_{52}, x_3 x_{52}) \\
E_2 &= R_{20} + \dots + R_{29} + R_{50} + R_{51} + R_{52} \in \pi_2(x_{16} x_{36}, x_{21} x_{36})
\end{aligned}
\]
Each of these domains obviously counts for the differential, and the compositions $E_1*D_{17} \in \pi_2(x_4 x_{52}, x_{23} x_{36})$ and $D_{18}*E_2 \in \pi_2 (x_3 x_{45}, x_{29} x_{52})$ are index-2 positive domains. The moduli spaces $\MM(E_1 * E_{17})$ and $\MM(D_{18} * E_2)$ are $1$-dimensional manifolds, so they each have an even number of ends.

By inspecting the list of all $1,013$ index-1 positive domains, it is easy to verify that the only other way to decompose $E_1*D_{17}$ into two such domains is as the composition of the indecomposable annulus $D_{10} \in \pi_2(x_4 x_{52}, x_{22} x_{52})$ and the boundary-punctured rectangle $R_7 + R_{36} \in \pi_2(x_{22} x_{52}, x_{23} x_{36})$, each of which admits a unique holomorphic representative. In order for $\MM(E_1 * D_{17})$, which is a $1$-dimensional manifold, to have an even number of ends, we see that $D_{17}$ must have an odd number of holomorphic representatives, regardless of the choice of complex structure, and thus must count for the differential.

For $D_{19}$, the situation is slightly more complicated. The composition $D_{18} * E_2$ can be split up in two other ways: (a) as the genus-1 domain $R_7 + \dots + R_{31} \in \pi_2(x_3 x_{45}, x_{30} x_{46})$, which has the same holomorphic geometry as $D_{17}$, composed with the bigon $R_{48} + \dots + R_{52}$; or (b) as the decomposable annulus $D_{13} \in \pi_2(x_3 x_{45}, x_5 x_{52})$ composed with the rectangle $R_9 + \dots + R_{24} \in \pi_2(x_5 x_{52}, x_{21}, x_{36})$. Since the moduli space of the composition has an even number of ends, and the composition in (a) provides an odd number of ends by the previous paragraph, it follows that
\[
\#\MM(D_{18}) + \#\MM(D_{13}) + 1  = 0 \pmod 2
\]
In other words, $D_{19}$ counts for the differential if and only if $D_{13}$ does not. Thus, if the complex structure is sufficiently stretched, then $\#\MM(D_{19}) = 1$ by Lemma \ref{lemma:decomposable}.
\end{proof}

\begin{figure}
\begin{center}
\labellist
 \pinlabel (a) at 10 214
 \pinlabel (b) at 192 214
 \pinlabel $x_3$ at 32 14
 \pinlabel $x_{23}$ at 32 210
 \pinlabel $x_{45}$ at 83 97
 \pinlabel $x_{46}$ at 83 123
 \pinlabel {{\color{red} $\alpha_1^L$}} at 123 79
 \pinlabel {{\color{darkred} $\alpha_2^L$}} at 120 36
 \pinlabel {{\color{blue} $\beta_1$}} at 27 147
 \pinlabel {{\color{darkgreen} $\alpha_2^R$}} at 160 137
 \pinlabel {{\color{darkblue} $\beta_2$}} at 79 144
 \pinlabel $\rho_3$ at 12 49
 \pinlabel $x_2$ at 237 20
 \pinlabel $x_{28}$ at 235 151
 \pinlabel $x_{37}$ at 316 62
 \pinlabel $x_{45}$ at 316 96
 \pinlabel {{\color{red} $\alpha_1^L$}} at 329 40
% \pinlabel {{\color{darkred} $\alpha_2^L$}} at 226 35
 \pinlabel {{\color{blue} $\beta_1$}} at 237 79
 \pinlabel {{\color{darkgreen} $\alpha_2^R$}} at 338 106
 \pinlabel {{\color{darkblue} $\beta_2$}} at 314 79
 \pinlabel $\rho_{23}$ at 217 173
\endlabellist
\includegraphics{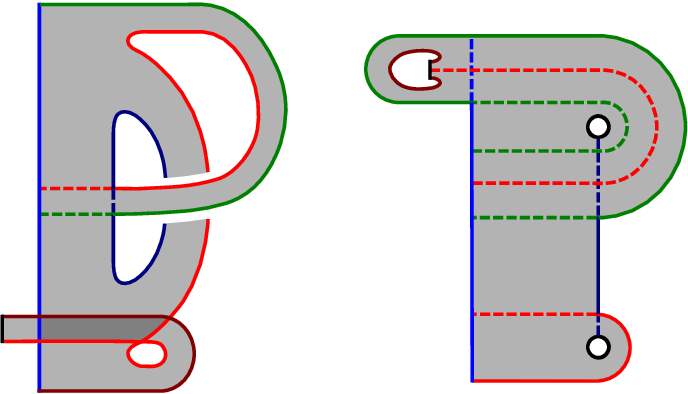}
\end{center}
\caption{The immersed genus-1 domains $D_{19}$ (a) and $D_{20}$ (b).} \label{fig:genus1immersed}
\end{figure}
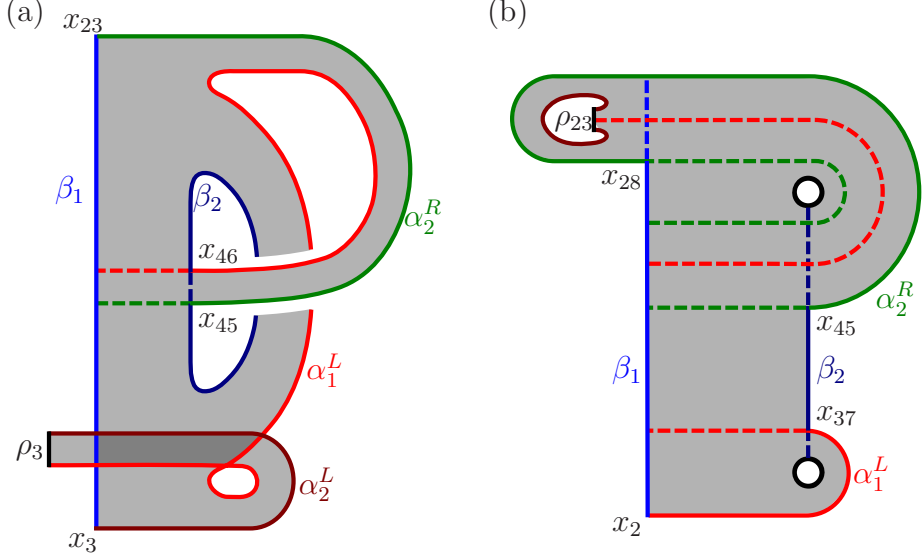

Finally, we consider the domains
\[
\begin{aligned}
D_{19} &= R_7 + 2R_8 + R_9 + R_{10} + \dots + R_{24} + R_{31} + R_{37}+ R_{48} \in \pi_2(x_3 x_{45}, x_{23} x_{46}) \\
D_{20} &= R_6 + 2R_7 + 2R_8 + R_9 + \dots + R_{31} + R_{36} + R_{37} + R_{48} + R_{49},
\end{aligned}
\]
which are represented by immersed genus-$1$ surfaces as shown in Figure \ref{fig:genus1immersed}. From the figure, we can see that $D_{19}$ is simply an immersed copy of a domain with the same geometry as $D_{17}$ (with a single Reeb chord), so it always counts for the differential. (Alternately, we may give an explicit proof by considering the composition of $D_{19}$ with the rectangle $R_{25} + \dots + R_{30}$ and arguing as above.) There are 8 domains of this type. On the other hand, $D_{20}$ does not have a cut that can go to the boundary, which would be needed to make it compatible with the sequence $(-\rho_2, -\rho_3)$, so it does not have any holomorphic representatives. There are two enlargements of $D_{20}$ with the same property. \\

Combining the results above, we have:
\begin{proposition} \label{prop:gooddomains}
If the complex structure on $\Sigma$ is sufficiently stretched, then $780$ of the $1,013$ positive, index-1 domains count for the differential on $\CFDD(\HH',0)$, as indicated in Table \ref{table:domains}.
\end{proposition}
Using this list, we may then record the differential on $\CFDD(\HH',0)$ as a
$245 \times 245$ matrix with entries in $\AA_\rho \otimes \AA_\sigma$, although
for obvious reasons we do not record this matrix here.

By counting the multiplicity at $w$ of each domain (whether or not the domain counts for the
differential), we can determine the relative Alexander gradings of all
of the generators using \eqref{eq:relalex}. We find that the generators of $\CFDD(\HH',0)$ all fall into
three consecutive gradings, which for now we arbitrarily declare to be $-1$,
$0$, and $1$. In the end, after we evaluate all tensor products, the symmetry
of $\CFK(S^3, D_{J,s}(K,t))$ will show that this was the correct choice. We do
not explicitly list all of the gradings here, however.

\subsection{Algebraic computation of \texorpdfstring{$\CFAA$}{CFAA}} \label{subsec:cfaa}

Using our \emph{Mathematica} package \verb"TorusAlgebra.nb", we may apply the edge cancellation algorithm explained in Subsection \ref{subsec:edge} to simplify $\CFDD(\HH',0)$, canceling only edges that preserve the filtration level. By abuse of notation, we denote the resulting bimodule by $\CFDD(\YY,B_3,0)$.

\begin{theorem} \label{thm:cfdd}
The type $DD$ structure $\CFDD(\YY,B_3,0)$ has a basis $\{y_1, \dots, y_{19}\}$
with the following properties:
\begin{enumerate}
\item The Alexander gradings of the basis elements are:
\[
A(y_i) = \begin{cases} -1 & i=1 \\ 0 & i = 2, \dots,10 \\ 1 & i=11, \dots, 19.
\end{cases}
\]
\item The associated idempotents in $\AA_\rho$ and $\AA_\sigma$ of the
generators are:
\[
\begin{array}{|c|c|c|} \hline
 & \iota_0^\rho & \iota_1^\rho \\ \hline
 \iota_0^\sigma & y_4, y_5, y_7, y_{11},y_{13}, y_{17}, y_{19} & y_8,y_{10},y_{14},y_{16} \\ \hline
 \iota_1^\sigma & y_3,y_6,y_{12},y_{18} & y_1,y_2,y_9,y_{15} \\ \hline
\end{array}
\]
\item The differential is given by \[\delta_1 (y_i) = \sum_{j=1}^{19} a_{ij}
\otimes y_j,\] where $(a_{ij})$ is the following matrix:
\begin{scriptsize}
\renewcommand\arraycolsep{3pt}
\[
\left(
\begin{array}{c|ccccccccc|ccccccccc}
 0 & 0 & 0 & 0 & 0 & 0 & 0 & 0 & 0 & 0 & 0 & 0 & 0 & 0 & 0 & 0 & 0 & 0 & 0 \\ \hline
 1 & 0 & 0 & 0 & 0 & 0 & 0 & \sigma_2 & 0 & 0 & 0 & 0 & 0 & 0 & 0 & 0 & 0 & 0 & 0 \\
 \rho_1 & 0 & 0 & 0 & 0 & 0 & 0 & 0 & 0 & 0 & 0 & 0 & 0 & 0 & 0 & 0 & 0 & 0 & 0 \\
 \rho_1 \sigma_{123} & \rho_1 \sigma_3 & \sigma_3 & 0 & 0 & 0 & 0 & 0 & 0 & 0 & 0 & 0 & 0 & 0 & 0 & 0 & 0 & 0 & 0 \\
 0 & \rho_3 \sigma_{123}+\rho_{123} \sigma_3 & 0 & 0 & 0 & \sigma_3 & 0 & 0 & \rho_3 \sigma_{123} & \rho_3 & 0 & 0 & 0 & 0 & 0 & 0 & 0 & 0 & 0 \\
 \rho_{123} & 0 & 0 & 0 & 0 & 0 & \sigma_2 & 0 & \rho_3 & 0 & 0 & 0 & 0 & 0 & 0 & 0 & 0 & 0 & 0 \\
 0 & 0 & 0 & 0 & 0 & 0 & 0 & \rho_3 & 0 & 0 & 0 & 0 & 0 & 0 & 0 & 0 & 0 & 0 & 0 \\
 0 & 0 & \rho_2 \sigma_1 & 0 & 0 & 0 & 0 & 0 & 0 & 0 & 0 & 0 & 0 & 0 & 0 & 0 & 0 & 0 & 0 \\
 0 & 0 & \rho_2 & 0 & 0 & 0 & 0 & \sigma_2 & 0 & 0 & 0 & 0 & 0 & 0 & 0 & 0 & 0 & 0 & 0 \\
 \sigma_{123} & 0 & 0 & \rho_2 & 0 & 0 & 0 & 0 & \sigma_3 & 0 & 0 & 0 & 0 & 0 & 0 & 0 & 0 & 0 & 0 \\ \hline
 0 & \rho_1 \sigma_1 & \sigma_1 & 0 & 0 & 0 & 0 & \rho_1 & \rho_1 \sigma_1 & 0 & 0 & 0 & 0 & 0 & 0 & 0 & 0 & 0 & 0 \\
 0 & \rho_1 & 1 & 0 & 0 & 0 & 0 & 0 & 0 & 0 & \sigma_2 & 0 & 0 & 0 & 0 & 0 & 0 & 0 & 0 \\
 0 & 0 & \sigma_{123} & 1 & 0 & 0 & 0 & 0 & 0 & 0 & 0 & \sigma_3 & 0 & 0 & 0 & 0 & 0 & 0 & 0 \\
 0 & 0 & 0 & 0 & 1 & 0 & 0 & 0 & 0 & 0 & 0 & 0 & 0 & 0 & \sigma_3 & 0 & 0 & 0 & \rho_3 \\
 0 & \rho_{123} & 0 & 0 & 0 & 1 & 0 & 0 & 0 & 0 & 0 & 0 & 0 & 0 & 0 & \sigma_2 & 0 & \rho_3 & 0 \\
 0 & 0 & 0 & 0 & 0 & 0 & 1 & \rho_{123} & 0 & 0 & 0 & 0 & 0 & 0 & 0 & 0 & \rho_3 & 0 & 0 \\
 0 & 0 & 0 & 0 & 0 & 0 & 0 & 1 & 0 & 0 & \rho_2 & 0 & 0 & 0 & 0 & 0 & 0 & 0 & 0 \\
 0 & 0 & 0 & 0 & 0 & 0 & 0 & 0 & 1 & 0 & 0 & \rho_2 & 0 & 0 & 0 & 0 & \sigma_2 & 0 & 0 \\
 0 & \sigma_{123} & 0 & 0 & 0 & 0 & 0 & 0 & \sigma_{123} & 1 & 0 & 0 & \rho_2 & 0 & 0 & 0 & 0 & \sigma_3 & 0
\end{array}
\right)
\]
\renewcommand\arraycolsep{5pt}
\end{scriptsize}
The block decomposition indicates the filtration levels.
\end{enumerate}
\end{theorem}

Finally, to compute $\CFAA(\YY,B_3,0)$, we use the $AA$ identity bimodule
described in Theorem \ref{thm:identityAA}:
\[
\CFAA(\YY,B_3,0) \simeq \CFAA(\I,0) \underset{\AA_\sigma}{\boxtimes}
(\CFAA(\I,0) \underset{\AA_\rho}{\boxtimes} \CFDD(\YY,B_3,0)).
\]
We evaluate this tensor product using our \emph{Mathematica} package. The
filtration on $\CFDD(\YY,B_3,0)$ induces a filtration on $\CFAA(\YY,B_3,0)$,
and we again use the edge cancellation procedure to reduce the number of
generators. For further details on the computation, see Appendix \ref{sec:appendix}.

\begin{theorem} \label{thm:cfaa}
The filtered $AA$-module $\CFAA(\YY,B_3,0)$ has a basis
\[
\{a_1, \dots, a_5, b_1, \dots, b_6, c_1, d_1, \dots, d_4, e_1, e_2, e_3\}
\]
with the following properties:
\begin{enumerate}
\item The Alexander gradings of the basis elements are:
\[
\begin{split}
A(c_1) &= -1 \\
A(a_i) = A(d_i) &= 0 \\
A(b_i) = A(e_i) &= 1. \\
\end{split}
\]
\item The associated idempotents in $\AA_\rho$ and $\AA_\sigma$ of the
generators are:
\[
\begin{array}{|c|c|c|} \hline
 & \iota_0^\rho & \iota_1^\rho \\ \hline
 \iota_0^\sigma & a_1,a_3,a_4, b_1,b_3,b_4,b_6  & d_1,d_3,e_1,e_3  \\ \hline
 \iota_1^\sigma & a_2,a_5, b_2,b_5     &  c_1,d_2,d_4,e_2  \\ \hline
\end{array}
\]
\item The $A_\infty$ multiplications are presented in the matrices that follow.
For $x,y \in \{a,b,c,d,e\}$, the entry in the $i\Th$ row and $j\Th$ column of
the matrix $M_{xy}$ records the multiplications taking $x_i$ to $y_j$, as
described in Subsection \ref{subsec:objects}. The matrices $M_{ab}$, $M_{cb}$,
$M_{cd}$, $M_{ce}$, $M_{db}$, and $M_{de}$ are necessarily zero because of the
Alexander grading.
\end{enumerate}
\end{theorem}

\begin{scriptsize}

\[
\text{{\normalsize $M_{aa}=$}} \left(
\begin{array}{ccccc}
 0 & \sigma_1 & \sigma_{12} & \rho_{12} & \sigma_{123} \rho_{12}+\sigma_1 \rho_3 \rho_2 \rho_{12} \\
 0 & 0 & \sigma_2 & 0 & \sigma_{23} \rho_{12}+\rho_{12} \\
 0 & 0 & 0 & 0 & \sigma_3 \rho_{12} \\
 0 & 0 & 0 & 0 & \sigma_1 \\
 0 & 0 & 0 & 0 & 0
\end{array}
\right)
\]

\[
\text{{\normalsize $M_{ac}=$}} \left(
\begin{array}{l}
 \begin{aligned}
  \sigma_{123} & \rho_{123}+\sigma_{123} \sigma_{23} \rho_{123}+\sigma_3 \sigma_2 \sigma_1 \rho_{123}+\sigma_1 \sigma_{23} \rho_3 \rho_2 \rho_{123} \\[-3pt]
  & +\sigma_{123} \sigma_{23} \rho_3 \rho_2 \rho_1+\sigma_3 \sigma_2 \sigma_1 \rho_3 \rho_2 \rho_1+\sigma_1 \sigma_{23} \rho_3 \rho_2 \rho_3 \rho_2 \rho_1 \\
 \end{aligned} \\
 \sigma_{23} \rho_{123}+\rho_3 \rho_2 \rho_1+\sigma_{23} \sigma_{23} \rho_{123}+\sigma_{23} \sigma_{23} \rho_3 \rho_2 \rho_1 \\
 \sigma_3 \sigma_{23} \rho_{123}+\sigma_3 \rho_3 \rho_2 \rho_1+\sigma_3 \sigma_{23} \rho_3 \rho_2 \rho_1 \\
 \sigma_{123} \rho_3+\sigma_3 \sigma_2 \sigma_1 \rho_3 \\
 \sigma_{23} \rho_3
\end{array}
\right)
\]

\[
\text{{\normalsize $M_{ad}=$}} \left(
\begin{array}{llll}
 \rho_1 & 0 &
  \begin{aligned}
   \sigma_{12} &\rho_{123}+ \sigma_{123} \sigma_2 \rho_{123}+\sigma_{12} \rho_3 \rho_2 \rho_1+\sigma_1 \sigma_2 \rho_3 \rho_2
   \rho_{123}\\[-3pt]
   &+\sigma_{123} \sigma_2 \rho_3 \rho_2 \rho_1+\sigma_1 \sigma_2 \rho_3 \rho_2 \rho_3 \rho_2 \rho_1
  \end{aligned}
  & \sigma_{123} \rho_1+\sigma_1 \rho_3 \rho_2 \rho_1 \\
 0 & \rho_1 & \sigma_2 \rho_{123}+\sigma_{23} \sigma_2 \rho_{123}+\sigma_2 \rho_3 \rho_2 \rho_1+\sigma_{23} \sigma_2 \rho_3 \rho_2 \rho_1 & \sigma_{23} \rho_1 \\
 0 & 0 & \sigma_3 \sigma_2 \rho_{123}+\sigma_3 \sigma_2 \rho_3 \rho_2 \rho_1+\rho_1 & \sigma_3 \rho_1 \\
 0 & 0 & \sigma_{12} \rho_3 & 0 \\
 0 & 0 & \sigma_2 \rho_3 & 0
\end{array}
\right)
\]

\[
\text{{\normalsize $M_{ba}=$}} \left(
\begin{array}{ccccc}
 1 & 0 & 0 & \sigma_3 \sigma_2 \rho_{123} \rho_2+\sigma_3 \sigma_2 \rho_3 \rho_2 \rho_1 \rho_2 & \sigma_{123} \rho_3 \rho_2 \rho_{12}+\sigma_1 \rho_3 \rho_2 \rho_3 \rho_2 \rho_{12} \\
 0 & 1 & 0 & 0 & \rho_3 \rho_2 \rho_{12}+\sigma_{23} \rho_3 \rho_2 \rho_{12} \\
 0 & 0 & 1 & 0 & \sigma_3 \rho_3 \rho_2 \rho_{12} \\
 0 & 0 & 0 & 1+\sigma_3 \sigma_2 \rho_3 \rho_2 & \sigma_{123} \\
 0 & 0 & 0 & 0 & 1+\sigma_{23} \\
 0 & 0 & 0 & 0 & \sigma_3
\end{array}
\right)
\]

\[
\text{{\normalsize $M_{bb}=$}} \left(
\begin{array}{cccccc}
 0 & \sigma_1 & \sigma_{12} & \rho_{12} & 0 & 0 \\
 0 & 0 & \sigma_2 & 0 & \rho_{12} & 0 \\
 0 & 0 & 0 & 0 & 0 & \rho_{12} \\
 0 & 0 & 0 & 0 & \sigma_1 & \sigma_{12} \\
 0 & 0 & 0 & 0 & 0 & \sigma_2 \\
 0 & 0 & 0 & 0 & 0 & 0
\end{array}
\right)
\]

\[
\text{{\normalsize $M_{bc}=$}} \left(
\begin{array}{l}
 \begin{aligned}
  \sigma_3 \sigma_2 & \sigma_{123} \rho_{123} \rho_{23} + \sigma_{123} \sigma_{23} \rho_3 \rho_2 \rho_{123}+\sigma_{123} \rho_3 \rho_2 \rho_3 \rho_2 \rho_1 + \sigma_1 \sigma_{23} \rho_3 \rho_2 \rho_3 \rho_2 \rho_{123} \\[-3pt]
  &+ \sigma_3 \sigma_2 \sigma_3 \sigma_2 \sigma_1 \rho_{123} \rho_{23} +\sigma_3 \sigma_2 \sigma_{123} \rho_3 \rho_2 \rho_1 \rho_{23} + \sigma_{123} \sigma_{23} \rho_3 \rho_2 \rho_3 \rho_2 \rho_1 \\[-3pt]
  &+ \sigma_1 \sigma_{23} \rho_3 \rho_2 \rho_3 \rho_2 \rho_3 \rho_2 \rho_1+\sigma_3 \sigma_2 \sigma_3 \sigma_2 \sigma_1 \rho_3 \rho_2 \rho_1 \rho_{23}
 \end{aligned} \\
 \sigma_{23} \rho_3 \rho_2 \rho_{123}+\rho_3 \rho_2 \rho_3 \rho_2 \rho_1+\sigma_{23} \sigma_{23} \rho_3 \rho_2 \rho_{123}+\sigma_{23} \sigma_{23} \rho_3 \rho_2 \rho_3 \rho_2 \rho_1 \\
 \sigma_3 \sigma_{23} \rho_3 \rho_2 \rho_{123}+\sigma_3 \rho_3 \rho_2 \rho_3 \rho_2 \rho_1+\sigma_3 \sigma_{23} \rho_3 \rho_2 \rho_3 \rho_2 \rho_1 \\
 \sigma_3 \sigma_2 \sigma_{123} \rho_3 \rho_{23}+\sigma_3 \sigma_2 \sigma_3 \sigma_2 \sigma_1 \rho_3 \rho_{23} \\
 0 \\
 0
\end{array}
\right)
\]

\renewcommand\arraycolsep{3pt}
\[
\text{{\normalsize $M_{bd}=$}} \left(
\begin{array}{llll}
 0 & \sigma_{123} \rho_{123} + \sigma_{123} \rho_3 \rho_2 \rho_1 &
   \begin{aligned}
     \sigma&_3 \sigma_2 \sigma_{12} \rho_{123} \rho_{23} + \sigma_{123} \sigma_2 \rho_3 \rho_2 \rho_{123}\\[-3pt]
     &+\sigma_1 \sigma_2 \rho_3 \rho_2 \rho_3 \rho_2 \rho_{123}+\sigma_3 \sigma_2 \sigma_{12} \rho_3 \rho_2 \rho_1 \rho_{23} \\[-3pt]
     &+\sigma_{123} \sigma_2 \rho_3 \rho_2 \rho_3 \rho_2 \rho_1 +\sigma_1 \sigma_2 \rho_3 \rho_2 \rho_3 \rho_2 \rho_3 \rho_2 \rho_1
   \end{aligned}
     &
   \begin{aligned}
     \sigma&_{123} \rho_{123} +\sigma_3 \sigma_2 \sigma_1 \rho_{123} \\[-3pt]
     & +\sigma_1 \rho_3 \rho_2 \rho_3 \rho_2 \rho_1\\[-3pt]
     & +\sigma_3 \sigma_2 \sigma_1 \rho_3 \rho_2 \rho_1
   \end{aligned} \\
 0 & \sigma_{23} \rho_{123} + \sigma_{23} \rho_3 \rho_2 \rho_1 &
  \begin{aligned}
    \sigma&_2 \rho_3 \rho_2 \rho_{123}+\sigma_{23} \sigma_2 \rho_3 \rho_2 \rho_{123}+\sigma_2 \rho_3 \rho_2 \rho_3 \rho_2 \rho_1 \\[-3pt]
    & +\sigma_{23} \sigma_2 \rho_3 \rho_2 \rho_3 \rho_2 \rho_1
  \end{aligned}
  & \sigma_{23} \rho_{123}+\rho_3 \rho_2 \rho_1 \\
 0 & \sigma_3 \rho_{123}+\sigma_3 \rho_3 \rho_2 \rho_1 &
  \begin{aligned}
    \rho&_3 \rho_2 \rho_1+\sigma_3 \sigma_2 \rho_3 \rho_2 \rho_{123} \\[-3pt]
    & +\sigma_3 \sigma_2 \rho_3 \rho_2 \rho_3 \rho_2 \rho_1+\rho_{123}
  \end{aligned}
  & \sigma_3 \rho_{123} \\
 0 & \sigma_{123} \rho_3 & \sigma_3 \sigma_2 \sigma_{12} \rho_3 \rho_{23} & \sigma_{123} \rho_3+\sigma_3 \sigma_2 \sigma_1 \rho_3 \\
 0 & \sigma_{23} \rho_3 & 0 & \sigma_{23} \rho_3 \\
 0 & \sigma_3 \rho_3 & \rho_3 & \sigma_3 \rho_3
\end{array}
\right)
\]
\renewcommand\arraycolsep{5pt}

\[
\text{{\normalsize $M_{be}=$}} \left(
\begin{array}{ccc}
 \rho_1 & 0 & 0 \\
 0 & \rho_1 & 0 \\
 0 & 0 & \rho_1 \\
 0 & 0 & 0 \\
 0 & 0 & 0 \\
 0 & 0 & 0
\end{array}
\right)
\]

\[
\text{{\normalsize $M_{cc}= (0)$}}
\]

\[
\text{{\normalsize $M_{da}=$}} \left(
\begin{array}{ccccc}
 0 & 0 & 0 & \rho_2 & 0 \\
 0 & 0 & 0 & 0 & \rho_2 \\
 0 & 0 & 0 & 0 & 0 \\
 0 & 0 & 0 & 0 & \rho_2
\end{array}
\right)
\]

\[
\text{{\normalsize $M_{dc}=$}} \left(
\begin{array}{c}
 \sigma_{123} \rho_{23}+\sigma_3 \sigma_2 \sigma_1 \rho_{23}+\sigma_{123} \\
 \sigma_{23} \rho_{23}+\sigma_{23} \\
 \sigma_3 \\
 1+\sigma_{23} \rho_{23}
\end{array}
\right)
\]

\[
\text{{\normalsize $M_{dd}=$}} \left(
\begin{array}{cccc}
 0 & \sigma_1 & \sigma_{12} \rho_{23}+\sigma_{12} & 0 \\
 0 & 0 & \sigma_2 \rho_{23}+\sigma_2 & 0 \\
 0 & 0 & 0 & 0 \\
 0 & 0 & \sigma_2 \rho_{23} & 0
\end{array}
\right)
\]

\[
\text{{\normalsize $M_{ea}=$}} \left(
\begin{array}{ccccc}
 0 & 0 & 0 & \sigma_3 \sigma_2 \rho_{23} \rho_2 & 0 \\
 0 & 0 & 0 & 0 & 0 \\
 0 & 0 & 0 & 0 & 0
\end{array}
\right)
\]

\[
\text{{\normalsize $M_{eb}=$}} \left(
\begin{array}{cccccc}
 0 & 0 & 0 & \rho_2 & 0 & 0 \\
 0 & 0 & 0 & 0 & \rho_2 & 0 \\
 0 & 0 & 0 & 0 & 0 & \rho_2
\end{array}
\right)
\]

\[
\text{{\normalsize $M_{ec}=$}} \left(
\begin{array}{c}
 \sigma_3 \sigma_2 \sigma_{123} \rho_{23} \rho_{23} \\
 0 \\
 0
\end{array}
\right)
\]

\[
\text{{\normalsize $M_{ed}=$}} \left(
\begin{array}{cccc}
 1 & \sigma_{123} \rho_{23} & \sigma_3 \sigma_2 \sigma_{12} \rho_{23} \rho_{23} & \sigma_{123} \rho_{23}+\sigma_3 \sigma_2 \sigma_1 \rho_{23}+\sigma_{123} \\
 0 & 1+\sigma_{23} \rho_{23} & 0 & \sigma_{23} \rho_{23}+\sigma_{23} \\
 0 & \sigma_3 \rho_{23} & 1+\rho_{23} & \sigma_3 \rho_{23}+\sigma_3
\end{array}
\right)
\]

\[
\text{{\normalsize $M_{ee}=$}} \left(
\begin{array}{ccc}
 0 & \sigma_1 & \sigma_{12} \\
 0 & 0 & \sigma_2 \\
 0 & 0 & 0
\end{array}
\right)
\]

\end{scriptsize}

Because we are ultimately interested in the tensor product of
$\CFAA(\YY,B_3,0)$ with $\CFD(\XX_J^s)$ and $\CFD(\XX_K^t)$, we may disregard
any higher multiplication that uses sequences of algebra elements that cannot
occur in these type $D$ structures. Specifically, by Proposition
\ref{prop:norho1rho2}, we may disregard any sequence containing $\rho_2
\rho_3$, $\rho_1 \rho_2$, $\rho_1 \rho_{23}$, $\sigma_2 \sigma_3$, $\sigma_1
\sigma_2$, or $\sigma_1 \sigma_{23}$. Accordingly, for the discussion that
follows, we may replace $M_{ac}$, $M_{ad}$, $M_{ba}$, $M_{bc}$, and $M_{bd}$
with the following:

\begin{scriptsize}

\[
\text{{\normalsize $M'_{ac}=$}} \left(
\begin{array}{c}
 \sigma_{123} \rho_{123}+\sigma_{123} \sigma_{23} \rho_{123}+\sigma_3 \sigma_2 \sigma_1 \rho_{123}+ \sigma_{123} \sigma_{23} \rho_3 \rho_2 \rho_1+\sigma_3 \sigma_2 \sigma_1 \rho_3 \rho_2 \rho_1  \\
 \sigma_{23} \rho_{123}+\rho_3 \rho_2 \rho_1+\sigma_{23} \sigma_{23} \rho_{123}+\sigma_{23} \sigma_{23} \rho_3 \rho_2 \rho_1 \\
 \sigma_3 \sigma_{23} \rho_{123}+\sigma_3 \rho_3 \rho_2 \rho_1+\sigma_3 \sigma_{23} \rho_3 \rho_2 \rho_1 \\
 \sigma_{123} \rho_3+\sigma_3 \sigma_2 \sigma_1 \rho_3 \\
 \sigma_{23} \rho_3
\end{array}
\right)
\]

\[
\text{{\normalsize $M'_{ad}=$}} \left(
\begin{array}{cccc}
 \rho_1 & 0 & \sigma_{12} \rho_{123}+\sigma_{123} \sigma_2 \rho_{123}+\sigma_{12} \rho_3 \rho_2 \rho_1+ \sigma_{123} \sigma_2 \rho_3 \rho_2 \rho_1 & \sigma_{123} \rho_1+\sigma_1 \rho_3 \rho_2 \rho_1 \\
 0 & \rho_1 & \sigma_2 \rho_{123}+\sigma_{23} \sigma_2 \rho_{123}+\sigma_2 \rho_3 \rho_2 \rho_1+\sigma_{23} \sigma_2 \rho_3 \rho_2 \rho_1 & \sigma_{23} \rho_1 \\
 0 & 0 & \sigma_3 \sigma_2 \rho_{123}+\sigma_3 \sigma_2 \rho_3 \rho_2 \rho_1+\rho_1 & \sigma_3 \rho_1 \\
 0 & 0 & \sigma_{12} \rho_3 & 0 \\
 0 & 0 & \sigma_2 \rho_3 & 0
\end{array}
\right)
\]

\[
\text{{\normalsize $M'_{ba}=$}} \left(
\begin{array}{ccccc}
 1 & 0 & 0 & \sigma_3 \sigma_2 \rho_{123} \rho_2 & \sigma_{123} \rho_3 \rho_2 \rho_{12} \\
 0 & 1 & 0 & 0 & \rho_3 \rho_2 \rho_{12}+\sigma_{23} \rho_3 \rho_2 \rho_{12} \\
 0 & 0 & 1 & 0 & \sigma_3 \rho_3 \rho_2 \rho_{12} \\
 0 & 0 & 0 & 1+\sigma_3 \sigma_2 \rho_3 \rho_2 & \sigma_{123} \\
 0 & 0 & 0 & 0 & 1+\sigma_{23} \\
 0 & 0 & 0 & 0 & \sigma_3
\end{array}
\right)
\]

\[
\text{{\normalsize $M'_{bc}=$}} \left(
\begin{array}{c}
 \sigma_3 \sigma_2 \sigma_{123} \rho_{123} \rho_{23} + \sigma_{123} \sigma_{23} \rho_3 \rho_2 \rho_{123}  \\
 \sigma_{23} \rho_3 \rho_2 \rho_{123}+ \sigma_{23} \sigma_{23} \rho_3 \rho_2 \rho_{123}  \\
 \sigma_3 \sigma_{23} \rho_3 \rho_2 \rho_{123}  \\
 \sigma_3 \sigma_2 \sigma_{123} \rho_3 \rho_{23}  \\
 0 \\
 0
\end{array}
\right)
\]

\[
\text{{\normalsize $M'_{bd}=$}} \left(
\begin{array}{llll}
 0 & \sigma_{123} \rho_{123}+\sigma_{123} \rho_3 \rho_2 \rho_1 &
     \sigma_3 \sigma_2 \sigma_{12} \rho_{123} \rho_{23} + \sigma_{123} \sigma_2 \rho_3 \rho_2 \rho_{123} &
   \begin{aligned}
     \sigma_{123} & \rho_{123} +\sigma_3 \sigma_2 \sigma_1 \rho_{123} \\[-3pt]
     & +\sigma_3 \sigma_2 \sigma_1 \rho_3 \rho_2 \rho_1
   \end{aligned} \\
 0 & \sigma_{23} \rho_{123}+\sigma_{23} \rho_3 \rho_2 \rho_1 & \sigma_2 \rho_3 \rho_2 \rho_{123}+\sigma_{23} \sigma_2 \rho_3 \rho_2 \rho_{123} & \sigma_{23} \rho_{123}+\rho_3 \rho_2 \rho_1 \\
 0 & \sigma_3 \rho_{123}+\sigma_3 \rho_3 \rho_2 \rho_1 & \rho_3 \rho_2 \rho_1+\sigma_3 \sigma_2 \rho_3 \rho_2 \rho_{123} + \rho_{123} & \sigma_3 \rho_{123} \\
 0 & \sigma_{123} \rho_3 & \sigma_3 \sigma_2 \sigma_{12} \rho_3 \rho_{23} & \sigma_{123} \rho_3+\sigma_3 \sigma_2 \sigma_1 \rho_3 \\
 0 & \sigma_{23} \rho_3 & 0 & \sigma_{23} \rho_3 \\
 0 & \sigma_3 \rho_3 & \rho_3 & \sigma_3 \rho_3
\end{array}
\right)
\]

\end{scriptsize}

\section{Evaluation of the tensor product} \label{sec:tensorproduct}
Using the computation of $\CFAA(\YY,B_3,0)$ given in the previous section, we
may now compute the double tensor product
\[
(\CFAA(\YY,B_3,0) \underset{\AA_\rho}{\boxtimes} \CFD(\XX_J^s))
\underset{\AA_\sigma}{\boxtimes} \CFD(\XX_K^t).
\]
In what follows, we evaluate the tensor product over $\AA_\rho$ and simplify
the resulting filtered type $A$ module before evaluating the tensor product
over $\AA_\sigma$. Then we use the edge cancellation algorithm to compute
$\tau(D_{J,s}(K,t))$. As a reminder, we restate the main theorem:
\[
\tau(D_{J,s}(K,t)) =
\begin{cases}
1 & s<2\tau(J) \text{ and } t<2 \tau(K) \\
-1 & s>2\tau(J) \text{ and } t>2 \tau(K) \\
0 & \text{otherwise}.
\end{cases}
\]

Notice that it suffices to consider only the cases where $s \le 2\tau(J)$,
since if $s > 2 \tau(J)$, the behavior of $\tau$ under mirroring implies:
\[
\begin{split}
\tau(D_{J,s}(K,t)) &= -\tau\left( \overline{D_{J,s}(K,t)} \right) \\
&= -\tau( D_{\bar J,-s}(\bar K, -t) ) \\
&= \begin{cases} -1 & -t < 2 \tau(\bar K) \\ 0 & -t \ge 2\tau(\bar K)
\end{cases}
\\
&= \begin{cases} -1 & t > -2 \tau(K) \\ 0 & t \le -2\tau(K) \end{cases}
\end{split}
\]

With only slightly more bookkeeping, we could also write down a formula for the
knot Floer homology groups $\HFK(D_{J,s}(K,t))$, but since we are primarily
interested in the value of $\tau$ and its applications to knot and link
concordance, we do not bother to do that here.
%\footnote{A simplified
%computation in which the knot $J$ is assumed to be the unknot can be found in a
%previous version of \cite{LevineBingWhitehead}, available online at
%\arxiv{0912.5222v1}.}

\subsection{Tensor product over \texorpdfstring{$\AA_\rho$}{A\textunderscore\textrho}}

Let $\VV$ denote the bordered solid torus obtained by gluing together $\YY$ and
$\XX_J^s$, and let $D_{J,s}$ denote the image of the knot $B_3$ in the union.
By the gluing theorem, $\CFA(\VV,D_{J,s}) \simeq \CFAA(\YY,B_3,0)
\underset{\AA_\rho}{\boxtimes} \CFD(\XX_J^s)$. We shall describe this tensor
product as a direct sum of subspaces corresponding to the stable and unstable
chains in $\CFD(\XX_J^s)$. This decomposition will not be a direct sum of
$\AA_\infty$ modules, but we will be able to keep track of the few
multiplications that do not respect the decomposition, and ultimately they will
not affect the computation of $\tau(D_{J,s}(K,t))$.

The generators of $\iota_1\CFD(\XX_J^s)$ all lie in the interiors of the
chains, so the corresponding generators of the tensor product can be grouped in
a natural way, but it is not obvious \emph{a priori} how to divide up the
generators coming from $\iota_0 \CFD(\XX_J^s)$. Consider the two specified
bases for $\iota_0 \CFD(\XX_J^s)$: $\{\eta_0, \dots, \eta_{2n}\}$ and $\{\xi_0,
\dots, \xi_{2n}\}$. Depending on the structure of the unstable chain, the
generators $\xi_i$ have outgoing arrows labeled $\rho_1$, $\rho_{12}$, or
$\rho_{123}$, while the $\eta_i$ have outgoing arrows labeled $\rho_3$ and
incoming arrows labeled $\rho_2$ or $\rho_{12}$. Accordingly, we should try to
pair the generators of $\CFAA(\YY,B_3,0)\iota_0$ with the $\xi_i$ or $\eta_i$
depending on whether they have outgoing $\rho_1$, $\rho_{12}$, and $\rho_{123}$ arrows or outgoing $\rho_3$s and incoming $\rho_2$s. If we consider only
the $\AA_\infty$ maps in $\CFAA(\YY,B_3,0)$ that use a single element of
$\AA_\rho$, we notice that each of the generators $a_1, \dots, a_5$ and $b_1,
\dots, b_6$ satisfies exactly one of these conditions. Specifically, define the
following subspaces of $\CFAA(\YY,B_3,0) \underset{\AA_\rho}{\boxtimes}
\CFD(\XX_J^s)$:
%Modified 3/4/12

\begin{equation} \label{eq:Pbases}
\begin{split}
P_{\ver}^j &= \gen{a_1,a_2,a_3,b_1,b_2,b_3} \boxtimes \gen{
\xi_{2j-1}, \xi_{2j}} \\
& \quad + \gen{c_1,d_1,d_2,d_3,d_4,e_1,e_2,e_3} \boxtimes
\gen{\kappa_i^j \mid 1 \le i \le k_j} \\
P_{\hor}^j &= \gen{a_4,a_5,b_4,b_5,b_6} \boxtimes \gen{
\eta_{2j-1}, \eta_{2j} } \\
& \quad + \gen{c_1,d_1,d_2,d_3,d_4,e_1,e_2,e_3} \boxtimes
\gen{\lambda_i^j \mid 1 \le i \le l_j} \\
P_{\unst} &= \gen{a_1,a_2,a_3,b_1,b_2,b_3} \boxtimes \gen{\xi_0} \\
& \quad + \gen{a_4,a_5,b_4,b_5,b_6} \boxtimes \gen{\eta_0} \\
& \quad + \gen{c_1,d_1,d_2,d_3,d_4,e_1,e_2,e_3} \boxtimes \gen{\lambda_i \mid 1
\le i \le r}.
\end{split}
\end{equation}
We thus obtain a direct sum decomposition of $\CFAA(\YY,B_3,0)
\underset{\AA_\rho}{\boxtimes} \CFD(\XX_J^s)$ as a vector space:
\begin{equation} \label{eq:Psplit}
\CFAA(\YY,B_3,0) \underset{\AA_\rho}{\boxtimes} \CFD(\XX_J^s) =
\bigoplus_{j=1}^n P_{\ver}^j \oplus \bigoplus_{j=1}^n P_{\hor}^j \oplus
P_{\unst}.
\end{equation}
By inspecting the matrices $M_{xy}$, we see that any $\AA_\infty$
multiplication on the tensor product that comes from a multiplication in
$\CFAA(\YY,B_3,0)$ that uses at most one element of $\AA_\rho$ preserves this
decomposition. These multiplications are illustrated in Figures \ref{fig:Pver}
through \ref{fig:Punst,s=2tau}. In these and subsequent figures, the dashed
arrows represent repeated sections. For instance, the dashed arrow from
$e_1\kappa^j_1$ to $d_2 \kappa^j_{k_j}$ in Figure \ref{fig:Pver} means that
there are multiplications $e_1 \kappa^j_i \xrightarrow{\sigma_{123}} d_2
\kappa^j_{i+1}$ for each $i=1, \dots, k_j-1$. The Alexander filtration is
indicated by horizontal position, increasing from left to right.

\begin{figure} \centering
\begin{tiny}
\[
\xy (50,0)*{a_1 \xi_{2j}}="a1xi2j"; (100,0)*{b_1 \xi_{2j}}="b1xi2j";
(50,-12.5)*{a_2 \xi_{2j}}="a2xi2j"; (100,-12.5)*{b_2 \xi_{2j}}="b2xi2j";
(50,-25)*{a_3 \xi_{2j}}="a3xi2j"; (100,-25)*{b_3 \xi_{2j}}="b3xi2j";
(50,-50)*{d_1 \kappa^j_1}="d1ka1"; (100,-50)*{e_1 \kappa^j_1}="e1ka1";
(50,-62.5)*{d_2 \kappa^j_1}="d2ka1"; (100,-62.5)*{e_2 \kappa^j_1}="e2ka1";
(50,-75)*{d_3 \kappa^j_1}="d3ka1"; (100,-75)*{e_3 \kappa^j_1}="e3ka1";
(0,-87.5)*{c_1 \kappa^j_1}="c1ka1";  (50,-87.5)*{d_4 \kappa^j_1}="d4ka1";
 (50,-112.5)*{d_1 \kappa^j_{k_j}}="d1kak";
 (100,-112.5)*{e_1\kappa^j_{k_j}}="e1kak";
 (50,-125)*{d_2 \kappa^j_{k_j}}="d2kak";
 (100,-125)*{e_2 \kappa^j_{k_j}}="e2kak";
 (50,-137.5)*{d_3 \kappa^j_{k_j}}="d3kak";
 (100,-137.5)*{e_3 \kappa^j_{k_j}}="e3kak";
 (0,-150)*{c_1 \kappa^j_{k_j}}="c1kak";
 (50,-150)*{d_4 \kappa^j_{k_j}}="d4kak";
 (50,-175)*{a_1 \xi_{2j-1}}="a1xi2j1"; (100,-175)*{b_1 \xi_{2j-1}}="b1xi2j1";
 (50,-187.5)*{a_2 \xi_{2j-1}}="a2xi2j1"; (100,-187.5)*{b_2 \xi_{2j-1}}="b2xi2j1";
 (50,-200)*{a_3 \xi_{2j-1}}="a3xi2j1"; (100,-200)*{b_3 \xi_{2j-1}}="b3xi2j1";
%internal to xi_{2j}
 {\ar|{\sigma_1} "a1xi2j";"a2xi2j"};
 {\ar|{\sigma_1} "b1xi2j";"b2xi2j"};
 {\ar|{\sigma_2} "a2xi2j";"a3xi2j"};
 {\ar|{\sigma_2} "b2xi2j";"b3xi2j"};
 {\ar "b1xi2j";"a1xi2j"};
 {\ar "b2xi2j";"a2xi2j"};
 {\ar "b3xi2j";"a3xi2j"};
 {\ar@/^1.5pc/|(.35){\sigma_{12}} "a1xi2j";"a3xi2j"};
 {\ar@/^1.5pc/|{\sigma_{12}} "b1xi2j";"b3xi2j"};
%xi_{2j} -rho_{123}-> kappa_1
 {\ar@/_2.5pc/|(0.45){\threeline{\sigma_{123}+}{\sigma_{123}\sigma_{23}+}{\sigma_3\sigma_2\sigma_1}} "a1xi2j";"c1ka1"};
 {\ar@/_2pc/|(0.5){\twoline{\sigma_{23}+}{\sigma_{23}\sigma_{23}}} "a2xi2j";"c1ka1"};
 {\ar@/_1.5pc/|{\sigma_{3}\sigma_{23}} "a3xi2j";"c1ka1"};
 {\ar@/_3.1pc/|(0.37){\twoline{\sigma_{12}+}{\sigma_{123}\sigma_2}} "a1xi2j";"d3ka1"};
 {\ar@/_2.25pc/|(0.41){\twoline{\sigma_{2}+}{\sigma_{23}\sigma_2}} "a2xi2j";"d3ka1"};
 {\ar@/_1.5pc/|(0.45){\sigma_{3}\sigma_{2}} "a3xi2j";"d3ka1"};
 {\ar@/_1.5pc/|(0.53){\sigma_{123}} "b1xi2j";"d2ka1"};
 {\ar@/_1.5pc/|{\sigma_{23}} "b2xi2j";"d2ka1"};
 {\ar@/_1.5pc/|(0.3){\sigma_{3}} "b3xi2j";"d2ka1"};
 {\ar@/_1pc/|(0.55){\twoline{\sigma_{123}+}{\sigma_3\sigma_2\sigma_1}} "b1xi2j";"d4ka1"};
 {\ar@/_0.5pc/|(0.4){\sigma_{23}} "b2xi2j";"d4ka1"};
 {\ar|(0.3){\sigma_{3}} "b3xi2j";"d4ka1"};
 {\ar@/_0.3pc/ "b3xi2j";"d3ka1"};
%internal to kappa_1
 {\ar|{\sigma_1} "d1ka1";"d2ka1"};
 {\ar|{\sigma_1} "e1ka1";"e2ka1"};
 {\ar|{\sigma_2} "d2ka1";"d3ka1"};
 {\ar|{\sigma_2} "e2ka1";"e3ka1"};
 {\ar@/^1.5pc/|(0.4){\sigma_{12}} "d1ka1";"d3ka1"};
 {\ar@/^1.5pc/|(0.35){\sigma_{12}} "e1ka1";"e3ka1"};
 {\ar "e1ka1";"d1ka1"};
 {\ar "e2ka1";"d2ka1"};
 {\ar "e3ka1";"d3ka1"};
 {\ar "d4ka1";"c1ka1"};
 {\ar|{\sigma_{123}} "d1ka1";"c1ka1"};
 {\ar|(0.55){\sigma_{23}} "d2ka1";"c1ka1"};
 {\ar|(0.6){\sigma_{3}} "d3ka1";"c1ka1"};
 {\ar|{\sigma_{123}} "e1ka1";"d4ka1"};
 {\ar|(0.35){\sigma_{23}} "e2ka1";"d4ka1"};
 {\ar|{\sigma_{3}} "e3ka1";"d4ka1"};
%internal to kappa_k
 {\ar|{\sigma_1} "d1kak";"d2kak"};
 {\ar|{\sigma_1} "e1kak";"e2kak"};
 {\ar|{\sigma_2} "d2kak";"d3kak"};
 {\ar|{\sigma_2} "e2kak";"e3kak"};
 {\ar@/^1.5pc/|(0.4){\sigma_{12}} "d1kak";"d3kak"};
 {\ar@/_1.5pc/|(0.35){\sigma_{12}} "e1kak";"e3kak"};
 {\ar "e1kak";"d1kak"};
 {\ar "e2kak";"d2kak"};
 {\ar "e3kak";"d3kak"};
 {\ar "d4kak";"c1kak"};
 {\ar|{\sigma_{123}} "d1kak";"c1kak"};
 {\ar|{\sigma_{23}} "d2kak";"c1kak"};
 {\ar|{\sigma_{3}} "d3kak";"c1kak"};
 {\ar|{\sigma_{123}} "e1kak";"d4kak"};
 {\ar|(0.35){\sigma_{23}} "e2kak";"d4kak"};
 {\ar|{\sigma_{3}} "e3kak";"d4kak"};
%kappa_1 - - rho_{23} - -> kappa_k
 {\ar@{-->}@/_2.5pc/|(0.5){\twoline{\sigma_{123}+}{\sigma_3\sigma_2\sigma_1}} "d1ka1";"c1kak"};
 {\ar@{-->}@/_2pc/|(0.55){\sigma_{23}} "d2ka1";"c1kak"};
 {\ar@{-->}|{\sigma_{23}} "d4ka1";"c1kak"};
 {\ar@{-->}@/_3pc/|(0.35){\sigma_{12}} "d1ka1";"d3kak"};
 {\ar@{-->}@/_2.25pc/|(0.39){\sigma_{2}} "d2ka1";"d3kak"};
 {\ar@{-->}@/_1.5pc/|(0.45){\sigma_{2}} "d4ka1";"d3kak"};
 {\ar@{-->}@/_1.5pc/|(0.53){\sigma_{123}} "e1ka1";"d2kak"};
 {\ar@{-->}@/_1.5pc/|{\sigma_{23}} "e2ka1";"d2kak"};
 {\ar@{-->}@/_1.5pc/|(0.3){\sigma_{3}} "e3ka1";"d2kak"};
 {\ar@{-->}@/_1pc/|(0.55){\twoline{\sigma_{123}+}{\sigma_3\sigma_2\sigma_1}} "e1ka1";"d4kak"};
 {\ar@{-->}@/_0.5pc/|(0.4){\sigma_{23}} "e2ka1";"d4kak"};
 {\ar@{-->}|(0.3){\sigma_{3}} "e3ka1";"d4kak"};
 {\ar@{-->}@/_0.3pc/ "e3ka1";"d3kak"};
%internal to xi_{2j-1}
 {\ar|{\sigma_1} "a1xi2j1";"a2xi2j1"};
 {\ar|{\sigma_1} "b1xi2j1";"b2xi2j1"};
 {\ar|{\sigma_2} "a2xi2j1";"a3xi2j1"};
 {\ar|{\sigma_2} "b2xi2j1";"b3xi2j1"};
 {\ar "b1xi2j1";"a1xi2j1"};
 {\ar "b2xi2j1";"a2xi2j1"};
 {\ar "b3xi2j1";"a3xi2j1"};
 {\ar@/^1.5pc/|(0.35){\sigma_{12}} "a1xi2j1";"a3xi2j1"};
 {\ar@/_1.5pc/|(0.35){\sigma_{12}} "b1xi2j1";"b3xi2j1"};
%xi_{2j-1} -rho_1-> kappa_k
 {\ar|{\sigma_{123}} "a1xi2j1";"d4kak"};
 {\ar@/_1.5pc/|(0.55){\sigma_{23}} "a2xi2j1";"d4kak"};
 {\ar@/_3pc/|(0.6){\sigma_{3}} "a3xi2j1";"d4kak"};
 {\ar@/^2pc/ "a1xi2j1";"d1kak"};
 {\ar@/^2pc/ "a2xi2j1";"d2kak"};
 {\ar@/^2pc/ "a3xi2j1";"d3kak"};
 {\ar@/_2pc/ "b1xi2j1";"e1kak"};
 {\ar@/_2pc/ "b2xi2j1";"e2kak"};
 {\ar@/_2pc/ "b3xi2j1";"e3kak"};
\endxy
\]
\end{tiny}
\caption{The subspace $P_{\ver}^j$, corresponding to a vertical stable chain
$\xi_{2j} \xrightarrow{\rho_{123}} \kappa^j_1 \xrightarrow{\rho_{23}} \cdots
\xrightarrow{\rho_{23}} \kappa^j_{k_j} \xleftarrow{\rho_1} \xi_{2j-1}$.}
\label{fig:Pver}
\end{figure}

\begin{figure}
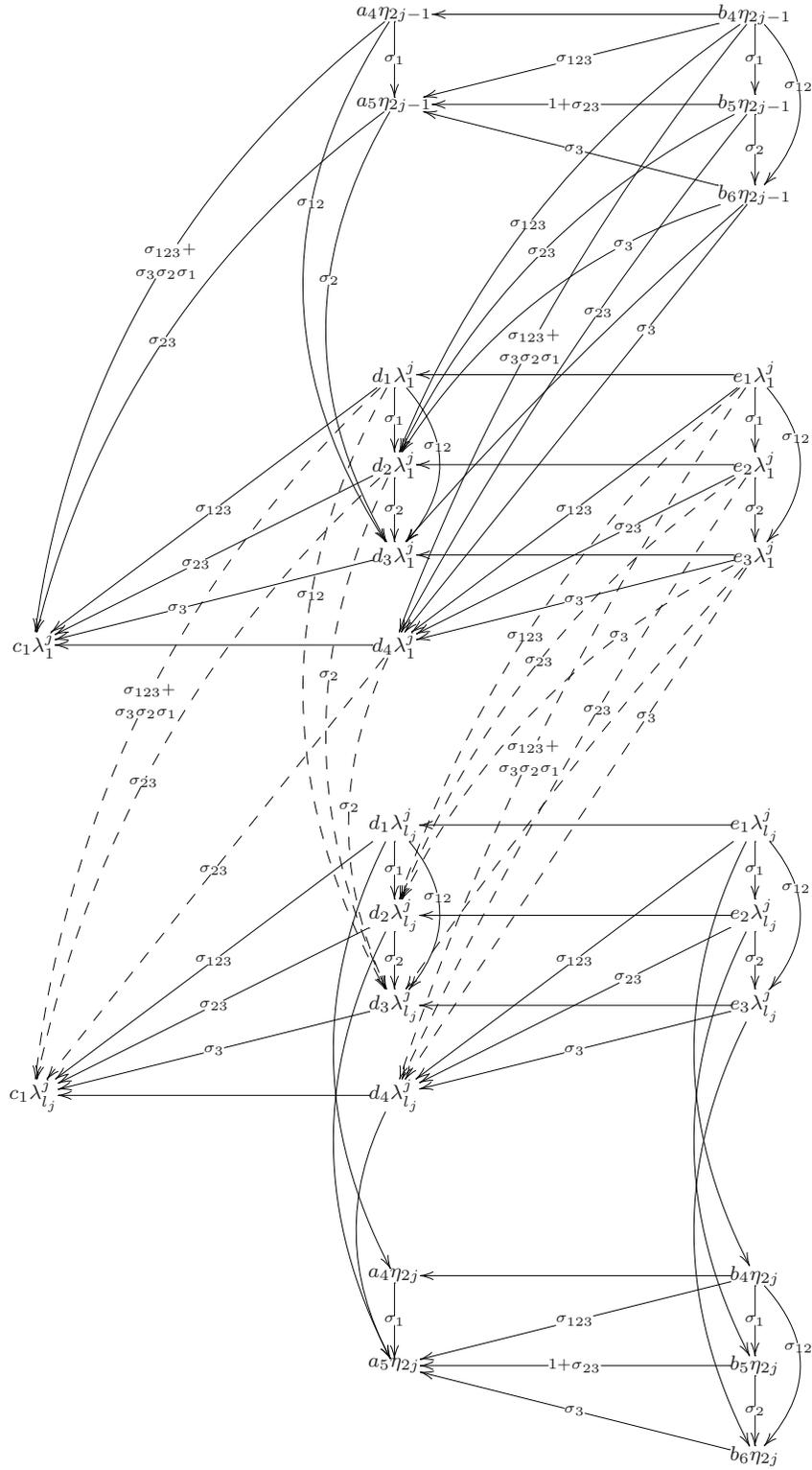
 \centering
\begin{tiny}
\[
\xy
 (50,0)*{a_4 \eta_{2j-1}}="a4eta2j1"; (100,0)*{b_4 \eta_{2j-1}}="b4eta2j1";
 (50,-12.5)*{a_5 \eta_{2j-1}}="a5eta2j1"; (100,-12.5)*{b_5 \eta_{2j-1}}="b5eta2j1";
 (100,-25)*{b_6 \eta_{2j-1}}="b6eta2j1";
 (50,-50)*{d_1 \lambda^j_1}="d1la1";
 (100,-50)*{e_1 \lambda^j_1}="e1la1";
 (50,-62.5)*{d_2 \lambda^j_1}="d2la1";
 (100,-62.5)*{e_2 \lambda^j_1}="e2la1";
 (50,-75)*{d_3 \lambda^j_1}="d3la1";
 (100,-75)*{e_3 \lambda^j_1}="e3la1";
 (0,-87.5)*{c_1 \lambda^j_1}="c1la1";
 (50,-87.5)*{d_4 \lambda^j_1}="d4la1";
 (50,-112.5)*{d_1 \lambda^j_{l_j}}="d1lal";
 (100,-112.5)*{e_1\lambda^j_{l_j}}="e1lal";
 (50,-125)*{d_2 \lambda^j_{l_j}}="d2lal";
 (100,-125)*{e_2 \lambda^j_{l_j}}="e2lal";
 (50,-137.5)*{d_3 \lambda^j_{l_j}}="d3lal";
 (100,-137.5)*{e_3 \lambda^j_{l_j}}="e3lal";
 (0,-150)*{c_1 \lambda^j_{l_j}}="c1lal";
 (50,-150)*{d_4 \lambda^j_{l_j}}="d4lal";
 (50,-175)*{a_4 \eta_{2j}}="a4eta2j"; (100,-175)*{b_4 \eta_{2j}}="b4eta2j";
 (50,-187.5)*{a_5 \eta_{2j}}="a5eta2j"; (100,-187.5)*{b_5 \eta_{2j}}="b5eta2j";
 (100,-200)*{b_6 \eta_{2j}}="b6eta2j";
%internal to eta_{2j-1}
 {\ar|{\sigma_1} "a4eta2j1";"a5eta2j1"};
 {\ar|{\sigma_1} "b4eta2j1";"b5eta2j1"};
 {\ar|{\sigma_2} "b5eta2j1";"b6eta2j1"};
 {\ar "b4eta2j1";"a4eta2j1"};
 {\ar|{\sigma_{123}} "b4eta2j1";"a5eta2j1"};
 {\ar|{1+\sigma_{23}} "b5eta2j1";"a5eta2j1"};
 {\ar|{\sigma_3} "b6eta2j1";"a5eta2j1"};
 {\ar@/^1.5pc/|(0.4){\sigma_{12}} "b4eta2j1";"b6eta2j1"};
%eta_{2j-1} -rho_{3}-> lambda_1
 {\ar@/_2.5pc/|(0.45){\twoline{\sigma_{123}+}{\sigma_3\sigma_2\sigma_1}} "a4eta2j1";"c1la1"};
 {\ar@/_2pc/|(0.5){\sigma_{23}} "a5eta2j1";"c1la1"};
 {\ar@/_3pc/|(0.35){\sigma_{12}} "a4eta2j1";"d3la1"};
 {\ar@/_2.25pc/|(0.39){\sigma_{2}} "a5eta2j1";"d3la1"};
 {\ar@/_1.5pc/|(0.53){\sigma_{123}} "b4eta2j1";"d2la1"};
 {\ar@/_1.5pc/|{\sigma_{23}} "b5eta2j1";"d2la1"};
 {\ar@/_1.5pc/|(0.3){\sigma_{3}} "b6eta2j1";"d2la1"};
 {\ar@/_1pc/|(0.55){\twoline{\sigma_{123}+}{\sigma_3\sigma_2\sigma_1}} "b4eta2j1";"d4la1"};
 {\ar@/_0.5pc/|(0.4){\sigma_{23}} "b5eta2j1";"d4la1"};
 {\ar|(0.3){\sigma_{3}} "b6eta2j1";"d4la1"};
 {\ar@/_0.3pc/ "b6eta2j1";"d3la1"};
%internal to lambda_1
 {\ar|{\sigma_1} "d1la1";"d2la1"};
 {\ar|{\sigma_1} "e1la1";"e2la1"};
 {\ar|{\sigma_2} "d2la1";"d3la1"};
 {\ar|{\sigma_2} "e2la1";"e3la1"};
 {\ar@/^1.5pc/|(0.4){\sigma_{12}} "d1la1";"d3la1"};
 {\ar@/^1.5pc/|(0.35){\sigma_{12}} "e1la1";"e3la1"};
 {\ar "e1la1";"d1la1"};
 {\ar "e2la1";"d2la1"};
 {\ar "e3la1";"d3la1"};
 {\ar "d4la1";"c1la1"};
 {\ar|{\sigma_{123}} "d1la1";"c1la1"};
 {\ar|(0.55){\sigma_{23}} "d2la1";"c1la1"};
 {\ar|(0.6){\sigma_{3}} "d3la1";"c1la1"};
 {\ar|{\sigma_{123}} "e1la1";"d4la1"};
 {\ar|(0.35){\sigma_{23}} "e2la1";"d4la1"};
 {\ar|{\sigma_{3}} "e3la1";"d4la1"};
%internal to lambda_k
 {\ar|{\sigma_1} "d1lal";"d2lal"};
 {\ar|{\sigma_1} "e1lal";"e2lal"};
 {\ar|{\sigma_2} "d2lal";"d3lal"};
 {\ar|{\sigma_2} "e2lal";"e3lal"};
 {\ar@/^1.5pc/|(0.4){\sigma_{12}} "d1lal";"d3lal"};
 {\ar@/^1.5pc/|(0.35){\sigma_{12}} "e1lal";"e3lal"};
 {\ar "e1lal";"d1lal"};
 {\ar "e2lal";"d2lal"};
 {\ar "e3lal";"d3lal"};
 {\ar "d4lal";"c1lal"};
 {\ar|{\sigma_{123}} "d1lal";"c1lal"};
 {\ar|{\sigma_{23}} "d2lal";"c1lal"};
 {\ar|{\sigma_{3}} "d3lal";"c1lal"};
 {\ar|{\sigma_{123}} "e1lal";"d4lal"};
 {\ar|(0.35){\sigma_{23}} "e2lal";"d4lal"};
 {\ar|{\sigma_{3}} "e3lal";"d4lal"};
%lambda_1 - - rho_{23} - -> lambda_k
 {\ar@{-->}@/_2.5pc/|(0.5){\twoline{\sigma_{123}+}{\sigma_3\sigma_2\sigma_1}} "d1la1";"c1lal"};
 {\ar@{-->}@/_2pc/|(0.55){\sigma_{23}} "d2la1";"c1lal"};
 {\ar@{-->}|{\sigma_{23}} "d4la1";"c1lal"};
 {\ar@{-->}@/_3pc/|(0.35){\sigma_{12}} "d1la1";"d3lal"};
 {\ar@{-->}@/_2.25pc/|(0.39){\sigma_{2}} "d2la1";"d3lal"};
 {\ar@{-->}@/_1.5pc/|(0.45){\sigma_{2}} "d4la1";"d3lal"};
 {\ar@{-->}@/_1.5pc/|(0.53){\sigma_{123}} "e1la1";"d2lal"};
 {\ar@{-->}@/_1.5pc/|{\sigma_{23}} "e2la1";"d2lal"};
 {\ar@{-->}@/_1.5pc/|(0.3){\sigma_{3}} "e3la1";"d2lal"};
 {\ar@{-->}@/_1pc/|(0.55){\twoline{\sigma_{123}+}{\sigma_3\sigma_2\sigma_1}} "e1la1";"d4lal"};
 {\ar@{-->}@/_0.5pc/|(0.4){\sigma_{23}} "e2la1";"d4lal"};
 {\ar@{-->}|(0.3){\sigma_{3}} "e3la1";"d4lal"};
 {\ar@{-->}@/_0.3pc/ "e3la1";"d3lal"};
%internal to eta_{2j}
 {\ar|{\sigma_1} "a4eta2j";"a5eta2j"};
 {\ar|{\sigma_1} "b4eta2j";"b5eta2j"};
 {\ar|{\sigma_2} "b5eta2j";"b6eta2j"};
 {\ar "b4eta2j";"a4eta2j"};
 {\ar|{\sigma_{123}} "b4eta2j";"a5eta2j"};
 {\ar|{1+\sigma_{23}} "b5eta2j";"a5eta2j"};
 {\ar|{\sigma_3} "b6eta2j";"a5eta2j"};
 {\ar@/^1.5pc/|(0.4){\sigma_{12}} "b4eta2j";"b6eta2j"};
%lambda_l -rho_2-> eta_{2j}
 {\ar@/_2pc/ "d1lal";"a4eta2j"};
 {\ar@/_2pc/ "d2lal";"a5eta2j"};
 {\ar@/_1.25pc/ "d4lal";"a5eta2j"};
 {\ar@/_2pc/ "e1lal";"b4eta2j"};
 {\ar@/_2pc/ "e2lal";"b5eta2j"};
 {\ar@/_2pc/ "e3lal";"b6eta2j"};
\endxy
\]
\end{tiny}
\caption{The subspace $P_{\hor}^j$, corresponding to a horizontal stable chain
$\eta_{2j-1} \xrightarrow{\rho_3} \lambda^j_1 \xrightarrow{\rho_{23}} \cdots
\xrightarrow{\rho_{23}} \lambda^j_{l_j} \xrightarrow{\rho_2} \eta_{2j}$.}
\label{fig:Phor}
\end{figure}

\begin{figure}
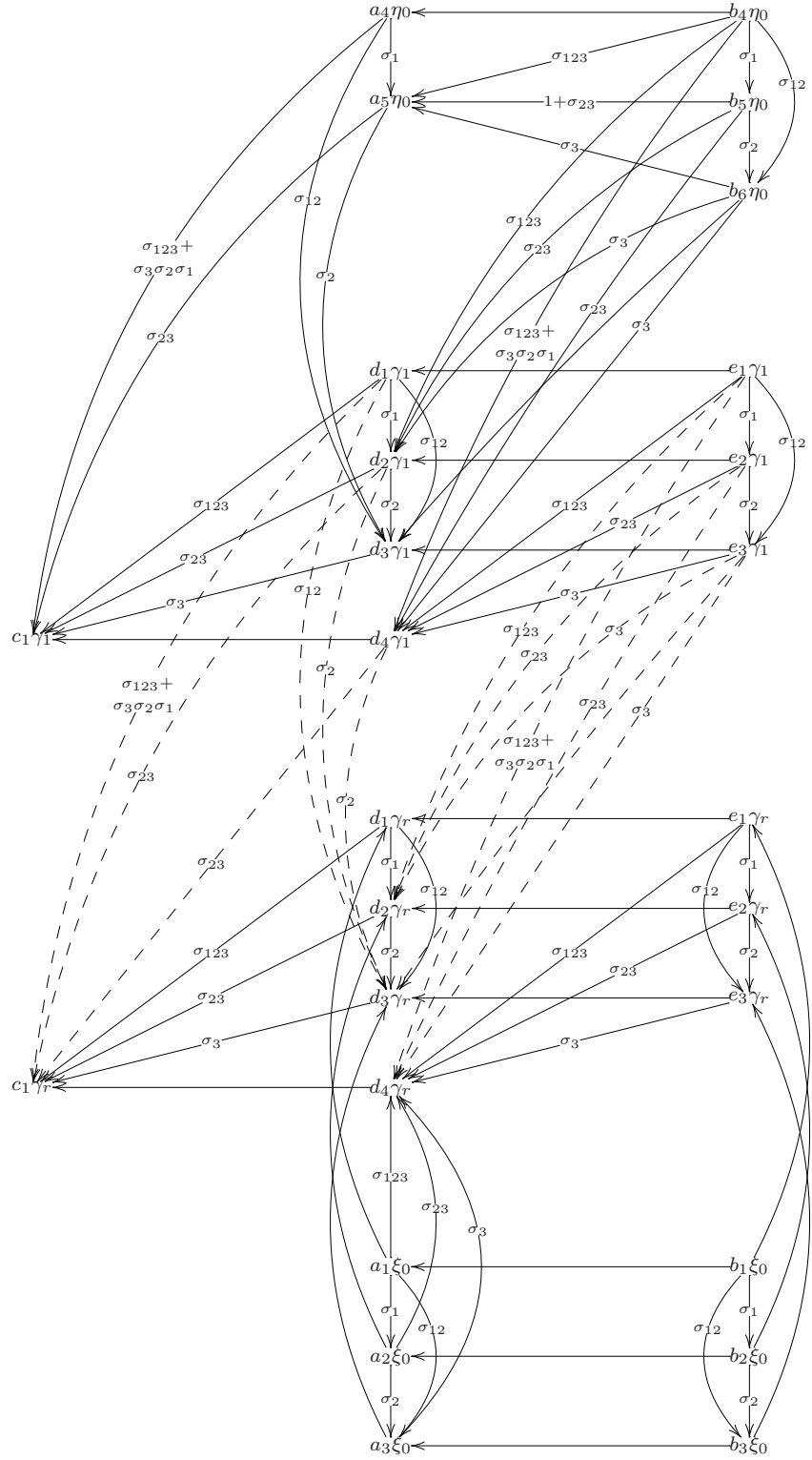
 \centering
\begin{tiny}
\[
\xy
 (50,0)*{a_4 \eta_0}="a4eta0"; (100,0)*{b_4 \eta_0}="b4eta0";
 (50,-12.5)*{a_5 \eta_0}="a5eta0"; (100,-12.5)*{b_5 \eta_0}="b5eta0";
 (100,-25)*{b_6 \eta_0}="b6eta0";
 (50,-50)*{d_1 \gamma_1}="d1ga1";
 (100,-50)*{e_1 \gamma_1}="e1ga1";
 (50,-62.5)*{d_2 \gamma_1}="d2ga1";
 (100,-62.5)*{e_2 \gamma_1}="e2ga1";
 (50,-75)*{d_3 \gamma_1}="d3ga1";
 (100,-75)*{e_3 \gamma_1}="e3ga1";
 (0,-87.5)*{c_1 \gamma_1}="c1ga1";
 (50,-87.5)*{d_4 \gamma_1}="d4ga1";
 (50,-112.5)*{d_1 \gamma_r}="d1gar";
 (100,-112.5)*{e_1\gamma_r}="e1gar";
 (50,-125)*{d_2 \gamma_r}="d2gar";
 (100,-125)*{e_2 \gamma_r}="e2gar";
 (50,-137.5)*{d_3 \gamma_r}="d3gar";
 (100,-137.5)*{e_3 \gamma_r}="e3gar";
 (0,-150)*{c_1 \gamma_r}="c1gar";
 (50,-150)*{d_4 \gamma_r}="d4gar";
 (50,-175)*{a_1 \xi_0}="a1xi0"; (100,-175)*{b_1 \xi_0}="b1xi0";
 (50,-187.5)*{a_2 \xi_0}="a2xi0"; (100,-187.5)*{b_2 \xi_0}="b2xi0";
 (50,-200)*{a_3 \xi_0}="a3xi0"; (100,-200)*{b_3 \xi_0}="b3xi0";
%internal to eta_0
 {\ar|{\sigma_1} "a4eta0";"a5eta0"};
 {\ar|{\sigma_1} "b4eta0";"b5eta0"};
 {\ar|{\sigma_2} "b5eta0";"b6eta0"};
 {\ar "b4eta0";"a4eta0"};
 {\ar|{\sigma_{123}} "b4eta0";"a5eta0"};
 {\ar|{1+\sigma_{23}} "b5eta0";"a5eta0"};
 {\ar|{\sigma_3} "b6eta0";"a5eta0"};
 {\ar@/^1.5pc/|(0.4){\sigma_{12}} "b4eta0";"b6eta0"};
%eta_0 -rho_{3}-> gamma_1
 {\ar@/_2.5pc/|(0.45){\twoline{\sigma_{123}+}{\sigma_3\sigma_2\sigma_1}} "a4eta0";"c1ga1"};
 {\ar@/_2pc/|(0.5){\sigma_{23}} "a5eta0";"c1ga1"};
 {\ar@/_3pc/|(0.35){\sigma_{12}} "a4eta0";"d3ga1"};
 {\ar@/_2.25pc/|(0.39){\sigma_{2}} "a5eta0";"d3ga1"};
 {\ar@/_1.5pc/|(0.53){\sigma_{123}} "b4eta0";"d2ga1"};
 {\ar@/_1.5pc/|{\sigma_{23}} "b5eta0";"d2ga1"};
 {\ar@/_1.5pc/|(0.3){\sigma_{3}} "b6eta0";"d2ga1"};
 {\ar@/_1pc/|(0.55){\twoline{\sigma_{123}+}{\sigma_3\sigma_2\sigma_1}} "b4eta0";"d4ga1"};
 {\ar@/_0.5pc/|(0.4){\sigma_{23}} "b5eta0";"d4ga1"};
 {\ar|(0.3){\sigma_{3}} "b6eta0";"d4ga1"};
 {\ar@/_0.3pc/ "b6eta0";"d3ga1"};
%internal to gamma_1
 {\ar|{\sigma_1} "d1ga1";"d2ga1"};
 {\ar|{\sigma_1} "e1ga1";"e2ga1"};
 {\ar|{\sigma_2} "d2ga1";"d3ga1"};
 {\ar|{\sigma_2} "e2ga1";"e3ga1"};
 {\ar@/^1.5pc/|(0.4){\sigma_{12}} "d1ga1";"d3ga1"};
 {\ar@/^1.5pc/|(0.4){\sigma_{12}} "e1ga1";"e3ga1"};
 {\ar "e1ga1";"d1ga1"};
 {\ar "e2ga1";"d2ga1"};
 {\ar "e3ga1";"d3ga1"};
 {\ar "d4ga1";"c1ga1"};
 {\ar|{\sigma_{123}} "d1ga1";"c1ga1"};
 {\ar|(0.55){\sigma_{23}} "d2ga1";"c1ga1"};
 {\ar|(0.6){\sigma_{3}} "d3ga1";"c1ga1"};
 {\ar|{\sigma_{123}} "e1ga1";"d4ga1"};
 {\ar|(0.35){\sigma_{23}} "e2ga1";"d4ga1"};
 {\ar|{\sigma_{3}} "e3ga1";"d4ga1"};
%internal to gamma_r
 {\ar|{\sigma_1} "d1gar";"d2gar"};
 {\ar|{\sigma_1} "e1gar";"e2gar"};
 {\ar|{\sigma_2} "d2gar";"d3gar"};
 {\ar|{\sigma_2} "e2gar";"e3gar"};
 {\ar@/^1.5pc/|(0.4){\sigma_{12}} "d1gar";"d3gar"};
 {\ar@/_1.5pc/|(0.4){\sigma_{12}} "e1gar";"e3gar"};
 {\ar "e1gar";"d1gar"};
 {\ar "e2gar";"d2gar"};
 {\ar "e3gar";"d3gar"};
 {\ar "d4gar";"c1gar"};
 {\ar|{\sigma_{123}} "d1gar";"c1gar"};
 {\ar|{\sigma_{23}} "d2gar";"c1gar"};
 {\ar|{\sigma_{3}} "d3gar";"c1gar"};
 {\ar|{\sigma_{123}} "e1gar";"d4gar"};
 {\ar|(0.35){\sigma_{23}} "e2gar";"d4gar"};
 {\ar|{\sigma_{3}} "e3gar";"d4gar"};
%gamma_1 - - rho_{23} - -> gamma_r
 {\ar@{-->}@/_2.5pc/|(0.5){\twoline{\sigma_{123}+}{\sigma_3\sigma_2\sigma_1}} "d1ga1";"c1gar"};
 {\ar@{-->}@/_2pc/|(0.55){\sigma_{23}} "d2ga1";"c1gar"};
 {\ar@{-->}|{\sigma_{23}} "d4ga1";"c1gar"};
 {\ar@{-->}@/_3pc/|(0.35){\sigma_{12}} "d1ga1";"d3gar"};
 {\ar@{-->}@/_2.25pc/|(0.39){\sigma_{2}} "d2ga1";"d3gar"};
 {\ar@{-->}@/_1.5pc/|(0.45){\sigma_{2}} "d4ga1";"d3gar"};
 {\ar@{-->}@/_1.5pc/|(0.53){\sigma_{123}} "e1ga1";"d2gar"};
 {\ar@{-->}@/_1.5pc/|{\sigma_{23}} "e2ga1";"d2gar"};
 {\ar@{-->}@/_1.5pc/|(0.3){\sigma_{3}} "e3ga1";"d2gar"};
 {\ar@{-->}@/_1pc/|(0.55){\twoline{\sigma_{123}+}{\sigma_3\sigma_2\sigma_1}} "e1ga1";"d4gar"};
 {\ar@{-->}@/_0.5pc/|(0.4){\sigma_{23}} "e2ga1";"d4gar"};
 {\ar@{-->}|(0.3){\sigma_{3}} "e3ga1";"d4gar"};
 {\ar@{-->}@/_0.3pc/ "e3ga1";"d3gar"};
%internal to xi_0
 {\ar|{\sigma_1} "a1xi0";"a2xi0"};
 {\ar|{\sigma_1} "b1xi0";"b2xi0"};
 {\ar|{\sigma_2} "a2xi0";"a3xi0"};
 {\ar|{\sigma_2} "b2xi0";"b3xi0"};
 {\ar "b1xi0";"a1xi0"};
 {\ar "b2xi0";"a2xi0"};
 {\ar "b3xi0";"a3xi0"};
 {\ar@/^1.5pc/|(0.35){\sigma_{12}} "a1xi0";"a3xi0"};
 {\ar@/_1.5pc/|(0.35){\sigma_{12}} "b1xi0";"b3xi0"};
%xi_0 -rho_1-> gamma_r
 {\ar|{\sigma_{123}} "a1xi0";"d4gar"};
 {\ar@/_1.5pc/|(0.55){\sigma_{23}} "a2xi0";"d4gar"};
 {\ar@/_3pc/|(0.6){\sigma_{3}} "a3xi0";"d4gar"};
 {\ar@/^2pc/ "a1xi0";"d1gar"};
 {\ar@/^2pc/ "a2xi0";"d2gar"};
 {\ar@/^2pc/ "a3xi0";"d3gar"};
 {\ar@/_2pc/ "b1xi0";"e1gar"};
 {\ar@/_2pc/ "b2xi0";"e2gar"};
 {\ar@/_2pc/ "b3xi0";"e3gar"};
\endxy
\]
\end{tiny}
\caption{The subspace $P_{\unst}$ when $s< 2\tau(J)$, corresponding to the
unstable chain $\eta_0 \xrightarrow{\rho_3} \gamma_1 \xrightarrow{\rho_{23}}
\cdots \xrightarrow{\rho_{23}} \gamma_s \xleftarrow{\rho_1} \xi_0$.}
\label{fig:Punst,s<2tau}
\end{figure}

\begin{figure}
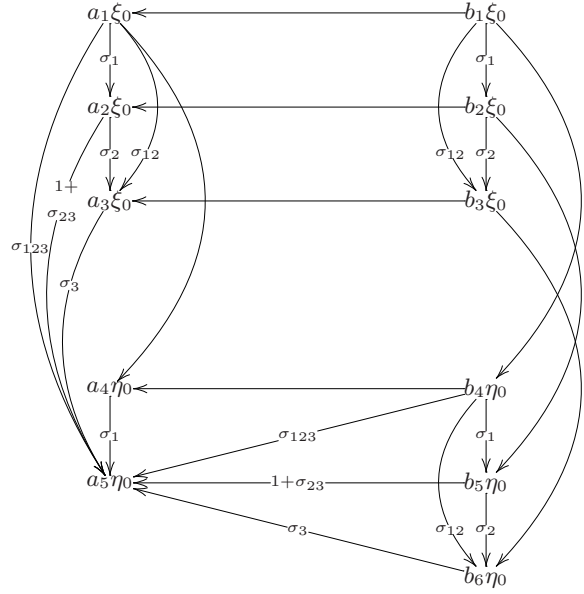
 \centering
\begin{scriptsize}
\[
\xy
 (50,0)*{a_1 \xi_0}="a1xi0"; (100,0)*{b_1 \xi_0}="b1xi0";
 (50,-12.5)*{a_2 \xi_0}="a2xi0"; (100,-12.5)*{b_2 \xi_0}="b2xi0";
 (50,-25)*{a_3 \xi_0}="a3xi0"; (100,-25)*{b_3 \xi_0}="b3xi0";
 (50,-50)*{a_4 \eta_0}="a4eta0"; (100,-50)*{b_4 \eta_0}="b4eta0";
 (50,-62.5)*{a_5 \eta_0}="a5eta0"; (100,-62.5)*{b_5 \eta_0}="b5eta0";
 (100,-75)*{b_6 \eta_0}="b6eta0";
%internal to xi_0
 {\ar|{\sigma_1} "a1xi0";"a2xi0"};
 {\ar|{\sigma_1} "b1xi0";"b2xi0"};
 {\ar|{\sigma_2} "a2xi0";"a3xi0"};
 {\ar|{\sigma_2} "b2xi0";"b3xi0"};
 {\ar "b1xi0";"a1xi0"};
 {\ar "b2xi0";"a2xi0"};
 {\ar "b3xi0";"a3xi0"};
 {\ar@/^1.5pc/|(.75){\sigma_{12}} "a1xi0";"a3xi0"};
 {\ar@/_1.5pc/|(.75){\sigma_{12}} "b1xi0";"b3xi0"};
%internal to eta_0
 {\ar|{\sigma_1} "a4eta0";"a5eta0"};
 {\ar|{\sigma_1} "b4eta0";"b5eta0"};
 {\ar|{\sigma_2} "b5eta0";"b6eta0"};
 {\ar "b4eta0";"a4eta0"};
 {\ar|{\sigma_{123}} "b4eta0";"a5eta0"};
 {\ar|{1+\sigma_{23}} "b5eta0";"a5eta0"};
 {\ar|{\sigma_3} "b6eta0";"a5eta0"};
 {\ar@/_1.5pc/|(.75){\sigma_{12}} "b4eta0";"b6eta0"};
%xi_0 -rho_12-> eta_0
 {\ar@/^3pc/ "a1xi0";"a4eta0"};
 {\ar@/_2.5pc/|{\sigma_{123}} "a1xi0";"a5eta0"};
 {\ar@/_2pc/|(.25){\twoline{1+}{\sigma_{23}}} "a2xi0";"a5eta0"};
 {\ar@/_1.5pc/|(.3){\sigma_3} "a3xi0";"a5eta0"};
 {\ar@/^3pc/ "b1xi0";"b4eta0"};
 {\ar@/^3pc/ "b2xi0";"b5eta0"};
 {\ar@/^3pc/ "b3xi0";"b6eta0"};
\endxy
\]
\end{scriptsize}
\caption{The subspace $P_{\unst}$ when $s= 2\tau(J)$, corresponding to the
unstable chain $\xi_0 \xrightarrow{\rho_{12}} \eta_0$.}
\label{fig:Punst,s=2tau}
\end{figure}

In addition, there are a few more multiplications that preserve the splitting,
coming from multiplications in $\CFAA(\YY,B_3,0)$ that use sequences like
$\rho_3 \rho_2$, $\rho_3 \rho_{23}$, or $\rho_{23} \rho_{23}$. These
multiplications are not shown in Figures \ref{fig:Pver} through
\ref{fig:Punst,s=2tau}. They are as follows:

\begin{itemize}
\item In $P_{\ver}^j$, when $k_j>1$, there are multiplications
\begin{equation} \label{eq:Pver-longer}
\begin{alignedat}{3}
b_1 \xi_{2j} & \xrightarrow{\sigma_3 \sigma_2 \sigma_{12}} d_3 \kappa^j_2
 &\qquad\quad
b_1 \xi_{2j} & \xrightarrow{\sigma_3 \sigma_2 \sigma_{123}} c_1 \kappa^j_2
 &\qquad& \\
e_1 \kappa^j_i & \xrightarrow{\sigma_3 \sigma_2 \sigma_{12}} d_3 \kappa^j_{i+2}
 & e_1 \kappa^j_i & \xrightarrow{\sigma_3 \sigma_2 \sigma_{123}} c_1 \kappa^j_{i+2}
 & (i&=1,\dots, k_j-2).
\end{alignedat}
\end{equation}

\item In $P_{\hor}^j$, when $l_j=1$, there is a multiplication $b_4\eta_{2j-1} \xrightarrow{\sigma_3\sigma_2} a_4
\eta_{2j}$. When $l_j>1$, there are multiplications
\begin{equation} \label{eq:Phor-longer}
\begin{alignedat}{3}
b_4 \eta_{2j-1} & \xrightarrow{\sigma_3 \sigma_2 \sigma_{12}} d_3 \lambda^j_2
 &\qquad
b_4 \eta_{2j-1} & \xrightarrow{\sigma_3 \sigma_2 \sigma_{123}} c_1 \lambda^j_2
 &\qquad&
 \\
e_1 \lambda^j_{l_j-1} & \xrightarrow{\sigma_3 \sigma_2} a_4 \eta_{2j} && && \\
 e_1 \lambda^j_i & \xrightarrow{\sigma_3 \sigma_2 \sigma_{12}} d_3\lambda^j_{i+2}  &
 e_1 \lambda^j_i & \xrightarrow{\sigma_3 \sigma_2 \sigma_{123}} c_1 \lambda^j_{i+2} &
 (i&=1,\dots, l_j-2). \\
\end{alignedat}
\end{equation}

\item In $P_{\unst}$ in the case when $s< 2\tau(J)-1$, there are multiplications
\begin{equation} \label{eq:Punst-longer}
\begin{alignedat}{3}
b_4 \eta_0 & \xrightarrow{\sigma_3 \sigma_2 \sigma_{12}} d_3 \gamma_2
 &\qquad\quad
b_4 \eta_0 & \xrightarrow{\sigma_3 \sigma_2 \sigma_{123}} c_1 \gamma_2
 &\qquad& \\
e_1 \gamma_i & \xrightarrow{\sigma_3 \sigma_2 \sigma_{12}} d_3
 \gamma_{i+2} &
e_1 \gamma_i & \xrightarrow{\sigma_3 \sigma_2 \sigma_{123}} c_1
\gamma_{i+2} & (i&=1,\dots, r-2). \\
\end{alignedat}
\end{equation}

\end{itemize}

Finally, we must consider the multiplications in the tensor product that do not
respect the splitting in \eqref{eq:Psplit}. These arise from sequences of
arrows in $\CFD(\XX_J^s)$ that involve multiple stable or unstable chains, and
they depend on the change-of-basis coefficients relating $\{\eta_0, \dots,
\eta_{2n}\}$ and $\{\xi_0, \dots, \xi_{2n}\}$.

For instance, if $\eta_{2j} = \xi_{2h}$ (where $j,h \in \{1, \dots,n\}$), then
$\CFD(\XX_J^s)$ contains a string of arrows of the form
\[
\eta_{2j-1} \xrightarrow{\rho_3} \lambda^j_1 \xrightarrow{\rho_{23}} \cdots
\xrightarrow{\rho_{23}} \lambda^j_{l_j} \xrightarrow{\rho_2} \eta_{2j}
\xrightarrow{\rho_{123}} \kappa^h_1 \xrightarrow{\rho_{23}} \cdots
\xrightarrow{\rho_{23}} \kappa^h_{k_h}.
\]
Any multiplication in $\CFAA(\YY,B_3,0)$ that uses a contiguous subsequence of
\[
\rho_3, \underbrace{\rho_{23}, \dots, \rho_{23}}_{l_j-1 \text{ times}}, \rho_2,
\rho_{123}, \underbrace{\rho_{23}, \dots, \rho_{23}}_{k_h-1 \text{ times}}
\]
contributes a nonzero multiplication in the tensor product that need not
respect the splitting. Similarly, if $\eta_{2j} = \xi_{2h-1}$, then the same is
true for contiguous subsequences of
\[
\rho_3, \underbrace{\rho_{23}, \dots, \rho_{23}}_{l_j-1 \text{ times}}, \rho_2,
\rho_1.
\]
Similar sequences may also occur near the unstable chain, where we take $\xi_0$
instead of $\xi_{2h-1}$ or $\xi_{2h}$. By Proposition \ref{prop:norho1rho2},
these are the only such sequences that occur. More generally, if the
coefficient of $\xi_p$ in $\eta_{2j}$ is nonzero, we obtain multiplications
that do not respect the splitting in \eqref{eq:Psplit}. We make this notion
more precise below.

By inspecting the matrices $M_{xy}$, we see that the only sequences of this
form that actually occur in $\CFAA(\YY,B_3,0)$ are $\rho_3 \rho_2 \rho_{123}$,
$\rho_3 \rho_2 \rho_{12}$, and $\rho_3 \rho_2 \rho_1$, which occur in the first
three rows of $M_{ac}$, $M_{ad}$, $M_{ba}$, $M_{bc}$, and $M_{bd}$.
Accordingly, the only multiplications that do not preserve the splitting arise
when there is a horizontal edge $\eta_{2j-1} \to \eta_{2j}$ of length $1$, and
they act on the elements $a_i \boxtimes \eta_{2j-1}$ and $b_i \boxtimes
\eta_{2j-1}$ $(i=1,2,3)$.

Notice that there are no multiplications into or out of any of the subspaces
$P_{\hor}^j$. Therefore, each $P_{\hor}^j$ is actually a direct summand of
$\CFA(\VV, D_{J,s})$ as an $\AA_\infty$ submodule, as is $P = \bigoplus_{j=1}^n
P_{\ver}^j \oplus P_{\unst}$. This implies that the tensor product $\CFA(\VV,
D_{J,s}) \boxtimes \CFD(\XX_K^t)$ (whose total homology, ignoring the
filtration, is $\HF(S^3) \cong \F$) will also split as a direct sum. We shall
see in Subsection \ref{subsec:Asigma} that the direct summand coming from $P$ has odd rank, which implies that it must have a nonzero contribution to the total homology. Therefore, each summand coming from $P_{\hor}^j$ is
acyclic and hence does not affect the computation of $\tau(D_{J,s}(K,t))$.
Thus, we shall henceforth ignore the submodules $P_{\hor}^j$.

It is preferable to describe all of the multiplications that do not respect the
splitting in terms of the bases specified in \eqref{eq:Pbases}. Recall that
$(x_{p,q})$ and $(y_{p,q})$ are the change-of-basis matrices, so that $\xi_p =
\sum_{q=0}^{2n} x_{p,q} \eta_q$ and $\eta_p = \sum_{q=0}^{2n} y_{p,q} \xi_q$.
Let $\mathfrak{j}$ denote the set $\{j \in \{1, \dots, n\} \mid l_j=1\}$. For
each $p \in \{0, \dots, 2n\}$ and $h \in \{1,\dots,n\}$, each $j \in
\mathfrak{j}$ for which $x_{p,2j-1}=1$ and $y_{2j,2h-1}=1$ contributes
multiplications (which we will specify shortly) from $a_i \xi_p$ and/or $b_i
\xi_p$ $(i=1,2,3)$ into $P_{\ver}^h$ via the sequence $\rho_3 \rho_2
\rho_{123}$. Of course, multiple values of $j$ may satisfy this criterion, but
they all contribute the same multiplications, so we really only care about the
count of such $j$ modulo $2$. That is, define $u_{p,h} = \sum_{j \in
\mathfrak{j}} x_{p,2j-1} y_{2j,2h-1}$; there are multiplications from $a_i
\xi_p$ and $b_i \xi_p$ into $P_{\ver}^h$ iff $u_{p,h}=1$.

Similarly, each $j$ for which $x_{p,2j-1}=1$ and $y_{2j,2h}$ $(h=1, \dots, n)$
contributes multiplications via $\rho_3 \rho_2 \rho_1$, so define $v_{p,h} =
\sum_{j \in \mathfrak{j}} x_{p,2j-1} y_{2j,2h}$. Finally, we set $w_p = \sum_{j
\in \mathfrak{j}} x_{p,2j-1} y_{2h,0}$; this determines whether there are
additional multiplications from $a_i \xi_p$ and $b_i \xi_p$ into the unstable
chain via $\rho_3\rho_2\rho_1$, $\rho_3\rho_2\rho_{12}$, or
$\rho_3\rho_2\rho_{123}$, according to whether $s < 2\tau(J)$, $s = 2\tau(J)$
or $s>2\tau(K)$, respectively (although we are ignoring the third case).

We now specify these multiplications:

\begin{figure}
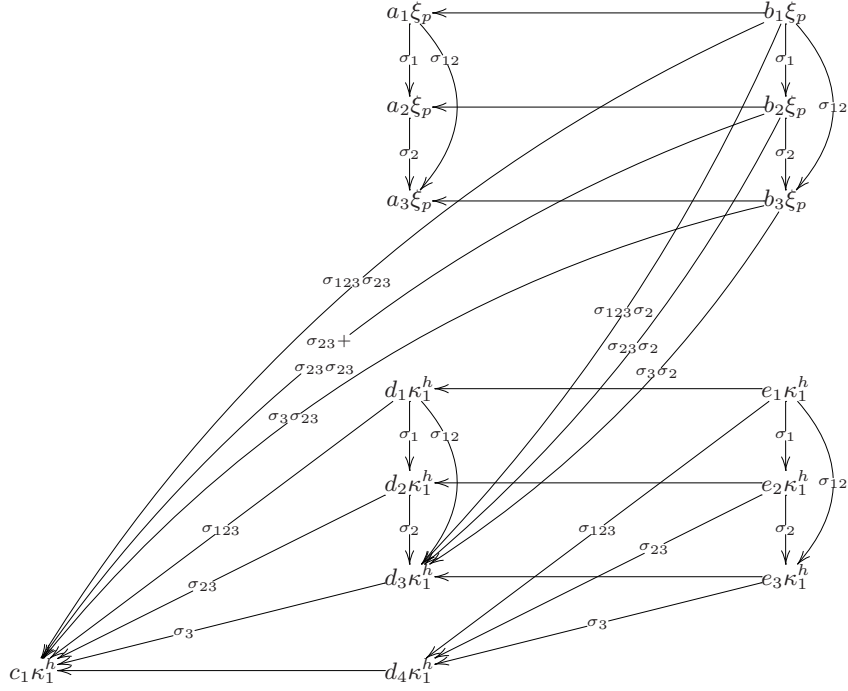
 \centering
\begin{scriptsize}
\[
\xy
 (50,0)*{a_1 \xi_p}="a1xip"; (100,0)*{b_1 \xi_p}="b1xip";
 (50,-12.5)*{a_2 \xi_p}="a2xip"; (100,-12.5)*{b_2 \xi_p}="b2xip";
 (50,-25)*{a_3 \xi_p}="a3xip"; (100,-25)*{b_3 \xi_p}="b3xip";
 (50,-50)*{d_1 \kappa^h_1}="d1ka1";
 (100,-50)*{e_1 \kappa^h_1}="e1ka1";
 (50,-62.5)*{d_2 \kappa^h_1}="d2ka1";
 (100,-62.5)*{e_2 \kappa^h_1}="e2ka1";
 (50,-75)*{d_3 \kappa^h_1}="d3ka1";
 (100,-75)*{e_3 \kappa^h_1}="e3ka1";
 (0,-87.5)*{c_1 \kappa^h_1}="c1ka1";
 (50,-87.5)*{d_4 \kappa^h_1}="d4ka1";
%internal to eta_{2j-1}
 {\ar|{\sigma_1} "a1xip";"a2xip"};
 {\ar|{\sigma_1} "b1xip";"b2xip"};
 {\ar|{\sigma_2} "a2xip";"a3xip"};
 {\ar|{\sigma_2} "b2xip";"b3xip"};
 {\ar "b1xip";"a1xip"};
 {\ar "b2xip";"a2xip"};
 {\ar "b3xip";"a3xip"};
 {\ar@/^1.5pc/|(.25){\sigma_{12}} "a1xip";"a3xip"};
 {\ar@/^1.5pc/|(.5){\sigma_{12}} "b1xip";"b3xip"};
%internal to kappa_1
 {\ar|{\sigma_1} "d1ka1";"d2ka1"};
 {\ar|{\sigma_1} "e1ka1";"e2ka1"};
 {\ar|{\sigma_2} "d2ka1";"d3ka1"};
 {\ar|{\sigma_2} "e2ka1";"e3ka1"};
 {\ar@/^1.5pc/|(0.25){\sigma_{12}} "d1ka1";"d3ka1"};
 {\ar@/^1.5pc/|(0.5){\sigma_{12}} "e1ka1";"e3ka1"};
 {\ar "e1ka1";"d1ka1"};
 {\ar "e2ka1";"d2ka1"};
 {\ar "e3ka1";"d3ka1"};
 {\ar "d4ka1";"c1ka1"};
 {\ar|{\sigma_{123}} "d1ka1";"c1ka1"};
 {\ar|(0.55){\sigma_{23}} "d2ka1";"c1ka1"};
 {\ar|(0.6){\sigma_{3}} "d3ka1";"c1ka1"};
 {\ar|{\sigma_{123}} "e1ka1";"d4ka1"};
 {\ar|(0.35){\sigma_{23}} "e2ka1";"d4ka1"};
 {\ar|{\sigma_{3}} "e3ka1";"d4ka1"};
%eta_{2j-1} -rho_3 rho_2 rho_123 -> kappa_1
 {\ar@/_2.5pc/|{\sigma_{123} \sigma_{23}} "b1xip";"c1ka1"};
 {\ar@/_2.5pc/|(.55){\twoline{\sigma_{23}+}{\sigma_{23}\sigma_{23}}} "b2xip";"c1ka1"};
 {\ar@/_2.5pc/|(.6){\sigma_3 \sigma_{23}} "b3xip";"c1ka1"};
 {\ar@/^1pc/|{\sigma_{123} \sigma_2} "b1xip";"d3ka1"};
 {\ar@/^1pc/|(.47){\sigma_{23} \sigma_2} "b2xip";"d3ka1"};
 {\ar@/^1pc/|(.4){\sigma_3 \sigma_2} "b3xip";"d3ka1"};
\endxy
\]
\end{scriptsize}
\caption{Multiplications coming from a sequence $\rho_3 \rho_2 \rho_{123}$ when
$u_{p,h}=1$.} \label{fig:rho3rho2rho123}
\end{figure}

\begin{figure}
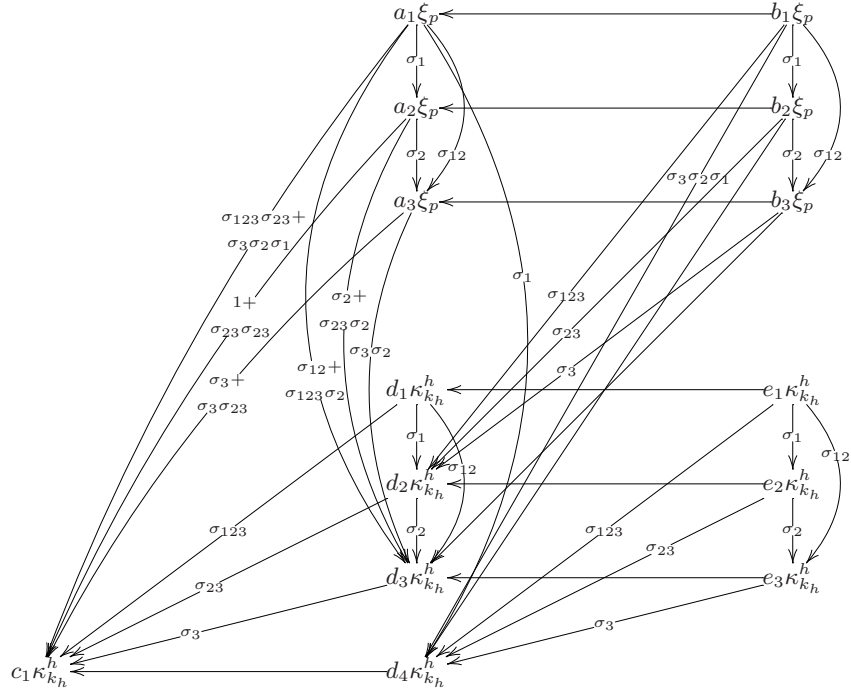
 \centering
\begin{scriptsize}
\[
\xy
 (50,0)*{a_1 \xi_p}="a1xip"; (100,0)*{b_1 \xi_p}="b1xip";
 (50,-12.5)*{a_2 \xi_p}="a2xip"; (100,-12.5)*{b_2 \xi_p}="b2xip";
 (50,-25)*{a_3 \xi_p}="a3xip"; (100,-25)*{b_3 \xi_p}="b3xip";
 (50,-50)*{d_1 \kappa^h_{k_h}}="d1kak";
 (100,-50)*{e_1 \kappa^h_{k_h}}="e1kak";
 (50,-62.5)*{d_2 \kappa^h_{k_h}}="d2kak";
 (100,-62.5)*{e_2 \kappa^h_{k_h}}="e2kak";
 (50,-75)*{d_3 \kappa^h_{k_h}}="d3kak";
 (100,-75)*{e_3 \kappa^h_{k_h}}="e3kak";
 (0,-87.5)*{c_1 \kappa^h_{k_h}}="c1kak";
 (50,-87.5)*{d_4 \kappa^h_{k_h}}="d4kak";
%internal to eta_{2j-1}
 {\ar|{\sigma_1} "a1xip";"a2xip"};
 {\ar|{\sigma_1} "b1xip";"b2xip"};
 {\ar|{\sigma_2} "a2xip";"a3xip"};
 {\ar|{\sigma_2} "b2xip";"b3xip"};
 {\ar "b1xip";"a1xip"};
 {\ar "b2xip";"a2xip"};
 {\ar "b3xip";"a3xip"};
 {\ar@/^1.5pc/|(.75){\sigma_{12}} "a1xip";"a3xip"};
 {\ar@/^1.5pc/|(.75){\sigma_{12}} "b1xip";"b3xip"};
%internal to kappa_1
 {\ar|{\sigma_1} "d1kak";"d2kak"};
 {\ar|{\sigma_1} "e1kak";"e2kak"};
 {\ar|{\sigma_2} "d2kak";"d3kak"};
 {\ar|{\sigma_2} "e2kak";"e3kak"};
 {\ar@/^1.5pc/|(0.42){\sigma_{12}} "d1kak";"d3kak"};
 {\ar@/^1.5pc/|(0.35){\sigma_{12}} "e1kak";"e3kak"};
 {\ar "e1kak";"d1kak"};
 {\ar "e2kak";"d2kak"};
 {\ar "e3kak";"d3kak"};
 {\ar "d4kak";"c1kak"};
 {\ar|{\sigma_{123}} "d1kak";"c1kak"};
 {\ar|(0.55){\sigma_{23}} "d2kak";"c1kak"};
 {\ar|(0.6){\sigma_{3}} "d3kak";"c1kak"};
 {\ar|{\sigma_{123}} "e1kak";"d4kak"};
 {\ar|(0.35){\sigma_{23}} "e2kak";"d4kak"};
 {\ar|{\sigma_{3}} "e3kak";"d4kak"};
%eta_{2j-1} -rho_3 rho_2 rho_1 -> kappa_1
 {\ar@/_1pc/|(.35){\twoline{\sigma_{123}\sigma_{23} +}{\sigma_3 \sigma_2 \sigma_1}} "a1xip";"c1kak"};
 {\ar@/_1pc/|(.4){\twoline{1 +}{\sigma_{23} \sigma_{23}}} "a2xip";"c1kak"};
 {\ar@/_1pc/|(.45){\twoline{\sigma_3+}{\sigma_3 \sigma_{23}}} "a3xip";"c1kak"};
 {\ar@/_3.5pc/|(.65){\twoline{\sigma_{12} +}{\sigma_{123} \sigma_2}} "a1xip";"d3kak"};
 {\ar@/_2.3pc/|(.43){\twoline{\sigma_{2} +}{\sigma_{23} \sigma_2}} "a2xip";"d3kak"};
 {\ar@/_1.5pc/|(.4){\sigma_3 \sigma_2} "a3xip";"d3kak"};
 {\ar@/^3.5pc/|(.4){\sigma_1} "a1xip";"d4kak"};
 {\ar|(.6){\sigma_{123}} "b1xip";"d2kak"};
 {\ar|(.6){\sigma_{23}} "b2xip";"d2kak"};
 {\ar|(.6){\sigma_3} "b3xip";"d2kak"};
 {\ar "b3xip";"d3kak"};
 {\ar|(.25){\sigma_3 \sigma_2 \sigma_1} "b1xip";"d4kak"};
 {\ar "b2xip";"d4kak"};
\endxy
\]
\end{scriptsize}
\caption[Multiplications coming from a sequence $\rho_3 \rho_2 \rho_1$ when
$v_{p,h}=1$.] {Multiplications coming from a sequence $\rho_3 \rho_2 \rho_1$
when $v_{p,h}=1$. If $w_p=1$ and $s<2\tau(K)$, we obtain the same
multiplications by replacing $\kappa^h_{k_h}$ by $\gamma_r$.}
\label{fig:rho3rho2rho1}
\end{figure}

\begin{itemize}
\item If $u_{p,h}=1$, the sequence $\rho_3 \rho_2 \rho_{123}$ provides the
multiplications shown in Figure \ref{fig:rho3rho2rho123}.

\item If $v_{p,h}=1$, the sequence $\rho_3 \rho_2 \rho_1$ provides the
multiplications shown in Figure \ref{fig:rho3rho2rho1}.

\item If $s<2\tau(J)$ and $w_p =1$, the sequence $\rho_3 \rho_2 \rho_1$ provides the multiplications shown
in Figure \ref{fig:rho3rho2rho1}, where we replace $\kappa^h_{k_h}$ by
$\gamma_r$.

\item Finally, if $s=2 \tau(K)$ and $w_p=1$,
the sequence $\rho_3 \rho_2 \rho_{12}$ provides the following multiplications:
\begin{equation} \label{eq:rho3rho2rho12}
\begin{split}
a_1 \xi_p & \xrightarrow{\sigma_1} a_5 \eta_0 \\
b_1 \xi_p & \xrightarrow{\sigma_{123}} a_5 \eta_0 \\
b_2 \xi_p & \xrightarrow{1+\sigma_{23}} a_5 \eta_0 \\
b_3 \xi_p & \xrightarrow{\sigma_3} a_5 \eta_0. \\
\end{split}
\end{equation}
\end{itemize}

%\subsection{Simplification of $\protect\CFA(\VV, D_{J,s})$}
\subsection{Simplification of \texorpdfstring{$\widehat{\mathrm{CFA}}(\VV,D_{J,s})$}{CFA(V,D\textunderscore\{J,s\})}}

Next, we may simplify $\CFA(\VV, D_{J,s})$ by canceling unmarked edges that
preserve the filtration level. In order to keep track of additional edges that
may appear, we must look carefully at the order of cancellation. As mentioned
above, we ignore the direct summands $P_{\hor}^j$. Define $P^0 = P_{\unst}$ and
$P^j = P_{\ver}^j$.

Assume first that $s<2 \tau(J)$.

For each $j \in \{1, \dots, n\}$, in $P^j$, we may cancel the differentials
$b_1 \xi_{2j-1} \to e_1 \kappa^j_{k_j}$, $b_2 \xi_{2j-1} \to e_2
\kappa^j_{k_j}$, $b_2 \xi_{2j-1} \to e_2 \kappa^j_{k_j}$, and $a_1 \xi_{2j-1}
\to d_1 \kappa^j_{k_j}$. Since the targets of those arrows do not lie at the
heads of any other arrows, no additional arrows are introduced. Similarly, in
$P^0$, cancel $b_1 \xi_0 \to e_1 \gamma_r$, $b_2 \xi_0 \to e_2 \gamma_r$, $b_2
\xi_0 \to e_2 \gamma_r$, and $a_1 \xi_0 \to d_1 \gamma_r$.

Next, we cancel the differentials $a_2 \xi_{2j-1} \to d_2 \kappa^j_{k_j}$ and
$a_2 \xi_0 \to d_2 \gamma_r$. Because of the edge $a_2 \xi_{2j-1}
\xrightarrow{\sigma_{23}} d_4 \kappa^j_{k_j}$, canceling $a_2 \xi_{2j-1} \to
d_2 \kappa^j_{k_j}$ introduces new multiplications:
\begin{equation} \label{eq:cancel1}
\begin{alignedat}{2}
e_1 \kappa^j_{k_j-1} & \xrightarrow{\sigma_{123} \sigma_{2}} a_3 \xi_{2j-1}
&\qquad
 e_1 \kappa^j_{k_j-1} & \xrightarrow{\sigma_{123} \sigma_{23}} d_4 \kappa^j_{k_j} \\
e_2 \kappa^j_{k_j-1} & \xrightarrow{\sigma_{23} \sigma_{2}} a_3 \xi_{2j-1} &
 e_2 \kappa^j_{k_j-1} & \xrightarrow{\sigma_{23} \sigma_{23}} d_4 \kappa^j_{k_j} \\
e_3 \kappa^j_{k_j-1} & \xrightarrow{\sigma_3 \sigma_{2}} a_3 \xi_{2j-1} &
 e_3 \kappa^j_{k_j-1} & \xrightarrow{\sigma_3 \sigma_{23}} d_4 \kappa^j_{k_j}.
\end{alignedat}
\end{equation}
(If $k_j=1$, then replace $e_i \kappa^j_{k_j-1}$ by $b_i \xi_{2j}$ in
\eqref{eq:cancel1}.) We shall examine the effects of these cancellations on the
edges that do not respect the splitting momentarily.

Next, because of the edge $a_3 \xi_{2j-1} \xrightarrow{\sigma_3} d_4
\kappa^j_{k_j}$, canceling $a_3 \xi_{2j-1} \to d_3 \kappa^j_{k_j}$ removes the
edge $e_3\kappa^j_{k_j-1} \xrightarrow{\sigma_3} d_4 \kappa^j_{k_j}$ and adds
edges
\begin{equation} \label{eq:cancel2}
\begin{split}
d_1 \kappa^j_{k_j-1} & \xrightarrow{\sigma_{12} \sigma_3} d_4 \kappa^j_{k_j} \\
d_2 \kappa^j_{k_j-1} & \xrightarrow{\sigma_2 \sigma_3} d_4 \kappa^j_{k_j} \\
d_4 \kappa^j_{k_j-1} & \xrightarrow{\sigma_2 \sigma_3} d_4 \kappa^j_{k_j} \\
e_1 \kappa^j_{k_j-2} & \xrightarrow{\sigma_3 \sigma_2 \sigma_{12} \sigma_3} d_4 \kappa^j_{k_j}. \\
\end{split}
\end{equation}
Because we will ultimately tensor with $\CFD(X_K^t)$, in which the sequences
$\sigma_2 \sigma_3$ and $\sigma_{12} \sigma_3$ do not appear, we may disregard
these four edges. We also eliminate the edge $e_3 \kappa^j_{k_j-1}
\xrightarrow{\sigma_3} d_4 \kappa^j_{k_j-1}$. The same thing occurs in $P^0$
when we cancel $a_3 \xi_0 \to d_3 \gamma_r$.

\begin{figure}
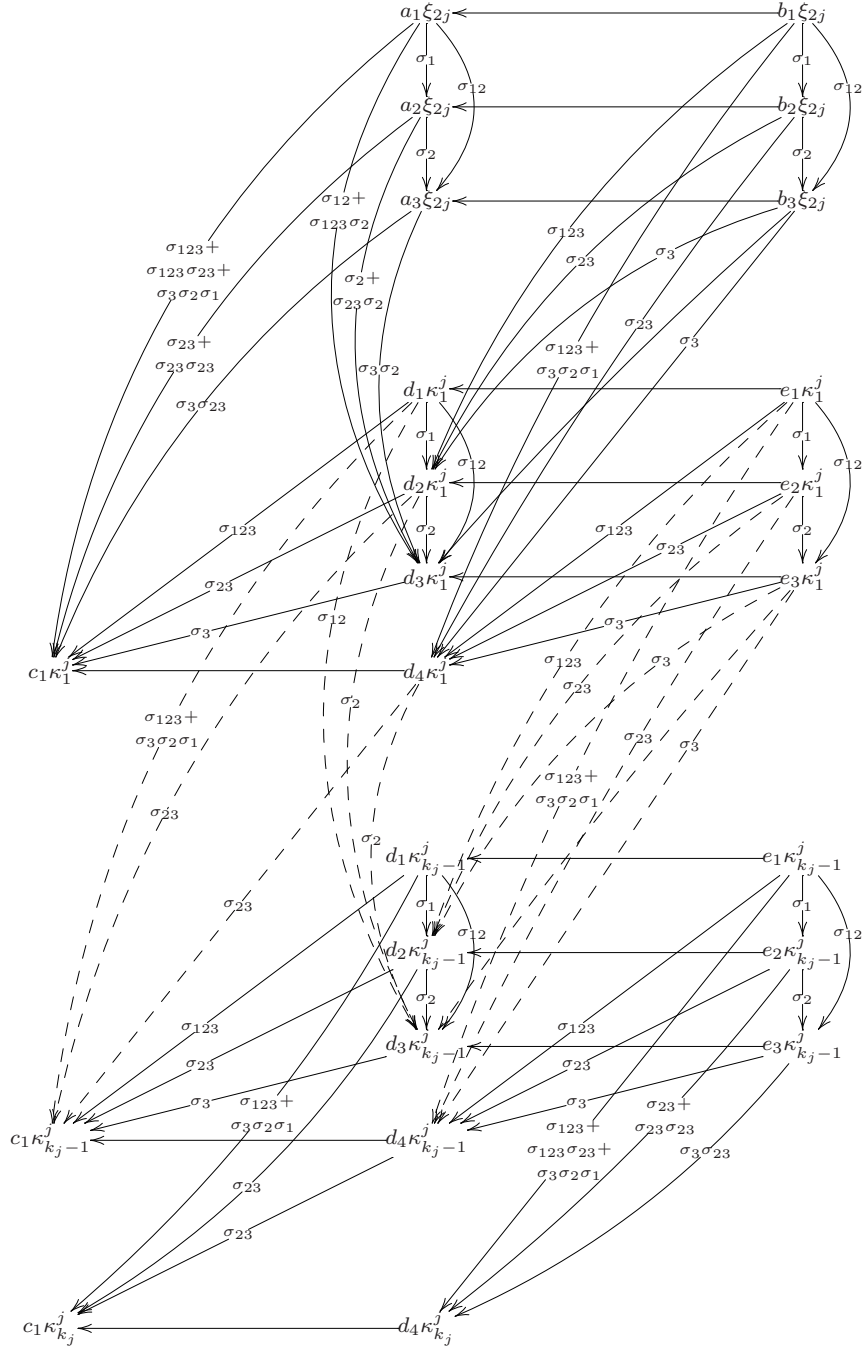
 \centering
\begin{tiny}
\[
\xy (50,0)*{a_1 \xi_{2j}}="a1xi2j"; (100,0)*{b_1 \xi_{2j}}="b1xi2j";
(50,-12.5)*{a_2 \xi_{2j}}="a2xi2j"; (100,-12.5)*{b_2 \xi_{2j}}="b2xi2j";
(50,-25)*{a_3 \xi_{2j}}="a3xi2j"; (100,-25)*{b_3 \xi_{2j}}="b3xi2j";
(50,-50)*{d_1 \kappa^j_1}="d1ka1"; (100,-50)*{e_1 \kappa^j_1}="e1ka1";
(50,-62.5)*{d_2 \kappa^j_1}="d2ka1"; (100,-62.5)*{e_2 \kappa^j_1}="e2ka1";
(50,-75)*{d_3 \kappa^j_1}="d3ka1"; (100,-75)*{e_3 \kappa^j_1}="e3ka1";
(0,-87.5)*{c_1 \kappa^j_1}="c1ka1";  (50,-87.5)*{d_4 \kappa^j_1}="d4ka1";
 (50,-112.5)*{d_1 \kappa^j_{k_j-1}}="d1kak1";
 (100,-112.5)*{e_1\kappa^j_{k_j-1}}="e1kak1";
 (50,-125)*{d_2 \kappa^j_{k_j-1}}="d2kak1";
 (100,-125)*{e_2 \kappa^j_{k_j-1}}="e2kak1";
 (50,-137.5)*{d_3 \kappa^j_{k_j-1}}="d3kak1";
 (100,-137.5)*{e_3 \kappa^j_{k_j-1}}="e3kak1";
 (0,-150)*{c_1 \kappa^j_{k_j-1}}="c1kak1";
 (50,-150)*{d_4 \kappa^j_{k_j-1}}="d4kak1";
 (0,-175)*{c_1 \kappa^j_{k_j}}="c1kak";
 (50,-175)*{d_4 \kappa^j_{k_j}}="d4kak";
%internal to xi_{2j}
 {\ar|{\sigma_1} "a1xi2j";"a2xi2j"};
 {\ar|{\sigma_1} "b1xi2j";"b2xi2j"};
 {\ar|{\sigma_2} "a2xi2j";"a3xi2j"};
 {\ar|{\sigma_2} "b2xi2j";"b3xi2j"};
 {\ar "b1xi2j";"a1xi2j"};
 {\ar "b2xi2j";"a2xi2j"};
 {\ar "b3xi2j";"a3xi2j"};
 {\ar@/^1.5pc/|(.4){\sigma_{12}} "a1xi2j";"a3xi2j"};
 {\ar@/^1.5pc/|(.4){\sigma_{12}} "b1xi2j";"b3xi2j"};
%xi_{2j} -rho_{123}-> kappa_1
 {\ar@/_2.5pc/|(0.45){\threeline{\sigma_{123}+}{\sigma_{123}\sigma_{23}+}{\sigma_3\sigma_2\sigma_1}} "a1xi2j";"c1ka1"};
 {\ar@/_2pc/|(0.5){\twoline{\sigma_{23}+}{\sigma_{23}\sigma_{23}}} "a2xi2j";"c1ka1"};
 {\ar@/_1.5pc/|{\sigma_{3}\sigma_{23}} "a3xi2j";"c1ka1"};
 {\ar@/_3pc/|(0.35){\twoline{\sigma_{12}+}{\sigma_{123}\sigma_2}} "a1xi2j";"d3ka1"};
 {\ar@/_2.25pc/|(0.39){\twoline{\sigma_{2}+}{\sigma_{23}\sigma_2}} "a2xi2j";"d3ka1"};
 {\ar@/_1.5pc/|(0.45){\sigma_{3}\sigma_{2}} "a3xi2j";"d3ka1"};
 {\ar@/_1.5pc/|(0.53){\sigma_{123}} "b1xi2j";"d2ka1"};
 {\ar@/_1.5pc/|{\sigma_{23}} "b2xi2j";"d2ka1"};
 {\ar@/_1.5pc/|(0.3){\sigma_{3}} "b3xi2j";"d2ka1"};
 {\ar@/_1pc/|(0.55){\twoline{\sigma_{123}+}{\sigma_3\sigma_2\sigma_1}} "b1xi2j";"d4ka1"};
 {\ar@/_0.5pc/|(0.4){\sigma_{23}} "b2xi2j";"d4ka1"};
 {\ar|(0.3){\sigma_{3}} "b3xi2j";"d4ka1"};
 {\ar@/_0.3pc/ "b3xi2j";"d3ka1"};
%internal to kappa_1
 {\ar|{\sigma_1} "d1ka1";"d2ka1"};
 {\ar|{\sigma_1} "e1ka1";"e2ka1"};
 {\ar|{\sigma_2} "d2ka1";"d3ka1"};
 {\ar|{\sigma_2} "e2ka1";"e3ka1"};
 {\ar@/^1.5pc/|(0.4){\sigma_{12}} "d1ka1";"d3ka1"};
 {\ar@/^1.5pc/|(0.4){\sigma_{12}} "e1ka1";"e3ka1"};
 {\ar "e1ka1";"d1ka1"};
 {\ar "e2ka1";"d2ka1"};
 {\ar "e3ka1";"d3ka1"};
 {\ar "d4ka1";"c1ka1"};
 {\ar|{\sigma_{123}} "d1ka1";"c1ka1"};
 {\ar|(0.55){\sigma_{23}} "d2ka1";"c1ka1"};
 {\ar|(0.6){\sigma_{3}} "d3ka1";"c1ka1"};
 {\ar|{\sigma_{123}} "e1ka1";"d4ka1"};
 {\ar|(0.35){\sigma_{23}} "e2ka1";"d4ka1"};
 {\ar|{\sigma_{3}} "e3ka1";"d4ka1"};
%internal to kappa_{k-1}
 {\ar|{\sigma_1} "d1kak1";"d2kak1"};
 {\ar|{\sigma_1} "e1kak1";"e2kak1"};
 {\ar|{\sigma_2} "d2kak1";"d3kak1"};
 {\ar|{\sigma_2} "e2kak1";"e3kak1"};
 {\ar@/^1.5pc/|(0.4){\sigma_{12}} "d1kak1";"d3kak1"};
 {\ar@/^1.5pc/|(0.4){\sigma_{12}} "e1kak1";"e3kak1"};
 {\ar "e1kak1";"d1kak1"};
 {\ar "e2kak1";"d2kak1"};
 {\ar "e3kak1";"d3kak1"};
 {\ar "d4kak1";"c1kak1"};
 {\ar|(0.6){\sigma_{123}} "d1kak1";"c1kak1"};
 {\ar|(0.6){\sigma_{23}} "d2kak1";"c1kak1"};
 {\ar|(0.6){\sigma_{3}} "d3kak1";"c1kak1"};
 {\ar|(0.6){\sigma_{123}} "e1kak1";"d4kak1"};
 {\ar|(0.6){\sigma_{23}} "e2kak1";"d4kak1"};
 {\ar|(0.6){\sigma_{3}} "e3kak1";"d4kak1"};
%kappa_1 - - rho_{23} - -> kappa_{k-1}
 {\ar@{-->}@/_2.5pc/|(0.5){\twoline{\sigma_{123}+}{\sigma_3\sigma_2\sigma_1}} "d1ka1";"c1kak1"};
 {\ar@{-->}@/_2pc/|(0.55){\sigma_{23}} "d2ka1";"c1kak1"};
 {\ar@{-->}|{\sigma_{23}} "d4ka1";"c1kak1"};
 {\ar@{-->}@/_3.25pc/|(0.35){\sigma_{12}} "d1ka1";"d3kak1"};
 {\ar@{-->}@/_2.5pc/|(0.39){\sigma_{2}} "d2ka1";"d3kak1"};
 {\ar@{-->}@/_1.75pc/|(0.45){\sigma_{2}} "d4ka1";"d3kak1"};
 {\ar@{-->}@/_1.5pc/|(0.53){\sigma_{123}} "e1ka1";"d2kak1"};
 {\ar@{-->}@/_1.5pc/|{\sigma_{23}} "e2ka1";"d2kak1"};
 {\ar@{-->}@/_1.5pc/|(0.3){\sigma_{3}} "e3ka1";"d2kak1"};
 {\ar@{-->}@/_1pc/|(0.55){\twoline{\sigma_{123}+}{\sigma_3\sigma_2\sigma_1}} "e1ka1";"d4kak1"};
 {\ar@{-->}@/_0.5pc/|(0.4){\sigma_{23}} "e2ka1";"d4kak1"};
 {\ar@{-->}|(0.3){\sigma_{3}} "e3ka1";"d4kak1"};
 {\ar@{-->}@/_0.3pc/ "e3ka1";"d3kak1"};
%\kappa_{k-1} --> \kappa_k
 {\ar@/^1pc/|{\twoline{\sigma_{123}+}{\sigma_3\sigma_2\sigma_1}} "d1kak1";"c1kak"};
 {\ar@/^1.25pc/|(0.55){\sigma_{23}} "d2kak1";"c1kak"};
 {\ar|{\sigma_{23}} "d4kak1";"c1kak"};
 {\ar|(0.62){\threeline{\sigma_{123}+} {\sigma_{123}\sigma_{23}+} {\sigma_3\sigma_2\sigma_1}} "e1kak1";"d4kak"};
 {\ar@/^.6pc/|(0.4){\twoline{\sigma_{23}+} {\sigma_{23}\sigma_{23}}} "e2kak1";"d4kak"};
 {\ar@/^1pc/|(0.3){\sigma_3\sigma_{23}} "e3kak1";"d4kak"};
% {\ar@/_3pc/|(0.35){\sigma_{12} \sigma_3} "d1kak1";"d4kak"};
% {\ar@/_2pc/|(0.39){\sigma_2 \sigma_3} "d2kak1";"d4kak"};
% {\ar|(0.45){\sigma_2 \sigma_3} "d4kak1";"d4kak"};
%internal to \kappa_k
 {\ar "d4kak";"c1kak"};
\endxy
\]
\end{tiny}
\caption{The subspace $Q^j$ $(j>0)$ obtained from $P_{\ver}^j$ by canceling
edges.} \label{fig:Qj}
\end{figure}

Let $Q^j$ denote the module resulting from $P^j$ after the cancellations just
described. The multiplications on $Q^j$ are shown in Figures \ref{fig:Qj} and
\ref{fig:Q0,s<2tau} and equations \eqref{eq:Pver-longer} and
\eqref{eq:Punst-longer}.

\begin{figure}
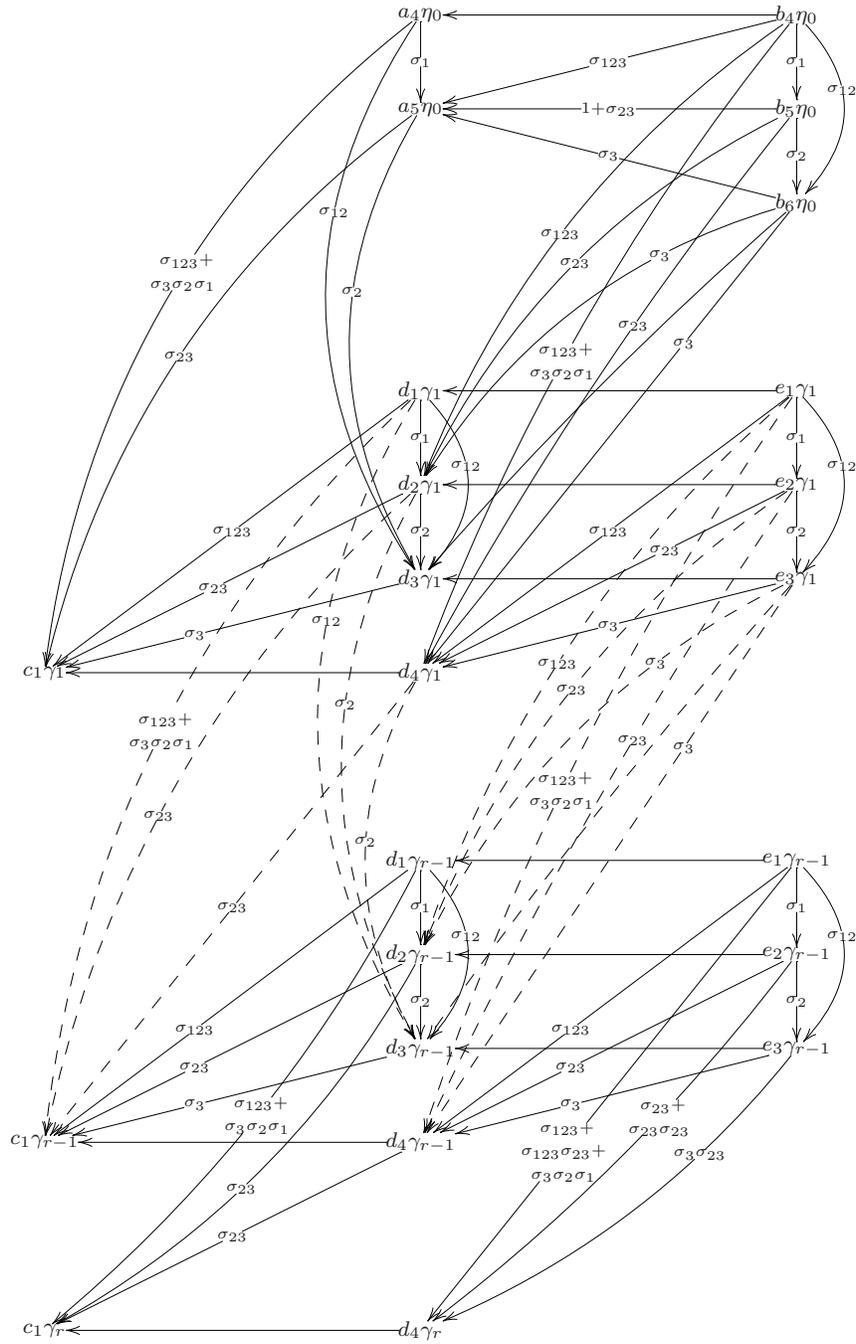
 \centering
\begin{tiny}
\[
\xy
 (50,0)*{a_4 \eta_0}="a4eta0"; (100,0)*{b_4 \eta_0}="b4eta0";
 (50,-12.5)*{a_5 \eta_0}="a5eta0"; (100,-12.5)*{b_5 \eta_0}="b5eta0";
 (100,-25)*{b_6 \eta_0}="b6eta0";
 (50,-50)*{d_1 \gamma_1}="d1ga1";
 (100,-50)*{e_1 \gamma_1}="e1ga1";
 (50,-62.5)*{d_2 \gamma_1}="d2ga1";
 (100,-62.5)*{e_2 \gamma_1}="e2ga1";
 (50,-75)*{d_3 \gamma_1}="d3ga1";
 (100,-75)*{e_3 \gamma_1}="e3ga1";
 (0,-87.5)*{c_1 \gamma_1}="c1ga1";
 (50,-87.5)*{d_4 \gamma_1}="d4ga1";
 (50,-112.5)*{d_1 \gamma_{r-1}}="d1gar1";
 (100,-112.5)*{e_1 \gamma_{r-1}}="e1gar1";
 (50,-125)*{d_2 \gamma_{r-1}}="d2gar1";
 (100,-125)*{e_2 \gamma_{r-1}}="e2gar1";
 (50,-137.5)*{d_3 \gamma_{r-1}}="d3gar1";
 (100,-137.5)*{e_3 \gamma_{r-1}}="e3gar1";
 (0,-150)*{c_1 \gamma_{r-1}}="c1gar1";
 (50,-150)*{d_4 \gamma_{r-1}}="d4gar1";
 (0,-175)*{c_1 \gamma_r}="c1gar";
 (50,-175)*{d_4 \gamma_r}="d4gar";
%internal to eta_0
 {\ar|{\sigma_1} "a4eta0";"a5eta0"};
 {\ar|{\sigma_1} "b4eta0";"b5eta0"};
 {\ar|{\sigma_2} "b5eta0";"b6eta0"};
 {\ar "b4eta0";"a4eta0"};
 {\ar|{\sigma_{123}} "b4eta0";"a5eta0"};
 {\ar|{1+\sigma_{23}} "b5eta0";"a5eta0"};
 {\ar|{\sigma_3} "b6eta0";"a5eta0"};
 {\ar@/^1.5pc/|(0.4){\sigma_{12}} "b4eta0";"b6eta0"};
%eta_0 -rho_{3}-> gamma_1
 {\ar@/_2.5pc/|(0.45){\twoline{\sigma_{123}+}{\sigma_3\sigma_2\sigma_1}} "a4eta0";"c1ga1"};
 {\ar@/_2pc/|(0.5){\sigma_{23}} "a5eta0";"c1ga1"};
 {\ar@/_3pc/|(0.35){\sigma_{12}} "a4eta0";"d3ga1"};
 {\ar@/_2.25pc/|(0.39){\sigma_{2}} "a5eta0";"d3ga1"};
 {\ar@/_1.5pc/|(0.53){\sigma_{123}} "b4eta0";"d2ga1"};
 {\ar@/_1.5pc/|{\sigma_{23}} "b5eta0";"d2ga1"};
 {\ar@/_1.5pc/|(0.3){\sigma_{3}} "b6eta0";"d2ga1"};
 {\ar@/_1pc/|(0.55){\twoline{\sigma_{123}+}{\sigma_3\sigma_2\sigma_1}} "b4eta0";"d4ga1"};
 {\ar@/_0.5pc/|(0.4){\sigma_{23}} "b5eta0";"d4ga1"};
 {\ar|(0.3){\sigma_{3}} "b6eta0";"d4ga1"};
 {\ar@/_0.3pc/ "b6eta0";"d3ga1"};
%internal to gamma_1
 {\ar|{\sigma_1} "d1ga1";"d2ga1"};
 {\ar|{\sigma_1} "e1ga1";"e2ga1"};
 {\ar|{\sigma_2} "d2ga1";"d3ga1"};
 {\ar|{\sigma_2} "e2ga1";"e3ga1"};
 {\ar@/^1.5pc/|(0.4){\sigma_{12}} "d1ga1";"d3ga1"};
 {\ar@/^1.5pc/|(0.4){\sigma_{12}} "e1ga1";"e3ga1"};
 {\ar "e1ga1";"d1ga1"};
 {\ar "e2ga1";"d2ga1"};
 {\ar "e3ga1";"d3ga1"};
 {\ar "d4ga1";"c1ga1"};
 {\ar|{\sigma_{123}} "d1ga1";"c1ga1"};
 {\ar|(0.55){\sigma_{23}} "d2ga1";"c1ga1"};
 {\ar|(0.6){\sigma_{3}} "d3ga1";"c1ga1"};
 {\ar|{\sigma_{123}} "e1ga1";"d4ga1"};
 {\ar|(0.35){\sigma_{23}} "e2ga1";"d4ga1"};
 {\ar|{\sigma_{3}} "e3ga1";"d4ga1"};
%internal to gamma_{r-1}
 {\ar|{\sigma_1} "d1gar1";"d2gar1"};
 {\ar|{\sigma_1} "e1gar1";"e2gar1"};
 {\ar|{\sigma_2} "d2gar1";"d3gar1"};
 {\ar|{\sigma_2} "e2gar1";"e3gar1"};
 {\ar@/^1.5pc/|(0.4){\sigma_{12}} "d1gar1";"d3gar1"};
 {\ar@/^1.5pc/|(0.4){\sigma_{12}} "e1gar1";"e3gar1"};
 {\ar "e1gar1";"d1gar1"};
 {\ar "e2gar1";"d2gar1"};
 {\ar "e3gar1";"d3gar1"};
 {\ar "d4gar1";"c1gar1"};
 {\ar|(0.6){\sigma_{123}} "d1gar1";"c1gar1"};
 {\ar|(0.6){\sigma_{23}} "d2gar1";"c1gar1"};
 {\ar|(0.6){\sigma_{3}} "d3gar1";"c1gar1"};
 {\ar|(0.6){\sigma_{123}} "e1gar1";"d4gar1"};
 {\ar|(0.6){\sigma_{23}} "e2gar1";"d4gar1"};
 {\ar|(0.6){\sigma_{3}} "e3gar1";"d4gar1"};
%gamma_1 - - rho_{23} - -> gamma_{r-1}
 {\ar@{-->}@/_2.5pc/|(0.5){\twoline{\sigma_{123}+}{\sigma_3\sigma_2\sigma_1}} "d1ga1";"c1gar1"};
 {\ar@{-->}@/_2pc/|(0.55){\sigma_{23}} "d2ga1";"c1gar1"};
 {\ar@{-->}|{\sigma_{23}} "d4ga1";"c1gar1"};
 {\ar@{-->}@/_3.25pc/|(0.35){\sigma_{12}} "d1ga1";"d3gar1"};
 {\ar@{-->}@/_2.5pc/|(0.39){\sigma_{2}} "d2ga1";"d3gar1"};
 {\ar@{-->}@/_1.75pc/|(0.45){\sigma_{2}} "d4ga1";"d3gar1"};
 {\ar@{-->}@/_1.5pc/|(0.53){\sigma_{123}} "e1ga1";"d2gar1"};
 {\ar@{-->}@/_1.5pc/|{\sigma_{23}} "e2ga1";"d2gar1"};
 {\ar@{-->}@/_1.5pc/|(0.3){\sigma_{3}} "e3ga1";"d2gar1"};
 {\ar@{-->}@/_1pc/|(0.55){\twoline{\sigma_{123}+}{\sigma_3\sigma_2\sigma_1}} "e1ga1";"d4gar1"};
 {\ar@{-->}@/_0.5pc/|(0.4){\sigma_{23}} "e2ga1";"d4gar1"};
 {\ar@{-->}|(0.3){\sigma_{3}} "e3ga1";"d4gar1"};
 {\ar@{-->}@/_0.3pc/ "e3ga1";"d3gar1"};
%\gamma_{r-1} --> \gamma_r
 {\ar@/^1pc/|{\twoline{\sigma_{123}+}{\sigma_3\sigma_2\sigma_1}} "d1gar1";"c1gar"};
 {\ar@/^1.25pc/|(0.55){\sigma_{23}} "d2gar1";"c1gar"};
 {\ar|{\sigma_{23}} "d4gar1";"c1gar"};
 {\ar|(0.62){\threeline{\sigma_{123}+} {\sigma_{123}\sigma_{23}+} {\sigma_3\sigma_2\sigma_1}} "e1gar1";"d4gar"};
 {\ar@/^.6pc/|(0.4){\twoline{\sigma_{23}+} {\sigma_{23}\sigma_{23}}} "e2gar1";"d4gar"};
 {\ar@/^1pc/|(0.3){\sigma_3\sigma_{23}} "e3gar1";"d4gar"};
%internal to \gamma_4
 {\ar "d4gar";"c1gar"};
\endxy
\]
\end{tiny}
\caption{The subspace $Q^0$ obtained from $P_{\unst}$ by canceling edges, when
$s<2\tau(J)$.} \label{fig:Q0,s<2tau}
\end{figure}

Now we keep track of what these cancellations do to the edges that do not
respect the splitting, as shown in Figures \ref{fig:rho3rho2rho123} and
\ref{fig:rho3rho2rho1}.

If $u_{p,j}=1$, then there are edges from $b_i \xi_p$ to $d_3 \kappa^j_1$, as
shown in Figure \ref{fig:rho3rho2rho123}. If $k_j=1$, then canceling $a_3
\xi_{2j-1} \to d_3 \kappa^j_1$ will introduce new multiplications coming from
$b_i \xi_p$, but all of these multiplications involve $\sigma_2\sigma_3$ or
$\sigma_{12}\sigma_3$ and may thus be disregarded. Also, when $p=2m+1$ or $p=0$
these edges are eliminated when we cancel $b_i \xi_{2m+1} \to e_i
\kappa^m_{k_m}$ or $b_i \xi_0 \to e_i \gamma_r$, respectively.

\begin{figure}
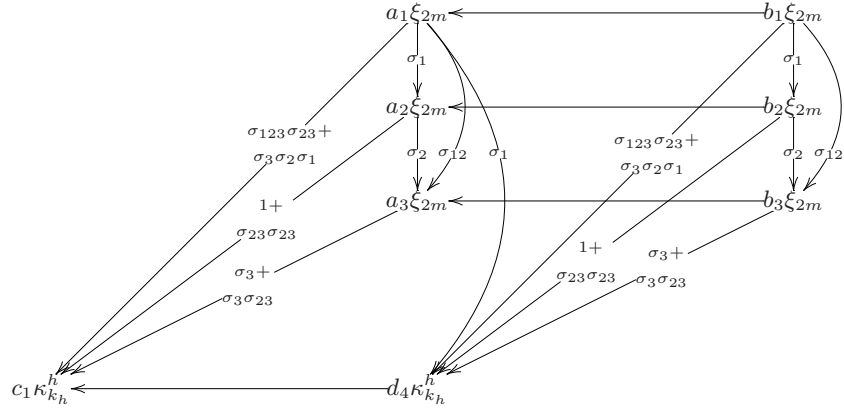
 \centering
\begin{scriptsize}
\[
\xy
 (50,0)*{a_1 \xi_{2m}}="a1xi2m"; (100,0)*{b_1 \xi_{2m}}="b1xi2m";
 (50,-12.5)*{a_2 \xi_{2m}}="a2xi2m"; (100,-12.5)*{b_2 \xi_{2m}}="b2xi2m";
 (50,-25)*{a_3 \xi_{2m}}="a3xi2m"; (100,-25)*{b_3 \xi_{2m}}="b3xi2m";
 (0,-50)*{c_1 \kappa^h_{k_h}}="c1kak";
 (50,-50)*{d_4 \kappa^h_{k_h}}="d4kak";
%internal to xi_2m
 {\ar|{\sigma_1} "a1xi2m";"a2xi2m"};
 {\ar|{\sigma_1} "b1xi2m";"b2xi2m"};
 {\ar|{\sigma_2} "a2xi2m";"a3xi2m"};
 {\ar|{\sigma_2} "b2xi2m";"b3xi2m"};
 {\ar "b1xi2m";"a1xi2m"};
 {\ar "b2xi2m";"a2xi2m"};
 {\ar "b3xi2m";"a3xi2m"};
 {\ar@/^1.5pc/|(.75){\sigma_{12}} "a1xi2m";"a3xi2m"};
 {\ar@/^1.5pc/|(.75){\sigma_{12}} "b1xi2m";"b3xi2m"};
%internal to kappa_1
 {\ar "d4kak";"c1kak"};
%xi_{2m} -rho_3 rho_2 rho_1 -> kappa_1
 {\ar|(.35){\twoline{\sigma_{123}\sigma_{23} +}{\sigma_3 \sigma_2 \sigma_1}} "a1xi2m";"c1kak"};
 {\ar|(.4){\twoline{1 +}{\sigma_{23} \sigma_{23}}} "a2xi2m";"c1kak"};
 {\ar|(.45){\twoline{\sigma_3+}{\sigma_3 \sigma_{23}}} "a3xi2m";"c1kak"};
 {\ar@/^2.75pc/|(.375){\sigma_1} "a1xi2m";"d4kak"};
 {\ar|(.375){\twoline{\sigma_{123}\sigma_{23} +}{\sigma_3 \sigma_2 \sigma_1}} "b1xi2m";"d4kak"};
 {\ar|(.55){\twoline{1+}{\sigma_{23} \sigma_{23}}}  "b2xi2m";"d4kak"};
 {\ar|(.35){\twoline{\sigma_3+}{\sigma_3 \sigma_{23}}}  "b3xi2m";"d4kak"};\endxy
\]
\end{scriptsize}
\caption{Reduced form of Figure \ref{fig:rho3rho2rho1} when $p=2m$, $m>0$.}
\label{fig:rho3rho2rho1red-peven}
\end{figure}

\begin{figure}
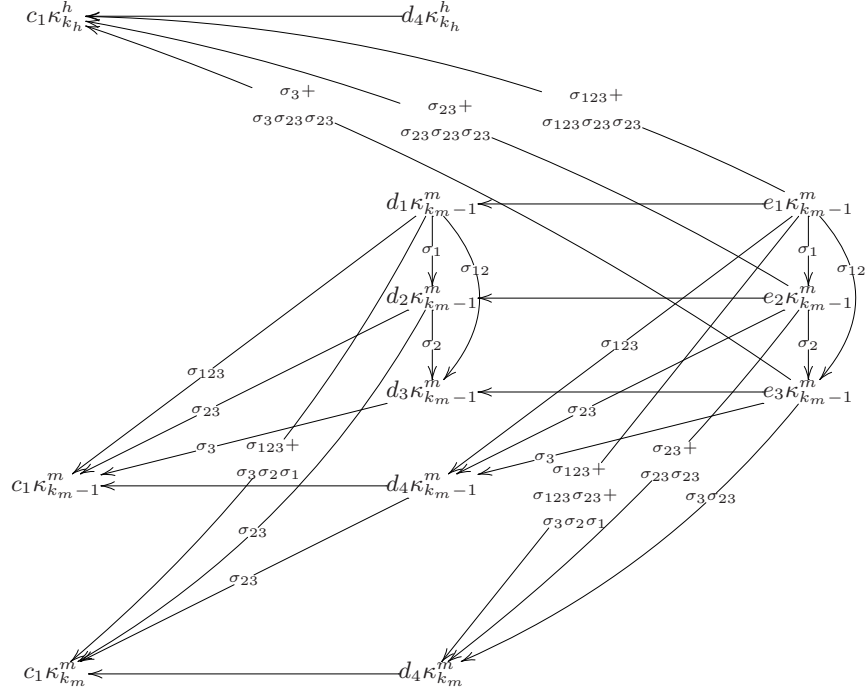
 \centering
\begin{scriptsize}
\[
\xy
 (50,0)*{d_1 \kappa^m_{k_m-1}}="d1kakm1";
 (100,0)*{e_1\kappa^m_{k_m-1}}="e1kakm1";
 (50,-12.5)*{d_2 \kappa^m_{k_m-1}}="d2kakm1";
 (100,-12.5)*{e_2 \kappa^m_{k_m-1}}="e2kakm1";
 (50,-25)*{d_3 \kappa^m_{k_m-1}}="d3kakm1";
 (100,-25)*{e_3 \kappa^m_{k_m-1}}="e3kakm1";
 (0,-37.5)*{c_1 \kappa^m_{k_m-1}}="c1kakm1";
 (50,-37.5)*{d_4 \kappa^m_{k_m-1}}="d4kakm1";
 (0,-62.5)*{c_1 \kappa^m_{k_m}}="c1kakm";
 (50,-62.5)*{d_4 \kappa^m_{k_m}}="d4kakm";
 (0,25)*{c_1 \kappa^h_{k_h}}="c1kakh";
 (50,25)*{d_4 \kappa^h_{k_h}}="d4kakh";
%internal to kappa_{k_m-1}
 {\ar|{\sigma_1} "d1kakm1";"d2kakm1"};
 {\ar|{\sigma_1} "e1kakm1";"e2kakm1"};
 {\ar|{\sigma_2} "d2kakm1";"d3kakm1"};
 {\ar|{\sigma_2} "e2kakm1";"e3kakm1"};
 {\ar@/^1.5pc/|(0.35){\sigma_{12}} "d1kakm1";"d3kakm1"};
 {\ar@/^1.5pc/|(0.35){\sigma_{12}} "e1kakm1";"e3kakm1"};
 {\ar "e1kakm1";"d1kakm1"};
 {\ar "e2kakm1";"d2kakm1"};
 {\ar "e3kakm1";"d3kakm1"};
 {\ar "d4kakm1";"c1kakm1"};
 {\ar|(0.6){\sigma_{123}} "d1kakm1";"c1kakm1"};
 {\ar|(0.6){\sigma_{23}} "d2kakm1";"c1kakm1"};
 {\ar|(0.6){\sigma_{3}} "d3kakm1";"c1kakm1"};
 {\ar|{\sigma_{123}} "e1kakm1";"d4kakm1"};
 {\ar|(0.6){\sigma_{23}} "e2kakm1";"d4kakm1"};
 {\ar|(0.7){\sigma_{3}} "e3kakm1";"d4kakm1"};
%\kappa_{k_m-1} --> \kappa_{k_m}
 {\ar@/^1pc/|{\twoline{\sigma_{123}+}{\sigma_3\sigma_2\sigma_1}} "d1kakm1";"c1kakm"};
 {\ar@/^1.25pc/|(0.55){\sigma_{23}} "d2kakm1";"c1kakm"};
 {\ar|{\sigma_{23}} "d4kakm1";"c1kakm"};
 {\ar|(0.62){\threeline{\sigma_{123}+} {\sigma_{123}\sigma_{23}+} {\sigma_3\sigma_2\sigma_1}} "e1kakm1";"d4kakm"};
 {\ar@/^.6pc/|(0.4){\twoline{\sigma_{23}+} {\sigma_{23}\sigma_{23}}} "e2kakm1";"d4kakm"};
 {\ar@/^1pc/|(0.3){\sigma_3\sigma_{23}} "e3kakm1";"d4kakm"};
%internal to \kappa_{k_m}
 {\ar "d4kakm";"c1kakm"};
%internal to \kappa_{k_h}
 {\ar "d4kakh";"c1kakh"};
%jumping
 {\ar@/_1.5pc/|(.3){\twoline{\sigma_{123} +}{\sigma_{123} \sigma_{23} \sigma_{23}}} "e1kakm1";"c1kakh"};
 {\ar@/_1.25pc/|(.5){\twoline{\sigma_{23} +}{\sigma_{23} \sigma_{23} \sigma_{23}}} "e2kakm1";"c1kakh"};
 {\ar@/_1pc/|(.7){\twoline{\sigma_{3} +}{\sigma_{3} \sigma_{23} \sigma_{23}}} "e3kakm1";"c1kakh"};
\endxy
\]
\end{scriptsize}
\caption[Reduced form of Figure \ref{fig:rho3rho2rho1} in the case where
$p=2m+1$ or $p=0$.] {Reduced form of Figure \ref{fig:rho3rho2rho1} in the case
where $p=2m+1$ (or $p=0$, replacing $\kappa^m_{k_m-1}$ by $\gamma_{r-1}$ and
$\kappa^m_{k_m}$ by $\gamma_r$).} \label{fig:rho3rho2rho1red-podd}
\end{figure}

If $v_{p,j}=1$, when we cancel $a_2 \xi_{2j-1} \to d_2 \kappa^j_{k_j}$, we
obtain multiplications
\begin{equation} \label{eq:cancel3}
\begin{alignedat}{2}
b_1 \xi_p & \xrightarrow{\sigma_{123} \sigma_{2}} a_3 \xi_{2j-1} &\qquad
 b_1 \xi_p & \xrightarrow{\sigma_{123} \sigma_{23}} d_4 \kappa^j_{k_j} \\
b_2 \xi_p & \xrightarrow{\sigma_{23} \sigma_{2}} a_3 \xi_{2j-1} &
 b_2 \xi_p & \xrightarrow{\sigma_{23} \sigma_{23}} d_4 \kappa^j_{k_j} \\
b_3 \xi_p & \xrightarrow{\sigma_3 \sigma_{23}} d_4 \kappa^j_{k_j} &
 b_3 \xi_p & \xrightarrow{\sigma_3 \sigma_{2}} a_3 \xi_{2j-1}
\end{alignedat}
\end{equation}
in addition to the ones already appearing in Figure \ref{fig:rho3rho2rho1}.
When we then cancel $a_3 \xi_{2j-1} \to d_3 \kappa^j_{k_j}$, we obtain new
multiplications:
\begin{equation} \label{eq:cancel4}
\begin{split}
a_1\xi_p & \xrightarrow{\sigma_{12}\sigma_3 + \sigma_{123}\sigma_2 \sigma_3}
d_4 \kappa^j_{k_j} \\
a_2\xi_p & \xrightarrow{\sigma_2\sigma_3 + \sigma_{23}\sigma_2 \sigma_3}
d_4 \kappa^j_{k_j} \\
a_3\xi_p & \xrightarrow{\sigma_3 \sigma_2 \sigma_3} d_4 \kappa^j_{k_j} \\
%b_1 \xi_p & \xrightarrow{\sigma_{123} \sigma_{2}} d_4 \xi_{2j-1} \\
%b_2 \xi_p & \xrightarrow{\sigma_{23} \sigma_{2}} d_4 \xi_{2j-1} \\
%b_3 \xi_p & \xrightarrow{\sigma_3 \sigma_{2}} d_4 \xi_{2j-1} \\
%These appear to have been a mistake, but I'll leave them commented in the text just to be careful.
b_3 \xi_p & \xrightarrow{\sigma_3} d_4 \kappa^j_{k_j} \\
\end{split}
\end{equation}
Most of these may be disregarded by Proposition \ref{prop:norho1rho2}. If
$p=2m$ for $m>0$, the resulting reduced form of Figure \ref{fig:rho3rho2rho1}
is shown in Figure \ref{fig:rho3rho2rho1red-peven}. On the other hand, if
$p=2m+1$, we also cancel the edges $a_i \xi_{2m+1} \to d_i \kappa^m_{k_m}$ and
$b_i \xi_{2m+1} \to e_i \kappa^m_{k_m}$, introducing the multiplications shown
in Figure \ref{fig:rho3rho2rho1red-podd}. Similarly, if $p=0$, we cancel the
edges $a_i \xi_{0} \to d_i \gamma_{r}$ and $b_i \xi_{0} \to e_i \gamma_{r}$,
introducing similar multiplications.

We now return to the case where $s=2\tau(J)$. In $P_{\unst}$, the edges $a_1
\xi_0 \to a_4 \eta_0$, $b_1 \xi_0 \to b_4 \eta_0$, $b_2 \xi_0 \to b_5 \eta_0$,
and $b_3 \xi_0 \to b_6 \eta_0$ cancel, and since their targets do not have any
other incoming edges, no new multiplications are introduced. The only three
remaining generators are $a_2 \xi_0$, $a_3 \xi_0$, and $a_5 \eta_0$, all in
filtration level $0$, with the following multiplications:
\begin{equation} \label{eq:Q0,s=2tau}
\xy
 (0,-7.5)*{a_2 \xi_0}="a2xi0";
 (15,7.5)*{a_3 \xi_0}="a3xi0";
 (30,-7.5)*{a_5 \eta_0}="a5eta0";
 {\ar|{\sigma_2} "a2xi0";"a3xi0"};
 {\ar|{1+\sigma_{23}} "a2xi0";"a5eta0"};
 {\ar|{\sigma_3} "a3xi0";"a5eta0"};
\endxy
\end{equation}
As above, $a_2 \xi_0$ and $a_3 \xi_0$ may have some outgoing edges, and $a_5
\eta_0$ may have some incoming ones. The rest of the argument goes through
unchanged.

\subsection{Tensor product over \texorpdfstring{$\AA_\sigma$}{A\textunderscore\textsigma}} \label{subsec:Asigma}

Let $Q = \bigoplus_{j=0}^n Q^j$, with multiplications as described in the
previous subsection. We consider the tensor product $Q
\underset{\AA_\sigma}{\boxtimes} \CFD(\XX_K^t)$. Again, the goal is to obtain a
decomposition of the tensor product according to the stable and unstable chains
in $\CFD(\XX_K^t)$.

It is convenient to give the generators of $Q^j$ new names. For $j=1, \dots, n$
and $i=1, \dots, k_j-1$, define:
\[
\begin{aligned}
A^j &= a_1 \xi_{2j} & A'^j &= b_1 \xi_{2j} &  E^j_i &= d_1 \kappa^j_i & E'^j_i &= e_1 \kappa^j_i \\
B^j &= a_2 \xi_{2j} & B'^j &= b_2 \xi_{2j} &  F^j_i &= d_2 \kappa^j_i &  F'^j_i &= e_2 \kappa^j_i \\
C^j &= a_3 \xi_{2j} & C'^j &= b_3 \xi_{2j} &  G^j_i &= d_3 \kappa^j_i & G'^j_i &= e_3 \kappa^j_i \\
D^j &= c_1 \kappa^j_{k_j} & D'^j &= d_4 \kappa^j_{k_j} & H^j_i &= c_1 \kappa^j_i & H'^j_i &= d_4 \kappa^j_i \\
\end{aligned}
\]
When $s<2\tau(J)$, for $i=1, \dots, r-1$, define:
\[
\begin{aligned}
A^0 &= a_4 \eta_0 & A'^0 &= b_4 \eta_0  & E^0_i &= d_1 \gamma_i & E'^0_i &= e_1 \gamma_i \\
B^0 &= a_5 \eta_0 & B'^0 &= b_5 \eta_0  & F^0_i &= d_2 \gamma_i &  F'^0_i &= e_2 \gamma_i \\
&& C'^0 &= b_6 \eta_0 & G^0_i &= d_3 \gamma_i & G'^0_i &= e_3 \gamma_i \\
D^0 &= c_1 \gamma_r & D'^0 &= d_4 \gamma_r  & H^0_i &= c_1 \gamma_i & H'^0_i &= d_4 \gamma_i \\
\end{aligned}
\]
Also, for notational convenience, define $k_0=r$.

We divide up the generators of the subspaces $Q^j$ by Alexander grading and
idempotent:

\begin{center}
\begin{tabular}{|c|c|c|c|} \hline
 & $A=-1$ & $A=0$ & $A=1$ \\ \hline
 $\iota_0^\sigma$ &  & $A^j, C^j, E^j_i, G^j_i$ & $A'^j, C'^j, E'^j_i, G'^j_i$ \\ \hline
 $\iota_1^\sigma$ & $D^j, H^j_i$ & $B^j, D'^j, F^j_i, H'^j_i$ & $B'^j, F'^j_i$ \\ \hline
\end{tabular}
\end{center}

In Figures \ref{fig:Qj} and \ref{fig:Q0,s<2tau}, notice that of the generators
in idempotent $\iota_0$, $A^j$, $A'^j$, $E^j_i$, and $E'^j_i$ have outgoing
edges labeled $\sigma_1$, $\sigma_{12}$, and $\sigma_{123}$, while $C^j$,
$C'^j$, $G^j_i$, and $G'^j_i$ have outgoing edges labeled $\sigma_3$ and
incoming edges labeled $\sigma_2$ and $\sigma_{12}$. Accordingly, it makes
sense to associate the former with the vertical chains and the latter with the
horizontal chains. That is, for each $J \in \{1, \dots, N\}$ and $j \in \{0,
\dots,n\}$, define:
\begin{equation} \label{eq:Zbases}
\begin{split} Z_{\ver}^{J,j} &= \gen{A^j, A'^j, E^j_i, E'^j_i}
\boxtimes \gen{
\Xi_{2J-1}, \Xi_{2J}} \\
& \quad + \gen{B^j, B'^j, D^j, D'^j, F^j_i, F'^j_i, H^j_i, H'^j_i} \boxtimes
\gen{\Kappa_I^J \mid 1 \le I \le K_J} \\
Z_{\hor}^{J,j} &= \gen{C^j, C'^j, G^j_i, G'^j_i} \boxtimes \gen{
\Eta_{2J-1}, \Eta_{2J} } \\
& \quad + \gen{B^j, B'^j, D^j, D'^j, F^j_i, F'^j_i, H^j_i, H'^j_i} \boxtimes
\gen{\Lambda_I^J \mid 1 \le I \le L_J} \\
Z_{\unst}^j &= \gen{A^j, A'^j, E^j_i, E'^j_i} \boxtimes \gen{\Xi_0} \\
& \quad + \gen{C^j, C'^j, G^j_i, G'^j_i} \boxtimes \gen{\Eta_0} \\
& \quad + \gen{B^j, B'^j, D^j, D'^j, F^j_i, F'^j_i, H^j_i, H'^j_i} \boxtimes
\gen{\Gamma_i \mid 1 \le I \le R}.
\end{split}
\end{equation}
Then, as a vector space,
\begin{equation} \label{eq:Zsplit}
Q \boxtimes \CFD(\XX_K^t) = \bigoplus_{\substack{ J=1,\dots,N \\ j=0,\dots,n}}
Z_{\ver}^{J,j} \, \oplus \, \bigoplus_{\substack{ J=1,\dots,N \\ j=0,\dots,n}}
Z_{\hor}^{J,j} \, \oplus \, \bigoplus_{j=0}^n Z_{\unst}^j.
\end{equation}
For fixed $J$, we write $Z_{\ver}^{J,*} = \bigoplus_{j=0}^n Z_{\ver}^{J,j}$,
and so on.

\begin{figure}
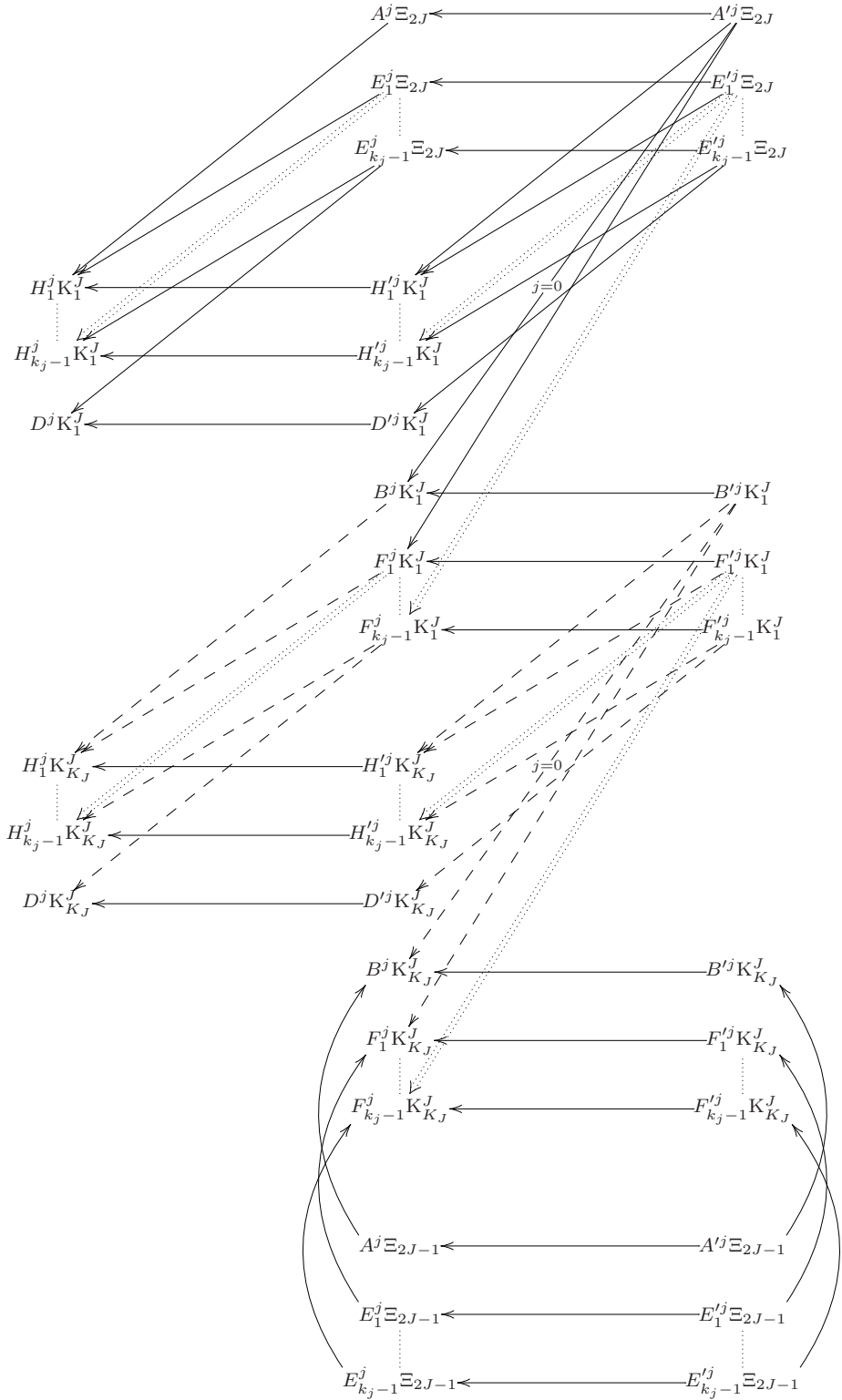
 \centering
\begin{tiny}
\[
\xy
 {(50,0)*{A^j \Xi_{2J}}="AXi2J"};
 {(100,0)*{A'^j \Xi_{2J}}="A'Xi2J"};
 {(50,-10)*{E^j_1 \Xi_{2J}}="E1Xi2J"};
 {(100,-10)*{E'^j_1 \Xi_{2J}}="E'1Xi2J"};
 {(50,-20)*{E^j_{k_j-1} \Xi_{2J}}="Ek1Xi2J"};
 {(100,-20)*{E'^j_{k_j-1} \Xi_{2J}}="E'k1Xi2J"};
 {(0,-40)*{H^j_1 \Kappa^J_1}="H1Ka1"};
 {(50,-40)*{H'^j_1 \Kappa^J_1}="H'1Ka1"};
 {(0,-50)*{H^j_{k_j-1} \Kappa^J_1}="Hk1Ka1"};
 {(50,-50)*{H'^j_{k_j-1} \Kappa^J_1}="H'k1Ka1"};
 {(0,-60)*{D^j \Kappa^J_1}="DKa1"};
 {(50,-60)*{D'^j \Kappa^J_1}="D'Ka1"};
 {(50,-70)*{B^j \Kappa^J_1}="BKa1"};
 {(100,-70)*{B'^j \Kappa^J_1}="B'Ka1"};
 {(50,-80)*{F^j_1 \Kappa^J_1}="F1Ka1"};
 {(100,-80)*{F'^j_1 \Kappa^J_1}="F'1Ka1"};
 {(50,-90)*{F^j_{k_j-1} \Kappa^J_1}="Fk1Ka1"};
 {(100,-90)*{F'^j_{k_j-1} \Kappa^J_1}="F'k1Ka1"};
 {(0,-110)*{H^j_1 \Kappa^J_{K_J}}="H1KaK"};
 {(50,-110)*{H'^j_1 \Kappa^J_{K_J}}="H'1KaK"};
 {(0,-120)*{H^j_{k_j-1} \Kappa^J_{K_J}}="Hk1KaK"};
 {(50,-120)*{H'^j_{k_j-1} \Kappa^J_{K_J}}="H'k1KaK"};
 {(0,-130)*{D^j \Kappa^J_{K_J}}="DKaK"};
 {(50,-130)*{D'^j \Kappa^J_{K_J}}="D'KaK"};
 {(50,-140)*{B^j \Kappa^J_{K_J}}="BKaK"};
 {(100,-140)*{B'^j \Kappa^J_{K_J}}="B'KaK"};
 {(50,-150)*{F^j_1 \Kappa^J_{K_J}}="F1KaK"};
 {(100,-150)*{F'^j_1 \Kappa^J_{K_J}}="F'1KaK"};
 {(50,-160)*{F^j_{k_j-1} \Kappa^J_{K_J}}="Fk1KaK"};
 {(100,-160)*{F'^j_{k_j-1} \Kappa^J_{K_J}}="F'k1KaK"};
 {(50,-180)*{A^j \Xi_{2J-1}}="AXi2J1"};
 {(100,-180)*{A'^j \Xi_{2J-1}}="A'Xi2J1"};
 {(50,-190)*{E^j_1 \Xi_{2J-1}}="E1Xi2J1"};
 {(100,-190)*{E'^j_1 \Xi_{2J-1}}="E'1Xi2J1"};
 {(50,-200)*{E^j_{k_j-1} \Xi_{2J-1}}="Ek1Xi2J1"};
 {(100,-200)*{E'^j_{k_j-1} \Xi_{2J-1}}="E'k1Xi2J1"};
%internal to Xi_2J
 {\ar "A'Xi2J"; "AXi2J"};
 {\ar "E'1Xi2J"; "E1Xi2J"};
 {\ar "E'k1Xi2J"; "Ek1Xi2J"};
 {\ar@{.} "E1Xi2J"; "Ek1Xi2J"};
 {\ar@{.} "E'1Xi2J"; "E'k1Xi2J"};
%internal to Kappa_1
 {\ar "B'Ka1"; "BKa1"};
 {\ar "F'1Ka1"; "F1Ka1"};
 {\ar "F'k1Ka1"; "Fk1Ka1"};
 {\ar "H'1Ka1"; "H1Ka1"};
 {\ar "H'k1Ka1"; "Hk1Ka1"};
 {\ar "D'Ka1"; "DKa1"};
 {\ar@{.} "F1Ka1"; "Fk1Ka1"};
 {\ar@{.} "F'1Ka1"; "F'k1Ka1"};
 {\ar@{.} "H1Ka1"; "Hk1Ka1"};
 {\ar@{.} "H'1Ka1"; "H'k1Ka1"};
%internal to Kappa_K
 {\ar "B'KaK"; "BKaK"};
 {\ar "F'1KaK"; "F1KaK"};
 {\ar "F'k1KaK"; "Fk1KaK"};
 {\ar "H'1KaK"; "H1KaK"};
 {\ar "H'k1KaK"; "Hk1KaK"};
 {\ar "D'KaK"; "DKaK"};
 {\ar@{.} "F1KaK"; "Fk1KaK"};
 {\ar@{.} "F'1KaK"; "F'k1KaK"};
 {\ar@{.} "H1KaK"; "Hk1KaK"};
 {\ar@{.} "H'1KaK"; "H'k1KaK"};
%internal to Xi_{2J-1}
 {\ar "A'Xi2J1"; "AXi2J1"};
 {\ar "E'1Xi2J1"; "E1Xi2J1"};
 {\ar "E'k1Xi2J1"; "Ek1Xi2J1"};
 {\ar@{.} "E1Xi2J1"; "Ek1Xi2J1"};
 {\ar@{.} "E'1Xi2J1"; "E'k1Xi2J1"};
%Xi_{2J} -sigma_{123}-> Kappa_1
 {\ar|(0.57){j=0} "A'Xi2J";"BKa1"};
 {\ar "A'Xi2J";"F1Ka1"};
 {\ar "A'Xi2J";"H'1Ka1"};
 {\ar "AXi2J";"H1Ka1"};
 {\ar "E'1Xi2J"; "H'1Ka1"};
 {\ar "E1Xi2J"; "H1Ka1"};
 {\ar "E'k1Xi2J"; "H'k1Ka1"};
 {\ar "Ek1Xi2J"; "Hk1Ka1"};
 {\ar@{:>} "E'1Xi2J"; "Fk1Ka1"};
 {\ar@{:>} "E'1Xi2J"; "H'k1Ka1"};
 {\ar@{:>} "E1Xi2J"; "Hk1Ka1"};
 {\ar "E'k1Xi2J"; "D'Ka1"};
 {\ar "Ek1Xi2J"; "DKa1"};
%Kappa_1 - -sigma_{23}- -> Kappa_K
 {\ar@{-->} "B'Ka1";"F1KaK"};
 {\ar@{-->} "B'Ka1";"H'1KaK"};
 {\ar@{-->} "BKa1";"H1KaK"};
 {\ar@{-->} "F'1Ka1"; "H'1KaK"};
 {\ar@{-->} "F1Ka1"; "H1KaK"};
 {\ar@{-->} "F'k1Ka1"; "H'k1KaK"};
 {\ar@{-->} "Fk1Ka1"; "Hk1KaK"};
 {\ar@{:>} "F'1Ka1"; "Fk1KaK"};
 {\ar@{:>} "F'1Ka1"; "H'k1KaK"};
 {\ar@{:>} "F1Ka1"; "Hk1KaK"};
 {\ar@{-->} "F'k1Ka1"; "D'KaK"};
 {\ar@{-->} "Fk1Ka1"; "DKaK"};
 {\ar@{-->}|(0.57){j=0} "B'Ka1"; "BKaK"};
%Xi_{2J-1} -rho_1-> Kappa_K
 {\ar@/^1.5pc/ "AXi2J1"+UL;"BKaK"+DL};
 {\ar@/_1.5pc/ "A'Xi2J1"+UR;"B'KaK"+DR};
 {\ar@/^1.5pc/ "E1Xi2J1"+UL;"F1KaK"+DL};
 {\ar@/_1.5pc/ "E'1Xi2J1"+UR;"F'1KaK"+DR};
 {\ar@/^1.5pc/ "Ek1Xi2J1"+UL;"Fk1KaK"+DL};
 {\ar@/_1.5pc/ "E'k1Xi2J1"+UR;"F'k1KaK"+DR};
\endxy
\]
\end{tiny}
\caption{The subspace $Z_{\ver}^{J,j}$, corresponding to a vertical stable
chain $\Xi_{2J} \xrightarrow{\sigma_{123}} \Kappa^J_1 \xrightarrow{\sigma_{23}}
\cdots \xrightarrow{\sigma_{23}} \Kappa^J_{K_J} \xleftarrow{\sigma_1}
\Xi_{2J-1}$.} \label{fig:Zver}
\end{figure}

\begin{figure}
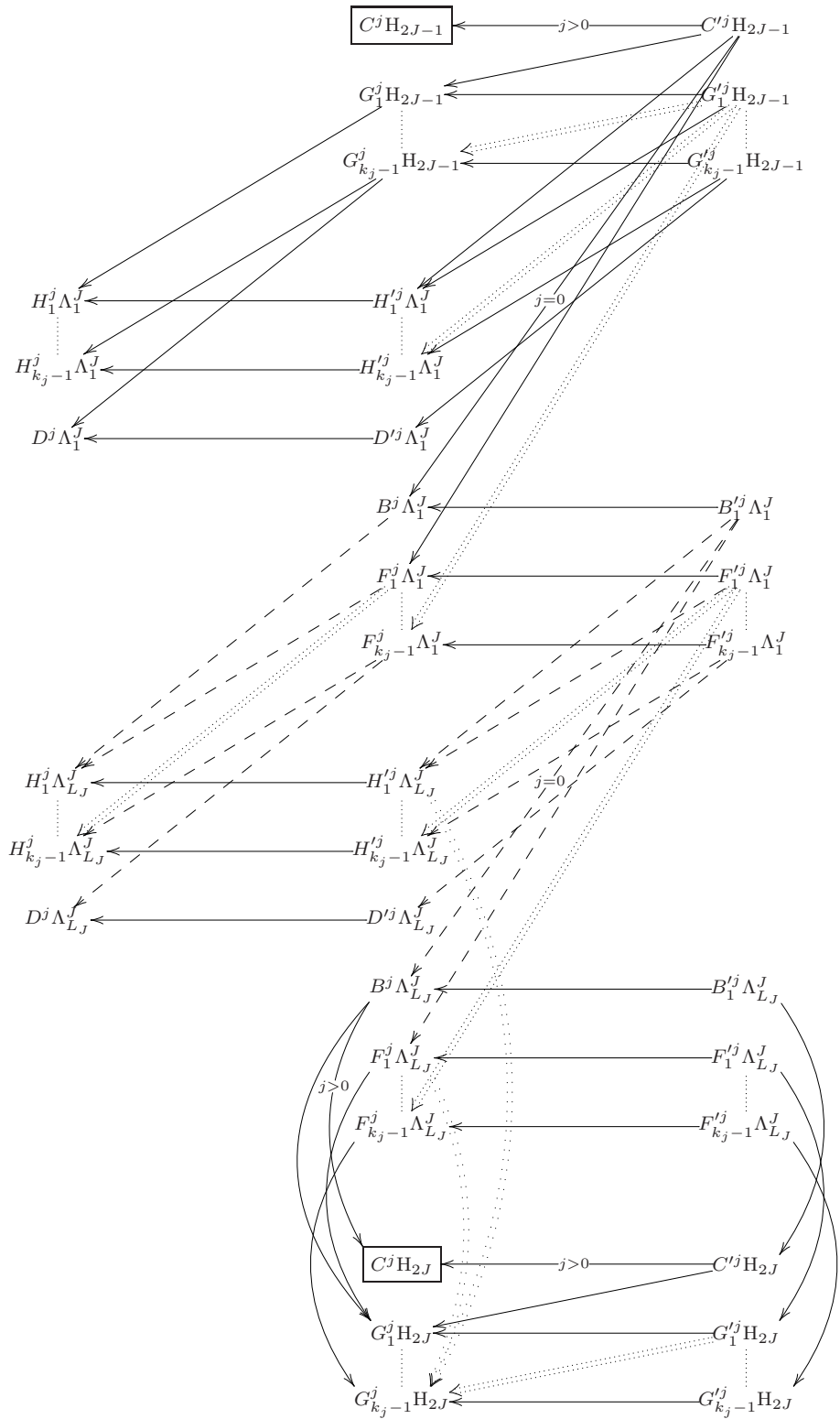
 \centering
\begin{tiny}
\[
\xy
 {(50,0)*{\fbox{$C^j \Eta_{2J-1}$}}="CEta2J1"};
 {(100,0)*{C'^j \Eta_{2J-1}}="C'Eta2J1"};
 {(50,-10)*{G^j_1 \Eta_{2J-1}}="G1Eta2J1"};
 {(100,-10)*{G'^j_1 \Eta_{2J-1}}="G'1Eta2J1"};
 {(50,-20)*{G^j_{k_j-1} \Eta_{2J-1}}="Gk1Eta2J1"};
 {(100,-20)*{G'^j_{k_j-1} \Eta_{2J-1}}="G'k1Eta2J1"};
 {(0,-40)*{H^j_1 \Lambda^J_1}="H1La1"};
 {(50,-40)*{H'^j_1 \Lambda^J_1}="H'1La1"};
 {(0,-50)*{H^j_{k_j-1} \Lambda^J_1}="Hk1La1"};
 {(50,-50)*{H'^j_{k_j-1} \Lambda^J_1}="H'k1La1"};
 {(0,-60)*{D^j \Lambda^J_1}="DLa1"};
 {(50,-60)*{D'^j \Lambda^J_1}="D'La1"};
 {(50,-70)*{B^j \Lambda^J_1}="BLa1"};
 {(100,-70)*{B'^j_1 \Lambda^J_1}="B'La1"};
 {(50,-80)*{F^j_1 \Lambda^J_1}="F1La1"};
 {(100,-80)*{F'^j_1 \Lambda^J_1}="F'1La1"};
 {(50,-90)*{F^j_{k_j-1} \Lambda^J_1}="Fk1La1"};
 {(100,-90)*{F'^j_{k_j-1} \Lambda^J_1}="F'k1La1"};
 {(0,-110)*{H^j_1 \Lambda^J_{L_J}}="H1LaL"};
 {(50,-110)*{H'^j_1 \Lambda^J_{L_J}}="H'1LaL"};
 {(0,-120)*{H^j_{k_j-1} \Lambda^J_{L_J}}="Hk1LaL"};
 {(50,-120)*{H'^j_{k_j-1} \Lambda^J_{L_J}}="H'k1LaL"};
 {(0,-130)*{D^j \Lambda^J_{L_J}}="DLaL"};
 {(50,-130)*{D'^j \Lambda^J_{L_J}}="D'LaL"};
 {(50,-140)*{B^j \Lambda^J_{L_J}}="BLaL"};
 {(100,-140)*{B'^j_1 \Lambda^J_{L_J}}="B'LaL"};
 {(50,-150)*{F^j_1 \Lambda^J_{L_J}}="F1LaL"};
 {(100,-150)*{F'^j_1 \Lambda^J_{L_J}}="F'1LaL"};
 {(50,-160)*{F^j_{k_j-1} \Lambda^J_{L_J}}="Fk1LaL"};
 {(100,-160)*{F'^j_{k_j-1} \Lambda^J_{L_J}}="F'k1LaL"};
 {(50,-180)*{\fbox{$C^j \Eta_{2J}$}}="CEta2J"};
 {(100,-180)*{C'^j \Eta_{2J}}="C'Eta2J"};
 {(50,-190)*{G^j_1 \Eta_{2J}}="G1Eta2J"};
 {(100,-190)*{G'^j_1 \Eta_{2J}}="G'1Eta2J"};
 {(50,-200)*{G^j_{k_j-1} \Eta_{2J}}="Gk1Eta2J"};
 {(100,-200)*{G'^j_{k_j-1} \Eta_{2J}}="G'k1Eta2J"};
%internal to Eta_2J-1
 {\ar "G'1Eta2J1"; "G1Eta2J1"};
 {\ar "G'k1Eta2J1"; "Gk1Eta2J1"};
 {\ar@{.} "G1Eta2J1"; "Gk1Eta2J1"};
 {\ar@{.} "G'1Eta2J1"; "G'k1Eta2J1"};
 {\ar|{j>0} "C'Eta2J1"; "CEta2J1"};
 {\ar "C'Eta2J1"; "G1Eta2J1"};
 {\ar@{:>} "G'1Eta2J1"; "Gk1Eta2J1"};
%internal to Eta_2J
 {\ar "G'1Eta2J"; "G1Eta2J"};
 {\ar "G'k1Eta2J"; "Gk1Eta2J"};
 {\ar@{.} "G1Eta2J"; "Gk1Eta2J"};
 {\ar@{.} "G'1Eta2J"; "G'k1Eta2J"};
 {\ar|{j>0} "C'Eta2J"; "CEta2J"};
 {\ar "C'Eta2J"; "G1Eta2J"};
 {\ar@{:>} "G'1Eta2J"; "Gk1Eta2J"};
%internal to Lambda_1
 {\ar "B'La1"; "BLa1"};
 {\ar "F'1La1"; "F1La1"};
 {\ar "F'k1La1"; "Fk1La1"};
 {\ar "H'1La1"; "H1La1"};
 {\ar "H'k1La1"; "Hk1La1"};
 {\ar "D'La1"; "DLa1"};
 {\ar@{.} "F1La1"; "Fk1La1"};
 {\ar@{.} "F'1La1"; "F'k1La1"};
 {\ar@{.} "H1La1"; "Hk1La1"};
 {\ar@{.} "H'1La1"; "H'k1La1"};
%internal to Lambda_L
 {\ar "B'LaL"; "BLaL"};
 {\ar "F'1LaL"; "F1LaL"};
 {\ar "F'k1LaL"; "Fk1LaL"};
 {\ar "H'1LaL"; "H1LaL"};
 {\ar "H'k1LaL"; "Hk1LaL"};
 {\ar "D'LaL"; "DLaL"};
 {\ar@{.} "F1LaL"; "Fk1LaL"};
 {\ar@{.} "F'1LaL"; "F'k1LaL"};
 {\ar@{.} "H1LaL"; "Hk1LaL"};
 {\ar@{.} "H'1LaL"; "H'k1LaL"};
%Eta_{2J-1} -sigma_{3}-> Lambda_1
 {\ar "C'Eta2J1"; "F1La1"};
 {\ar "C'Eta2J1"; "H'1La1"};
 {\ar "G1Eta2J1"; "H1La1"};
 {\ar "G'1Eta2J1"; "H'1La1"};
 {\ar "Gk1Eta2J1"; "Hk1La1"};
 {\ar "G'k1Eta2J1"; "H'k1La1"};
 {\ar@{:>} "G'1Eta2J1"; "Fk1La1"};
 {\ar@{:>} "G'1Eta2J1"; "H'k1La1"};
 {\ar "G'k1Eta2J1"; "D'La1"};
 {\ar "Gk1Eta2J1"; "DLa1"};
 {\ar|(0.57){j=0} "C'Eta2J1"; "BLa1"};
%Lambda_1 - -sigma_{23}- -> Lambda_L
 {\ar@{-->} "B'La1";"F1LaL"};
 {\ar@{-->} "B'La1";"H'1LaL"};
 {\ar@{-->} "BLa1";"H1LaL"};
 {\ar@{-->} "F'1La1"; "H'1LaL"};
 {\ar@{-->} "F1La1"; "H1LaL"};
 {\ar@{-->} "F'k1La1"; "H'k1LaL"};
 {\ar@{-->} "Fk1La1"; "Hk1LaL"};
 {\ar@{:>} "F'1La1"; "Fk1LaL"};
 {\ar@{:>} "F'1La1"; "H'k1LaL"};
 {\ar@{:>} "F1La1"; "Hk1LaL"};
 {\ar@{-->} "F'k1La1"; "D'LaL"};
 {\ar@{-->} "Fk1La1"; "DLaL"};
 {\ar@{-->}|(0.57){j=0} "B'La1"; "BLaL"};
%Lambda_L -sigma_2-> Eta_{2J}
 {\ar@/_1.25pc/|(0.35){j>0} "BLaL"+DL; "CEta2J"+UL};
 {\ar@/^1.5pc/ "B'LaL"+DR; "C'Eta2J"+UR};
 {\ar@/_2.5pc/ "BLaL"+DL; "G1Eta2J"+UL};
 {\ar@/_1.5pc/ "F1LaL"+DL; "G1Eta2J"+UL};
 {\ar@/^1.5pc/ "F'1LaL"+DR; "G'1Eta2J"+UR};
 {\ar@/_1.5pc/ "Fk1LaL"+DL; "Gk1Eta2J"+UL};
 {\ar@/^1.5pc/ "F'k1LaL"+DR; "G'k1Eta2J"+UR};
 {\ar@{:>}@<2.5ex>@/^1.5pc/ "F1LaL"; "Gk1Eta2J"};
 {\ar@{:>}@<2.5ex>@/^3pc/ "H'1LaL"; "Gk1Eta2J"};
\endxy
\]
\end{tiny}
\caption{The subspace $Z_{\hor}^{J,j}$, corresponding to a horizontal stable
chain $\Eta_{2J-1} \xrightarrow{\sigma_3} \Lambda^J_1 \xrightarrow{\sigma_{23}}
\cdots \xrightarrow{\sigma_{23}} \Lambda^J_{L_J} \xrightarrow{\sigma_2}
\Eta_{2J}$.} \label{fig:Zhor}
\end{figure}

As before, it is easy to verify that the differentials on the tensor product
coming from $m_1$ and $m_2$ multiplications in Figures \ref{fig:Qj} and
\ref{fig:Q0,s<2tau} respect the splitting \eqref{eq:Zsplit}. These
differentials are illustrated in Figures \ref{fig:Zver} through
\ref{fig:Zunst-tequal}. Note that we obtain slightly different differentials
depending on whether $j=0$ or $j>0$. The double-dotted arrows correspond to the
dashed arrows in Figures \ref{fig:Pver} through \ref{fig:Punst,s=2tau}: for
instance, in Figure \ref{fig:Zver}, the double-dotted arrow from $E'^j_1
\Xi_{2J}$ to $H'^j_{k_j-1} \Kappa^J_1$ really means that there are
differentials $E'^j_i \Xi_{2J} \to H'^j_{i+1} \Kappa^J_1$ for $i = 1, \dots,
k_j-2$.

\begin{figure} \centering
\begin{tiny}
\[
\xy
 {(50,0)*{\fbox{$C^j \Eta_0$}}="CEta0"};
 {(100,0)*{C'^j \Eta_0}="C'Eta0"};
 {(50,-10)*{G^j_1 \Eta_0}="G1Eta0"};
 {(100,-10)*{G'^j_1 \Eta_0}="G'1Eta0"};
 {(50,-20)*{G^j_{k_j-1} \Eta_0}="Gk1Eta0"};
 {(100,-20)*{G'^j_{k_j-1} \Eta_0}="G'k1Eta0"};
 {(0,-40)*{H^j_1 \Gamma_1}="H1Ga1"};
 {(50,-40)*{H'^j_1 \Gamma_1}="H'1Ga1"};
 {(0,-50)*{H^j_{k_j-1} \Gamma_1}="Hk1Ga1"};
 {(50,-50)*{H'^j_{k_j-1} \Gamma_1}="H'k1Ga1"};
 {(0,-60)*{D^j \Gamma_1}="DGa1"};
 {(50,-60)*{D'^j \Gamma_1}="D'Ga1"};
 {(50,-70)*{B^j \Gamma_1}="BGa1"};
 {(100,-70)*{B'^j_1 \Gamma_1}="B'Ga1"};
 {(50,-80)*{F^j_1 \Gamma_1}="F1Ga1"};
 {(100,-80)*{F'^j_1 \Gamma_1}="F'1Ga1"};
 {(50,-90)*{F^j_{k_j-1} \Gamma_1}="Fk1Ga1"};
 {(100,-90)*{F'^j_{k_j-1} \Gamma_1}="F'k1Ga1"};
 {(0,-110)*{H^j_1 \Gamma_R}="H1GaR"};
 {(50,-110)*{H'^j_1 \Gamma_R}="H'1GaR"};
 {(0,-120)*{H^j_{k_j-1} \Gamma_R}="Hk1GaR"};
 {(50,-120)*{H'^j_{k_j-1} \Gamma_R}="H'k1GaR"};
 {(0,-130)*{D^j \Gamma_R}="DGaR"};
 {(50,-130)*{D'^j \Gamma_R}="D'GaR"};
 {(50,-140)*{B^j \Gamma_R}="BGaR"};
 {(100,-140)*{B'^j_1 \Gamma_R}="B'GaR"};
 {(50,-150)*{F^j_1 \Gamma_R}="F1GaR"};
 {(100,-150)*{F'^j_1 \Gamma_R}="F'1GaR"};
 {(50,-160)*{F^j_{k_j-1} \Gamma_R}="Fk1GaR"};
 {(100,-160)*{F'^j_{k_j-1} \Gamma_R}="F'k1GaR"};
 {(50,-180)*{A^j \Xi_0}="AXi0"};
 {(100,-180)*{A'^j \Xi_0}="A'Xi0"};
 {(50,-190)*{E^j_1 \Xi_0}="E1Xi0"};
 {(100,-190)*{E'^j_1 \Xi_0}="E'1Xi0"};
 {(50,-200)*{E^j_{k_j-1} \Xi_0}="Ek1Xi0"};
 {(100,-200)*{E'^j_{k_j-1} \Xi_0}="E'k1Xi0"};
%internal to Eta_0
 {\ar "G'1Eta0"; "G1Eta0"};
 {\ar "G'k1Eta0"; "Gk1Eta0"};
 {\ar@{.} "G1Eta0"; "Gk1Eta0"};
 {\ar@{.} "G'1Eta0"; "G'k1Eta0"};
 {\ar|{j>0} "C'Eta0"; "CEta0"};
 {\ar "C'Eta0"; "G1Eta0"};
 {\ar@{:>} "G'1Eta0"; "Gk1Eta0"};
%internal to Xi_0
 {\ar "A'Xi0"; "AXi0"};
 {\ar "E'1Xi0"; "E1Xi0"};
 {\ar "E'k1Xi0"; "Ek1Xi0"};
 {\ar@{.} "E1Xi0"; "Ek1Xi0"};
 {\ar@{.} "E'1Xi0"; "E'k1Xi0"};
%internal to Gamma_1
 {\ar "B'Ga1"; "BGa1"};
 {\ar "F'1Ga1"; "F1Ga1"};
 {\ar "F'k1Ga1"; "Fk1Ga1"};
 {\ar "H'1Ga1"; "H1Ga1"};
 {\ar "H'k1Ga1"; "Hk1Ga1"};
 {\ar "D'Ga1"; "DGa1"};
 {\ar@{.} "F1Ga1"; "Fk1Ga1"};
 {\ar@{.} "F'1Ga1"; "F'k1Ga1"};
 {\ar@{.} "H1Ga1"; "Hk1Ga1"};
 {\ar@{.} "H'1Ga1"; "H'k1Ga1"};
%internal to Gamma_L
 {\ar "B'GaR"; "BGaR"};
 {\ar "F'1GaR"; "F1GaR"};
 {\ar "F'k1GaR"; "Fk1GaR"};
 {\ar "H'1GaR"; "H1GaR"};
 {\ar "H'k1GaR"; "Hk1GaR"};
 {\ar "D'GaR"; "DGaR"};
 {\ar@{.} "F1GaR"; "Fk1GaR"};
 {\ar@{.} "F'1GaR"; "F'k1GaR"};
 {\ar@{.} "H1GaR"; "Hk1GaR"};
 {\ar@{.} "H'1GaR"; "H'k1GaR"};
%Eta_0 -sigma_{3}-> Gamma_1
 {\ar "C'Eta0"; "F1Ga1"};
 {\ar "C'Eta0"; "H'1Ga1"};
 {\ar "G1Eta0"; "H1Ga1"};
 {\ar "G'1Eta0"; "H'1Ga1"};
 {\ar "Gk1Eta0"; "Hk1Ga1"};
 {\ar "G'k1Eta0"; "H'k1Ga1"};
 {\ar@{:>} "G'1Eta0"; "Fk1Ga1"};
 {\ar@{:>} "G'1Eta0"; "H'k1Ga1"};
 {\ar "G'k1Eta0"; "D'Ga1"};
 {\ar "Gk1Eta0"; "DGa1"};
 {\ar|(0.57){j=0} "C'Eta0"; "BGa1"};
%Gamma_1 - -sigma_{23}- -> Gamma_L
 {\ar@{-->}|(0.57){j=0} "B'Ga1";"BGaR"};
 {\ar@{-->} "B'Ga1";"F1GaR"};
 {\ar@{-->} "B'Ga1";"H'1GaR"};
 {\ar@{-->} "BGa1";"H1GaR"};
 {\ar@{-->} "F'1Ga1"; "H'1GaR"};
 {\ar@{-->} "F1Ga1"; "H1GaR"};
 {\ar@{-->} "F'k1Ga1"; "H'k1GaR"};
 {\ar@{-->} "Fk1Ga1"; "Hk1GaR"};
 {\ar@{:>} "F'1Ga1"; "Fk1GaR"};
 {\ar@{:>} "F'1Ga1"; "H'k1GaR"};
 {\ar@{:>} "F1Ga1"; "Hk1GaR"};
 {\ar@{-->} "F'k1Ga1"; "D'GaR"};
 {\ar@{-->} "Fk1Ga1"; "DGaR"};
%Xi_0 -rho_1-> Gamma_K
 {\ar@/^1.5pc/ "AXi0"+UL;"BGaR"+DL};
 {\ar@/_1.5pc/ "A'Xi0"+UR;"B'GaR"+DR};
 {\ar@/^1.5pc/ "E1Xi0"+UL;"F1GaR"+DL};
 {\ar@/_1.5pc/ "E'1Xi0"+UR;"F'1GaR"+DR};
 {\ar@/^1.5pc/ "Ek1Xi0"+UL;"Fk1GaR"+DL};
 {\ar@/_1.5pc/ "E'k1Xi0"+UR;"F'k1GaR"+DR};
\endxy
\]
\end{tiny}
\caption{The subspace $Z_{\unst}^j$ when $t < 2\tau(K)$, corresponding to the
unstable chain $\Eta_0 \xrightarrow{\sigma_3} \Gamma_1
\xrightarrow{\sigma_{23}} \cdots \xrightarrow{\sigma_{23}} \Gamma_R
\xleftarrow{\sigma_1} \Xi_0$.} \label{fig:Zunst-tless}
\end{figure}

\begin{figure} \centering
\begin{tiny}
\[
\xy
 {(50,0)*{A^j \Xi_0}="AXi0"};
 {(100,0)*{A'^j \Xi_0}="A'Xi0"};
 {(50,-10)*{E^j_1 \Xi_0}="E1Xi0"};
 {(100,-10)*{E'^j_1 \Xi_0}="E'1Xi0"};
 {(50,-20)*{E^j_{k_j-1} \Xi_0}="Ek1Xi0"};
 {(100,-20)*{E'^j_{k_j-1} \Xi_0}="E'k1Xi0"};
 {(0,-40)*{H^j_1 \Gamma_1}="H1Ga1"};
 {(50,-40)*{H'^j_1 \Gamma_1}="H'1Ga1"};
 {(0,-50)*{H^j_{k_j-1} \Gamma_1}="Hk1Ga1"};
 {(50,-50)*{H'^j_{k_j-1} \Gamma_1}="H'k1Ga1"};
 {(0,-60)*{D^j \Gamma_1}="DGa1"};
 {(50,-60)*{D'^j \Gamma_1}="D'Ga1"};
 {(50,-70)*{B^j \Gamma_1}="BGa1"};
 {(100,-70)*{B'^j \Gamma_1}="B'Ga1"};
 {(50,-80)*{F^j_1 \Gamma_1}="F1Ga1"};
 {(100,-80)*{F'^j_1 \Gamma_1}="F'1Ga1"};
 {(50,-90)*{F^j_{k_j-1} \Gamma_1}="Fk1Ga1"};
 {(100,-90)*{F'^j_{k_j-1} \Gamma_1}="F'k1Ga1"};
 {(0,-110)*{H^j_1 \Gamma_R}="H1GaR"};
 {(50,-110)*{H'^j_1 \Gamma_R}="H'1GaR"};
 {(0,-120)*{H^j_{k_j-1} \Gamma_R}="Hk1GaR"};
 {(50,-120)*{H'^j_{k_j-1} \Gamma_R}="H'k1GaR"};
 {(0,-130)*{D^j \Gamma_R}="DGaR"};
 {(50,-130)*{D'^j \Gamma_R}="D'GaR"};
 {(50,-140)*{B^j \Gamma_R}="BGaR"};
 {(100,-140)*{B'^j \Gamma_R}="B'GaR"};
 {(50,-150)*{F^j_1 \Gamma_R}="F1GaR"};
 {(100,-150)*{F'^j_1 \Gamma_R}="F'1GaR"};
 {(50,-160)*{F^j_{k_j-1} \Gamma_R}="Fk1GaR"};
 {(100,-160)*{F'^j_{k_j-1} \Gamma_R}="F'k1GaR"};
 {(50,-180)*{\fbox{$C^j \Eta_0$}}="CEta0"};
 {(100,-180)*{C'^j \Eta_0}="C'Eta0"};
 {(50,-190)*{G^j_1 \Eta_0}="G1Eta0"};
 {(100,-190)*{G'^j_1 \Eta_0}="G'1Eta0"};
 {(50,-200)*{G^j_{k_j-1} \Eta_0}="Gk1Eta0"};
 {(100,-200)*{G'^j_{k_j-1} \Eta_0}="G'k1Eta0"};
%internal to Xi_2J
 {\ar "A'Xi0"; "AXi0"};
 {\ar "E'1Xi0"; "E1Xi0"};
 {\ar "E'k1Xi0"; "Ek1Xi0"};
 {\ar@{.} "E1Xi0"; "Ek1Xi0"};
 {\ar@{.} "E'1Xi0"; "E'k1Xi0"};
%internal to Kappa_1
 {\ar "B'Ga1"; "BGa1"};
 {\ar "F'1Ga1"; "F1Ga1"};
 {\ar "F'k1Ga1"; "Fk1Ga1"};
 {\ar "H'1Ga1"; "H1Ga1"};
 {\ar "H'k1Ga1"; "Hk1Ga1"};
 {\ar "D'Ga1"; "DGa1"};
 {\ar@{.} "F1Ga1"; "Fk1Ga1"};
 {\ar@{.} "F'1Ga1"; "F'k1Ga1"};
 {\ar@{.} "H1Ga1"; "Hk1Ga1"};
 {\ar@{.} "H'1Ga1"; "H'k1Ga1"};
%internal to Kappa_K
 {\ar "B'GaR"; "BGaR"};
 {\ar "F'1GaR"; "F1GaR"};
 {\ar "F'k1GaR"; "Fk1GaR"};
 {\ar "H'1GaR"; "H1GaR"};
 {\ar "H'k1GaR"; "Hk1GaR"};
 {\ar "D'GaR"; "DGaR"};
 {\ar@{.} "F1GaR"; "Fk1GaR"};
 {\ar@{.} "F'1GaR"; "F'k1GaR"};
 {\ar@{.} "H1GaR"; "Hk1GaR"};
 {\ar@{.} "H'1GaR"; "H'k1GaR"};
%internal to Eta_0
 {\ar "G'1Eta0"; "G1Eta0"};
 {\ar "G'k1Eta0"; "Gk1Eta0"};
 {\ar@{.} "G1Eta0"; "Gk1Eta0"};
 {\ar@{.} "G'1Eta0"; "G'k1Eta0"};
 {\ar|{j>0} "C'Eta0"; "CEta0"};
 {\ar "C'Eta0"; "G1Eta0"};
 {\ar@{:>} "G'1Eta0"; "Gk1Eta0"};
%Xi_0 -sigma_{123}-> Gamma_1
 {\ar|(0.57){j=0} "A'Xi0";"BGa1"};
 {\ar "A'Xi0";"F1Ga1"};
 {\ar "A'Xi0";"H'1Ga1"};
 {\ar "AXi0";"H1Ga1"};
 {\ar "E'1Xi0"; "H'1Ga1"};
 {\ar "E1Xi0"; "H1Ga1"};
 {\ar "E'k1Xi0"; "H'k1Ga1"};
 {\ar "Ek1Xi0"; "Hk1Ga1"};
 {\ar@{:>} "E'1Xi0"; "Fk1Ga1"};
 {\ar@{:>} "E'1Xi0"; "H'k1Ga1"};
 {\ar@{:>} "E1Xi0"; "Hk1Ga1"};
 {\ar "E'k1Xi0"; "D'Ga1"};
 {\ar "Ek1Xi0"; "DGa1"};
%Gamma_1 - -sigma_{23}- -> Gamma_R
 {\ar@{-->}|(0.57){j=0} "B'Ga1";"BGaR"};
 {\ar@{-->} "B'Ga1";"F1GaR"};
 {\ar@{-->} "B'Ga1";"H'1GaR"};
 {\ar@{-->} "BGa1";"H1GaR"};
 {\ar@{-->} "F'1Ga1"; "H'1GaR"};
 {\ar@{-->} "F1Ga1"; "H1GaR"};
 {\ar@{-->} "F'k1Ga1"; "H'k1GaR"};
 {\ar@{-->} "Fk1Ga1"; "Hk1GaR"};
 {\ar@{:>} "F'1Ga1"; "Fk1GaR"};
 {\ar@{:>} "F'1Ga1"; "H'k1GaR"};
 {\ar@{:>} "F1Ga1"; "Hk1GaR"};
 {\ar@{-->} "F'k1Ga1"; "D'GaR"};
 {\ar@{-->} "Fk1Ga1"; "DGaR"};
%Lambda_L -sigma_2-> Eta_0
 {\ar@/_1.25pc/|(0.35){j>0} "BGaR"+DL; "CEta0"+UL};
 {\ar@/^1.5pc/ "B'GaR"+DR; "C'Eta0"+UR};
 {\ar@/_2.5pc/ "BGaR"+DL; "G1Eta0"+UL};
 {\ar@/_1.5pc/ "F1GaR"+DL; "G1Eta0"+UL};
 {\ar@/^1.5pc/ "F'1GaR"+DR; "G'1Eta0"+UR};
 {\ar@/_1.5pc/ "Fk1GaR"+DL; "Gk1Eta0"+UL};
 {\ar@/^1.5pc/ "F'k1GaR"+DR; "G'k1Eta0"+UR};
 {\ar@{:>}@<2.5ex>@/^1.5pc/ "F1GaR"; "Gk1Eta0"};
 {\ar@{:>}@<2.5ex>@/^3pc/ "H'1GaR"; "Gk1Eta0"};
\endxy
\]
\end{tiny}
\caption{The subspace $Z_{\unst}^j$ when $t > 2\tau(K)$, corresponding to the
unstable chain $\Xi_0 \xrightarrow{\sigma_{123}} \Gamma_1
\xrightarrow{\sigma_{23}} \cdots \xrightarrow{\sigma_{23}} \Gamma_R
\xrightarrow{\sigma_2} \Eta_0$.} \label{fig:Zunst-tgreater}
\end{figure}

\begin{figure}
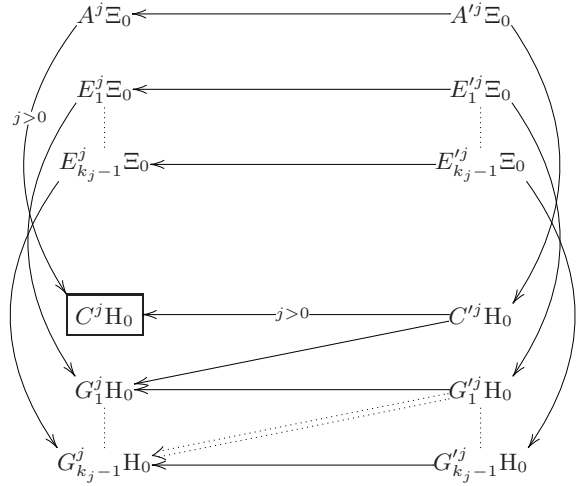
 \centering
\begin{scriptsize}
\[
\xy
 {(0,0)*{}="blank"};
 {(50,0)*{A^j \Xi_0}="AXi0"};
 {(100,0)*{A'^j \Xi_0}="A'Xi0"};
 {(50,-10)*{E^j_1 \Xi_0}="E1Xi0"};
 {(100,-10)*{E'^j_1 \Xi_0}="E'1Xi0"};
 {(50,-20)*{E^j_{k_j-1} \Xi_0}="Ek1Xi0"};
 {(100,-20)*{E'^j_{k_j-1} \Xi_0}="E'k1Xi0"};
 {(50,-40)*{\fbox{$C^j \Eta_0$}}="CEta0"};
 {(100,-40)*{C'^j \Eta_0}="C'Eta0"};
 {(50,-50)*{G^j_1 \Eta_0}="G1Eta0"};
 {(100,-50)*{G'^j_1 \Eta_0}="G'1Eta0"};
 {(50,-60)*{G^j_{k_j-1} \Eta_0}="Gk1Eta0"};
 {(100,-60)*{G'^j_{k_j-1} \Eta_0}="G'k1Eta0"};
%internal to Xi_0
 {\ar "A'Xi0"; "AXi0"};
 {\ar "E'1Xi0"; "E1Xi0"};
 {\ar "E'k1Xi0"; "Ek1Xi0"};
 {\ar@{.} "E1Xi0"; "Ek1Xi0"};
 {\ar@{.} "E'1Xi0"; "E'k1Xi0"};
%internal to Eta_0
 {\ar "G'1Eta0"; "G1Eta0"};
 {\ar "G'k1Eta0"; "Gk1Eta0"};
 {\ar@{.} "G1Eta0"; "Gk1Eta0"};
 {\ar@{.} "G'1Eta0"; "G'k1Eta0"};
 {\ar|{j>0} "C'Eta0"; "CEta0"};
 {\ar "C'Eta0"; "G1Eta0"};
 {\ar@{:>} "G'1Eta0"; "Gk1Eta0"};
%Xi_0 -sigma_{12}-> Eta_0
 {\ar@/_1.5pc/|(0.35){j>0} "AXi0"+DL; "CEta0"+UL};
 {\ar@/^1.5pc/ "A'Xi0"+DR; "C'Eta0"+UR};
 {\ar@/_1.5pc/ "E1Xi0"+DL; "G1Eta0"+UL};
 {\ar@/^1.5pc/ "E'1Xi0"+DR; "G'1Eta0"+UR};
 {\ar@/_1.5pc/ "Ek1Xi0"+DL; "Gk1Eta0"+UL};
 {\ar@/^1.5pc/ "E'k1Xi0"+DR; "G'k1Eta0"+UR};
\endxy
\]
\end{scriptsize}
\caption{The subspace $Z_{\unst}^j$ when $t = 2\tau(K)$, corresponding to the
unstable chain $\Xi_0 \xrightarrow{\sigma_{12}} \Eta_0$.}
\label{fig:Zunst-tequal}
\end{figure}

Next, we must consider the differentials coming from the remaining
multiplications on $Q$. First, we look at differentials that respect the
splitting. If $k_j>1$, the relevant multiplications on $Q^j$ are:
\begin{alignat*}{4}
 A^j & \xrightarrow{\sigma_{123} \sigma_2} G^j_1 &\qquad\quad
 B^j & \xrightarrow{\sigma_{23} \sigma_2} G^j_1 &\qquad\quad
 C^j & \xrightarrow{\sigma_3 \sigma_2} G^j_1 &\qquad&
 \text{if } j>0 \\
 A'^j & \xrightarrow{\sigma_{123} \sigma_{23}} D'^j &
 B'^j & \xrightarrow{\sigma_{23} \sigma_{23}} D'^j &
 C'^j & \xrightarrow{\sigma_3 \sigma_{23}} D'^j &&
 \text{if } k_j=1 \\
 E'^j_{k_j-1} & \xrightarrow{\sigma_{123} \sigma_{23}} D'^j &
 F'^j_{k_j-1} & \xrightarrow{\sigma_{23} \sigma_{23}} D'^j &
 G'^j_{k_j-1} & \xrightarrow{\sigma_3 \sigma_{23}} D'^j &&
 \text{if } k_j>1.
\end{alignat*}
Therefore:
\begin{itemize}
\item In $Z^{J,j}_{\ver}$, if $K_J>1$, there are differentials
$E'^j_{k_{j-1}} \Xi_{2J} \to D'^j \Kappa^J_2$ and $F'^j_{k_{j-1}} \Kappa^J_I
\to D'^j \Kappa^J_{I+2}$.

\item In $Z^{J,j}_{\hor}$, if $K_J>1$, there are differentials
$G'^j_{k_{j-1}} \Eta_{2J-1} \to D'^j \Lambda^J_2$ and $F'^j_{k_{j-1}}
\Lambda^J_I \to D'^j \Lambda^J_{I+2}$. Additionally, when $j>0$, there are
differentials $C^j \Eta_{2J-1} \to G^j_1 \Eta_{2J}$ if $K_J=1$, and $B^j
\Lambda^J_{K_J-1} \to G^j_1 \Eta_{2J}$ if $K_J>1$.

\item In $Z^j_{\unst}$, if $t<2\tau(K)-1$, there are differentials
$G'^j_{k_{j-1}} \Eta_0 \to D'^j \Gamma_2$ and $F'^j_{k_{j-1}} \Gamma_I \to D'^j
\Gamma_{I+2}$. If $t=2\tau(K)+1$, there are differentials $A^j \Xi_0 \to G^j_1
\Eta_0$ for $j>0$. If $t>2 \tau(K)+1$, there are differentials $E'^j_{k_{j-1}}
\Xi_0 \to D'^j \Gamma_2$ and $F'^j_{k_{j-1}} \Gamma_I \to D'^j \Gamma_{I+2}$
for all $j$, and $B^j \Gamma_{R-1} \to G^j_1 \Eta_0$ for $j>0$.
\end{itemize}

Next, we may have some differentials that preserve the decomposition
\[
\bigoplus_J Z_{\ver}^{J,*} \oplus \bigoplus_J Z_{\hor}^{J,*} \oplus Z_{\unst}^*
\]
but which come from the multiplications on $Q$ that do not preserve the
splitting $Q = \bigoplus_{j=0}^n Q^j$, shown in Figures
\ref{fig:rho3rho2rho123}, \ref{fig:rho3rho2rho1red-peven}, and
\ref{fig:rho3rho2rho1red-podd}. The resulting differentials are shown in Table
\ref{table:uph,vph}. In each line that involves expressions like $\Kappa^J_I$,
$\Lambda^J_I$, and $\Gamma_I$, we assume that $K_J$, $L_J$, or $R$ is
sufficiently large for the indices to make sense and that $I$ ranges over
appropriate bounds. The symbol $^*$ denotes both primed and unprimed symbols;
thus, for instance, the notation $A^{*j} \Xi_{2J} \to D^{*h} \Kappa^J_2$ means
that there are differentials $A^j \Xi_{2J} \to D^h \Kappa^J_2$ and $A'^j
\Xi_{2J} \to D'^h \Kappa^J_2$. Additionally, note that if $k_h=1$, then we
replace $H^h_1$ by $D^h$ where it appears; if $k_j=1$, we replace
$E'^j_{k_j-1}$, $F'^j_{k_j-1}$, and $G'^j_{k_j-1}$ by $A'^j$, $B'^j$, and
$C'^j$, respectively.

Notice that almost all of the differentials in Table \ref{table:uph,vph} drop
the filtration level by a nonzero amount. The two exceptions are $A^j
\Xi_{2J-1} \to D'^h \Kappa^J_{K_J}$ and $A^j \Xi_0 \to D'^h \Gamma_R$ in the
second column.

\begin{table}
\renewcommand\arraycolsep{4pt}
\begin{small}
\[
\begin{array}{|c|c|c|c|} \hline
& u_{2j,h}=1, \ j,h>0 & v_{2j,h}=1, \ j>0 & v_{2j-1,h}=1 \text{ or } w_h=1 \\
\hline Z_{\ver}^J &
 \begin{aligned}
  A'^j \Xi_{2J} &\to H^h_1 \Kappa^J_2 \\
  B'^j \Kappa^J_I &\to H^h_1 \Kappa^J_{I+1}  \\
  B'^j \Kappa^J_I &\to H^h_1 \Kappa^J_{I+2} \\
 \end{aligned}
&
 \begin{aligned}
  A^{*j} \Xi_{2J} &\to D^{*h} \Kappa^J_2 \\
  B^{*j} \Kappa^J_I &\to D^{*h} \Kappa^J_I  \\
  B^{*j} \Kappa^J_I &\to D^{*h} \Kappa^J_{I+2} \\
  A^j \Xi_{2J-1} &\to D'^h \Kappa^J_{K_J} \\
 \end{aligned}
&
 \begin{aligned}
  E'^j_{k_j-1} \Xi_{2J} &\to D^h \Kappa^J_1 \\
  E'^j_{k_j-1} \Xi_{2J} &\to D^h \Kappa^J_3 \\
  F'^j_{k_j-1} \Kappa^J_I &\to D^h \Kappa^J_{I+1} \\
  F'^j_{k_j-1} \Kappa^J_I &\to D^h \Kappa^J_{I+3} \\
 \end{aligned}
\\ \hline
Z_{\hor}^J &
 \begin{aligned}
  \text{If } L_J=1: & \\
  C'^j \Eta_{2J-1} & \to G^h_1 \Eta_{2J-1} \\
  \text{If } L_J>1: & \\
  C'^j \Eta_{2J-1} &\to H^h_1 \Lambda^J_2 \\
  B'^j \Lambda^J_I &\to H^h_1 \Lambda^J_{I+1}  \\
  B'^j \Lambda^J_I &\to H^h_1 \Lambda^J_{I+2}  \\
  B'^j \Lambda^J_{K_J-1} &\to G^h_1 \Eta_{2J-1} \\
 \end{aligned}
&
 \begin{aligned}
  C^{*j} \Eta_{2J-1} &\to D^{*J} \Lambda^J_1 \\
  C^{*j} \Eta_{2J-1} &\to D^{*J} \Lambda^J_2 \\
  B^{*j} \Lambda^J_I &\to D^{*h} \Lambda^J_I  \\
  B^{*j} \Lambda^J_I &\to D^{*h} \Lambda^J_{I+2} \\
 \end{aligned}
&
 \begin{aligned}
  G'^j_{k_j-1} \Eta_{2J-1} &\to D^h \Lambda^J_1 \\
  G'^j_{k_j-1} \Eta_{2J-1} &\to D^h \Lambda^J_3 \\
  F'^j_{k_j-1} \Lambda^J_I &\to D^h \Lambda^J_{I+1} \\
  F'^j_{k_j-1} \Lambda^J_I &\to D^h \Lambda^J_{I+3} \\
 \end{aligned}
\\ \hline
 \begin{array}{c}
  Z_{\unst}, \\ t<2\tau(K)
 \end{array}
&
 \begin{aligned}
  C'^j \Eta_0 &\to H^h_1 \Gamma_2 \\
  B'^j \Gamma_I &\to H^h_1 \Gamma_{I+1}  \\
  B'^j \Gamma_I &\to H^h_1 \Gamma_{I+2}  \\
 \end{aligned}
&
 \begin{aligned}
  C^{*j} \Eta_0 &\to D^{*J} \Gamma_1 \\
  C^{*j} \Eta_0 &\to D^{*J} \Gamma_2  \\
  B^{*j} \Gamma_I &\to D^{*h} \Gamma_I  \\
  B^{*j} \Gamma_I &\to D^{*h} \Gamma_{I+2} \\
  A^j \Xi_0 &\to D'^h \Gamma_R \\
 \end{aligned}
&
 \begin{aligned}
  G'^j_{k_j-1} \Eta_0 &\to D^h \Gamma_1 \\
  G'^j_{k_j-1} \Eta_0 &\to D^h \Gamma_3  \\
  F'^j_{k_j-1} \Gamma_I &\to D^h \Gamma_{I+1} \\
  F'^j_{k_j-1} \Gamma_I &\to D^h \Gamma_{I+3} \\
 \end{aligned}
\\ \hline
 \begin{array}{c}
  Z_{\unst}, \\ t>2\tau(K)
 \end{array}
&
 \begin{aligned}
  \text{If } R=1: & \\
  A'^j \Xi_0 & \to G^h_1 \Eta_0 \\
  \text{If } R>1: & \\
  A'^j \Xi_0 &\to H^h_1 \Gamma_2 \\
  B'^j \Gamma_I &\to H^h_1 \Gamma_{I+1}  \\
  B'^j \Gamma_I &\to H^h_1 \Gamma_{I+2} \\
  B'^j \Gamma_{R-1} &\to G^h_1 \Eta_0
 \end{aligned}
&
 \begin{aligned}
  A^{*j} \Xi_0 &\to D^{*h} \Gamma_2 \\
  B^{*j} \Gamma_I &\to D^{*h} \Gamma_I  \\
  B^{*j} \Gamma_I &\to D^{*h} \Gamma_{I+2} \\
 \end{aligned}
&
 \begin{aligned}
  E'^j_{k_j-1} \Xi_0 &\to D^h \Gamma_1 \\
  E'^j_{k_j-1} \Xi_0 &\to D^h \Gamma_3 \\
  F'^j_{k_j-1} \Gamma_I &\to D^h \Gamma_{I+1} \\
  F'^j_{k_j-1} \Gamma_I &\to D^h \Gamma_{I+3} \\
 \end{aligned}
\\ \hline
\end{array}
\]
\end{small}
\renewcommand\arraycolsep{5pt}
\caption{Differentials arising from the multiplications in Figures
\ref{fig:rho3rho2rho123}, \ref{fig:rho3rho2rho1red-peven}, and
\ref{fig:rho3rho2rho1red-podd}.} \label{table:uph,vph}
\end{table}

Finally, we must look at differentials that do not respect the splitting at
all. Notice that the sequence $\sigma_3 \sigma_2 \sigma_1$ occurs several times
in Figures \ref{fig:Qj} and \ref{fig:Q0,s<2tau}, and the sequences $\sigma_3
\sigma_2 \sigma_{12}$ and $\sigma_3 \sigma_2 \sigma_{123}$ occur in Equations
\eqref{eq:Pver-longer} and \eqref{eq:Punst-longer}, and these are the only such
sequences that appear. More precisely, in $Q^j$ with $k_j>1$, we have the
following multiplications:
\begin{equation} \label{eq:sigma3sigma2sigma1}
\begin{alignedat}{3}
 A^j & \xrightarrow{\sigma_3 \sigma_2 \sigma_1} H^j_1 &\qquad\quad A'^j & \xrightarrow{\sigma_3 \sigma_2 \sigma_1} H'^j_1 &\qquad& \\
 E^j_i & \xrightarrow{\sigma_3 \sigma_2 \sigma_1} H^j_{i+1} & E'^j_i & \xrightarrow{\sigma_3 \sigma_2 \sigma_1} H'^j_{i+1} & (i&=1, \dots, k_j-2)  \\
 E^j_{k_j-1} & \xrightarrow{\sigma_3 \sigma_2 \sigma_1} D^j & E'^j_{k_j-1} & \xrightarrow{\sigma_3 \sigma_2 \sigma_1} D'^j && \\
 A'^j &\xrightarrow{\sigma_3 \sigma_2 \sigma_{12}} G^j_2 & A'^j &\xrightarrow{\sigma_3 \sigma_2 \sigma_{123}} H^j_2 && \\
 E'^j_i &\xrightarrow{\sigma_3 \sigma_2 \sigma_{12}} G^j_{i+2} & E'^j_i &\xrightarrow{\sigma_3 \sigma_2 \sigma_{123}} H^j_{i+2} & (i&=1, \dots, k_j-3) \\
 && E'^j_{k_j-2} &\xrightarrow{\sigma_3 \sigma_2 \sigma_{123}} D^j && \\
\end{alignedat}
\end{equation}
If $k_j=1$, then we simply have $A^j \xrightarrow{\sigma_3 \sigma_2 \sigma_1}
D^j$ and $A'^j \xrightarrow{\sigma_3 \sigma_2 \sigma_1} D'^j$. Finally, from
Figure \ref{fig:rho3rho2rho1red-peven}, if $v_{2j,h}=1$, then there are
multiplications $A^j \xrightarrow{\sigma_3 \sigma_2 \sigma_1} D^h$ and $A'^j
\xrightarrow{\sigma_3 \sigma_2 \sigma_1} D'^h$.

Notice that all of these multiplications come out of $A^j$, $A'^j$, $E^j_i$, or
$E'^j_i$, all of which are paired with $\{\Xi_0, \dots,  \Xi_{2N}\}$ rather
than $\{\Eta_0, \dots, \Eta_{2N}\}$ in \eqref{eq:Zbases}. It follows that each
group $Z^{J,*}_{\hor}$ is actually a direct summand as a chain complex. We
shall see that the generator of the total homology comes from $\bigoplus_J
Z_{\ver}^{J,*} \oplus Z^*_{\unst}$, so we may ignore each of the $Z^{J,*}_{\hor}$ summands as before.
Furthermore, if we define $U_{P,M}$, $V_{P,M}$, and $W_P$ analogously to
$u_{p,h}$, $v_{p,m}$, and $w_p$ above, then we obtain differentials from $A^j
\Xi_P$, $A'^j\Xi_P$, $E^j_i\Xi_P$, and/or $E'^j_i\Xi_P$ to elements of
$Z^M_{\ver}$ and $Z_{\unst}$ whenever $U_{P,M}$, $V_{P,M}$, or $W_P$ is
nonzero. Specifically:

\begin{itemize}
\item If $V_{P,M} = 1$, then there are differentials
\begin{equation} \label{eq:VPM}
\begin{alignedat}{3}
A^j \Xi_P & \to H^j_1 \Kappa^M_{K_M} &\qquad\quad
A'^j \Xi_P & \to H'^j_1 \Kappa^M_{K_M} &\qquad& \\
E^j_i \Xi_P & \to H^j_{i+1} \Kappa^M_{K_M} &
E'^j_i \Xi_P & \to H'^j_{i+1} \Kappa^M_{K_M} & (i&=1, \dots, k_j-2) \\
E^j_{k_j-1} \Xi_P & \to D^j \Kappa^M_{K_M} &
E'^j_{k_j-1} \Xi_P & \to D'^j \Kappa^M_{K_M} && \\
\end{alignedat}
\end{equation}
if $k_j>1$, and $A^j \Xi_P \to D^j \Kappa^M_{K_M}$ and $A'^j \Xi_P \to D'^j
\Kappa^M_{K_M}$ if $k_j=1$. Also, if $v_{2j,h}=1$, then there are differentials
$A^j \Xi_P \to D^h \Kappa^M_{K_M}$ and $A'^j \Xi_P \to D'^h \Kappa^M_{K_M}$.

Similarly, if $W_P=1$ and $t<2\tau(K)$, then we obtain similar differentials
going into $Z_{\unst}$, replacing $\Kappa^M_{K_M}$ by $\Gamma_R$.

\item If $U_{P,M}=1$, then there are differentials
\begin{equation} \label{eq:UPM}
\begin{split}
A'^j \Xi_P &\to H^j_2 \Kappa^M_1 \\
E'^j_i \Xi_P &\to H^j_{i+2} \Kappa^M_1 \ (i=1, \dots, k_j-3) \\
E'^j_{k_j-2} \Xi_P &\to D^j \Kappa^M_1. \\
\end{split}
\end{equation}
Similarly, if $W_P=1$ and $t>2\tau(K)$, then we obtain similar differentials
going into $Z_{\unst}$, replacing $\Kappa^M_1$ by $\Gamma_1$.

\item Finally, if $W_P=1$ and $t=2\tau(K)$, there are differentials
\begin{equation} \label{eq:WP}
\begin{split}
A'^j \Xi_P &\to G^j_2 \Eta_0 \\
E'^j_i \Xi_P &\to G^j_{i+2} \Eta_0 \ (i=1, \dots, k_j-3) \\
\end{split}
\end{equation}
\end{itemize}

\subsection{Computation of \texorpdfstring{$\tau(D_{J,s}(K,t))$}
{\texttau(D\textunderscore\{J,s\}(K,t)) }}

We now describe the edge cancellations that occur in each of the pieces. Recall
that we must cancel edges in increasing order of the amount by which they drop
filtration level. We shall see that a single generator survives. The filtration
level of this generator, by definition, is $\tau(D_{J,s}(K,t))$.

We start by canceling the filtration-preserving edges in $Z_{\ver}^{J,j}$. Note
that there are are no other edges into $B'^j \Kappa^J_{K_J}$ or $F'^j_i
\Kappa^J_{K_J}$, so eliminating the edges coming from these does not introduce
any new edges. If $V_{2J-1},M = 1$, or if $W_{2J-1}=1$ and $t<2\tau(K)$, then
canceling the edges $A^j \Xi_{2J-1} \to B^j \Kappa^J_{K_J}$ and $E^j_i \to
F^j_k \Kappa^J_{K_J}$ introduces some new edges, which all reduce filtration
level by $2$. Note also that the filtration-preserving edges $A^j \Xi_{2j-1}
\to D'^h \Kappa^J_{K_J}$ $(j>0)$ in Table \ref{table:uph,vph} are eliminated,
since $B^j \Kappa^J_{K_j}$ has no other incoming edges when $j>0$.

In $Z_{\unst}^j$, when $t<2\tau(K)$, we perform the same cancellations as in
$Z_{\ver}^{J,j}$, \emph{mutatis mutandis}. When $t>2 \tau(K)$, there are $2k_j$
filtration-preserving edges to cancel when $j>0$ (namely, $B^{*j} \Gamma^R \to
C^{*j} \Eta_0$ and $F^{*j}_i \Gamma^R \to G^{*j}_i \Eta_0$ for $i=1, \dots,
k_j-1$), but only $2k_0-1$ such edges in $Z_{\unst}^0$, since the generator
$C_0 \Eta_0$ does not exist. Thus, the generator $B^0 \Gamma^R$ survives after
these cancellations. Also, note that canceling $B^j \Gamma^R \to C^j \Eta_0$
and $F^j_i \Gamma^R \to G^j_i \Eta_0$ may introduce some new differentials
using the arrows in Table \ref{table:uph,vph}, but they all filtration level by
$2$.

When $t=2\tau(K)$, the only generator in $Z^0_{\unst}$ that survives is $A^0
\Xi_0$. Notice, however, that by \eqref{eq:VPM}, there is a differential $A^0
\Xi_0 \to H^0_1 \Kappa^M_{K_M}$ for any $M$ with $V_{0,M}=1$. All the
generators of $Z^j_{\unst}$ for $j>0$ are canceled.

We have now canceled all edges that preserve the filtration level, so we now
begin canceling differentials that drop filtration level by $1$. Specifically,
starting at the top of Figure \ref{fig:Zver} and working down, we cancel every edge of the form $X' \to X$, where $X'$ and $X$ are two generators in the same row (e.g., $A'^j \Xi_{2J} \to A^j \Xi_{2J}$). We use the following key
observations:
\begin{itemize}
\item If $X$ is in filtration level $0$ and $X'$ is in level $1$, then $X$ has no
other incoming edges, since by induction we have already eliminated everything
above $X$ and $X'$, and Table \ref{table:uph,vph} and Equations \eqref{eq:VPM}
and \eqref{eq:UPM} contain no differentials that go into $A^j \Xi_{2J}$, $E^j_i
\Xi_{2J}$, $B^j \Kappa^J_I$, or $F^j_i \Kappa^J_I$ from elsewhere.

\item If $X$ is in filtration level $-1$ and $X'$ is in level $0$, then $X'$
has no other outgoing edges, since Table \ref{table:uph,vph} and Equations
\eqref{eq:VPM} and \eqref{eq:UPM} contain no differentials that go out of
$H'^j_i \Kappa^J_I$ or $D'^j \Kappa^J_I$.
\end{itemize}
Thus, we can completely cancel $Z_{\ver}^{J,j}$.

If $t=2\tau(K)$, we have now eliminated all generators except $A^0\Xi_0$, which
is in filtration level $0$, so $\tau(D_{J_s}(K,t)) = 0$ when $s<2\tau(J)$ and
$t=2\tau(K)$.

If $t>2\tau(K)$, we proceed with $Z_{\unst}^j$ just as with $Z_{\ver}^{J,j}$.
When $j>0$, all generators in $Z_{\unst}^j$ cancel; when $j=0$, the one surviving generator
is $B^0 \Gamma_R$, which is in filtration level $0$. Thus, $\tau(D_{J_s}(K,t))
= 0$ when $s<2\tau(J)$ and $t>2\tau(K)$.

If $t<2\tau(K)$, when $j>0$, we start by canceling $C'^j \Eta_0 \to C^j \Eta_0$
and proceeding downward in Figure \ref{fig:Zunst-tless}, as before, eliminating
all generators. When $j=0$, we start by canceling $G'^j_1 \Eta_0 \to G^j_1
\Eta_0$ and proceed downward, and we thus see that the only surviving generator
is $C'^j \Eta_0$, which is in filtration level $1$. Thus, $\tau(D_{J_s}(K,t)) =
1$ when $s<2\tau(J)$ and $t<2\tau(K)$.

Finally, we must return to the case where $s=2\tau(J)$. Recall that $Q_0$ in
this case consists of three generators, all in filtration level $0$, as in
\eqref{eq:Q0,s=2tau}. For $j>0$, the definitions of $Z^{J,j}_{\ver}$,
$Z^{J,j}_{\hor}$, and $Z^j_{\unst}$ go through the same way, and we see again
that all of the resulting generators eventually cancel. It follows that the
surviving generator must be in filtration level $0$, so $\tau(D_{J,s}(K,t))=0$
whenever $s=2\tau(K)$. \qed

%\end{document}

\section{Other results regarding \texorpdfstring{$D_{J,s}(K,t)$}{D\textunderscore\{J,s\}(K,t)}} \label{sec:doubling-topology}
Prior to Hedden's complete computation of $\HFK$ and $\tau$ of all twisted
Whitehead doubles \cite{HeddenWhitehead}, Livingston and Naik
\cite{LivingstonNaikDoubled} used the formal properties of $\tau$ to understand
the asymptotic behavior of $\tau$ for large values of the twisting parameter. They proved:
\begin{theorem} \label{thm:LivingstonNaik}
Suppose $\nu$ is any homomorphism from the smooth knot concordance group to
$\Z$ with the properties that $\abs{\nu(K)} \le g_4(K)$ and $\nu(T_{p,q}) =
(p-1)(q-1)/2$, where $p,q>0$ and $T_{p,q}$ denotes the $(p,q)$ torus knot. Then
for any knot $K$, there exists $t_\nu(K) \in \Z$ such that
\[
\nu(Wh_+(K,t))= \begin{cases} 1 & t\le t_\nu(K) \\ 0 & t > t_\nu(K) \end{cases}
\]
and $TB(K) \le t_\nu(K) < -TB(-K)$ (where $TB(K)$ denotes the maximal
Thurston-Bennequin number of $K$).
\end{theorem}
Two invariants satisfying the hypotheses of Theorem \ref{thm:LivingstonNaik}
are $\tau(K)$ and $-s(K)/2$, a renormalization of Rasmussen's concordance
invariant $s(K)$ \cite{RasmussenGenus}. Around the same time, Hedden and Ording
\cite{HeddenOrding} proved that these two invariants are not equal by showing
that $\tau(Wh_+(T_{2,3},2)) = 0$ while $s(Wh_+(T_{2,3},2)) = -2$, disproving a
conjecture of Rasmussen. Later, Hedden \cite{HeddenWhitehead} showed that
$t_\tau(K) = 2\tau(K)-1$ for any knot $K$. Finding a general formula for the
$s$ invariant of Whitehead doubles remains an open question.

We may extend the techniques of Livingston and Naik to study knots of the form
$D_{J,s}(K,t)$ as well.
\begin{proposition} \label{prop:nu}
Let $\nu$ be an invariant satisfying the hypotheses of Theorem
\ref{thm:LivingstonNaik}, and fix knots $J$ and $K$.
\begin{enumerate}
\item If $s \le TB(J)$ and $t \le TB(K)$, then $\nu(D_{J,s}(K,t))=1$. If $s \ge
-TB(-J)$ and $t \ge -TB(-K)$, then $\nu(D_{J,s}(K,t))=-1$.
\item For fixed $s$ (resp.~$t$), the function $t \mapsto \nu(D_{J,s}(K,t))$
(resp.~$s \mapsto \nu(D_{J,s}(K,t))$) is non-increasing and has as its image
either $\{-1,0\}$, $\{0\}$, or $\{0,1\}$.
\end{enumerate}
\end{proposition}

\begin{proof}
The proof is very similar to that of \cite[Theorem 2]{LivingstonNaikDoubled}.

Let $A(J,s)$ denote an annulus in $S^3$, embedded along $J$ with framing $s$, and define $A(K,t)$ analogously. We may obtain a Seifert surface for $D_{J,s}(K,t)$ as a plumbing $A(J,s)*A(K,t)$. By results of Rudolph \cite{RudolphAnnuli, RudolphPlumbing}, when $s \le TB(J)$ and $t \le TB(K)$, the annuli $A(J,s)$ and $A(K,t)$ are \emph{quasipositive surfaces}, so $A(J,s)*A(K,t)$ is also a quasipositive surface. Thus, $D_{J,s}(K,t)$ is a strongly quasipositive knot with genus $1$, and hence $\nu(D_{J,s}(K,t))=1$ \cite[Theorem 4]{LivingstonComputations}. Mirroring gives the second half of (1).

The non-increasing statement in (2) follows from the fact that $D_{J,s}(K,t)$
is obtained from $D_{J,s-1}(K,t)$ or $D_{J,s}(K,t-1)$ by changing a positive
crossing to a negative crossing, which can only preserve or decrease $\nu$
\cite[Corollary 3]{LivingstonComputations}. Also, since $D_{J,s}(K,t)$ is related to
$D_{J,s'}(K,t)$ or $D_{J,s}(K,t')$ (for any $s'$ or $t'$) by changing the number of twists in a band of a Seifert surface, each of the two functions can assume at most two values, either $\{-1,0\}$ or $\{0,1\}$ \cite[Corollary 5]{LivingstonNaikDoubled}. Finally, we rule out the possibility that either of
the functions in (1) is constant and nonzero. Suppose, without loss of
generality, that $\nu(D_{J,s}(K,t))=1$ for a fixed $s$ and all $t$. In
particular, $\nu(D_{J,s}(K,-TB(-K)))=1$. On the other hand,
$\nu(D_{J,-TB(-J)}(K,-TB(-K)))=-1$, which contradicts the fact that the image
of the function $s \mapsto \nu(D_{J,s}(K,-TB(-K)))$ contains at most two
consecutive integers.
\end{proof}

It follows that if $s \le TB(J)$ and $t \ge -TB(-K)$, or $s \ge -TB(-J)$ and $t
\le TB(K)$, then $\nu(D_{J,s}(K,t))=0$. Thus, for large absolute values of $s$
and $t$, $\nu(D_{J,s}(K,t)) = \tau(D_{J,s}(K,t))$. On the other hand, the
behavior of $\nu(D_{J,s}(K,t)$ for small $s$ and $t$ (specifically, when
$TB(J)<s<-TB(-J)$ or $TB(K) < t < -TB(-K)$) may be more complicated than the
simple behavior of $\tau$ given by Theorem \ref{thm:taudjskt}.

In another direction, we may also look for instances when $D_{J,s}(K,t)$ is
actually smoothly slice. The following proposition generalizes Casson's
argument \cite[page~227]{KauffmanOnKnots} that the $p(p+1)$-twisted positive
Whitehead double of the $(p,p+1)$ torus knot is smoothly slice. For an oriented
knot $K$ and relatively prime integers $p,q$, let $C_{p,q}(K)$ denote the
$(p,q)$-cable of $K$. (Note that $C_{p,q}(K)^r = C_{-p,-q}(K) = C_{p,q}(K^r)$
and $\overline{C_{p,q}(K)} = C_{p,-q}(\bar K)$.)

\begin{figure}
\centering \psfrag{J,s}[cc][cc]{{\small $J,s$}} \psfrag{K,t}[cc][cc]{{\small
$K,t$}} \psfrag{gp}{$\gamma_p$}
\includegraphics{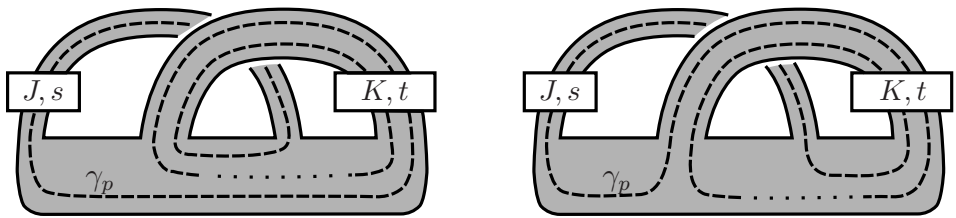}
\caption[The Seifert surface $F$ with the curve $\gamma_p$.] {The Seifert
surface $F$ with the curve $\gamma_p$, in the cases where $p<0$ (left) and
$p>0$ (right).} \label{fig:seifert}
\end{figure}

\begin{proposition} \label{prop:smoothlyslice}
Let $K$ be any knot, and let $p,t \in \Z$. If $J$ is any knot that is smoothly
concordant to $-C_{p,pt\pm 1}(K)$, then $D_{J,-p(pt \pm 1)}(K,t)$ is smoothly
slice.
\end{proposition}

\begin{proof}
Let $F$ be the Seifert surface for $D_{J,s}(K,t)$ shown in Figure
\ref{fig:seifert}, and let $\gamma_p$ be a curve that winds once around the
band tied into $J$ and $p$ times around the band tied into $K$, as indicated.
The knot type of $\gamma_p$ is $C_{p, pt+1}(K)$, and the surface framing on
$\gamma_p$ is $s + p+ p^2t$. Thus, if $J$ is smoothly concordant to
$-C_{p,pt+1}(K)$ and $s = -p(pt+1)$, we may surger $F$ along $\gamma_p$ in
$D^4$ along a smooth slice disk for $J \# {C_{p,pt+1}(K)}$, resulting in a
smooth slice disk for $D_{J,s}(K,t)$.

If we reverse the crossing between the two bands of $F$, we obtain the result
with the opposite signs.
\end{proof}

Proposition \ref{prop:smoothlyslice} is quite interesting in light of work of Hom \cite{HomTau}, who found a general formula for the $\tau$ invariant of all cable knots in terms of $p$, $q$, $\tau(K)$, and an invariant $\epsilon(K) \in \{-1, 0,1\}$ that depends solely on the knot Floer complex of
$K$. She proved:
\begin{theorem} \label{thm:Hom}
Let $K$ be a knot, and let $p>0$. Then:
\begin{itemize}
\item If $\epsilon(K) =1$, then $\tau(C_{p,q}(K)) = p \tau(K) + \frac12
(p-1)(q-1)$ for all $q$.
\item If $\epsilon(K) =-1$, then $\tau(C_{p,q}(K)) = p \tau(K) + \frac12
(p-1)(q+1)$ for all $q$.
\item If $\epsilon(K)=0$, then $\tau(K)=0$, and
\[
\tau(C_{p,q}(K)) = \begin{cases} \frac12(p-1)(q+1)  & q<0 \\
\frac12(p-1)(q-1) & q>0. \end{cases}
\]
\end{itemize}
\end{theorem}

We may use Theorem \ref{thm:Hom} to compute the value of $\tau$ for the cable
knots appearing in Proposition \ref{prop:smoothlyslice}, where we take $t =
2\tau(K)$.

\begin{corollary} \label{cor:boundarycases}
For any knot $K$, if either $\epsilon(K)\ge 0$ and $p>0$, or $\epsilon(K) \le
0$ and $p<0$, there exists a knot $J$ such that $D_{J, 2\tau(J)-p}(K,
2\tau(K))$ is smoothly slice, while $\tau(D_{J, 2\tau(J)-p}(K, 2\tau(K)-
\frac{p}{\abs{p}})) \ne 0$.
\end{corollary}

\begin{proof}
Suppose that $\epsilon(K)=1$ and $p> 0$. Set $J = -C_{p,2p\tau(K)+1}(K)$, so
that:
\[
\begin{split}
2\tau(J) - p &= - 2\tau(C_{p,2p\tau(K)+1} (K)) -p \\
&=  -2p\tau(K) - (p-1)(2p\tau(K)) - p \\
&= -2p^2\tau(K)-p \\
&= -p(2p\tau(K)+1).
\end{split}
\]
By Proposition \ref{prop:smoothlyslice}, $D_{J, 2\tau(J)-p} (K, 2\tau(K))$ is
smoothly slice. On the other hand, $\tau(D_{J,2\tau(J)-p}(K, 2\tau(K)-1) = 1$
by Theorem \ref{thm:taudjskt}. The case where $\epsilon(K)=-1$ and $p<0$
follows by mirroring, since $\epsilon(\bar K) = - \epsilon(K)$. Finally, if
$\epsilon(K)=0$, we set $J = -C_{p, 1}(K)$ if $p>0$ and $J = -C_{-p,-1}$ if
$p<0$.
\end{proof}

Theorem \ref{thm:taudjskt} says that the set $\{(s,t) \in \Z^2 \mid
D_{J,s}(K,t)=0\}$ always has the same shape for any $J$ and $K$, up to
translation: the union of the second and fourth quadrants of the $\Z^2$
lattice, including both axes. Corollary \ref{cor:boundarycases} implies that
any point on the boundary of this region may be realized by a smoothly slice
knot $D_{J,s}(K,t)$ for suitable choices of $J$ and $K$.

Finally, recall that the main idea of the proof of Theorem \ref{thm:taudjskt}
is that only the form of the unstable chains in $\CFD(\XX_J^s)$ and
$\CFD(\XX_K^t)$ matters for the computation of $\tau(D_{J,s}(K,t))$. Petkova
\cite{PetkovaCables} and Hom \cite{HomTau} have observed similar behavior in
using bordered Heegaard Floer homology to compute $\tau(C_{p,q}(K))$. The
invariant $\epsilon(K)$ defined by Hom describes the structure of the part of
$\CFD(\XX_K^t)$ ``near'' the unstable chain. Specifically, when we take
vertically and horizontally reduced bases $\{\tilde\xi_0, \dots,
\tilde\xi_{2n}\}$ and $\{\tilde\eta_0, \dots, \tilde\eta_{2n}\}$ for
$\operatorname{CFK}^-(K)$, we may arrange that $\tilde\xi_0 = \tilde\eta_i$ for
some $i$. The cases $\epsilon(K) = 1$, $\epsilon(K)=-1$, and $\epsilon(K)=0$
correspond, respectively, to whether $i$ is even and positive, odd and
positive, or zero. Within each case, Hom showed that only the form of the
unstable chain matters for computing $\tau(C_{p,q}(K))$. It is an interesting
question whether the behavior of $\tau$ for more general classes of satellite
knots can be described in this way.

\appendix
\section{Notes on the computation} \label{sec:appendix}
This section provides further details about the \emph{Mathematica} notebooks used for the computations in Section \ref{sec:heegaard}. The computation makes use of two packages that are designed to be useful for Heegaard Floer homology computations: \verb"HeegaardDiagram.nb", which is used to find the index-1 domains in a Heegaard diagram, and \verb"TorusAlgebra.nb", which provides algebraic tools for working with $\AA_\infty$-modules and type-$D$ structures over the torus algebra. In the hope that these tools will be of use to other researchers in the future, Sections \ref{subsec:HeegaardDiagram.nb} and \ref{subsec:TorusAlgebra.nb} provide brief user's guides. The computation of $\CFAA(\YY,B_3,0)$ is found in the notebook \verb"Borromean.nb", which is described in Section \ref{subsec:Borromean.nb}. All of these notebooks are available as ancillary materials in this article's arXiv folder: \url{http://arxiv.org/src/1008.3349/anc}.

\subsection{\texorpdfstring{\texttt{HeegaardDiagram.nb}}{HeegaardDiagram.nb}} \label{subsec:HeegaardDiagram.nb}

The file \verb"HeegaardDiagram.nb" contains functions for finding all of the positive domains of index $1$ in a Heegaard diagram in order to compute $\CF$ of a closed $3$-manifold or $\CFD$ of a bordered manifold with all boundary components of genus $1$. This is fundamentally a problem of solving systems of linear equations. Obviously, the program does not determine whether or not a given domain supports holomorphic representatives, but it generates a list of domains whose moduli spaces that can then be checked by hand (or with computer assistance, as in the present setting).

\subsubsection*{Preliminary input}

The basic input for the program consists of the following data. We label all of the intersection points between the $\alpha$ and $\beta$ curves $x_1, \dots, x_n$, and the regions of the diagram $R_1, \dots, R_m$, so that the basepointed region is $R_1$. Assume the Heegaard diagram has $k$ boundary components and that the genus is $g$. (Typically $k \in \{0, 1, 2\}$, but in principle we can take $k$ to be arbitrary.)

To input the Heegaard diagram, one must specify the following data:
\begin{itemize}
\item \verb"pointdata" is an $n \times 4$ array recording which regions are incident to the intersection points. Specifically, if we draw a neighborhood of $x_i$ such that the $\alpha$ curve is the horizontal axis, the $\beta$ curve is the vertical axis, and the four quadrants are $R_a$, $R_b$, $R_c$, and $R_d$ (starting with the upper-right quadrant and going counterclockwise), then the $i\Th$ entry of \verb"pointdata" is the list \verb"{a,b,c,d}".

\item \verb"euler" is a list of the Euler measures of the regions $R_1, \dots, R_m$.

\item \verb"boundary" is a $k \times 3$ array that records which regions abut each boundary component of the diagram. For a closed diagram, this is simply the empty list. For a bordered diagram, if the three regions (other than $R_1$) adjacent to the $i\Th$ boundary component are $R_a$, $R_b$, and $R_c$ --- adjacent to the $\rho_1, \rho_2, \rho_3$ arcs, respectively, following the labeling convention for $\CFD$ --- then the $i\Th$ entry of \verb"boundary" is \verb"{a,b,c}".

\item \verb"alphaarcs" is a $k \times 2$ array of lists recording which intersection points are on each of the $\alpha$ arcs. Specifically, if $\alpha_1^i$ and $\alpha_2^i$ are the arcs that abut the $i\Th$ boundary component, where $\alpha_1^i$ meets the boundary between the basepoint and $\rho_1$, and between $\rho_2$ and $\rho_3$, the $(i,1)\Th$ (resp.~$(i,2)\Th$) entry of \verb"alphaarcs" is the list of the indices of the intersection points on $\alpha_1^i$ (resp.~$\alpha_2^i$).

\item \verb"alphacircles" is a list of length $g-k$ whose $i\Th$ entry is the list of the indices of the intersection points on the $i\Th$ $\alpha$ circle.

\item \verb"beta" is a list of length $g$ whose $i\Th$ entry is the list of the indices of the intersection points on the $g\Th$ $\beta$ circle.
\end{itemize}

The first command to execute is \verb"Initialize[]", which initializes the values of several other variables that are used throughout the computation.

A generator $(x_{i_1}, \dots, x_{i_g})$, where $x_i \in \beta_i$, is represented by the ordered $g$-tuple of indices $(i_1, \dots, i_g)$. Note that the points should be written in the same order as the $\beta$ circles. The list of all generators should be stored as \verb"generators", which \verb"Initialize[]" does automatically (using the function \verb"FindGenerators[]"), but one can also define such a list manually.

\subsubsection*{Finding domains}

The following functions are used in finding domains:

\begin{itemize}
\item \verb"PositiveDomain[from, to, constraints]" takes three arguments: two generators \verb"from" and \verb"to", given in the format described above, and a list \verb"constraints" indicating the constraints imposed on the multiplicities of certain regions. The latter is a list of pairs $(a_i,b_i)$, where each pair corresponds to requiring the region $R_{a_i}$ to have multiplicity $b_i$. We always require that there are enough constraints so that there is at most one solution to the linear equations \eqref{eq:pointconditions}. (If there are no provincial periodic domains in the diagram, then constraining the multiplicities of each of the boundary regions is sufficient; compare Lemma \ref{lemma:uniquesolution}.) \verb"PositiveDomain[from, to, constraints]" returns
    a list of length $m$ consisting of the multiplicities of each of the regions in the unique positive domain from \verb"from" to \verb"to" if one exists, and \verb"{}" otherwise.

\item \verb"Index1Domain[from, to, constraints]" likewise returns the unique index-1 positive domain from \verb"from" to \verb"to" if one exists, and \verb"{}" otherwise. \emph{Note:} Because of the way the index of a domain in a bordered diagram is computed (see page \pageref{pageref:iota}), this function only works properly if the multiplicity of each boundary region is either $0$ or $1$. By Proposition \ref{prop:chordsallowed}, we need only consider such domains in order to compute $\CFD$.

\item \verb"CoeffsToList[domain]" takes as its argument a positive domain \verb"domain" in the format output by \verb"PositiveDomain" and \verb"Index1Domain" as above, and it outputs the list of which regions have nonzero multiplicity, with repetitions for multiplicities greater than $1$. This format is more convenient for inspecting domains manually, especially when the number of regions is large. The function \verb"ListToCoeffs[domain]" reverses this process.

\item \verb"FindIndex1Domains[constraints]" takes as its argument a list of constraints as above, and it outputs a list of triples \verb"{from, to, domain}" consisting all of the index-1 domains satisfying the given constraints. Specifically, it applies the function \verb"Index1Domain" to every pair of generators (taken from the list \verb"generators"), and applies \verb"CoeffsToList" to each of the outputs. Because this involves solving a large system of linear equations for each pair, it can take a long time to run.

\item \verb"FastPositiveDomain", \verb"FastIndex1Domain", and \verb"FastFindIndex1Domains" are more efficient versions of the functions above. The function \verb"Initialize[]" generates a list of domains connecting each pair of consecutive generators in \verb"generators" and finds a basis for the group of periodic domains. \begin{verbatim}
    FastPositiveDomain[i, j, constraints]
    \end{verbatim}
    finds a positive domain, if one exists, from the \verb"i"$\Th$ entry of \verb"generators" to the \verb"j" entry by a two-step process: first, it finds a ``test domain'' by adding together entries in the preloaded list of domains, and then it solves a system of linear equations to determine what linear combination of periodic domains, if any, can be added to the test domain to give a domain satisfying the needed constraints. This system typically involves far fewer variables than the used in \verb"PositiveDomain".

    The other two functions work analogously. \verb"FastFindIndex1Domains" saves additional time by considering only pairs of generators that occupy the appropriate $\alpha$ arcs to be compatible with the given constraints, as per Proposition \ref{prop:chordsallowed}; this feature can be disabled with the option setting \verb"TypeDOnly -> False".

    \emph{Important note:} At present, \verb"FastFindIndex1Domains" requires that all generators in \verb"generators" represent the same spin$^c$ structure --- i.e., that there is a domain connecting any two generators. If there are multiple spin$^c$ structures, they should each be handled separately, building the list of generators manually each time. The author plans to address this issue in a future version of the program.

\item \verb"AllIndex1Domains[]" finds all of the index-1 domains in the diagram. (This is the only function that a typical user needs to call.) Specifically, it cycles through all possible sets of constraints where each boundary region has multiplicity either $0$ and $1$ and calls \verb"FastFindIndex1Domains" for each set. To impose constraints on regions other than the boundary, \verb"AllIndex1Domains" can take an extra argument \verb"extraregions", a list of tuples $(a_i,b_i,c_i)$, where each one corresponds to letting the multiplicity of region $R_{a_i}$ range from $b_i$ to $c_i$; the default value of \verb"extraregions" is \verb"{}". This option should be used when the diagram contains periodic domains, although some thought is needed to determine the appropriate bounds. The option setting \verb"MonitorProgress->True" provides a progress indicator that indicates the time elapsed and which set of constraints is being considered.
\end{itemize}

\subsection{\texorpdfstring{\texttt{TorusAlgebra.nb}}{TorusAlgebra.nb}} \label{subsec:TorusAlgebra.nb}

The file \verb"TorusAlgebra.nb" contains functions used for computations with $\AA_\infty$ modules, type-$D$ structures, and bimodules over the torus algebra. It includes an implementation of the edge reduction algorithm described in Subsection \ref{subsec:edge}, and it can compute the box tensor product. Although gradings are not discussed in this paper, the package also contains some functionality for working with the non-abelian grading on bordered Floer homology, with certain caveats described below. It makes considerable use of \emph{Mathematica}'s capabilities for pattern-matching and symbolic manipulation.

\subsubsection*{Algebra basics}

The six Reeb elements in $\AA(T^2)$ are represented by the symbols \verb"rho[1]", \verb"rho[2]", \verb"rho[3]", \verb"rho[1,2]", \verb"rho[2,3]", and \verb"rho[1,2,3]". Here \verb"rho" can be any function initialized with the
command
\begin{verbatim}
rho[w___Integer, x_Integer, y_Integer, z___Integer] :=  0 /; x + 1 != y;
\end{verbatim}
This guarantees that an expression containing non-consecutive indices (e.g.
\verb"rho[1,3]") automatically becomes $0$.  The package automatically
initializes the commands \verb"rho", \verb"sigma", \verb"tau", and \verb"phi", along with the actual corresponding Greek letters.

Tensor products of algebra are elements denoted using the $\otimes$ symbol,
which in Mathematica is typed \emph{Esc}~\verb"c *"~\emph{Esc} or
\verb"\[CircleTimes]". One may also use the \verb"CircleTimes" command; thus,
\begin{verbatim}
CircleTimes[rho[3],rho[2]]
\end{verbatim}
is the same as \verb"rho[3]"$\otimes$\verb"rho[2]".

We list the basic commands for working with the algebra, but the typical user
does not need them.

\begin{itemize}
\item \verb"CollectVars" sorts the entries of a tensor product of algebra
elements according to which function they involve. Thus
\verb"CollectVars[rho[2]"$\otimes$\verb"sigma[1]"$\otimes$\verb"rho[1]]"
returns \verb"rho[2]"$\otimes$\verb"rho[1]"$\otimes$\verb"sigma[1]". It does
not change the order of factors using the same function.

\item \verb"Group" takes two arguments: a tensor product of algebra elements
and a list of which algebras are treated as type $A$. The other algebras are
treated as type $D$, meaning that multiplications are evaluated. For instance:
\begin{itemize}
\item 
The expression
\begin{flushleft}
\verb"Group[rho[1]"$\otimes$\verb"rho[2]"$\otimes$\verb"sigma[1], {rho}]"
\end{flushleft}
returns \verb"rho[1]"$\otimes$\verb"rho[2]"$\otimes$\verb"sigma[1]".

\item The expressions
\begin{flushleft}
\verb"Group[rho[1]"$\otimes$\verb"rho[2]"$\otimes$\verb"sigma[1], {sigma}]"
\end{flushleft}
and
\begin{flushleft}
\verb"Group[rho[1]"$\otimes$\verb"rho[2]"$\otimes$\verb"sigma[1], {}]"
\end{flushleft}
both return \verb"rho[1, 2]"$\otimes$\verb"sigma[1]", since $\rho_1 \rho_2 =
\rho_{12}$ in $\AA$.

\item \verb"Group[rho[2]"$\otimes$\verb"rho[1]"$\otimes$\verb"sigma[1], {sigma}]"
returns \verb"0", since $\rho_2 \rho_1 = 0$ in $\AA$.
\end{itemize}

\item \verb"SelectPart" and \verb"SelectRemaining" both take as arguments a
tensor product of algebra elements and the name of one of the algebras.
\verb"SelectPart" returns the factors that use the specified algebra;
\verb"SelectRemaining" returns the factors that do not use that algebra. Thus,
\begin{flushleft}
\verb"SelectPart[rho[1]"$\otimes$\verb"rho[2]"$\otimes$\verb"sigma[1], rho]"
\end{flushleft}
and
\begin{flushleft}
\verb"SelectRemaining[rho[1]"$\otimes$\verb"rho[2]"$\otimes$\verb"sigma[1], sigma]"
\end{flushleft}
both return \verb"rho[1]"$\otimes$\verb"rho[2]".
\end{itemize}

The commands \verb"Group", \verb"SelectPart", and \verb"SelectRemaining" all
distribute over addition in their first arguments.

\subsubsection*{Modules}

A module, bimodule, or even multimodule is represented by a list such as the following, which represents $\CFA$ of a solid torus with a particular framing:
\begin{verbatim}
SolidTorusA = {{{rho, 0}}, {{1}, {2}, {1}},
 {{0, rho[1], 1 + rho[1, 2]}, {0, 0, rho[2]}, {0, 0, 0}},
 {{1/2, {1, 0}}},
 {{0, {0, 0}}, {-1/2, {1/2, -1/2}}, {-1, {0, 0}}}
};
\end{verbatim}
The data for a module are as follows:

\begin{itemize}
\item The first entry records the algebras that act and whether the action is type
$D$ or type $A$. The convention is that \verb"0" means type $A$ and \verb"1" means type $D$.
Thus, the entry \verb"{{rho,0}}" means that we have an $A$ module over a single
copy of the algebra, in which elements are denoted \verb"rho[1]", etc. Likewise
\verb"{{rho,1}, {sigma,1}}" would signify a type $DD$ bimodule, and
\verb"{{rho,0}, {sigma,1}}" would signify a type $AD$ bimodule in which the
type $A$ action uses the \verb"rho" elements and the type $D$ action uses the
\verb"sigma" elements. The order in which the algebras are listed matters
throughout.

\item The second entry records the idempotents of the generators. The convention is
that \verb"1" corresponds to $\iota_0$ and \verb"2" corresponds to $\iota_1$. (This unfortunate convention is needed because \emph{Mathematica} indexes lists from $0$ rather than from $1$.) In the example, the first and third generators are in idempotent $\iota_0$ and the second is in $\iota_1$. For bimodules, we record the idempotents for each action in order.

\item The third entry is a matrix recording the differential or the $\AA_\infty$ multiplications. The convention is that the $(i,j)\Th$ entry records differentials from $x_i$ to $x_j$ (where $\{x_1, \dots, x_n\}$ is a basis). Thus, in the example shown, we have:
\begin{align*}
m_1(x_1) &= x_3 \\
m_2(x_1,\rho_1) &= x_2 \\
m_2(x_1,\rho_{12}) &= x_3 \\
m_2(x_2,\rho_2) &=x_3.
\end{align*}
For higher multiplications in a type $A$ structure, we would represent $m_4(x_i, \rho_3, \rho_2, \rho_1) = x_j$, for instance, with a \verb"rho[3]"$\otimes$\verb"rho[2]"$\otimes$\verb"rho[1]" term in the $(i,j)\Th$ entry.

When working in a bimodule, we must list the elements of the two algebras in the order in which they occur in the first entry.

\item The fourth and fifth entries are used when working with gradings; they may be omitted when gradings are not used. We work with the reduced grading group $G = \{(m;a,b) \mid m,a,b \in \frac12\Z)$; elements are written \verb"{m, {a, b}}". For bimodules, the gradings take the form \verb"{m, {a, b}, {c, d}}". The fourth entry denotes the gradings of a basis for the group of periodic domains, and the fifth entry is a list of coset representatives for all of the generators.

    \emph{Note:} At present, gradings do not work properly for mixed bimodules (type $AD$ or $DA$), but they do work with single modules and with $AA$ or $DD$ bimodules.

\item Extra gradings may be given as additional entries in the module. This is how we generally encode filtrations, such as the filtration on bordered Floer homology induced by a knot.
\end{itemize}

The function \verb"DSquared" can be used to verify that the differential on a type-$D$ module satisfies $\partial^2=0$. Note that for large matrices, this function can take a long time to run. At present, the package does not include an analogous function for verifying the $\AA_\infty$ relations on a type-$A$ module.

The function \verb"CheckGradings" checks whether the gradings on a module are
consistent with the differential, returning either \verb"True" or \verb"False".
It works for single modules and $AA$ and $DD$ bimodules. It also has a
\verb"FullForm" option which, if enabled, lists every single nontrivial
differential and whether or not it is consistent with the gradings.

\subsubsection*{Evaluating tensor products}

The command for evaluating tensor products is \verb"TensorProduct". It takes
four arguments: The first and third are the names of the $A$ and $D$ modules,
respectively. The second argument is a number indicating which algebra on the
first module is being used for the tensor product, and the fourth argument does
the same for the second module.

Thus, for example, if we have two modules of the form
\begin{verbatim}
AABimodule = {{{rho, 0}, {sigma, 0}}, ... }
DModule = {{{rho, 1}, ... }
\end{verbatim}
we compute their tensor product with the command
\begin{verbatim}
TensorProduct[AABimodule, 1, DModule, 1]
\end{verbatim}
resulting in an $A$ module of the form \verb"{{{sigma, 0}}, ... }". By default,
the matrix in the new module will be given as a \verb"SparseArray" object,
unless the \verb"SparseArrayForm" setting is set to \verb"False".

The variable names of the two algebras over which we are tensoring should match each other, and there should be no overlap in the remaining variable names. In the preceding example, to use the second algebra structure on \verb"AABimodule" rather than the first, one might enter
\begin{verbatim}
TensorProduct[AABimodule /. sigma->tau, 2, DModule /. rho->tau, 1]
\end{verbatim}
to change the variable names before evaluating the tensor product. Note, however, that this substitution does not work if the differentials are presented as \verb"SparseArray" objects.

If either of the modules in a tensor product does not contain grading information, include the option setting \verb"Graded -> False".

The following options allow the user to monitor the progress of the computation:
\begin{itemize}
\item \verb"MonitorProgress -> True" gives running updates of which steps are being performed.
\item \verb"ListGradings -> True" explains the computations of the gradings
of the generators in the tensor product.
\item \verb"ListEdges -> True" explains all of the differentials that occur in the
tensor product.
\end{itemize}
Each of these options is set to \verb"False" by default.

If either module has extra gradings or filtrations, you can have them extend to
the tensor product using the \verb"AFiltrations" or \verb"DFiltrations"
options. If the modules are graded, a filtration is typically the sixth entry,
so we might write something like
\begin{verbatim}
TensorProduct[AABimodule, 1, DModule, 1, DFiltrations -> {6}]
\end{verbatim}
to extend a filtration on \verb"DModule" to the tensor product. If filtrations
on both the $A$ and $D$ modules are used, those coming from the $A$ module are
listed first.

\subsubsection*{Reducing modules}

The \verb"ReduceModule" command is used to simplify a module by looking for
entries in the matrix that equal \verb"1". The basic syntax is
\verb"ReduceModule[Module]", where \verb"Module" is the name of the module.

If \verb"Module" does not contain grading information, the option \verb"Graded"
should be set to \verb"False".

If \verb"Module" has a filtration, you can keep track of the filtration as
edges are canceled in increasing order of the amount by which they drop
filtration level. If the filtration information is in the sixth entry in
\verb"Module", for instance, the command is
\begin{verbatim}
ReduceModule[Module, Filtrations -> {6}]
\end{verbatim}
The output is of the form \verb"{M, {{M1, ..., Mk}}}", where \verb"M" is the
module that results from performing all simplifications, and \verb"M1", \dots, \verb"Mk" are the pages of the spectral sequence associated to the filtration. Usually, one is interested in the module \verb"M1", which is
filtered chain homotopy equivalent to \verb"M"; one can obtain this directly
using the command
\begin{verbatim}
ReduceModule[Module, Filtrations -> {6}] [[2, 1, 1]]
\end{verbatim}

With multiple filtrations --- e.g., when computing link Floer homology
--- \verb"ReduceModule" cancels edges in the order of the amount by which
they drop the total filtration level, while keeping track of each individual
filtration. The syntax is, e.g.,
\begin{verbatim}
ReduceModule[Module, Filtrations -> {6, 7}]
\end{verbatim}

The option setting \verb"ListEdges -> True" causes \verb"ReduceModule" to
display every single cancellation that is performed. The option
\verb"MonitorProgress" is an integer, 0 by default. When it is set to a
positive number $n$, \verb"ListEdges" displays an update (with the elapsed
time) after every $n$ cancellations.

\subsubsection*{Examples}

Several useful bordered Heegaard Floer modules are built into in the \verb"TorusAlgebra.nb" package (including a few not discussed in this paper). All of these are given with gradings except where otherwise noted.

\begin{itemize}
\item \verb"SolidTorusA" is the simple $A$ module for a solid torus with a
particular framing, as described above.

\item \verb"IdentityAA" is the $AA$ identity bimodule.

\item \verb"LHTrefoil0D" and \verb"RHTrefoil0D" are the $D$ modules for the
complement of the two trefoils, taken with the $0$-framing.

\item \verb"BorromeanDD" and \verb"BorromeanAA" are the filtered $DD$ and $AA$ bimodules for
the complement of two components of the Borromean rings, with the filtration
induced by the third component, as computed in this paper. (The Maslov gradings on \verb"BorromeanDD" are computed in \verb"Borromean.nb"; those on \verb"BorromeanAA" were computed by hand.)

\item \verb"FramingSwitchDA" and \verb"FramingSwitchAA" are the bimodules for
the diffeomorphism of the torus that takes each slope to its perpendicular
slope. For instance, the following would produce the $D$ module for $\infty$
framing on the left-handed trefoil complement:
\begin{verbatim}
TensorProduct[FramingSwitchDA, 2, LHTrefoil0D /. rho -> sigma, 1]
\end{verbatim}
The module \verb"FramingSwitchDA" currently does not contain grading
information.

\item \verb"SplitBasepointDD" is a bimodule used for computing the link Floer homology of two-component links obtained by taking a knot in each piece of a bordered decomposition. Specifically, given nulhomologous knots $K_1 \subset Y_1$ and $K_2 \subset Y_2$, one would often like to be able to compute the link Floer homology of $K_1 \cup K_2$ in $Y_1 \cup Y_2$. However, $\CFA(Y_1,K_1) \boxtimes \CFD(Y_2,K_2)$ is only an invariant of the bouquet of circles obtained by connecting $K_1$ and $K_2$. The link Floer complex of the actual link is given by
\[
\CFL(Y_1 \cup Y_2, K_1 \cup K_2) \simeq \CFA(Y_1,K_1) \boxtimes (\CFA(Y_2,K_2)
\boxtimes \mathtt{SplitBasepointDD}),
\]
with the $\Z \times \Z$ filtration induced from the filtrations on
$\CFA(Y_1,K_1)$ and $\CFA(Y_2,K_2)$. (The author is grateful to Rumen Zarev for describing this construction.)

\item \verb"T42DD" is the $DD$ bimodule for the exterior of the $(4,-2)$
torus link with the $-2$ framing on both components. If we attach a copy of the
$0$-framed exterior of a knot $K$ to both boundary components, we obtain the
branched double cover of $Wh_+(K)$.
\end{itemize}

\subsection{\texorpdfstring{\texttt{Borromean.nb}}{Borromean.nb}} \label{subsec:Borromean.nb}

We begin by loading all of the functions in \verb"HeegaardDiagram.nb" and \verb"TorusAlgebra.nb" packages, encoding the data of the Heegaard diagram $\HH'$ (Figure \ref{fig:heegaard2}) as described above, and running \verb"Initialize[]". We then run \verb"AllIndex1Domains[]" and store the complete list as \verb"m1list" (m for Maslov). This part of the computation takes approximately $7.5$ minutes to run on the author's Lenovo X220 laptop, so the output is included as a separate cell for rapid pre-loading. (In contrast, finding the domains using \verb"FindIndex1Domains[]" rather than \verb"FastFindIndex1Domains[]", as described in Section \ref{subsec:HeegaardDiagram.nb}, takes over an hour and a half.)

The next step is to partition \verb"m1list" into different families of domains that share the same holomorphic geometry, saving each as a separate list and then deleting that list from \verb"m1list" to ensure no repetitions. For example, as noted in Section \ref{subsec:domains}, any index-1 domain that does not use the regions $R_2$, $R_4$, $R_7$, or $R_8$ (which are the only regions with negative Euler measure) counts for the differential. We form a list \verb"nicediagram" consisting of all such domains and then delete them from \verb"m1list" using the following commands:
\begin{verbatim}
nicediagram = Select[m1list, Intersection[#[[3]], {2,4,7,8}] == {} &];
m1list = Complement[m1list, nicediagram];
\end{verbatim}
Likewise, all instances of the domains $D_1$ and $D_2$ (which count for multiple pairs of generators) can be identified as follows:
\begin{verbatim}
rho23cutbigons =  Select[m1list, #[[3]] ==
  Sort[Join[{7,8,36,37}, Range[19,30], Range[49,52]]] &];
sigma12cutbigons = Select[m1list, #[[3]] ==
  {4,11,17,20,24,25,29,32,35,39,43,47,50} &];
m1list = Complement[m1list, rho23cutbigons, sigma12cutbigons];
\end{verbatim}
Other families of domains with the same geometry can be identified using the command \verb"SubsetQ[list1, list2]" (part of the \verb"HeegaardDiagram.nb" package), which determines whether or not \verb"list2" is a subset of \verb"list1". For instance, we may find one family of domains with the same geometry as $D_3$ using the command
\begin{verbatim}
Select[m1list, SubsetQ[#[[3]], {7,8,36,37}] &&
   SubsetQ[{36,37,7,8,48,49,31,30,18,19,38,10,5,3,13,41,15}, #[[3]]] &]
\end{verbatim}
(Compare the definition of the list \verb"rho23cutrectangles" in \verb"Borromean.nb", which includes these and other domains.) After all of these families have been defined and deleted from \verb"m1list", we verify that the resulting \verb"m1list" is empty, confirming that all of the domains have been classified. The reader can easily verify that every family matches up with the domains described in Section \ref{subsec:heegaarddiagram} and that the counts given in Table \ref{table:domains} are correct.

The domains that count for the differential (as per Proposition \ref{prop:gooddomains}) are stored in the list \verb"gooddomains". We have included a verification that every composition of two domains in \verb"gooddomains" either cancels against another such composition or yields a pair of algebra elements that multiply to $0$ (e.g. $\rho_2 \rho_1$), which proves that the condition \eqref{eq:d-relation} holds for $\CFDD(\HH,0)$.

The module $\CFDD(\HH,0)$ is stored as \verb"bigborromeandd" (using the conventions of Section \ref{subsec:TorusAlgebra.nb}). Although we have not discussed Maslov gradings anywhere in this paper, we can easily compute the relative grading between any two generators directly from the differential, so we include the gradings for future applications. The command \verb"CheckGradings[bigborromeandd]" verifies that the gradings are computed correctly.

We may apply the edge-reduction algorithm to $\CFDD(\HH,0)$ using the \verb"ReduceModule" command, as described above. The resulting module, with 19 generators, is stored as \verb"borromeandd". This step takes under a second to execute. For purely aesthetic reasons, we permute the basis elements to obtain the version described in Theorem \ref{thm:cfdd}. The matrix in the statement of that theorem is taken directly from the output of \emph{Mathematica}. Note that \verb"DSquared[borromeandd]" yields the zero matrix, as expected.

We then use the \verb"TensorProduct" and \verb"ReduceModule" comands to compute $\CFAD(\YY,B_3,0)$ and $\CFAA(\YY,B_3,0)$, which are respectively stored as \verb"borromeanad" and \verb"borromeanaa". Because the grading functions in \verb"TorusAlgebra.nb" do not work for mixed bimodules, we disable the \verb"Graded" option for each function. The matrices in Theorem \ref{thm:cfaa} are taken directly from the output of this computation. 

%
%\section{Poetic conclusion} \label{sec:conclusion}
%\begin{verse}
%Our goal is one whose application's nice \\
%For smooth four-manifold topology: \\
%To tell if certain knots and links are slice \\
%With bordered Heegaard Floer homology.
%
%We seek concordance data that detect \\
%Some links obtained by Whitehead doublings, \\
%As well as knots we get when we infect \\
%Along two of the three Borromean rings.
%
%Some lengthy work with bordered Floer then proves \\
%How $\tau$ for satellites like these is found. \\
%We see, by this result and cov'ring moves, \\
%That smooth slice disks our links can never bound.
%
%The theorem's proved, the dissertation's done, \\
%But all the work ahead has just begun.
%\end{verse}

\bibliography{bibliography}

\providecommand{\bysame}{\leavevmode\hbox to3em{\hrulefill}\thinspace}
\providecommand{\MR}{\relax\ifhmode\unskip\space\fi MR }
\providecommand{\MRhref}[2]{%
  \href{http://www.ams.org/mathscinet-getitem?mr=#1}{#2}
}
\providecommand{\href}[2]{#2}
\begin{thebibliography}{10}

\bibitem{Bizaca}
{\v{Z}}arko Bi{\v{z}}aca, \emph{An explicit family of exotic {C}asson handles},
  Proc. Amer. Math. Soc. \textbf{123} (1995), no.~4, 1297--1302.

\bibitem{FreedmanWhitehead3}
Michael~H. Freedman, \emph{{${\rm Whitehead}_3$} is a ``slice'' link}, Invent.
  Math. \textbf{94} (1988), no.~1, 175--182.

\bibitem{FreedmanQuinn}
Michael~H. Freedman and Frank Quinn, \emph{Topology of 4-manifolds}, Princeton
  Mathematical Series, vol.~39, Princeton University Press, Princeton, NJ,
  1990.

\bibitem{HeddenWhitehead}
Matthew Hedden, \emph{Knot {F}loer homology of {W}hitehead doubles}, Geom.
  Topol. \textbf{11} (2007), 2277--2338.

\bibitem{HeddenOrding}
Matthew Hedden and Philip Ording, \emph{The {O}zsv\'ath-{S}zab\'o and
  {R}asmussen concordance invariants are not equal}, Amer. J. Math.
  \textbf{130} (2008), no.~2, 441--453.

\bibitem{HomTau}
Jennifer Hom, \emph{Bordered {H}eegaard {F}loer homology and the tau-invariant
  of cable knots}, J. Topol. (2013), to appear, \arxiv{1202.1463}.

\bibitem{KauffmanOnKnots}
Louis~H. Kauffman, \emph{On knots}, Annals of Mathematics Studies, vol. 115,
  Princeton University Press, Princeton, NJ, 1987.

\bibitem{KirbyList}
Rob Kirby (ed.), \emph{Problems in low-dimensional topology}, AMS/IP Stud. Adv.
  Math., vol.~2, Amer. Math. Soc., Providence, RI, 1997.

\bibitem{LevineBingWhitehead}
Adam~S. Levine, \emph{Slicing mixed {B}ing-{W}hitehead doubles}, J. Topol.
  \textbf{5} (2012), no.~3, 713--726.

\bibitem{LOTBordered}
Robert Lipshitz, Peter Ozsv{\'a}th, and Dylan Thurston, \emph{Bordered
  {H}eegaard {F}loer homology}, \arxiv{0810.0687}, 2009.

\bibitem{LOTBimodules}
\bysame, \emph{Bimodules in bordered {H}eegaard {F}loer homology},
  \arxiv{1003.0598}, 2010.

\bibitem{LivingstonComputations}
Charles Livingston, \emph{Computations of the {O}zsv\'ath-{S}zab\'o knot
  concordance invariant}, Geom. Topol. \textbf{8} (2004), 735--742
  (electronic).

\bibitem{LivingstonNaikDoubled}
Charles Livingston and Swatee Naik, \emph{{O}zsv\'ath-{S}zab\'o and {R}asmussen
  invariants of doubled knots}, Algebr. Geom. Topol. \textbf{6} (2006),
  651--657 (electronic).

\bibitem{OSz4Genus}
Peter Ozsv{\'a}th and Zolt{\'a}n Szab{\'o}, \emph{Knot {F}loer homology and the
  four-ball genus}, Geom. Topol. \textbf{7} (2003), 615--639 (electronic).

\bibitem{OSzKnot}
\bysame, \emph{Holomorphic disks and knot invariants}, Adv. Math. \textbf{186}
  (2004), no.~1, 58--116.

\bibitem{PetkovaCables}
Ina Petkova, \emph{Cables of thin knots and bordered {H}eegaard {F}loer
  homology}, Quantum Topol. \textbf{4} (2013), no.~4, 377--409.

\bibitem{RasmussenGenus}
Jacob Rasmussen, \emph{Khovanov homology and the slice genus}, Invent. Math.
  \textbf{182} (2010), no.~2, 419--447.

\bibitem{RasmussenThesis}
Jacob~A. Rasmussen, \emph{Floer homology and knot complements}, Ph.D. thesis,
  Harvard University, 2003, \arxiv{math/0509499}.

\bibitem{RudolphAnnuli}
Lee Rudolph, \emph{Quasipositive annuli. ({C}onstructions of quasipositive
  knots and links. {IV})}, J. Knot Theory Ramifications \textbf{1} (1992),
  no.~4, 451--466.

\bibitem{RudolphObstruction}
\bysame, \emph{Quasipositivity as an obstruction to sliceness}, Bull. Amer.
  Math. Soc. (N.S.) \textbf{29} (1993), no.~1, 51--59.

\bibitem{RudolphPlumbing}
\bysame, \emph{Quasipositive plumbing (constructions of quasipositive knots and
  links. {V})}, Proc. Amer. Math. Soc. \textbf{126} (1998), no.~1, 257--267.

\bibitem{SarkarWang}
Sucharit Sarkar and Jiajun Wang, \emph{An algorithm for computing some
  {H}eegaard {F}loer homologies}, Ann. of Math. (2) \textbf{171} (2010), no.~2,
  1213--1236.

\end{thebibliography}
\bibliographystyle{amsplain}

\end{document}